\documentclass[a4paper,11pt,reqno,final]{amsart}
\newcommand{\choarticle}{1}
\newcommand{\chofrench}{0}
\newcommand{\choano}{0}
\newcommand{\choleg}{1}

\newcommand{\nopiecesix}{1}

\ifcase \chofrench
\usepackage{articleanglaiscr}
\or
\usepackage{articlefrancaiscr}
\fi
\usepackage{fmtcount}

\newcommand{\longueurl}{21.8}
\usepackage{amssymb}

\iflanguage{french}{%
\newcommand{\slw}{chemin autobouclant}

\newcommand{\slws}{chemins autobouclant}
\newcommand{\Slws}{Chemins autobouclant}
}{%
\newcommand{\slw}{perfectly looping walk}

\newcommand{\slws}{perfectly looping walks}
\newcommand{\Slws}{Perfectly looping walks}
}

\bibliography{enumeration_construction_circuit,enumeration_construction_circuit_new}

\newcommand{\pieceu}{1}
\newcommand{\pieced}{2}
\newcommand{\piecet}{3}
\newcommand{\pieceq}{4}
\newcommand{\piecec}{5}
\newcommand{\piececb}{6}
\newcommand{\pieces}{7}
\newcommand{\piecesb}{8}
\ifcase \nopiecesix
\newcommand{\nomtypma}{8}
\or
\newcommand{\nomtypma}{6}
\fi

\begin{document}
\ifcase \chofrench
\selectlanguage{english}
\or
\selectlanguage{french}
\fi

%%%%%%%%%%%%%%%%%%%%%%%%%%%%
% TITLE (for amsart class)
%%%%%%%%%%%%%%%%%%%%%%%%%%%%

%%%%%%%%% Begin of title

\iflanguage{french}{%
\title%
{Construction et énumération de circuits aptes à guider un véhicule miniature}
}%
{%
\title%
{Construction and enumeration of circuits capable of guiding a miniature vehicle}
}

\ifcase \choano
\author{Jérôme BASTIEN}
\or
\author{x Y}
\fi

\date{\today}

\ifcase \choano
\address{Centre de Recherche et d'Innovation sur le Sport\\
   U.F.R.S.T.A.P.S.\\
   Université Claude Bernard - Lyon 1\\
   27-29, Bd du 11 Novembre 1918\\
   69622 Villeurbanne Cedex\\
France}
\or
\address{ZZ\\
France}
\fi

\ifcase \choano
\email{\href{mailto:jerome.bastien@univ-lyon1.fr}{\nolinkurl{jerome.bastien@univ-lyon1.fr}}}
\or
\email{\href{mailto:x.y@zz.fr}{\nolinkurl{x.y@zz.fr}}}
\fi

\ifcase \choarticle
\or

\iflanguage{french}{%
\keywords{\Slws, Chemins fermés, Circuits de jeux, Combinatoire, \'Enumération exacte et asymptotique}
}%
{%
\keywords{\Slws, Closed paths, Toy tracks, Combinatorics, Exact and asymptotic enumeration}
}

\subjclass[2010]{Primary 82B41;
Secondary    
05A15
05A16}

\fi

\newcommand{\resumefrancais}{%
Contrairement aux circuits traditionnels de jeux pour trains miniatures, un système breveté  
permet de créer un grand nombre de circuits, avec un nombre minimum de rails, dont les boucles
se referment toujours bien.
Ces circuits ressemblent fortement aux traditionnels polygones autoévitants, dont l'énumération explicite n'est pas encore
résolue pour un nombre quelconque de carrés ; cependant, de nombreuses différences les opposent, dont notamment
le fait que  les contraintes géométriques sont différentes
de celles des polygones autoévitants.
Nous présentons la méthodologie permettant de construire et d'énumérer tous les circuits possibles contenant un nombre
 de pièces donné. Pour de petits nombres de pièces, plusieurs variantes sont proposées 
permettant de prendre ou non en compte le fait que les circuits obtenus peuvent être identiques à une isométrie près.
Pour des nombres de pièces plus importants, seule une estimation sera proposée.
Dans ce dernier cas, une construction aléatoire de circuit est aussi donnée. 
Nous donnerons quelques pistes de généralisations
à d'autres problèmes similaires, mais plus variés.%
}

\newcommand{\resumeanglais}{%
In contrast to traditional toy tracks, a patented  system allows the creation of a large number of tracks with a minimal number of pieces, and whose loops always close properly.
These circuits strongly resemble traditional self-avoiding polygons (whose explicit enumeration has not yet been resolved for an arbitrary number of squares) yet there are numerous differences, notably the fact that the geometric constraints are different  than those of self-avoiding polygons.
We present the methodology allowing the construction and enumeration of all of the possible tracks containing a given number of pieces. For small numbers of pieces, several variants are proposed which allow the consideration or not of the fact that the obtained circuits are identical up to a given isometry.
For greater numbers of pieces, only an estimation will be offered.
In the latter case, 
 a randomly construction of  circuits is also given. %
%a randomly constructed circuit is also given.
We will give some routes for generalizations for similar problems.%
}

\begin{abstract}
\iflanguage{french}{%
\selectlanguage{french}
\resumefrancais
\vspace{0.5 cm}

\noindent
\textsc{Abstract}.
\selectlanguage{english}
\resumeanglais
}%
{%
\selectlanguage{english}
\resumeanglais
}
\end{abstract}

%%%%%%%%% End of title

\ifcase \chofrench
\selectlanguage{english}
\or
\selectlanguage{french}
\fi

\maketitle

%%%%%%%%%%%%%%%%%%%%%%%%%%%%%%%%%%%%%%%%%%%%%%%%%%%%%%%%%%%%%
%%%%%%%%%%%%%%%%%%%%%%%%%%%%%%%%%%%%%%%%%%%%%%%%%%%%%%%%%%%%%
\section{Introduction}
\label{intro}

\iflanguage{french}{%
Les circuits pour enfants existent depuis longtemps et permettent de
transporter aussi bien des petits trains de bois que des voitures ou des trains miniatures.
Par exemple, on pourra se référer aux marques déposées
\textit{Brio} \textregistered,
\textit{Scalextric}  \textregistered\ ou 
\textit{Jouef} \textregistered.
Il existe aujourd'hui des circuits formés d'un très grand nombre de pièces 
différentes, c'est-à-dire supérieur à dix. Ce grand nombre de pièces est intéressant 
car il permet de fabriquer 
%un grand nombre de circuits 
différents
circuits  à partir du même 
ensemble de pièces. En revanche, à cause de ce grand nombre de pièces différentes, il 
existe aussi de nombreuses situations où il n'est pas possible de raccorder 
simplement les deux extrémités du circuit pour le refermer sur lui-même. Dans 
certains cas, cela n'est tout simplement pas possible. Pour remédier à cet 
inconvénient, l'une des solutions consiste à rajouter des pièces de guidage 
supplémentaires, permettant de raccorder les extrémités du circuit qui, auparavant, 
n'étaient pas raccordables. Cela conduit donc à augmenter le nombre de pièces 
différentes et crée de nouvelles situations où il n'est pas possible de raccorder les 
deux extrémités du circuit. 
Souvent, les circuits proposés présentent un jeu, plus ou moins important,
qui permet de créer de grands circuits. Ce jeu, d'origine géométrique, est 
pris en compte dans les pièces constituant ces circuits. Par exemple, les pièces de connexions (tenon/mortaise)
des rails 
\textit{Brio} \textregistered\   laissent les rails se déplacer très légèrement entre eux (<<Vario system>>, voir \url{http://www.woodenrailway.info/track/trackmath.html}).
Les rails de trains de modélisme peuvent être légèrement déformés afin de boucler le circuit.
L'accumulation de ces jeux permet certes de boucler le circuit construit, 
néanmoins  ce jeu rend souvent difficile  la fermeture de la boucle et, si celle-ci se referme,
il est aussi possible que la discontinuité produite fasse dérailler les véhicules miniatures
qui empruntent ces circuits.
\`A l'opposé, il existe des circuits formés de très peu de pièces différentes. Le 
nombre de circuits fermés réalisables avec un ensemble donné de ces pièces est 
alors très faible. Par exemple, il n'existe souvent qu'une seule combinaison des 
pièces de cet ensemble permettant de refermer le circuit sur lui-même.

Le système breveté%
%\footnote{jeu commercialisé : voir \url{http://easyloop.toys/}}
\footnote{Une commercialisation de ce jeu a été tentée  : voir \url{http://easyloop.toys/}
La marque \textit{Easyloop} \textregistered\ a été déposée.}
vise à pallier ces inconvénients en proposant 
un circuit ou un ensemble de pièces de guidage qui:
\begin{itemize}
\item
permet la réalisation d'un grand nombre de circuits fermés à partir du même ensemble de pièces,
\item
utilise un minimum de pièces de guidage différentes,
\item
garantit qu'il est toujours possible de refermer simplement le circuit sur lui-même.
\end{itemize}
Le jeu entre les pièces de ce système sera rigoureusement nul, contrairement aux systèmes traditionnels, 
permettant un emboîtement parfait des boucles du circuit. La fabrication
prévoira toutefois un  jeu très faible permettant la connexion entre elles des pièces du circuit grâce à des couples de tenons/mortaises.
Typiquement, ce système concerne le domaine des circuits de trains pour enfants. Il peut concerner également 
le domaine des circuits de petites voitures ou similaires.

La construction des pièces de ces circuits, qui a donné lieu à un brevet 
\ifcase \choano
\cite{brevetJB,pctJB},
\or
\cite{brevetJBano,pctJBano},
\fi
n'est pas l'objet de cet article, même si celle-ci est sommairement rappelée en section \ref{principe}.
Notons cependant que, d'un point de vue pédagogique et didactique, la construction de ces circuits a fait l'objet de divers exposés à différents publics,
du grand public jusqu'à un séminaire <<~détente~>> niveau L3 à Doctorat, en passant par  des élèves de lycées
\ifcase \choano
(voir \cite{bastienbrevetforum,bastienbrevetmmi})%
\or
(voir \cite{bastienbrevetforumano,bastienbrevetmmiano})%
\fi
 ; ces circuits font appel à différentes notions de géométrie  couvrant le programme du collège et du lycée 
jusqu'à la licence (théorème de Pythagore, tangente, cercle, parabole, courbes de Bézier,  
rayon de courbure, pavage, dénombrement), qui peuvent être évoquées de façon différenciée et adaptée au public concerné.
Une collaboration est envisagée avec Nicolas Pelay, de l'association  Plaisir Maths (\url{http://www.plaisir-maths.fr/}) et nous allons tenter ensemble de promouvoir le jeu 
sur son versant pédagogique et didactique.

%%%%%%%%%%%%%%%%%
% Macrolocale
\newcommand{\questionindustrielle}{%
<<~Est-ce qu'il est possible de dénombrer
tous les circuits réalisables à partir d'un nombre de pièces donné~?~>>}
%%%%%%%%%%%%%%%%%%

La motivation initiale des travaux présentés ici était une question posée par un industriel :
\questionindustrielle\
L'objectif de cet article est d'essayer de répondre à cette question.
Nous présenterons en section \ref{enumeration}
la méthodologie permettant de construire et d'énumérer tous les circuits possibles contenant un nombre
de pièces donné. 
Pour des nombre de pièces plus important, seule une estimation sera proposée (section \ref{estimation}).
Dans ce dernier cas, une construction aléatoire de circuit est aussi donnée (section \ref{constructionaleat}).
Nous donnerons quelques pistes de généralisations  en section \ref{generalisation}.

On pourra  consulter la page web en français 
dédiée\footnote{En cas de défaillance de \url{http://utbmjb.chez-alice.fr/}, utiliser le site miroir  
\ifcase \choano
\url{http://jerome.bastien2.free.fr/}
\or
\url{http://x.y.free.fr/}
\fi
et, dans toutes les url qui suivent remplacer \url{http://utbmjb.chez-alice.fr/} par 
\ifcase \choano
\url{http://jerome.bastien2.free.fr/}%
\or
\url{http://x.y.free.fr/}%
\fi
.} :

\noindent
\url{http://utbmjb.chez-alice.fr/recherche/brevet_rail/detail_brevet_rails.html}

Tous les algorithmes présentés dans cet article ont été implémentés informatiquement
et ont permis aussi bien de déterminer les différents circuits présentés ici que de les tracer à l'échelle.
Quatre exécutables (distribués uniquement pour Windows) et une documentation en français  permettant d'installer 
les librairies graphiques, de créer des circuits de façon manuelle ou aléatoire
et de les dessiner 
 sont disponibles sur 
Internet aux url :
\label{execatalogueweb}

\noindent
\url{http://utbmjb.chez-alice.fr/recherche/brevet_rail/MCRInstaller.exe}

\noindent
\url{http://utbmjb.chez-alice.fr/recherche/brevet_rail/creecircuit.exe}

\noindent
\url{http://utbmjb.chez-alice.fr/recherche/brevet_rail/creecircuitaleat.exe}

\noindent
\url{http://utbmjb.chez-alice.fr/recherche/brevet_rail/dessinecircuit.exe}

\noindent
\url{http://utbmjb.chez-alice.fr/recherche/brevet_rail/mode_emploi_rail_demo.pdf}

Un catalogue a été créé de façon totalement automatique (à part les quelques circuits déterminés manuellement), regroupant les différentes méthodes 
présentées ici et est disponible à l'url suivante : 

\noindent
\url{http://utbmjb.chez-alice.fr/recherche/brevet_rail/catalogue.pdf}

On pourra aussi consulter la totalité des %1102 
$1102$
premiers circuits réalisables (le nombre total de pièces utilisées varie de $4$  à $11$ et le nombre maximal de pièces par type est $4$) sur

\noindent
\url{http://utbmjb.chez-alice.fr/recherche/brevet_rail/catalogue_exhaustif_11rails.pdf}%
}%
{%
Children's tracks have existed for a long time, and allow the transportation of small wooden trains, as well as cars or model trains.
One may refer, for example, to the trademarks
\textit{Brio} \textregistered,
\textit{Scalextric}  \textregistered\ or
\textit{Jouef} \textregistered.
Today there are tracks formed of a very large number of different pieces, i.e. larger than ten. This large number of pieces is interesting, as it allows the production of 
%a large number of 
different circuits from the same set of pieces. On the other hand, due to this large number of different pieces, there are also numerous situations where it is not possible to simply connect the two extremities so as to close the circuit. In certain cases this is simply not possible. To overcome this inconvenience, one solution consists of including additional guide pieces, allowing the connection of the track's extremities which were previously unable to connect. This therefore leads to an increase in the number of different pieces, and creates new situations where it is not possible to connect the two extremities of the track. Often, the circuits offered demonstrate some play, to a greater or lesser extent, which allows the creation of large circuits. This play, of a geometric origin, is taken into account in the pieces constituting these circuits. For example, the (mortise and tenon) connecting parts of \textit{Brio} \textregistered\ track pieces allow them to move very slightly with respect to each other (the ``Vario system''; see \url{http://www.woodenrailway.info/track/trackmath.html}). Model train tracks can be slightly deformed in order to close the circuit. The accumulation of this play indeed allows the closure of the constructed circuit, however the play often renders it difficult to close the circuit and, if it does close, it is also possible that the resulting discontinuity derails the miniature vehicles which use these circuits. There also exist, on the opposite scale, circuits formed of very few different pieces. The number of closed circuits achievable with a given set of pieces is then very low. For example, there often exists only a single combination of the pieces of this set allowing the circuit to close.

The patented  system% 
%\footnote{commercial game: see \url{http://easyloop.toys/}} 
%\footnote{jeu commercialisé : voir \url{http://easyloop.toys/}}
\footnote{A commercialization of this game was attempted: see \url{http://easyloop.toys/}
The trademark \textit{Easyloop}\textregistered\ was registered.}
aims to overcome these drawbacks by offering a circuit or a set of guiding pieces which
\begin{itemize}
\item
allows the realization of a large number of closed circuits from the same set of pieces,
\item 
uses a minimum of different guide pieces,
\item
guarantees that it is always possible to simply close the circuit.
\end{itemize}
The play between the pieces of this system will be strictly zero, in contrast with traditional systems, allowing a perfect fit of circuit loops. The manufacturing will nonetheless provide a very small play, allowing the pieces to be connected to each other by mortise and tenon joints. Typically, this system concerns the domain of train tracks for children, though it may equally concern that of circuits for small cars and the like.

The construction of the pieces of the  circuits, which has led to a patent 
\ifcase \choano
\cite{brevetJB,pctJB},
\or
\cite{brevetJBano,pctJBano},
\fi
is not the subject of this article, though it is briefly recalled in section~\ref{principe}. We note however that, from a pedagogical and didactic point of view, the construction of these circuits has been the subject of various presentations to different audiences, from the general public, to high school students, to a informal seminar for final-year undergraduates to doctoral students 
\ifcase \choano
(see~\cite{bastienbrevetforum,bastienbrevetmmi}). 
\or
(see~\cite{bastienbrevetforumano,bastienbrevetmmiano}). 
\fi
These circuits employ various notions of geometry, spanning the curricula of middle and high school up to undergraduate studies (Pythagoras' theorem, tangents, circles, parabolas, Bézier curves, radii of curvature, tessellation and enumerative combinatorics), which may be brought up in a distinguished way, and adapted to the relevant public.
A collaboration is planned with Nicolas Pelay from the Plaisir Maths association (\url{http://www.plaisir-maths.fr/}), and we will try together to promote the game on it's pedagogical and didactic aspect.

%%%%%%%%%%%%%%%%%
% Macrolocale
% A traduire une seule fois !!!!
\newcommand{\questionindustrielle}{%
<<~Is it possible to tally all of the circuits which can be realized from a given number of pieces?~>>}
%%%%%%%%%%%%%%%%%%

The initial motivation of the work presented here was a question asked by a manufacturer:
\questionindustrielle\
The objective of this article is to attempt to respond to this question.
We will present in Section~\ref{enumeration}
the methodology allowing the construction and enumeration of all of the possible circuits containing a given number of pieces.
For greater numbers of pieces, only an estimation will be offered (Section~\ref{estimation}).
In the latter case, 
a randomly construction of circuits
%a randomly constructed circuit 
is also given (Section \ref{constructionaleat}).
We will give some routes for generalizations in Section~\ref{generalisation}.

One may refer to the dedicated web page in French\footnote{In case \url{http://utbmjb.chez-alice.fr/} is down, use the mirror 
\ifcase \choano
\url{http://jerome.bastien2.free.fr/}
\or
\url{http://x.y.free.fr/}
\fi
and, in all of the following URLs, replace \url{http://utbmjb.chez-alice.fr/} with 
\ifcase \choano
\url{http://jerome.bastien2.free.fr/}%
\or
\url{http://x.y.free.fr/}%
\fi
.} :

\noindent
\url{http://utbmjb.chez-alice.fr/recherche/brevet_rail/detail_brevet_rails.html}

All of the algorithms presented in this article have been implemented computationally and have allowed the determination of the different circuits presented, as well as their to-scale depiction. Four executables (distributed for Windows only) and a documentation in French allowing the installation of graphical libraries, the creation of circuits in a manual or random manner, and the drawing of the circuits are available on the internet at the URLs:
\label{execatalogueweb}

\noindent
\url{http://utbmjb.chez-alice.fr/recherche/brevet_rail/MCRInstaller.exe}

\noindent
\url{http://utbmjb.chez-alice.fr/recherche/brevet_rail/creecircuit.exe}

\noindent
\url{http://utbmjb.chez-alice.fr/recherche/brevet_rail/creecircuitaleat.exe}

\noindent
\url{http://utbmjb.chez-alice.fr/recherche/brevet_rail/dessinecircuit.exe}

\noindent
\url{http://utbmjb.chez-alice.fr/recherche/brevet_rail/mode_emploi_rail_demo.pdf}

A catalogue has been created in a totally automatic manner (apart from the few manually determined circuits), grouping together the different methods presented here, and is available at the following URL: 

\noindent
\url{http://utbmjb.chez-alice.fr/recherche/brevet_rail/catalogue.pdf}

One can also look up all of the first  %1102 
$1102$
 feasible circuits (the total number of pieces used varies from $4$ to $11$
and the maximal number of each kind of piece is $4$)
at

\noindent
\url{http://utbmjb.chez-alice.fr/recherche/brevet_rail/catalogue_exhaustif_11rails.pdf}%
}

\iflanguage{french}{%

%%%%%%%%%%%%%%%%%%%%%%%%%%%%%%%%%%%%%%%%%%%%%%%%%%%%%%%%%%%%%
%%%%%%%%%%%%%%%%%%%%%%%%%%%%%%%%%%%%%%%%%%%%%%%%%%%%%%%%%%%%%
\section{Principe du système breveté}
\label{principe}

%%%%%%%%%%%%%%%%%%%%%%%%%%%%%%%%%%%%%%%%%%%%%%%%%%%%%%%%%%%%%
\subsection{Construction des courbes de base}
\label{principecourbe}

Le principe de ce système est de définir  un chemin $\Gamma$ de $\Er^2$, de
classe ${\mathcal{C}}^1$, ce qui assurera la continuité entre deux pièces successives du circuit
ainsi que leur bon assemblage. 
Soit $N$  un entier naturel quelconque non nul.
Deux idées fondamentales sont utilisées :
\begin{itemize}
\item
On considère un ensemble de carrés ${\mathcal{C}}_i$, $1 \leq i\leq N$ appartenant tous à un pavage du plan par des carrés,
de coté défini par 
\begin{equation}
\label{longueurlunit}
L_0=1.
\end{equation}
On supposera donc, sans perte de 
généralité, que les coordonnées des centres des  carrés ${\mathcal{C}}_i$ sont entières.
%On considère un pavage du plan, formé par des carrés ${\mathcal{C}}_i$, $1 \leq i\leq N$, dont on supposera, sans perte de 
%généralité, que les coordonnées des centres sont entières.
Chaque carré contient une partie du chemin $\Gamma$ et  l'intersection d'un carré  ${\mathcal{C}}_i$
avec $\Gamma$ est notée $\Gamma_i$.
\item
Pour chacun des carré ${\mathcal{C}}_i$, la courbe $\Gamma_i$
doit vérifier les contraintes suivantes :  
\begin{itemize}
\item
elle est contenue dans le carré ${\mathcal{C}}_i$, 
\item
elle débute sur un sommet du carré ou un milieu d'un côté du carré en un point $A_i$ et se 
termine sur un autre sommet du carré ou un milieu d'un autre côté du carré en un 
point $B_i$, 
\item
elle est tangente en $A_i$ et en $B_i$ aux droites reliant respectivement le centre du carré aux 
points $A_i$ et $B_i$.
\end{itemize}
\end{itemize}
Ainsi, le chemin  $\Gamma$ sera défini comme la  réunion des courbes ${(\Gamma_i)}_{1\leq i \leq N}$.
Pour $1\leq i\leq N-1$,
chacun des carrés 
${\mathcal{C}}_i$ doit avoir en commun un unique sommet ou un unique côté avec le  carré voisin
${\mathcal{C}}_{i+1}$. Si $i=N$,
la même règle s'applique pour les carrés ${\mathcal{C}}_1$ et ${\mathcal{C}}_N$.
On  peut donc définir le chemin $\Gamma$, à partir des centre ${(c_i)}_{1\leq i \leq N}$
des carrés ${\mathcal{C}}_i$.
à coordonnées
entières.
Ce problème  s'apparente donc très fortement à la recherche de
chemins autoévitants, décrits dans 
\cite{MR2883859,MR1197356} dans le cas plan,
mais surtout au cas particulier correspondant au cas où l'origine et l'extrémités sont identiques,
c'est-à-dire, le cas des polygones autoévitants, décrits dans 
\cite{%
MR2883859,%
MR1985492,%
MR1718791,%
Guttmann2012,%
Guttmann2012b,%
MR2902304}.
Cinq différences essentielles opposent les circuits de jeux aux polygones autoévitants.
%hormis quatre différences essentielles.
D'une part, dans 
\cite{%
MR2883859,%
MR1985492,%
MR1718791,%
Guttmann2012,%
Guttmann2012b,%
MR2902304},
si les carrés  doivent être nécessairement deux à deux distincts,
le système  \textit{Easyloop} autorise deux carrés non successifs à se confondre ; nous reviendrons sur 
ce point en section \ref{contrainteslocales}. D'autre part, dans \cite{MR2883859,MR1985492},
deux carrés successifs ne peuvent avoir en commun qu'un côté, contrairement au système 
 \textit{Easyloop}.
De plus,  des contraintes supplémentaires dues au nombre de pièces disponibles seront à prendre en compte dans le système \textit{Easyloop}.
Il faudra ne conserver que les circuits, différents à une isométrie près. Voir section \ref{isometrie}.
Enfin, le nombre de pièces utilisées dans les polygones  autoévitants est nécessairement pair ; dans le cas impair, aucun polygone
n'existe, ce qui n'est pas le cas des circuits.

\begin{remark}
\label{newrem01}
Dans 
\cite{MR2104301}, il est mentionné une généralisation des chemins autoévitants,
évoluant, cette fois-ci, dans un pavage du plan formé de triangles
isocèles rectangles (voir figure 1 de cette référence), chacun d'eux correspondant à la moitié d'un carré. 
Nos circuits sont plus proches de ce type de chemin  autoévitants,
puisque, cette fois-ci, il est pris en compte la possibilité que deux
carrés voisins aient un commun un sommet.
Dans cette même référence, il est aussi évoqué le pavage de Kagomé, correspondant
à une tuile octogonale et une tuile carré (voir figure 1 de cette référence). Là
encore, plus de degré de liberté sont offerts, mais la géométrie proposée n'est pas exactement celle de nos circuits.
\end{remark}

Il reste maintenant à définir la géométrie de chacune des courbes $\Gamma_i$. Fixons $i\in \{1,...,N\}$.
Appelons ${\mathcal{H}}_i$, l'ensemble des huit points formés par les quatre milieux et les quatre sommets du carré
${\mathcal{C}}_{i}$. Pour avoir un nombre élevé de circuits, nous chercherons toutes les courbes
possibles correspondant à tous les choix possibles de couples de points distincts $A_i$ et $B_i$ dans ${\mathcal{H}}_i$, 
ce qui représente \textit{a priori} $C^2_8=28$ cas. Cependant, le carré possède un groupe d'isométries $\mathcal{S}$
le laissant invariant de cardinal 8, ce qui réduit à 6 le nombre de courbes possibles. 
On définit 6 type de courbes de la façon suivante :
\begin{itemize}
        \item un premier type regroupant uniquement les courbes pour 
lesquelles les points $A_i$ et $B_i$ sont les milieux de deux côtés opposés du carré,
        \item un deuxième type regroupant uniquement les courbes pour 
lesquelles les points $A_i$ et $B_i$ sont les milieux de deux côtés adjacents du carré, 
        \item un troisième type regroupant uniquement les courbes pour 
lesquelles les points $A_i$ et $B_i$ sont deux sommets diagonalement opposés du carré,
        \item un quatrième type regroupant uniquement les courbes pour 
lesquelles les points $A_i$ et $B_i$ sont deux sommets immédiatement consécutifs du carré,
        \item un cinquième type regroupant uniquement les courbes pour 
lesquelles les points $A_i$ et $B_i$ sont un milieu d'un côté et un sommet du côté opposé du 
carré,
        \item un sixième type regroupant uniquement les courbes pour 
lesquelles les points $A_i$ et B sont un milieu d'un côté et un sommet du même côté du 
carré.
\end{itemize}
De telles contraintes ne définissent pas encore totalement les courbes, mais nous les cherchons maintenant
dans l'ensemble des segments de droites ou des arcs de cercle, comme   dans le monde du jouet.
Sur la figure \ref{fig02tot}, les six types de courbes ont été présentés.
Les premier et troisième types contiennent uniquement des segments de 
droites de longueurs respectives $1$ et $\sqrt 2$ (figures \ref{fig02a} et \ref{fig02c}). Les deuxième et quatrième types contiennent uniquement des 
quarts de cercles, de rayons respectifs $1/2$ et $\sqrt 2/2$ (figures \ref{fig02b} et \ref{fig02d}).
Enfin, pour les deux derniers types, il n'existe pas d'arc de cercle. On cherche donc 
une solution par exemple sous la forme d'une parabole définie par deux points $A_i$ et $B_i$ et les deux tangentes associées.
On peut, soit déterminer l'unique parabole ainsi définie, soit, ce qui revient au même, déterminer
l'unique  courbe de Bézier de degré deux, qui est alors 
définie par les points de contrôle suivants : le point $A_i$, le centre du carré $c_i$ et le point $B_i$
 (figures \ref{fig02e} et \ref{fig02f}).
\ifcase \choano
On pourra consulter \cite{MR1928533,martinholweck13,PerrinBezier,LebosseHemery,bastienbrevetmmi}.
\or
On pourra consulter \cite{MR1928533,martinholweck13,PerrinBezier,LebosseHemery,bastienbrevetmmiano}.
\fi
}%
{%

%%%%%%%%%%%%%%%%%%%%%%%%%%%%%%%%%%%%%%%%%%%%%%%%%%%%%%%%%%%%%
%%%%%%%%%%%%%%%%%%%%%%%%%%%%%%%%%%%%%%%%%%%%%%%%%%%%%%%%%%%%%
\section{Principles of the patented system}
\label{principe}

%%%%%%%%%%%%%%%%%%%%%%%%%%%%%%%%%%%%%%%%%%%%%%%%%%%%%%%%%%%%%
\subsection{Construction of the basic curves}
\label{principecourbe}

The principle of this system is to define a path $\Gamma$ in $\Er^2$, and of class ${\mathcal{C}}^1$, which ensures continuity between two successive pieces of the circuit, as well as their good fit.
Let $N$ be any non-zero natural number.
Two fundamental ideas are used:
\begin{itemize}
\item
We consider a set of squares ${\mathcal{C}}_i$, $1 \leq i\leq N$ each belonging to a square tiling of the plane. 
The side of each square is 
defined by 
\begin{equation}
\label{longueurlunit}
L_0=1.
\end{equation}
We will then assume, without loss of generality, that the coordinates of the centers of the squares ${\mathcal{C}}_i$ are integers.
% We consider a tiling of the plane, formed by squares ${\mathcal{C}}_i$, $1 \leq i\leq N$, for which we assume, without loss of generality, that the 
% coordinates of the centers are integers.
Each square contains a part of the path $\Gamma$, and the intersection of a square ${\mathcal{C}}_i$
with $\Gamma$ is denoted $\Gamma_i$.
\item
For each of the squares ${\mathcal{C}}_i$, the curve $\Gamma_i$ must satisfy the following constraints:
\begin{itemize}
\item
it is contained within the square ${\mathcal{C}}_i$, 
\item
it begins on one vertex of the square, or in the middle of one side of the square, at a point $A_i$, and ends on another vertex of the square, or in the middle of another side, at a point $B_i$, 
\item
it is tangent at $A_i$ and at $B_i$ to the straight lines connecting respectively the center of the square to the points $A_i$ and $B_i$.
\end{itemize}
\end{itemize}
Thus, the path $\Gamma$ will be defined as the union of the curves ${(\Gamma_i)}_{1\leq i \leq N}$. For $1\leq i\leq N-1$, each of the squares ${\mathcal{C}}_i$ must have a unique vertex or side in common with the neighboring square ${\mathcal{C}}_{i+1}$. If $i=N$, then the same rule applies for the squares ${\mathcal{C}}_1$ and ${\mathcal{C}}_N$. One may hence define the path $\Gamma$, from the centers ${(c_i)}_{1\leq i \leq N}$
of the squares ${\mathcal{C}}_i$ with integer coordinates. 
This problem is therefore very similar to the research into self-avoiding walks, described in \cite{MR2883859,MR1197356}, in the planar case, 
and also in the particular case where the origin and the end are identical,
i.e. the case of the self-avoiding polygons, described in 
\cite{%
MR2883859,%
MR1985492,%
MR1718791,%
Guttmann2012,%
Guttmann2012b,%
MR2902304}.
%except for four essential differences. 
Five essential differences distinguish the game's circuits from self-avoiding polygons.
On the one hand, in 
\cite{%
MR2883859,%
MR1985492,%
MR1718791,%
Guttmann2012,%
Guttmann2012b,%
MR2902304},
while the squares must necessarily be distinct, the \textit{Easyloop} system allows two non-successive squares to be confounded; we will return to this point in Section~\ref{contrainteslocales}. On the other hand, in \cite{MR2883859,MR1985492}, two successive squares may only have one side in common, in contrast with the \textit{Easyloop} system. Furthermore, some additional constraints due to the number of available pieces are to be considered in the \textit{Easyloop} system.It will only be necessary to keep circuits which are different up to an isometry. See section \ref{isometrie}.
Finally, the number of pieces used in self-avoiding polygons is necessarily even; in the case of an odd number of pieces, no polygon exists, which is not the case for the circuits.

\begin{remark}
\label{newrem01}
In \cite{MR2104301}, a generalization of self-avoiding walks is mentioned. These evolve in a tiling of the plane formed of isosceles right triangles (see Figure 1 of this reference), each corresponding to half of a square. Our circuits are closer to this type of self-avoiding walk since, here, the possibility that two neighboring squares have a common vertex is taken into account.
In the same reference, the Kagom\'e lattice is also mentioned, corresponding to an octagonal tile and a square tile (see Figure~1 of the reference). There again, more degrees of freedom are offered, but the proposed geometry is not exactly those of our circuits.
\end{remark}

It now remains to define the geometry of each of the curves $\Gamma_i$. Let us fix $i\in \{1,...,N\}$. We call ${\mathcal{H}}_i$, the set of eight points formed by the four middles and the four vertices of the square ${\mathcal{C}}_{i}$. To have a high number of circuits, we seek all of the possible curves corresponding to all of the possible choices of pairs of distinct points $A_i$ and $B_i$ in ${\mathcal{H}}_i$, which represents, \textit{a priori}, $C^2_8=28$ cases. However, the square possesses a group of isometries $\mathcal{S}$ leaving it invariant, of %order 
cardinal 8, which reduces the number of possible curves to 6. 
We define 6 types of curve in the following way:
\begin{itemize}
        \item a first type, grouping together only the curves for which the points $A_i$ and $B_i$ are the middles of two opposite sides of the square,
        \item a second type, grouping together only the curves for which the points $A_i$ and $B_i$ are the middles of two adjacent sides of the square, 
        \item a third type, grouping together only the curves for which the points $A_i$ and $B_i$ are two diagonally opposite vertices of the square,
        \item a fourth type, grouping together only the curves for which the points $A_i$ and $B_i$ are two immediately consecutive vertices of the square,
        \item a fifth type, grouping together only the curves for which the points $A_i$ and $B_i$ are the middle of one side and a vertex of the opposite side of the square,
        \item a sixth type, grouping together only the curves for which the points $A_i$ and $B_i$ are the middle of one side and a vertex of the same side of the square.
\end{itemize}
Such constraints still do not totally define the curves, but we now seek them in the set of line segments or circular arcs, as in the world of the toy.
In Figure~\ref{fig02tot}, the six types of curves are presented.
The first and third types contain only line segments of respective lengths $1$ and $\sqrt 2$ (Figures \ref{fig02a} and \ref{fig02c}). The second and fourth types contain only quarter-circles, with respective radii $1/2$ and $\sqrt 2/2$ (Figures \ref{fig02b} and \ref{fig02d}). Finally, for the last two types, no circular arcs exist. We therefore seek a solution, for example, in the form of a parabola defined by two points $A_i$ and $B_i$ and the two associated tangents.
On may either determine the unique parabola thus defined, or equivalently, determine the unique Bézier curve of order two, which is then defined by the following control points: the point $A_i$, the center of the square $c_i$ and the point $B_i$ (Figures \ref{fig02e} and \ref{fig02f}). 
\ifcase \choano
One may refer to~\cite{MR1928533,martinholweck13,PerrinBezier,LebosseHemery,bastienbrevetmmi}.
\or
One may refer to~\cite{MR1928533,martinholweck13,PerrinBezier,LebosseHemery,bastienbrevetmmiano}.
\fi
}

% extrait de presentation_finale.tex
% figures obtenues grâce à 
%:\programmes\matlab\circuit_rail\trace_forme_base

\begin{figure}
\centering
%%% sous figure 1
\subfigure[\label{fig02a}\iflanguage{french}{Forme}{Form} 1]
{\epsfig{file=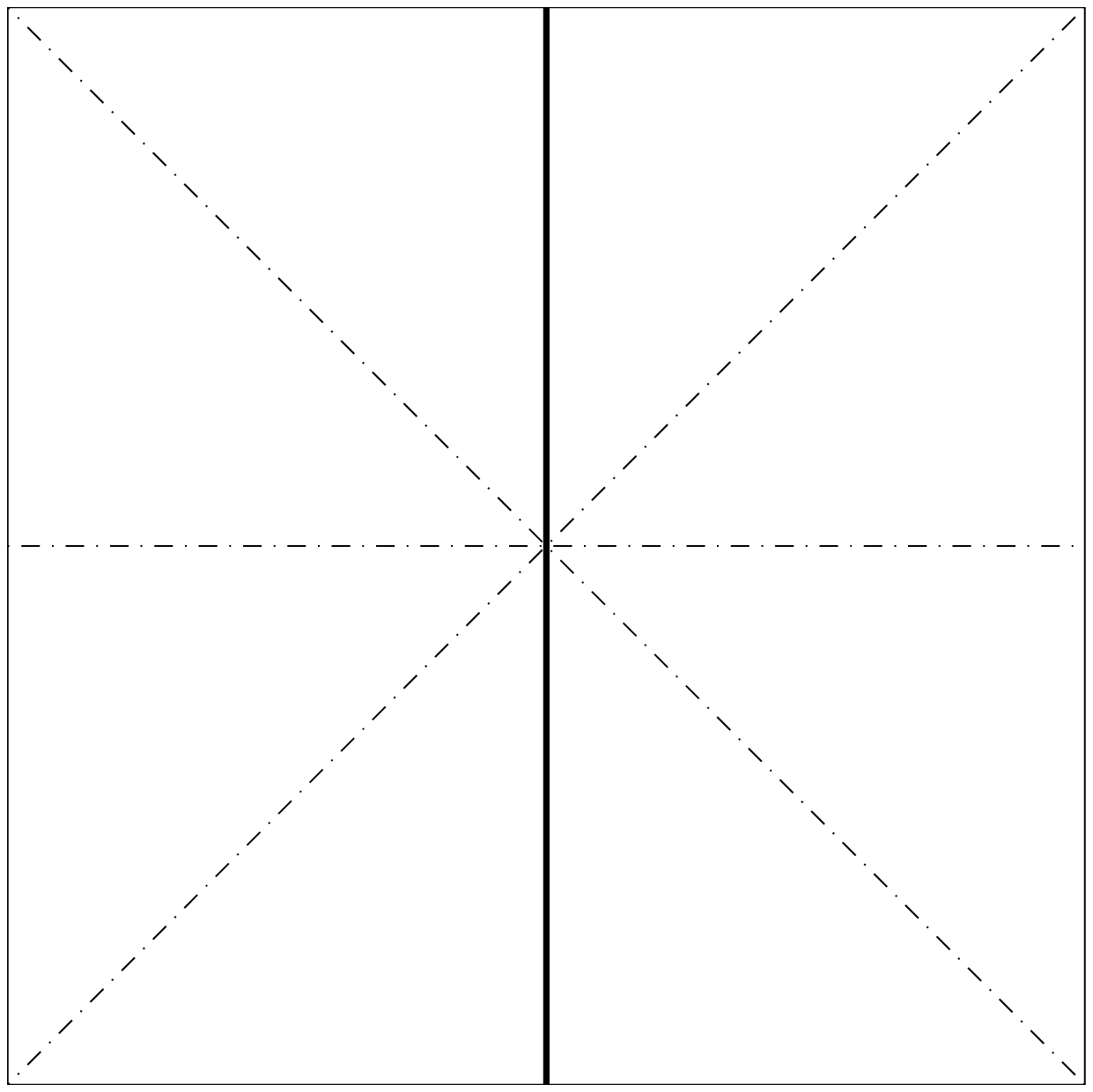, width=3.6 cm}}
\qquad
%%% sous figure 2
\subfigure[\label{fig02b}\iflanguage{french}{Forme}{Form} 2]
{\epsfig{file=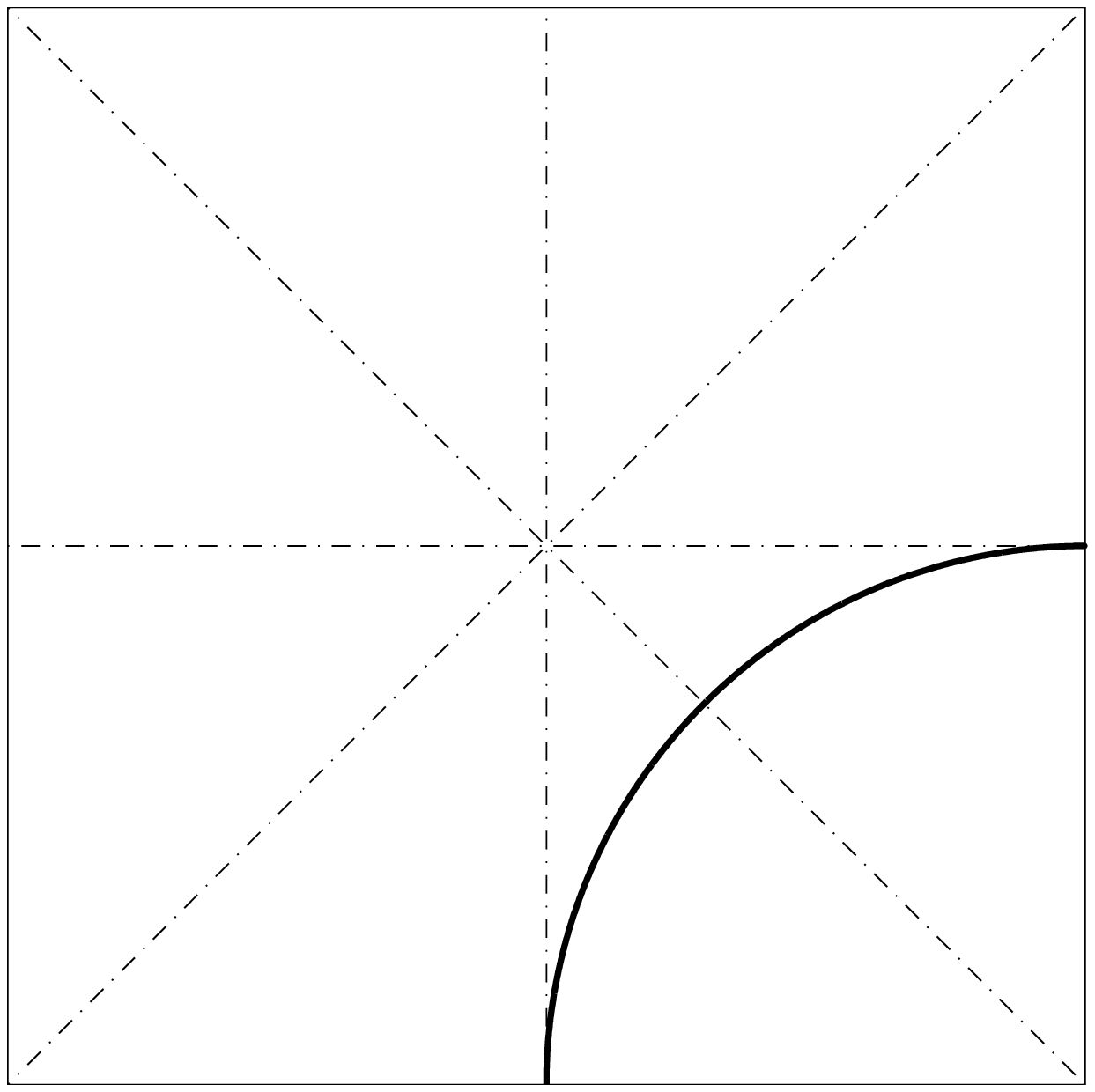, width=3.6 cm}}
\qquad
%%% sous figure 3
\subfigure[\label{fig02c}\iflanguage{french}{Forme}{Form} 3]
{\epsfig{file=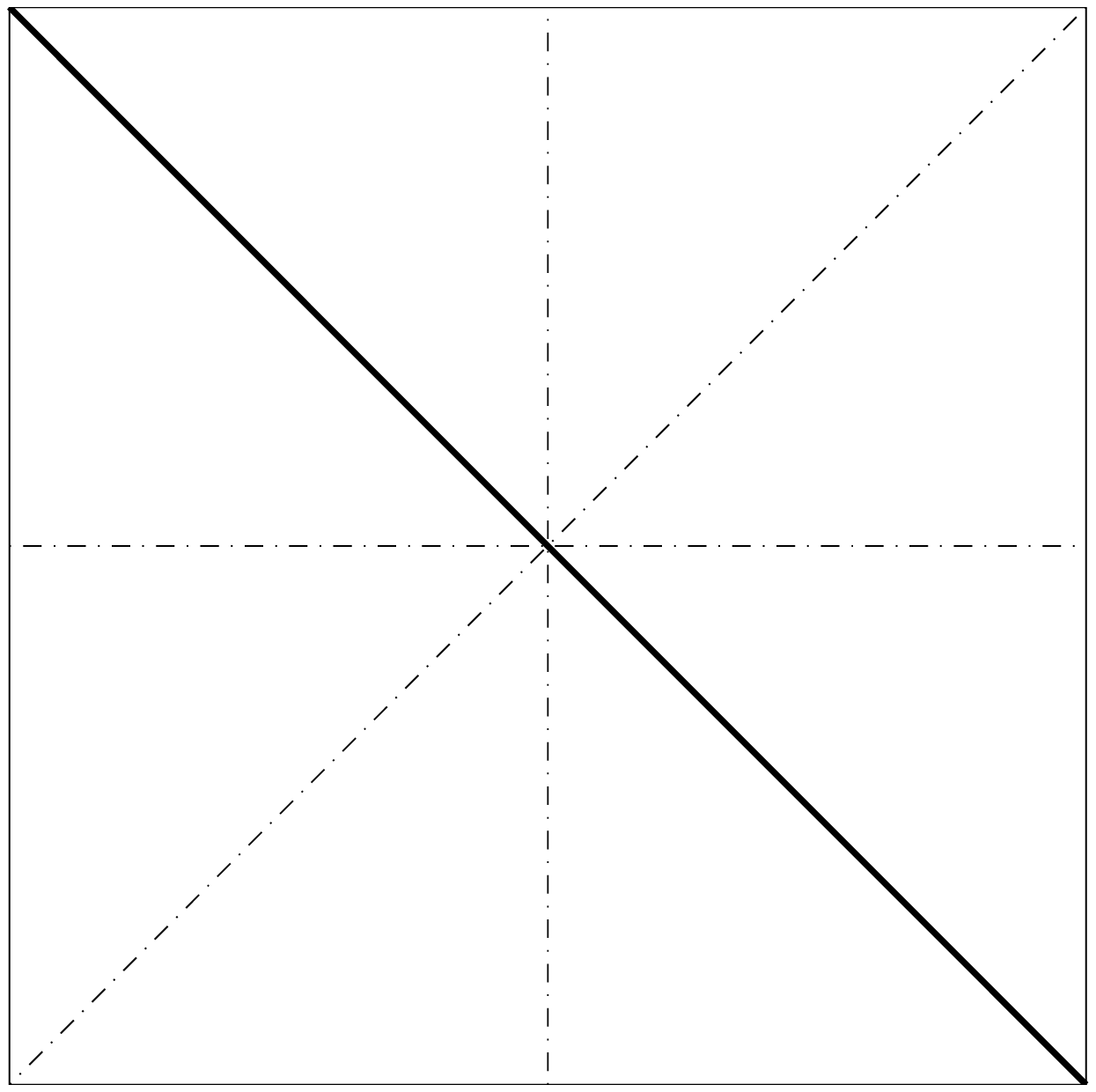, width=3.6 cm}}
\qquad
%%% sous figure 4
\subfigure[\label{fig02d}\iflanguage{french}{Forme}{Form} 4]
{\epsfig{file=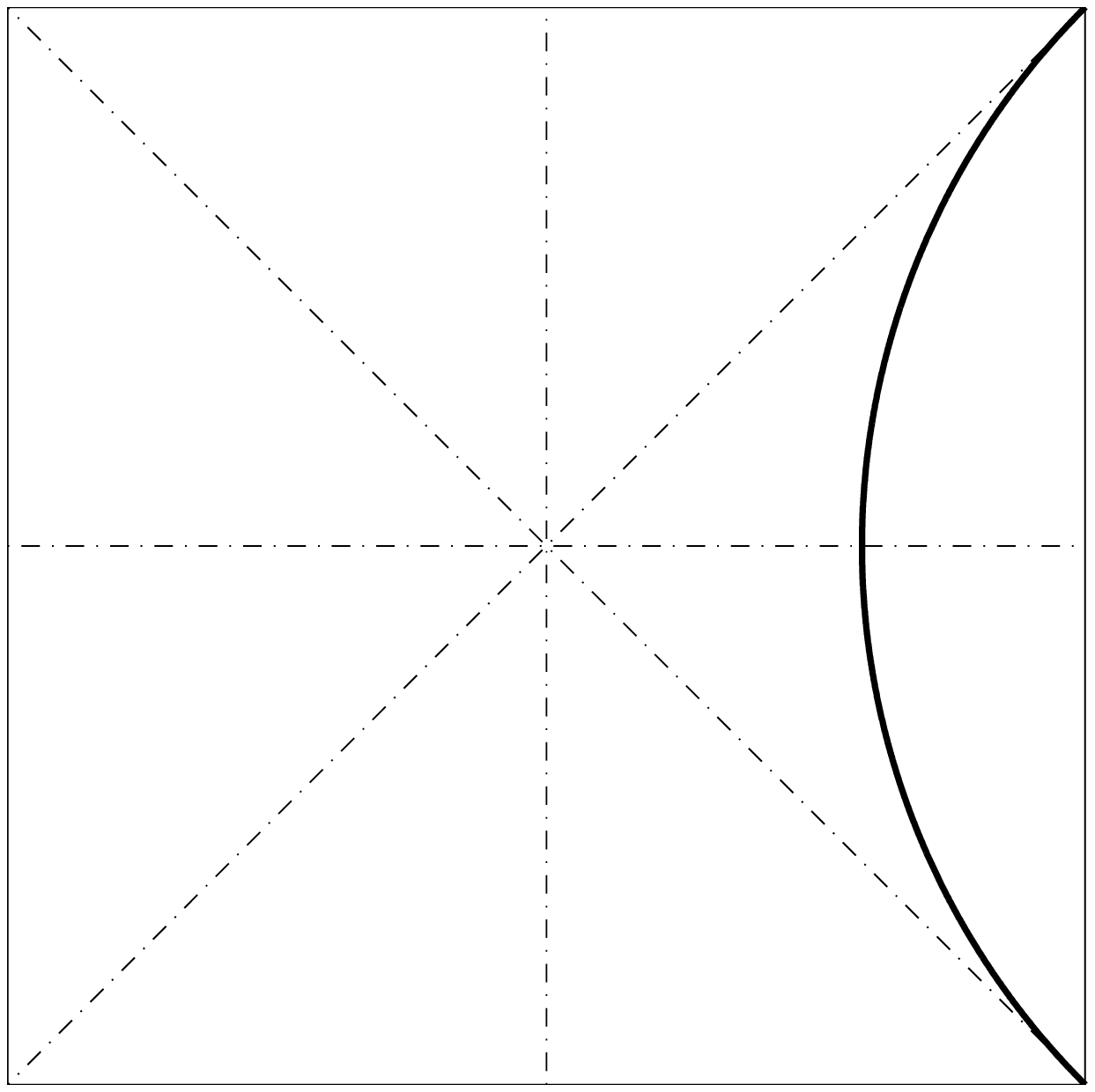, width=3.6 cm}}
\qquad
%%% sous figure 5
\subfigure[\label{fig02e}\iflanguage{french}{Forme}{Form} 5]
{\epsfig{file=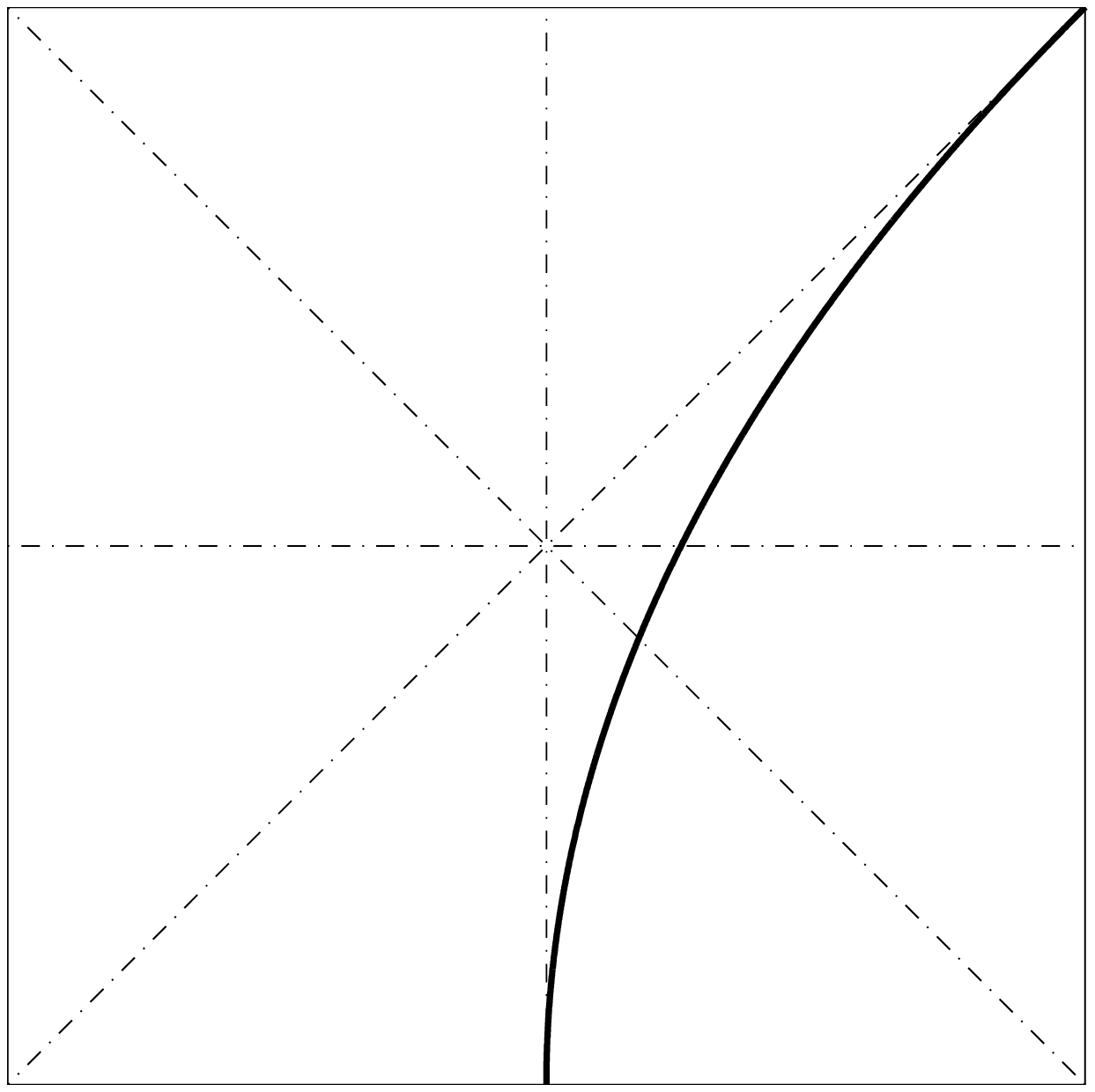, width=3.6 cm}}
\qquad
%%% sous figure 6
\subfigure[\label{fig02f}\iflanguage{french}{Forme}{Form} 5']
{\epsfig{file=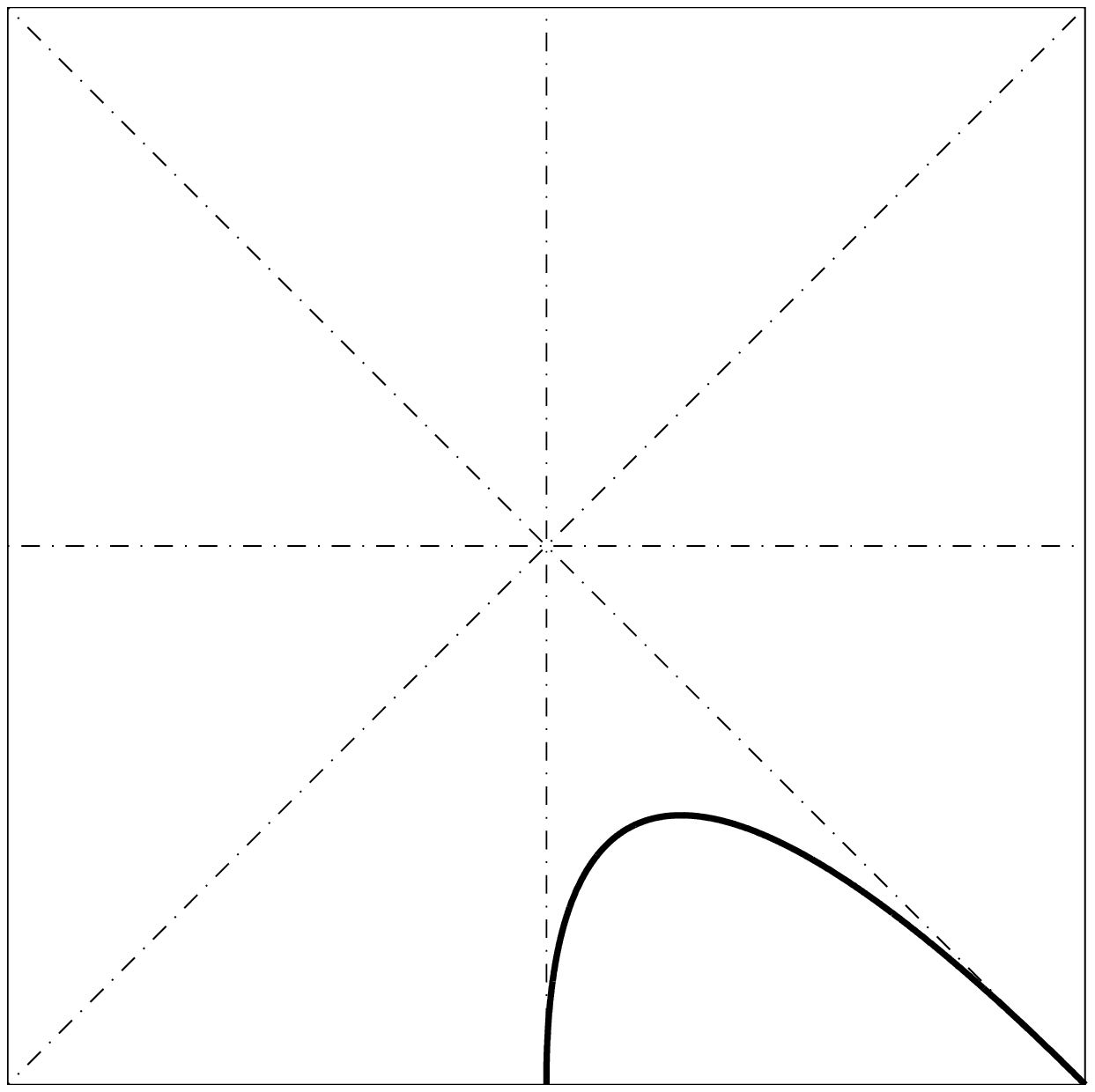, width=3.6 cm}}
\caption{\label{fig02tot}\iflanguage{french}{Les six formes de base.}{The six basic forms.}}
\end{figure}

% extrait de presentation_finale.tex
% figures obtenues grâce à 
%:\programmes\matlab\circuit_rail\trace_forme_base
\begin{figure}
\centering
%%% sous figure 1
\subfigure[\label{fig1000a}\iflanguage{french}{avec 5 pièces de base}{with 5 basic pieces}]
{\epsfig{file=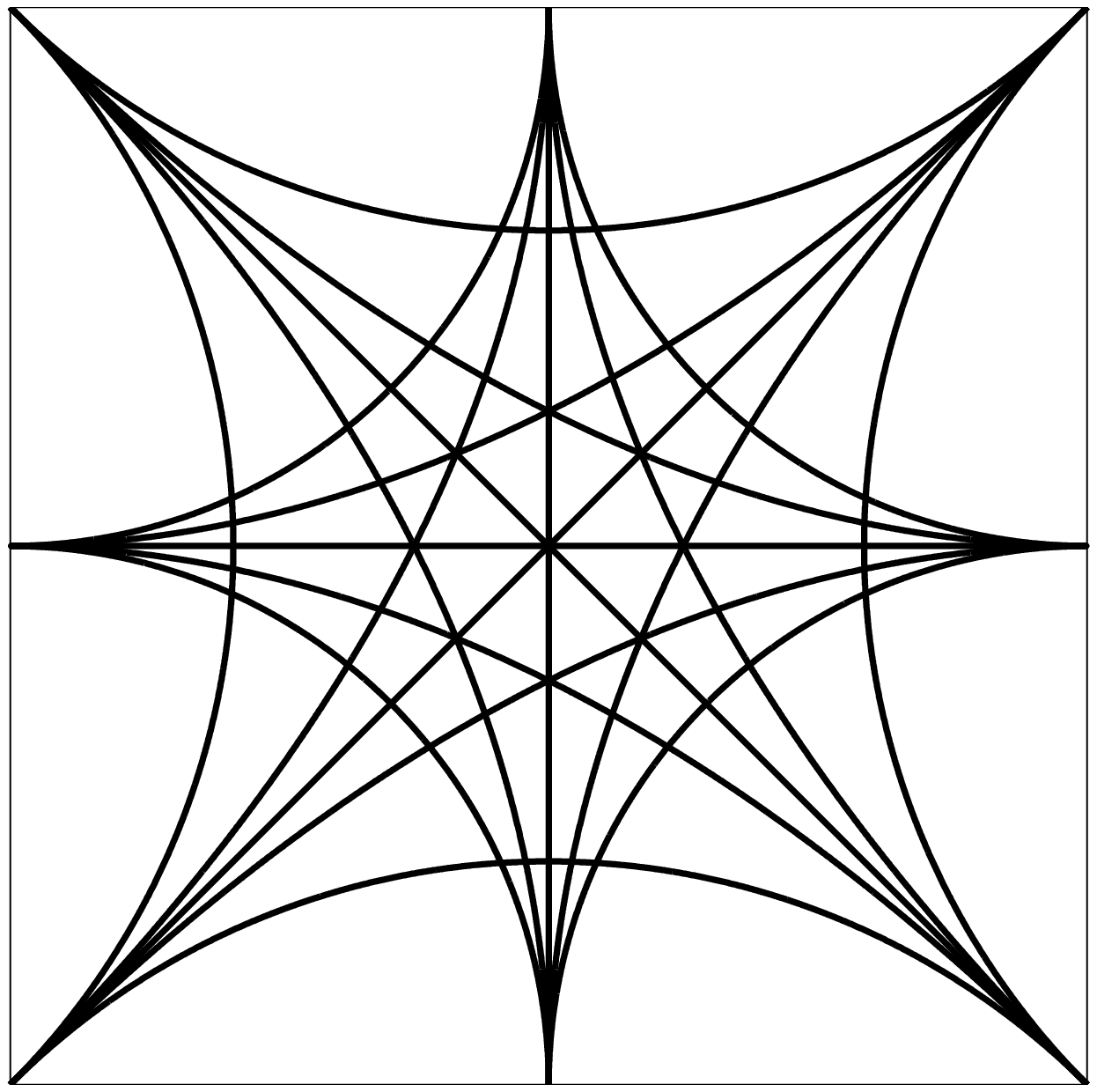, width=6cm}}
\qquad
%%% sous figure 2
\subfigure[\label{fig1000b}\iflanguage{french}{avec 6 pièces de base}{with 6 basic pieces}]
{\epsfig{file=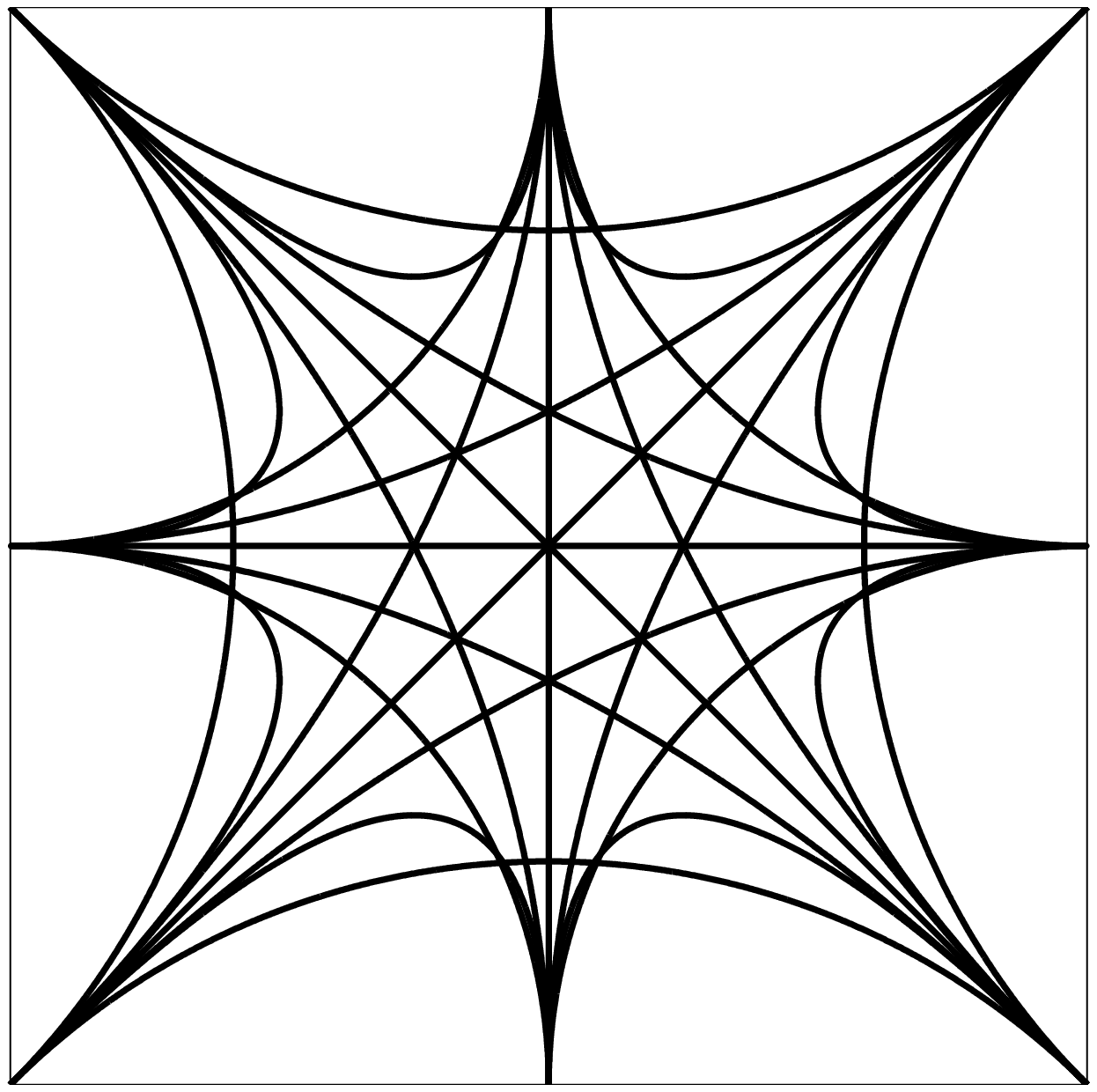, width=6cm}}
\caption{\label{fig1000}\iflanguage{french}{Ensemble des trajets possibles.}{The set of possible paths.}}
\end{figure}

\iflanguage{french}{%
Si fait opérer le groupe des 8 isométries $\mathcal{S}$ sur les 6 courbes de la figure 
\ref{fig02tot}, on obtient bien les 28 courbes possibles de la figure \ref{fig1000b}.
\ifcase \nopiecesix
\or
Le sixième type de la figure \ref{fig02f}
sera éliminé par la suite, puisque la pièce correspondante possède un rayon de courbure minimal trop petit
pour que des véhicules miniatures puissent y rouler (voir section \ref{principepiece}), 
ce qui réduit à 20 le nombre de trajets possibles (voir la figure \ref{fig1000a}).
Dans ce cas, la règle <<~toute courbe relie deux points distincts quelconques de ${\mathcal{H}}_i$~>>  est à remplacer
par  <<~toute courbe relie deux points quelconques distincts et non voisins  de ${\mathcal{H}}_i$~>>, ce qui, on le verra par la suite, offre
tout de même un grand nombre de circuits.
\fi

La courbe $\Gamma$ est de classe ${\mathcal{C}}^1$ ; en effet,
chacune des courbes $\Gamma_i$ est de classe ${\mathcal{C}}^{\infty}$.  En outre, la réunion 
de toutes ces courbes sera de classe ${\mathcal{C}}^1$ ; par construction, en effet, aux points de raccordement, qui 
ne peuvent être que des sommets ou des milieux des côtés de carrés, les courbes sont continues (puisqu'elles passent par 
les mêmes points de début et de fin) et à dérivée continue, puisque les tangentes coïncident.

Les courbes $\Gamma$ obtenues sont de classe ${\mathcal{C}}^1$
mais non de classe ${\mathcal{C}}^2$, à  cause de la discontinuité du rayon de courbure, contrairement aux réseaux réels
ferrés et de routes. Notons cependant que cette discontinuité est aussi présente dans les systèmes traditionnels déjà existants, constitués de
pièces rectilignes et circulaires, de rayons de courbure différents. Sur le plan mécanique,
cela engendre, pour le véhicule miniature qui emprunte ce circuit,  une discontinuité de l'accélération normale (à vitesse continue)
et de l'angle de braquage des roues. Ces contraintes, importantes pour les véhicules réels, interviennent directement sur le confort du passager et 
sur l'usure induite des pièces mécaniques, mais ne sont pas prises en compte dans le domaine des jouets. En effet, 
les masses et les vitesses des véhicules sont très faibles, donc les chocs dus aux discontinuité de l'accélération normale sont négligeables. 
De plus, la notion de confort du voyageur n'a pas de sens 
\ifcase \choarticle
ici\footnote{sauf si l'on considère 
des Schtroumpfs ou des Playmobils qui monteront sur les véhicules miniatures, mais qui ont toujours un grand sourire.}.
\or
ici.
\fi
Enfin, les roues des véhicules peuvent subir une discontinuité de l'angle de braquage 
parce qu'elles présentent un léger jeu par rapport au châssis. 
}%
{%

Acting on the 6 curves in Figure~\ref{fig02tot} with the group of 8 isometries $\mathcal{S}$, one indeed obtains the 28 possible curves of Figure~\ref{fig1000b}.
\ifcase \nopiecesix
\or
The sixth type in Figure~\ref{fig02f}
will be eliminated in the following, since the corresponding piece has a radius of curvature which is too small for the miniature vehicles to be able 
to ride there (see section \ref{principepiece}), 
which reduces the number of possible paths to 20 (see Figure~\ref{fig1000a}).
In this case, the rule <<~every curve linking any two distinct points in ${\mathcal{H}}_i$~>> is to be replaced by <<~every curve linking any two distinct non-neighboring points in ${\mathcal{H}}_i$~>>, which, as we will see in the following, nevertheless offers a large number of circuits.
\fi

The curve $\Gamma$ is of %the
class ${\mathcal{C}}^1$; %in fact, 
indeed, each of the curves $\Gamma_i$ is of class ${\mathcal{C}}^{\infty}$. Furthermore, the union of all of these curves will be of class ${\mathcal{C}}^1$. By construction, indeed, 
at the connecting points, which can only be vertices or middles of the sides of squares, the curves are continuous (since they pass by the same start and end points) and have a continuous derivative, since the tangents coincide.

The curves $\Gamma$ obtained are of class ${\mathcal{C}}^1$, but not of class ${\mathcal{C}}^2$, due to the discontinuity of the radius of curvature, in contrast with real rail and road networks. We note however that this discontinuity is also present in the existing traditional systems, constituted of straight-line and circular pieces of different radii of curvature. On a mechanical level, this generates, for the miniature vehicle which takes the circuit, a discontinuity of the normal acceleration (at constant velocity) and of the steering angle of the wheels. These constraints, significant for real vehicles, directly affect the comfort of the passenger and the wear induced on the mechanical parts, but are not taken into account in the domain of games. Indeed, the masses and the velocities of the vehicles are very low, and therefore the shocks due to the discontinuity of the normal acceleration are negligible.
Moreover, the notion of the comfort of the passenger has no meaning 
\ifcase \choarticle
here\footnote{except if one considers the Smurfs or the Playmobils which will be mounted on the miniature vehicles, but which always have a large smile.}. 
\or
here.
\fi
Finally, the wheels of the vehicles may be subjected to a discontinuity of the steering angle since they exhibit a slight play with respect to the chassis.%
}

\iflanguage{french}{%

%%%%%%%%%%%%%%%%%%%%%%%%%%%%%%%%%%%%%%%%%%%%%%%%%%%%%%%%%%%%%
\subsection{Construction des pièces}
\label{principepiece}

Une fois le chemin $\Gamma$ construit, il reste à définir les différents types de rails constituant le circuit.
Chacun des type de rails sera donc défini à partir de l'un des 
\ifcase \nopiecesix
six
\or
cinq 
\fi
types de courbes précédemment définis. 

\begin{itemize}
\item
Ces courbes constituent la ligne médiane de chacun des types de pièces.

\item
Les passages des roues
sont définis comme deux courbes à distance constante de cette courbe médiane, c'est-à-dire :
chaque point de l'une de ces deux courbes se trouve sur une droite perpendiculaire à la tangente 
à la ligne médiane à une distance constante de la courbe médiane. Les bords de rails sont définis de la même façon.

\item
La section transversale de la pièce est définie de façon classique%
\footnote{voir par exemple \url{http://pw1.netcom.com/~thoog/hnr/hnr.html}}.

\item
%%%%%%%%%%%%%%%%%%%%%%%%%%%%%%%%%%%%%%%%%%%%%%%%%
%%%\input{pb_rayon_courbure_trop_petit}
%%%%%%%%%%%%%%%%%%%%%%%%%%%%%%%%%%%%%%%%%%%%%%%%%

\item
Enfin, les connecteurs sont des tenons/mortaises, conçus 
de telle sorte que chaque rails possède un tenon et une mortaise. 
\end{itemize}

Les pièces de types  1 à 4 sont symétriques : elles possèdent un plan de symétrie perpendiculaire à la courbe médiane, 
et puisque la section est elle-même symétrique, il suffit donc de construire un seul type de pièces pour ces quatre types.
En revanche, la pièce de type 5 
\ifcase \nopiecesix
(ainsi que celle de type 6)
\or
\fi
n'est pas symétrique et les deux extrémités n'ont donc pas le même rôle. Il a donc fallu
construire deux pièces différentes où les prises mâles et femelles sont inversées pour pouvoir réaliser le circuit.%
}%
{%

%%%%%%%%%%%%%%%%%%%%%%%%%%%%%%%%%%%%%%%%%%%%%%%%%%%%%%%%%%%%%
\subsection{Construction of the pieces}
\label{principepiece}

Once the path $\Gamma$ has been constructed, it remains to define the different types of tracks constituting the circuit. Each of the types of track will be defined from one of the 
\ifcase \nopiecesix
six
\or
five
\fi
types of curves defined above.

\begin{itemize}
\item
These curves constitute the midline of each of the types of piece.

\item
The wheel passages are defined as two curves at a constant distance from this middle curve, i.e: each point from one of these two curves is found on a straight line perpendicular to the tangent of the midline at constant distance from the middle curve. The edges of the tracks are defined in the same way.

\item
The cross section of the piece is defined in a conventional way%
\footnote{See, for example \url{http://pw1.netcom.com/~thoog/hnr/hnr.html}}.

\item
%%%%%%%%%%%%%%%%%%%%%%%%%%%%%%%%%%%%%%%%%%%%%%%%%%%%%%%%%%%%
%\input{pb_rayon_courbure_trop_petit}
% fichier tex crée par MaTeXBuild02 le 29-Feb-2016 14:58:54
% à compiler avec 
% MaTeXBuild02('pb_rayon_courbure_trop_petit',0)
% après enumeration_construction_circuit_new
% Partiellement adapté de courbes_bezier_5_6.matex
\iflanguage{french}{%
Notons que dans \cite{brevetJB}, nous avons spécifié que la demi-largeur $e/2$ du rail devait être inférieure au rayon de courbure de la courbe médiane, 
pour éviter qu'au point considéré, la courbe construite à égale distance de la courbe médiane, ne présente un point stationnaire avec
un changement de sens du vecteur unitaire tangent.
Pour la parabole de la figure \ref{fig02e}, le rayon minimal est égal à $R_{\min}=1$, tandis que pour 
la parabole de la figure \ref{fig02f}, le rayon minimal est égal à $R_{\min}=\frac{\sqrt{5}}{25}$. Ainsi, si l'on considère les courbes des figures 
\ref{fig02a} à \ref{fig02e}, le rayon de courbure minimal est donc égal à 
\begin{equation*}
R_{\min}=\frac{1}{2},
\end{equation*}
tandis que si l'on considère
les courbes des figures 
\ref{fig02a} à \ref{fig02f}, le rayon de courbure minimal est donc égal à 
\begin{equation*}
R_{\min}=\frac{\sqrt{5}}{25}.
\end{equation*}
Dans le premier cas, la largeur $e$ du rail est donc inférieure strictement à $e_0$ donné par 
\begin{equation}
\label{maxlar5}
e_0=1,
\end{equation}
tandis que dans le second cas, $e_0$ est donné par 
\begin{equation}
\label{maxlar6}
e_0=\frac{2\, \sqrt{5}}{25}\approx   0.17889.
\end{equation}
Le choix d'une section standard, compatible avec les véhicules miniatures de type \text{Brio} \textregistered, correspond à 
\begin{equation}
\label{dminminexaeval}
e=0.18349.
%0.18349.
\end{equation}
Ainsi, pour les courbes des figures 
\ref{fig02a} à \ref{fig02e}, ce choix de largeur permet d'avoir aucun point stationnaire.
En revanche, pour les courbes des figures 
\ref{fig02a} à \ref{fig02f}, ce choix de largeur fait apparaître un point stationnaire, comme le montre 
les pièces  \pieces\ et \piecesb\ de 
la figure \ref{numerotation_piece}.%
}{%
Note that in \cite{brevetJB}, we specified that the half-width $e/2$ of the rail must be less than the radius of curvature of the midline in order to avoid that, at the considered point, the curve constructed at equal distance from the midline does not feature a stationary point with a change in the direction of the unit tangent vector.
For the parabola in Figure \ref{fig02e}, the minimal radius is $R_{\min}=1$, while for the parabola in Figure \ref{fig02f}, the minimal radius is $R_{\min}=\frac{\sqrt{5}}{25}$. Thus, if we consider the curves in Figures\ref{fig02a} to \ref{fig02e}, the minimal radius of curvature is then
\begin{equation*}
R_{\min}=\frac{1}{2},
\end{equation*}
whereas if we consider the curves in Figures \ref{fig02a} to \ref{fig02f}, the minimal radius of curvature is then 
\begin{equation*}
R_{\min}=\frac{\sqrt{5}}{25}.
\end{equation*}
In the first case, the width $e$ of the rail is therefore strictly less than $e_0$, given by 
\begin{equation}
\label{maxlar5}
e_0=1,
\end{equation}
while in the second case, $e_0$ is given by
\begin{equation}
\label{maxlar6}
e_0=\frac{2\, \sqrt{5}}{25}\approx   0.17889.
\end{equation}
The choice of a standard cross-section, compatible with \text{Brio}\textregistered-type vehicles, corresponds to
\begin{equation}
\label{dminminexaeval}
e=0.18349.
%0.18349.
\end{equation}
Thus, for the curves in Figures \ref{fig02a} to \ref{fig02e}, this choice of width allows us to have no stationary points. On the other hand, for the curves in Figures \ref{fig02a} to \ref{fig02f}, this choice of width presents a stationary point, as exhibited by pieces \pieces\ and \piecesb\ in Figure \ref{numerotation_piece}.%
}
%%%%%%%%%%%%%%%%%%%%%%%%%%%%%%%%%%%%%%%%%%%%%%%%%%%%%%%%%%%%

\item
Finally, the connectors are mortise and tenon joints, designed such that each track possesses one mortise and one tenon.
\end{itemize}

Piece-types 1 to 4 are symmetric: they possess a symmetry plane perpendicular to the middle curve, and since the cross section is itself symmetric, it is therefore sufficient to construct a single type of piece for each of these four types. On the other hand, the type 5 piece 
\ifcase \nopiecesix
(and that of type 6)
\or
\fi
isn't symmetric, and the two extremities therefore do not play the same role. In order to realize the circuit, it was therefore necessary to construct two different pieces where the male and female connectors are inverted.%
}%
\ifcase \choleg
\begin{figure}[h] 
\begin{center} 
{\epsfig{file=01290010.eps, width=15 cm}}
\end{center} 
\caption{\label{deux_type_piecec}\iflanguage{french}{Les deux pièces correspondant à la ligne médiane de type 5.}{The two pieces corresponding to the type 5 midline.}} 
\end{figure}
\or
\fi
\iflanguage{french}{%
\ifcase \choleg
Voir la figure \ref{deux_type_piecec}.
\or
\fi

Les extrémités des pièces pouvant être soit des milieux de côté, soit des sommets, 
il a fallu marquer   sur les rails les extrémités correspondant à des sommets,
\ifcase \choleg
ce qui a été réalisé grâce à une pastille jaune visible sur la figure \ref{deux_type_piecec}.
\or
ce qui a été réalisé grâce à une pastille jaune.
\fi 
Cette pastille figure aussi sur tous les plans qui seront présentés dans ce document.
L'enfant qui joue à assembler les pièces aura donc cette unique règle à respecter :
<<~n'assembler des pièces entre elles que si les deux extrémités des deux pièces contiguës
ont la même nature (absence ou présence simultanée de pastilles)~>>. Cette règle est la seule règle à respecter pour pouvoir faire des circuits
qui se rebouclent bien !

\ifcase \choleg
Enfin des prototypes ont été fabriqués, du type de ceux de la figure \ref{deux_type_piecec}.
\or
Enfin des prototypes ont été fabriqués.
\fi
Les carrés théoriques  présentés de côté donnés par \eqref{longueurlunit},
ont tous été multipliés par une longueur de référence donnée par 
\begin{equation}
\label{longueurl}
L=\longueurl \, \text{cm},
\end{equation}
qui est égale au côté du carré de base constituant le pavage réel.%
}%
{%
\ifcase \choleg
See Figure~\ref{deux_type_piecec}.
\or
\fi

The extremities of the pieces may be either %
%squares 
middles of side or vertices, and it was necessary to mark on the rails the extremities corresponding to vertices,
\ifcase \choleg
which was done using a yellow dot, visible in Figure~\ref{deux_type_piecec}. 
\or
which was done using a yellow dot. 
\fi
This dot also appears on all of the track designs which will be presented in this document. The child which plays at assembling the pieces will therefore have this single rule to obey: <<~only put pieces together if the two extremities of two contiguous pieces have the same nature  (simultaneous absence or presence of dots)~>>. This rule is the only one to be obeyed in order to be able to make circuits which loop properly!

\ifcase \choleg
Finally some prototypes were manufactured, of the type of those in Figure~\ref{deux_type_piecec}. 
\or
Finally some prototypes were manufactured.
\fi
The theoretical squares shown with side of length 
given by \eqref{longueurlunit}
have all been multiplied by a reference length given by
\begin{equation}
\label{longueurl}
L=\longueurl \, \text{cm},
\end{equation}
which is equal to the side of the basic square constituting the real tiling.
}

\begin{figure}[h] 
\begin{center} 
\epsfig{file=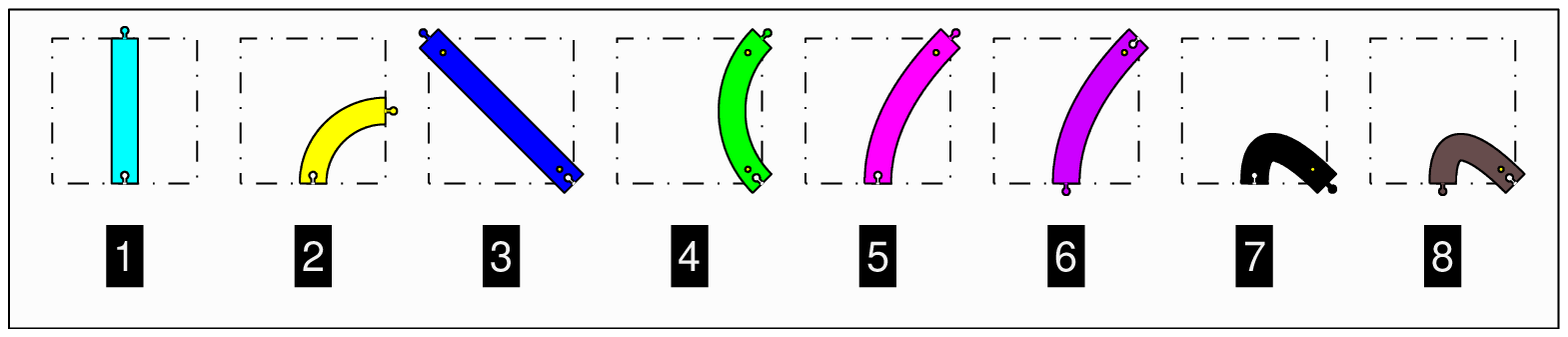, width=17 cm} 
\end{center} 
\caption{\label{numerotation_piece}\iflanguage{french}{La numérotation des pièces.}{The numbering of the parts.}} 
\end{figure}

\iflanguage{french}{%
Les pièces sont désormais numérotées comme l'indique la figure \ref{numerotation_piece}.

\ifcase \nopiecesix
\or
Rappelons que les pièces 
\pieces\ et \piecesb\
présentes en théorie dans le circuit, n'ont pas été réalisées en pratique, puisque trop incurvées.
Les calculs présentés par la suite peuvent tout à fait prendre en compte ces deux types de pièces, mais pour simplifier
nous supposerons désormais par la suite que  seuls les six premiers types de pièces sont utilisés.
Cependant, les programmes et les algorithmes décrits peuvent éventuellement prévoir la présence de ces deux pièces.
\fi%
}%
{%

The pieces are henceforth numbered as indicated in Figure~\ref{numerotation_piece}.

\ifcase \nopiecesix
\or
Recall that pieces \pieces\ and \piecesb\, present in theory in the circuit, were not produced in practice, since they are too curved. The calculations presented in the following can certainly take into account these two types of piece, but for simplicity we will henceforth assume that only the first six types of piece are used. However, the programs and algorithms described can also provide for the presence of these two pieces.
\fi%
}

%%%%%%%%%%%%%%%%%%%%%%%%%%%%%%%%%%%%%%%%%%%%%%%%%%%%%%%%%%%%%%
%\input{simulations_circuit/simulation400}
% fichier tex crée par MaTeXBuild02 le 03-Sep-2015 06:03:53
% A compiler avec MaTeXBuild02('simulation400',0)
% après le fichier 'enumeration_construction_circuit.matex'
% Sortie statique, pas très pertinent avec matex !!

\begin{figure}[h] 
\begin{center} 
\epsfig{file=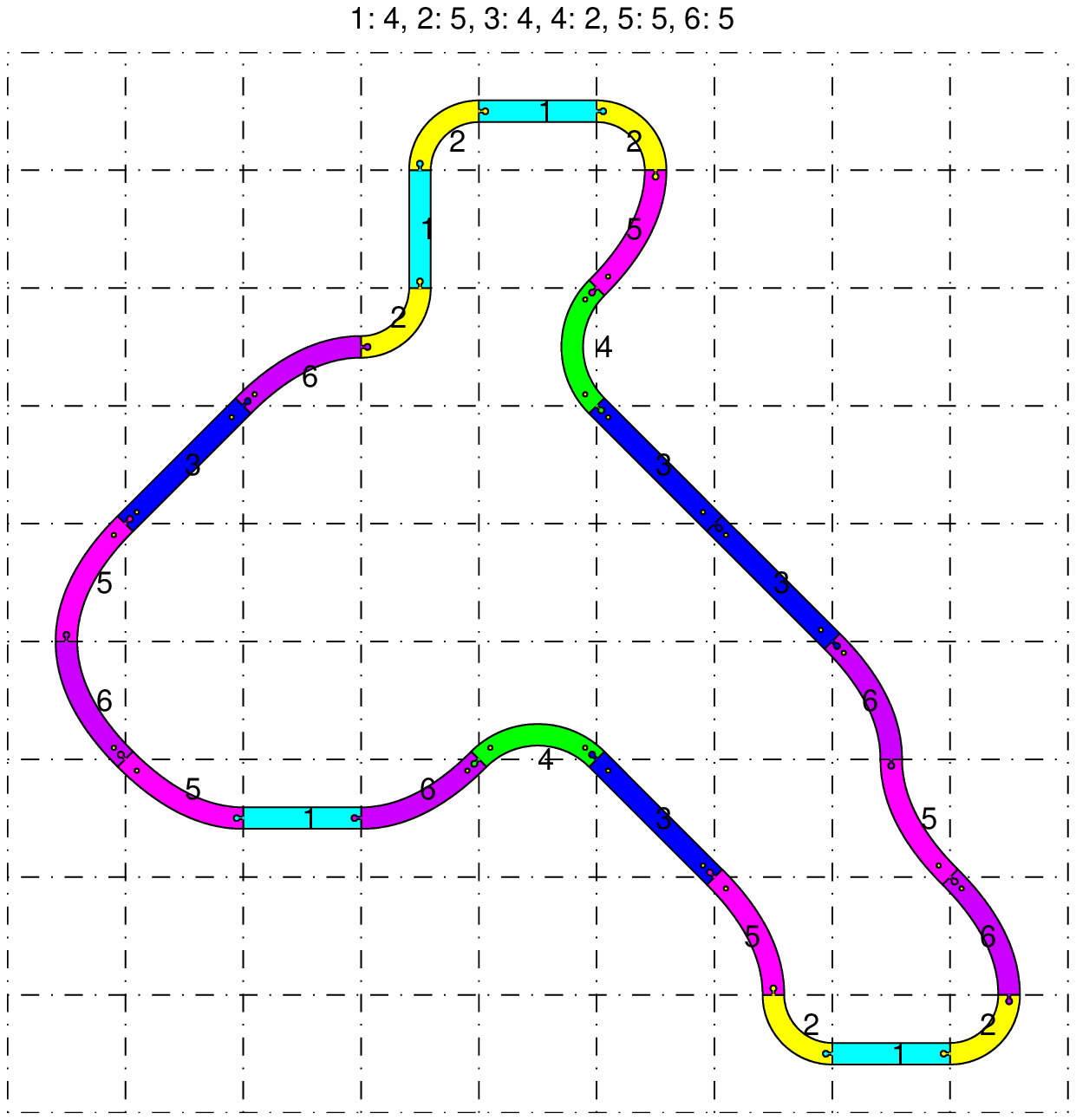, width=11cm}
\end{center} 
\caption{\label{exemple_circuit_plan}\iflanguage{french}{Un exemple de plan de circuit.}{An example track design.}} 
\end{figure}
%%%%%%%%%%%%%%%%%%%%%%%%%%%%%%%%%%%%%%%%%%%%%%%%%%%%%%%%%%%%%%

\ifcase \choleg
\begin{figure}[h] 
\begin{center} 
\epsfig{file=01300011.eps, width=15cm}
\end{center} 
\caption{\label{exemple_circuit_photo}\iflanguage{french}{Un exemple de circuit créé correspondant au plan de la figure \ref{exemple_circuit_plan}.}{An example track created corresponding to the design in Figure~\ref{exemple_circuit_plan}.}} 
\end{figure}
\or
\fi

\ifcase \choleg
\iflanguage{french}{Donnons à titre d'exemple un circuit réellement fabriqué sur les figures \ref{exemple_circuit_plan} et \ref{exemple_circuit_photo}.}{We give an actually manufactured circuit as an example in Figures~\ref{exemple_circuit_plan} and \ref{exemple_circuit_photo}.}
\or
\iflanguage{french}{Donnons à titre d'exemple un circuit réellement fabriqué sur la figure \ref{exemple_circuit_plan}.}{We give an actually manufactured circuit as an example in Figure~\ref{exemple_circuit_plan}.}
\fi

\ifcase \choarticle

\ifcase \nopiecesix

\input{simulations_circuit/simulation405}

\iflanguage{french}{Donnons aussi un circuit contenant des pièces de types   \pieces\ and  \piecesb\ en  figure \ref{exemple_circuit_plan_piecesix}.}{We give also a circuit with 
types \pieces\ and  \piecesb\ pieces in Figure~\ref{exemple_circuit_plan_piecesix}.}

\or
\fi

\or
\fi

%%%%%%%%%%%%%%%%%%%%%%%%%%%%%%%%%%%%%%%%%%%%%%%%%%%
%\input{simulations_circuit/simulation_new_2000}
% fichier tex crée par MaTeXBuild02 le 29-Feb-2016 16:28:06
 % A compiler avec MaTeXBuild02('simulation_new_2000',0)
% après le fichier 'enumeration_construction_circuit_new.matex'

\begin{figure}
\centering
%%% sous figure 1
\subfigure[\label{plancircuit3}\iflanguage{french}{Représentation par plan de circuit}{Circuit map representation}]
{\epsfig{file=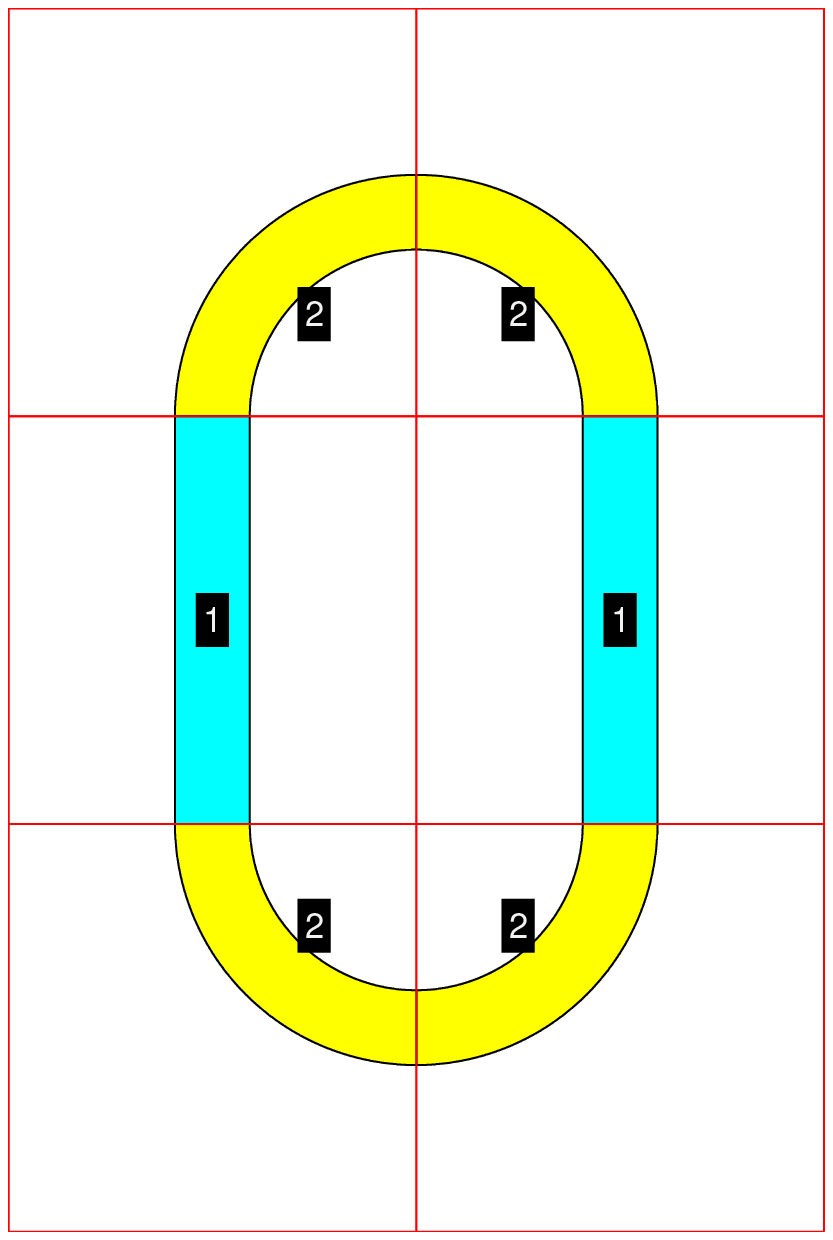, width=5  cm}}
%\qquad
%%% sous figure 2
\subfigure[\label{plancircuit3saw}\iflanguage{french}{Polygone auto-évitant}{Self-avoiding polygon}]
{\epsfig{file=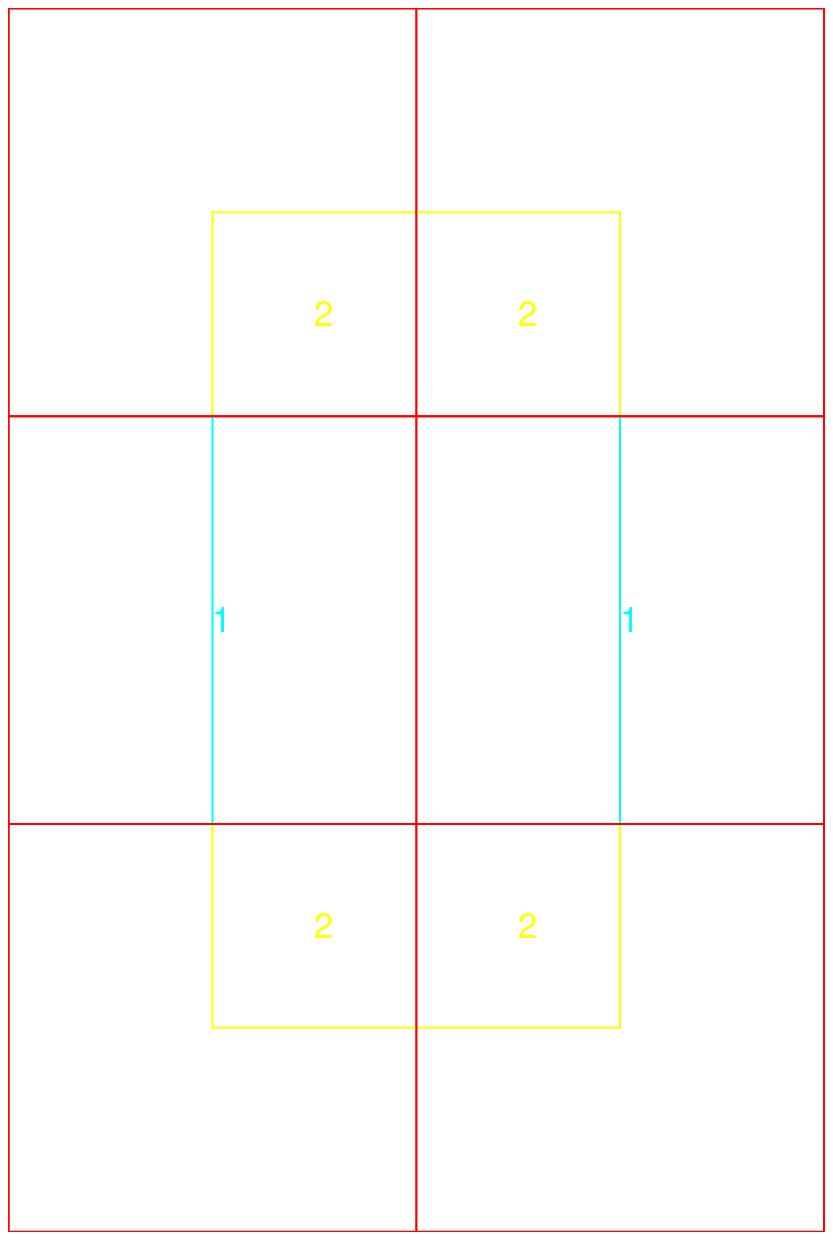, width=5  cm}}
\caption{\label{plancircuit3tot}\iflanguage{french}{Un exemple de circuit}{An example circuit}.}
\end{figure}

\iflanguage{french}{%
Sur la figure \ref{plancircuit3tot}, on a représenté un plan de circuit,
vu comme un plan de jeu (figure \ref{plancircuit3})  ou comme 
un polygone auto-évitant (figure \ref{plancircuit3saw}), obtenu en reliant les centres $c_i$.
Dans ce cas, le circuit est un polygone auto-évitant au sens classique du terme.}%
{%
In Figure \ref{plancircuit3tot}, we have represented a track design, seen as a game map (Figure \ref{plancircuit3}) or as a self-avoiding polygon (Figure \ref{plancircuit3saw}),
obtained by by joining the centers $c_i$.
In this case, the circuit is a self-avoiding polygon in the classical sense of the term.%
}

\begin{figure}
\centering
%%% sous figure 1
\subfigure[\label{plan24max04}\iflanguage{french}{Représentation par plan de circuit}{Circuit map representation}]
{\epsfig{file=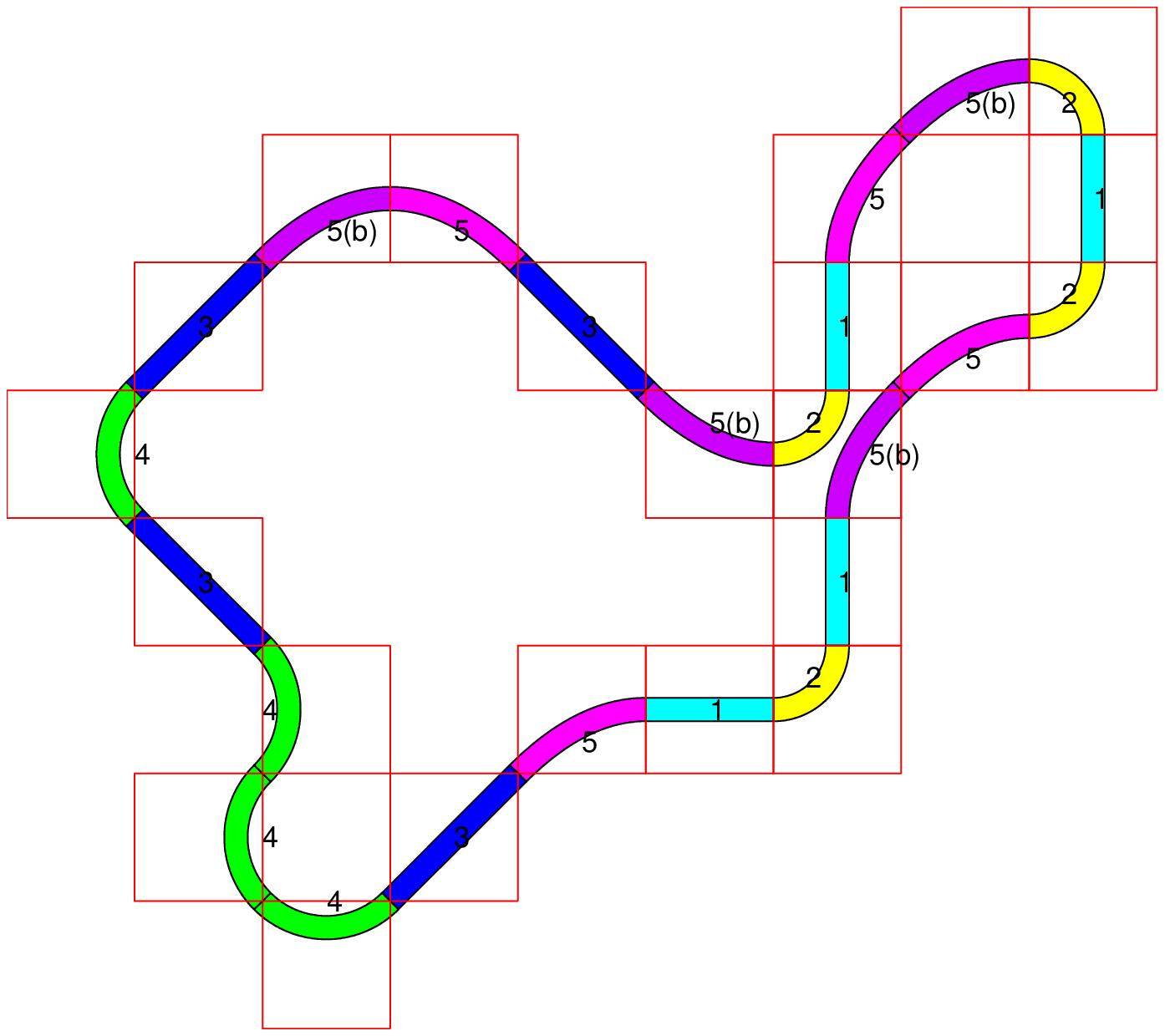, width=11  cm}}
%\qquad
%%% sous figure 2
\subfigure[\label{plan24max04saw}\iflanguage{french}{Polygone auto-évitant}{Self-avoiding polygon}]
{\epsfig{file=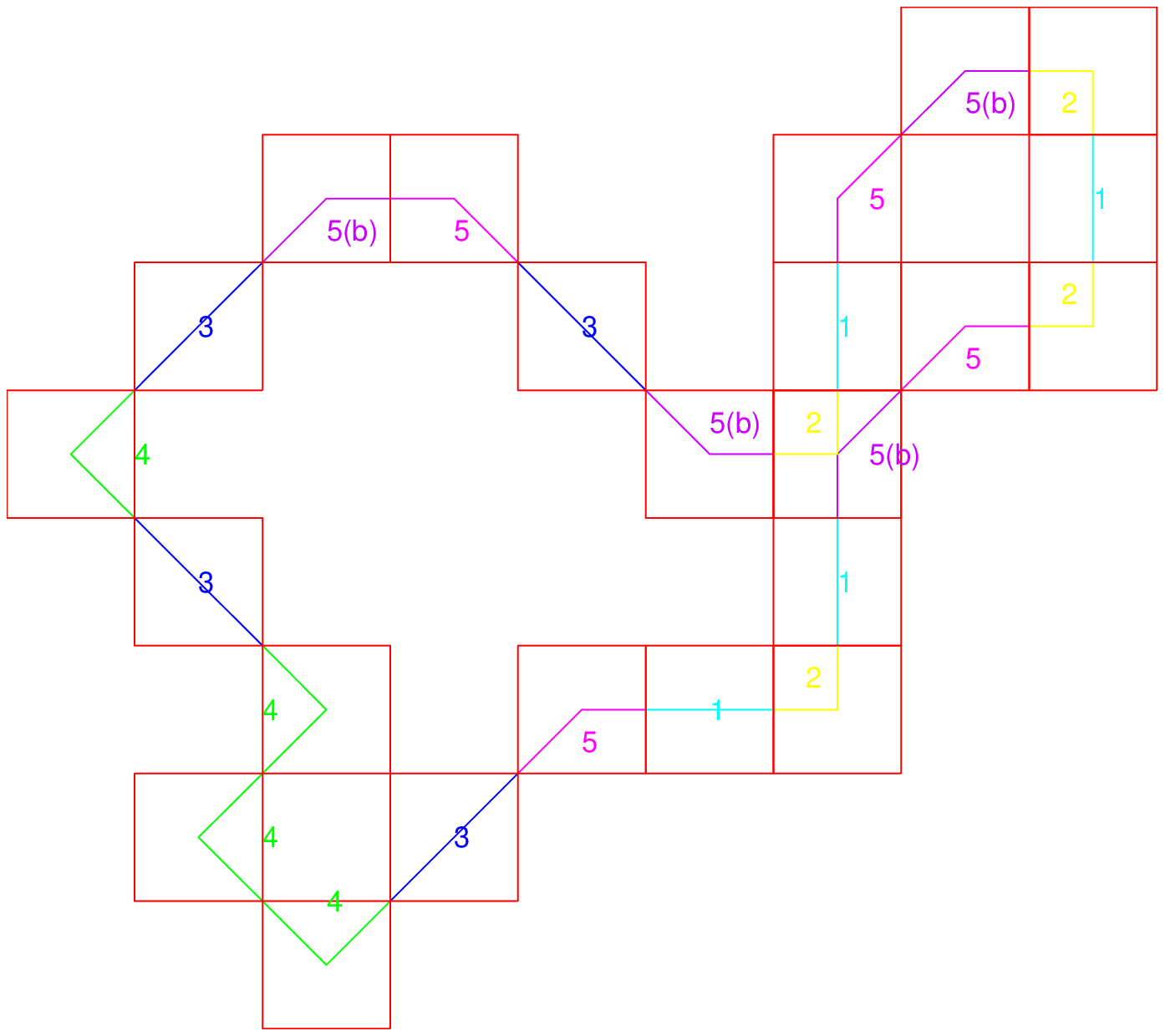, width=11  cm}}
\caption{\label{plan24max04tot}\iflanguage{french}{Un exemple de circuit}{An example circuit}.}
\end{figure}

\iflanguage{french}{%
Sur la figure \ref{plan24max04tot}, on a représenté un plan de circuit,
vu comme un plan de jeu (figure \ref{plan24max04})  ou comme 
un polygone auto-évitant (figure \ref{plan24max04saw}). Cette dernière figure met en évidence les différences essentielles qui 
existent entre nos polygones auto-évitants et ceux de la littérature : 
les nôtres autorisent des carrés successifs à avoir un commun un sommet et deux carrés différents peuvent être occupés par deux
portions différentes de trajectoires.}%
{%
In Figure \ref{plan24max04tot}, we have represented a track design seen as a game map (Figure \ref{plan24max04}) or as a self-avoiding polygon (Figure \ref{plan24max04saw}). The latter figure highlights the essential differences which exist between our self-avoiding polygons and those in the literature; ours allow successive squares to have a vertex in common, and two different squares may be occupied by two different trajectory parts.%
}
%%%%%%%%%%%%%%%%%%%%%%%%%%%%%%%%%%%%%%%%%%%%%%%%%%%

\iflanguage{french}{%

%%%%%%%%%%%%%%%%%%%%%%%%%%%%%%%%%%%%%%%%%%%%%%%%%%%%%%%%%%%%%
%%%%%%%%%%%%%%%%%%%%%%%%%%%%%%%%%%%%%%%%%%%%%%%%%%%%%%%%%%%%%
\section{Construction et énumération des circuits}
\label{enumeration}

%%%%%%%%%%%%%%%%%%%%%%%%%%%%%%%%%%%%%%%%%%%%%%%%%%%%%%%%%%%%%
\subsection{Position du problème}
\label{enumerationpositionprobleme}

On rappelle la question qui intéresse l'industriel :
\questionindustrielle\
Pour $1\leq j\leq \nomtypma$, on note $N_j\in \En\cup \{+\infty\}$ les nombres de pièces
\pieceu\ à 
\ifcase \nopiecesix
{\piecesb} disponibles
\or
{\piececb} disponibles
\fi
 et $N$ le nombre total de pièces utilisées et on cherche tous les circuits
qui se referment  contenant en tout exactement $N$ pièces et tels que, pour chaque type de pièce, le nombre de pièces utilisées
soit inférieur à $N_j$.
Nécessairement, on a 
\begin{equation}
\label{eq01}
N\leq \sum_{j=1} ^\nomtypma N_j.
\end{equation}
Le cas $N_j=+\infty$ correspond au cas où le type de pièces concerné n'est \textit{a priori }pas limité. Cependant,
le nombre de pièces de ce type  est nécessairement 
inférieur à $N$.

Un circuit est totalement déterminé par les $N$ centres ${(c_i)}_{1 \leq i \leq N}$ des carrés. Ces centres
étant donnés, il est donc possible, en prenant les milieux de deux centres successifs, de déterminer, pour chaque 
carré, les coordonnées des points $A_i$ et $B_i$, qui correspondent au début et au fin de la courbe $\Gamma_i$ dans le carré
${\mathcal{C}}_i$. Le fait de considérer des tenons/mortaises  oriente le circuit et il est nécessaire de considérer cette orientation 
pour les pièces 
\ifcase \nopiecesix
{\piecec}, {\piececb}, {\pieces} et {\piecesb}. 
\or
{\piecec} et {\piececb}. 
\fi
Choisissons une orientation des tenons/mortaises (ce qui revient à 
orienter  le circuit)  de la façon suivante : si l'on parcourt le circuit dans l'ordre croissant 
des indices de carrés, 1, 2 
, ..., $N$, alors dans chaque carré, la première extrémité de la pièce correspond à une prise femelle, la seconde 
correspondant  à une prise mâle.
Notons 
\ifcase \nopiecesix
$p_i\in\{\pieceu,...,\piecesb\}$,
\or
$p_i\in\{\pieceu,...,\piececb\}$,
\fi
le type de la pièce concernée dans le carré ${\mathcal{C}}_{i}$.
On notera respectivement $A_i$ et $B_i$ (éléments de ${\mathcal{H}}_i$), le début de la courbe $\Gamma_i$, 
correspondant à la prise femelle  et la fin de la courbe  $\Gamma_i$, 
correspondant donc à la prise mâle.
Le numéro $p_i$ dépend donc uniquement des points $A_i$ et $B_i$.
Par exemple, si ces deux points sont deux sommets successifs, la pièce est de type \pieceq.

Une autre façon de voir cela est de remarquer que chaque pièce est totalement déterminée par la position relative du carré
contenant la pièce précédente et de celui contenant la pièce suivante
ainsi que la nature des points $A_i$ et $B_i$ (c'est-à-dire, être un sommet ou un milieu).
\`A cet effet, pour $i\in\{2,...,N-1\}$,  considérons l'angle
\begin{equation} 
\label{defalpha}
\alpha_i=\widehat{\left(\overrightarrow{c_{i-1}c_i},\overrightarrow{c_ic_{i+1}}\right)}\in [0,2\pi[.
\end{equation} %
}%
{%

%%%%%%%%%%%%%%%%%%%%%%%%%%%%%%%%%%%%%%%%%%%%%%%%%%%%%%%%%%%%%
%%%%%%%%%%%%%%%%%%%%%%%%%%%%%%%%%%%%%%%%%%%%%%%%%%%%%%%%%%%%%
\section{Construction and enumeration of the circuits}
\label{enumeration}

%%%%%%%%%%%%%%%%%%%%%%%%%%%%%%%%%%%%%%%%%%%%%%%%%%%%%%%%%%%%%
\subsection{Posing the problem}
\label{enumerationpositionprobleme}

Recall the question which is of interest to the manufacturer:
\questionindustrielle\
For $1\leq j\leq \nomtypma$, we denote by $N_j\in \En\cup \{+\infty\}$ the number of pieces \pieceu\ to 
\ifcase \nopiecesix
{\piecesb} available,
\or
{\piececb} available, 
\fi
 and by $N$ the total number of pieces used, and we seek all of the circuits which close containing exactly $N$ pieces in all, and such that, for each type of piece, the number of pieces used in less than $N_j$.
We necessarily have
\begin{equation}
\label{eq01}
N\leq \sum_{j=1} ^\nomtypma N_j.
\end{equation}
The case $N_j=+\infty$ corresponds to the case where the type of pieces concerned is not \textit{a priori} limited. However, the number of pieces of this type is necessarily less than $N$.

A circuit is totally determined by the $N$ centers ${(c_i)}_{1 \leq i \leq N}$ of the squares. These centers being given, it is therefore possible to determine, by taking the middles of two successive centers, the coordinates of the points $A_i$ and $B_i$ for each square, which correspond to the start and end of the curve $\Gamma_i$ in square ${\mathcal{C}}_i$. Consideration of the mortises and tenons orients the circuit, and it is necessary to consider this orientation for pieces 
\ifcase \nopiecesix
{\piecec}, {\piececb}, {\pieces} and {\piecesb}. 
\or
{\piecec} and {\piececb}. 
\fi
Let us choose an orientation of the mortises and tenons (which amounts to orienting the circuit) in the following way: if the circuit is traversed in increasing order of the square indices, $1,2, \ldots ,$ then in each square, the first extremity of the piece corresponds to the female connector and the second corresponds to the male one.
We denote by
\ifcase \nopiecesix
$p_i\in\{\pieceu,...,\piecesb\}$,
\or
$p_i\in\{\pieceu,...,\piececb\}$,
\fi
the type of piece concerned in square ${\mathcal{C}}_{i}$.
We will write respectively $A_i$ and $B_i$ (elements of ${\mathcal{H}}_i$), for the start of the curve $\Gamma_i$, corresponding to the female connector, and the end of the curve $\Gamma_i$, corresponding to the male connector. The number  $p_i$ therefore depends only on the points $A_i$ and $B_i$. For example, if these two points are two successive vertices, the piece is of type \pieceq.

Another way to see this is to notice that each piece is totally determined by the relative position of the square containing the previous piece and the one containing the following piece, as well as the nature of the points $A_i$ and $B_i$ (that is, being a vertex or middle). To this end, 
for $i\in\{2,...,N-1\}$, 
we consider the angle
\begin{equation} 
\label{defalpha}
\alpha_i=\widehat{\left(\overrightarrow{c_{i-1}c_i},\overrightarrow{c_ic_{i+1}}\right)}\in [0,2\pi[.
\end{equation} %
}%
\begin{figure}[h] 
\psfrag{cm1}{${\mathcal{C}}_{i-1}$}
\psfrag{cn}{${\mathcal{C}}_{i}$}
\psfrag{cpu}{${\mathcal{C}}_{i+1}$}
\psfrag{Mm1}{$c_{i-1}$}
\psfrag{Mn}{$c_i$}
\psfrag{Mpu}{$c_{i+1}$}
\psfrag{a}{$\alpha_{i}$}
\begin{center} 
\epsfig{file=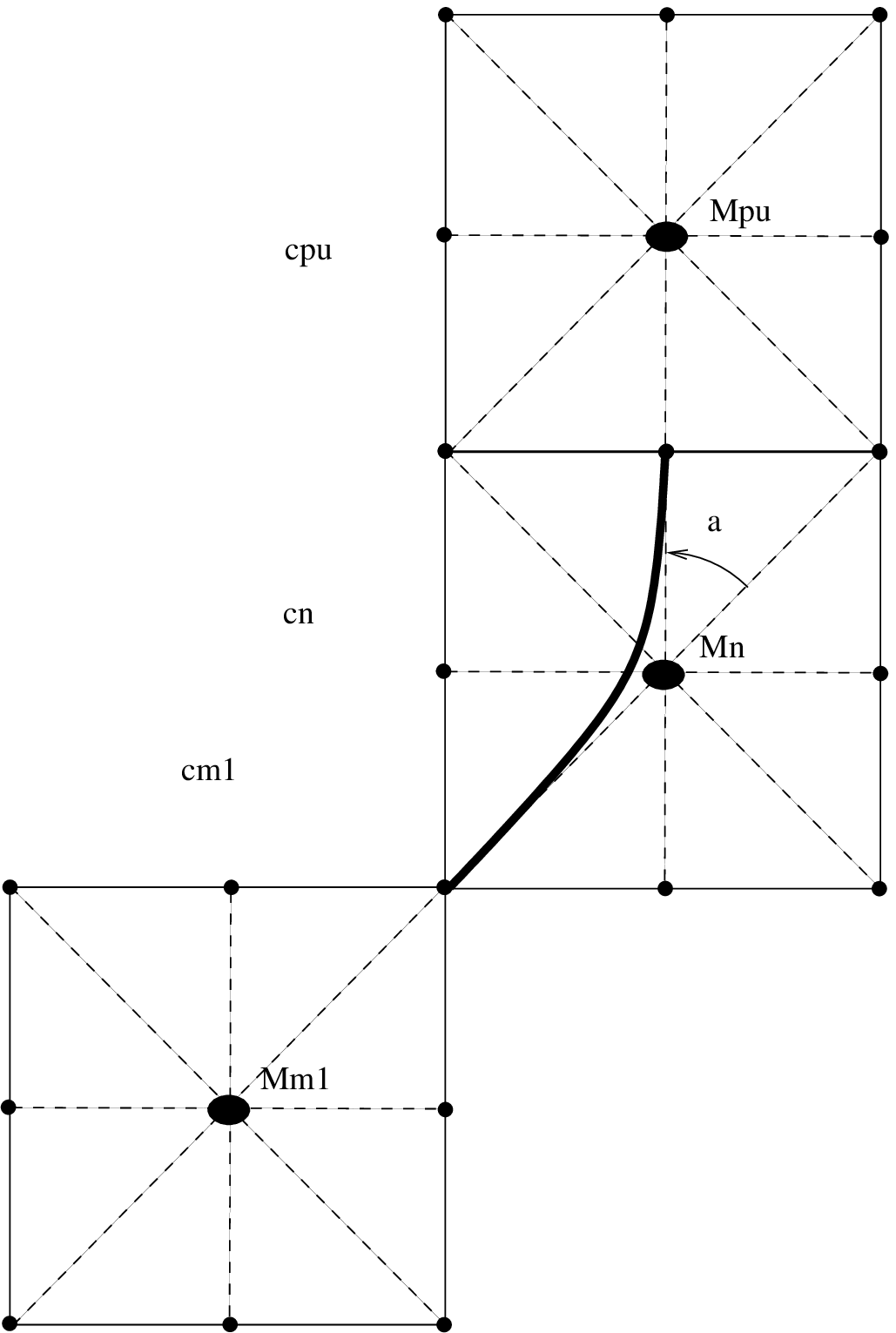, width=6 cm} 
\end{center} 
\caption{\label{piece_deux_carre}\iflanguage{french}%
%{Chaque pièce est définie par les deux carrés l'encadrant.}{Each 
%piece is defined by the two squares framing it}} 
%\caption{\label{piece_deux_carre}\iflanguage{french}%
{Chaque pièce est définie par les deux carrés la limitant et le carré 
auquel elle appartient.}%
{Each piece is defined by the two squares delimiting it and the square to which it belongs.}%
}
\end{figure}
\iflanguage{french}{%
Voir la figure \ref{piece_deux_carre}.
Si on pose 
\begin{equation}
\label{defKalpha}
\forall i\in\{2,...,N-1\},\quad 
\alpha_i=\frac{k_i \pi}{4},
\end{equation}
les seules valeurs possibles de $k_i$ décrivent l'ensemble
\ifcase \nopiecesix
$\{0,1,2,3,5,6,7\}.$
\or
$\{0,1,2,6,7\}.$
\fi
\begin{table}[h]
\begin{center}
\begin{tabular}{|l|l|l|l|}
\hline
$A_i$ & $B_i$ &
$k_i$ &
type de pièce 
\\
\hline
\hline
milieu &milieu & $0$& \pieceu
\\
\hline
milieu &milieu & $2$ ou $6$ & \pieced
\\
\hline
sommet &sommet & $0$& \piecet
\\
\hline
sommet &sommet & $2$ ou $6$ & \pieceq
\\
\hline
milieu &sommet & $1$ ou $7$ & \piecec
\\
\hline
sommet& milieu & $1$ ou $7$ & \piececb
\\
\hline
\ifcase \nopiecesix
milieu& sommet & $3$ ou $5$ & \pieces
\\
\hline
sommet& milieu & $3$ ou $5$ & \piecesb
\\
\hline
\or
\fi
\end{tabular}
\vspace{1 cm}
\caption{\label{typepiece}Numéros  $p_i$ des pièces en fonction de $A_i$, $B_i$ et $k_i$}
\end{center}
\end{table}
Les relations entre ces éléments sont donnés dans le tableau \ref{typepiece}.
Sur la figure \ref{piece_deux_carre}, l'exemple montre le cas d'une pièce 
correspondant à $k_i=1$, donc de type \piececb.

\begin{remark}
\label{remarkalpha}
Dans le cas où $k_i\not =0$, la valeur de $k_i$ permet de déterminer le sens du virage de la pièce en question (droite ou gauche).
Pour 
\ifcase \nopiecesix
$k_i\in\{1,2,3\}$,
\or
$k_i=1$ ou $k_i=2$, 
\fi
 la pièce tourne vers la gauche et pour
\ifcase \nopiecesix
$k_i\in\{5,6,7\}$,
\or
$k_i=6$ ou $k_i=7$, 
\fi
 elle tourne vers la droite. 
On peut donc  aussi adjoindre un signe aux numéros des pièces non rectilignes. 
 Seule la valeur absolue est importante pour le comptage des différents types de pièces, mais nous verrons plus loin
qu'il peut être nécessaire de conserver ce signe pour orienter le sens dans lequel tourne cette pièce. 
\end{remark}%
}%
{%
See Figure~\ref{piece_deux_carre}.
If we let
\begin{equation}
\label{defKalpha}
\forall i\in\{2,...,N-1\},\quad 
\alpha_i=\frac{k_i \pi}{4},
\end{equation}
then the sole possible values of $k_i$ describe the set
\ifcase \nopiecesix
$\{0,1,2,3,5,6,7\}.$
\or
$\{0,1,2,6,7\}.$
\fi
%$\{0,1,2,6,7\}.$
\begin{table}[h]
\begin{center}
\begin{tabular}{|l|l|l|l|}
\hline
$A_i$ & $B_i$ &
$k_i$ &
type of piece
\\
\hline
\hline
middle & middle & $0$& \pieceu
\\
\hline
middle & middle & $2$ or $6$ & \pieced
\\
\hline
vertex & vertex & $0$& \piecet
\\
\hline
vertex & vertex & $2$ or $6$ & \pieceq
\\
\hline
middle & vertex & $1$ or $7$ & \piecec
\\
\hline
vertex & middle & $1$ or $7$ & \piececb
\\
\hline
\ifcase \nopiecesix
middle&vertex  & $3$ or $5$ & \pieces
\\
\hline
vertex& middle & $3$ or $5$ & \piecesb
\\
\hline
\or
\fi
\end{tabular}
\vspace{1 cm}
\caption{\label{typepiece}Piece numbers $p_i$ as a function of $A_i$, $B_i$ and $k_i$}
\end{center}
\end{table}
The relations between these elements are given in Table~\ref{typepiece}. In Figure~\ref{piece_deux_carre}, the example shows the case of a piece corresponding to $k_i=1$, that is, of type \piececb.

\begin{remark}
\label{remarkalpha}
In the case where $k_i\not =0$, the value of $k_i$ allows the determination of the turning direction of the piece in question (right or left).
For 
\ifcase \nopiecesix
$k_i\in\{1,2,3\}$,
\or
$k_i=1$ or $k_i=2$, 
\fi
the piece turns towards the left, and for 
\ifcase \nopiecesix
$k_i\in\{5,6,7\}$,
\or
$k_i=6$ or $k_i=7$, 
\fi
it turns towards the right. One may therefore also associate a sign to the numbers of the curved pieces. Only the absolute value is important for counting of the different types of piece, but we will see later that it may be necessary to keep this sign in order to orient the direction in which this piece turns.
\end{remark}%
}

\iflanguage{french}{%

%%%%%%%%%%%%%%%%%%%%%%%%%%%%%%%%%%%%%%%%%%%%%%%%%%%%%%%%%%%%%
\subsection{Description de tous les circuits}
\label{enumerationtouscircuit}

Dans un premier temps, décrivons la recherche de tous les circuits
à $N$ pièces
dont le centre du premier carré est arbitrairement égal à l'origine et 
le centre du dernier est donné par $(x,y) \in \mathbb{Z}^2$.

Le deuxième carré, nécessairement voisin immédiat de l'origine, peut être 
donc choisi parmi 8 carrés possibles. Pour chacun de ces choix, 
on peut choisir librement les valeurs de $k_i$, $2\leq i \leq N-1$ dans 
\ifcase \nopiecesix
$\{0,1,2,3,5,6,7\},$
\or
$\{0,1,2,6,7\},$
\fi
ce qui fixe les valeurs de $c_i$,  $1\leq i \leq N$,
ainsi que les valeurs des angles $\alpha_i$, $2\leq i \leq N-1$.
On se donne, de plus, le premier point $A_1$ de la première courbe (dans ${\mathcal{H}}_1$)
et le dernier point $B_N$ de la première courbe (dans ${\mathcal{H}}_N$).
Le point $B_2$ est connu ; on en déduit  le numéro de pièce $p_1$.
De même, $p_N$ est connu. 
Pour tout $i$, $2\leq i\leq N-1$, 
les natures de tous les points $A_i$ et $B_i$ et la valeur de $k_i$ sont connues, dont
on déduit la valeur de $p_i$, grâce au tableau \ref{typepiece}.
Sur tous les circuits ainsi définis, on ne conservera que ceux correspondant à $c_N=(x,y)$.

Ainsi, en faisant varier un certain nombre de paramètres indépendants, on est capable
d'énumérer tous les circuits, de façon géométrique (détermination des $c_i$) et constitutive
(détermination des $p_i$) allant de l'origine à un point donné.

Si l'on cherche maintenant tous les  circuits qui se rebouclent, on considérera de même que
le centre du premier carré est arbitrairement égal à l'origine.
Le dernier carré $c_N$ ne pourra être qu'un des 8 voisins du premier. Par symétrie et rotation, on pourra se contenter 
de choisir $c_N\in\{(1,0),(1,1)\}$. Pour chacun des choix de $c_N$, on applique ce que l'on a vu précédemment pour
déterminer tous les circuits allant de l'origine à $c_N$. Dans ce cas, les sommets $A_1$ et $B_N$ sont nécessairement connus
et égaux, puisqu'ils seront nécessairement le sommet ou le milieu commun du premier et du dernier carré.
Notons qu'on peut poser, de façon analogue à  \eqref{defalpha} et \eqref{defKalpha},
\begin{subequations}
\label{defalphabis}
\begin{align} 
&\alpha_1=\widehat{\left(\overrightarrow{c_{N}c_1},\overrightarrow{c_1c_{2}}\right)}\in [0,2\pi[,\\
&\alpha_N=\widehat{\left(\overrightarrow{c_{N-1}c_N},\overrightarrow{c_Nc_{1}}\right)}\in [0,2\pi[,
\end{align} 
\end{subequations}
et 
\begin{equation}
\label{defKalphabis}
\forall i\in\{1,N\},\quad 
\alpha_i=\frac{k_i \pi}{4}.
\end{equation}
Nous conserverons uniquement les circuits tels que 
$k_1$ et $k_N$ appartiennent à l'ensemble
\ifcase \nopiecesix
$\{0,1,2,3,5,6,7\}.$
\or
$\{0,1,2,6,7\}.$
\fi
On en déduit ainsi la valeur des $N$ entiers ${(p_i)}_{1\leq i\leq N}$.
Enfin, 
%éventuellement, 
sur tous ces circuits, on ne conservera que ceux dont le nombre total de chaque pièce de chaque type est inférieur à $N_j$.

Cette méthode, fondée sur un balayage de paramètres, obtenus comme produit cartésien d'ensembles finis,
est très coûteuse en temps  et est vite limitée pour des valeurs de $N$ trop grandes (nous reviendrons sur ce point en section 
\ref{limitationinfo}). Dans 
\cite{%
MR2065628, 
MR2883859,%
MR1985492,%
MR1718791,%
Guttmann2012,%
Guttmann2012b,%
MR2902304},
il existe des technique beaucoup plus subtiles et parallélisables pour dénombrer tous 
les chemins ou les polygones autoévitants. Cependant, 
nous avons déjà signalé dans l'introduction, les différences essentielles entre nos circuits
et les chemins ou polygones autoévitants, ce qui peut rendre l'utilisation des méthodes de \cite{MR2065628} 
inopérante ici.
De plus, les trajets géométriques de ces circuits se doivent d'être déterminés pour vérifier la contrainte
portant sur les $N_j$. Ces circuits devront aussi subir des éliminations 
pour prendre en compte la répétition de circuits isométriques (voir section \ref{isometrie}) 
ainsi que la vérification d'une contrainte locale, qui n'est pas celle des chemins ou polygones autoévitants (voir section \ref{contrainteslocales}).
Enfin, il nous semble important, outre de compter tous les circuits, d'en présenter aussi la totalité, 
pour les petites valeurs de $N$ tout du moins.

L'inconvénient de cette méthode est de passer, d'abord par la détermination de 
tous les circuits possibles, ce représente un très grand nombre de circuits,
puis de ne conserver que ceux qui se rebouclent. 
Une autre méthode serait de déterminer tous 
les chemins géométriques se rebouclant, soit par produit cartésien, soit de façon récursive et parallèle. 
Sur tous ces chemins, on ne retient que ceux qui correspondent à un circuit.
Cette méthode semble être aussi longue que celle que l'on a montrée précédemment et présente tout autant
l'inconvénient d'avoir à déterminer un grand nombre d'éléments pour n'en conserver qu'une partie.

}%
{%

%%%%%%%%%%%%%%%%%%%%%%%%%%%%%%%%%%%%%%%%%%%%%%%%%%%%%%%%%%%%%
\subsection{Description of all circuits}
\label{enumerationtouscircuit}

Firstly, we describe the search for all circuits with $N$ pieces for which the center of the first piece is arbitrarily equal to the origin, and center of the last is given by $(x,y) \in \mathbb{Z}^2$.

The second square, which is necessarily neighboring the origin, may be therefore chosen among 8 possible squares. For each of these choices, one may choose freely the values of $k_i$, $2\leq i \leq N-1$ from
\ifcase \nopiecesix
$\{0,1,2,3,5,6,7\},$
\or
$\{0,1,2,6,7\},$
\fi
which fixes the values of $c_i$,  $1\leq i \leq N$, as well as the values of the angles $\alpha_i$, $2\leq i \leq N-1$. We are given, moreover, the first point $A_1$ of the first curve (in ${\mathcal{H}}_1$) and the last point $B_N$ of the last curve (in ${\mathcal{H}}_N$).
The point $B_2$ is known; from this we deduce the number of piece $p_1$. Likewise, $p_N$ is known. For all $i$, $2\leq i\leq N-1$, the natures of all of the points $A_i$ and $B_i$ and the value of $k_i$ are known, from which we deduce the value of $p_i$ using Table~\ref{typepiece}. Of all of the circuits thus defined, we will keep only those corresponding to $c_N=(x,y)$.

Thus, by varying a certain number of independent parameters, we are capable of enumerating all of the circuits, in a geometric (the determination of the $c_i$) and constitutive (determination of the $p_i$) way, going from the origin to a given point.

If one now seeks all of the circuits which form a loop, one will similarly consider the center of the first square to be arbitrarily equal to the origin. The last square $c_N$ can only be one of the 8 neighbors of the first one. By symmetry and rotation, one may simply choose $c_N\in\{(1,0),(1,1)\}$. For each choice of $c_N$, we apply that which we have seen above to determine all of the circuits going from the origin to $c_N$. In this case, the vertices $A_1$ and $B_N$ are necessarily known and equal, since they will necessarily be the vertex or the common middle of the first and last square. 
Note that one could also set, in a similar manner to \eqref{defalpha} and \eqref{defKalpha},
\begin{subequations}
\label{defalphabis}
\begin{align} 
&\alpha_1=\widehat{\left(\overrightarrow{c_{N}c_1},\overrightarrow{c_1c_{2}}\right)}\in [0,2\pi[,\\
&\alpha_N=\widehat{\left(\overrightarrow{c_{N-1}c_N},\overrightarrow{c_Nc_{1}}\right)}\in [0,2\pi[,
\end{align} 
\end{subequations}
and
\begin{equation}
\label{defKalphabis}
\forall i\in\{1,N\},\quad 
\alpha_i=\frac{k_i \pi}{4}.
\end{equation}
We will keep only the circuits such that $k_1$ and $k_N$ belong to the set
\ifcase \nopiecesix
$\{0,1,2,3,5,6,7\}.$
\or
$\{0,1,2,6,7\}.$
\fi
We hence deduce the values of the $N$ integers ${(p_i)}_{1\leq i\leq N}$.
Finally, of all of these circuits, we will keep only those for which the total number of each piece of each type is less than $N_j$.

This method, based on a parameter sweep, obtained as the Cartesian product of finite sets, is very costly in time, and is quickly limited for values of $N$ which are too large (we will return to this point in Section~\ref{limitationinfo}). In 
\cite{%
MR2065628, 
MR2883859,%
MR1985492,%
MR1718791,%
Guttmann2012,%
Guttmann2012b,%
MR2902304},
there exist some much more subtle and parallelizable techniques for the enumeration of all self-avoiding or polygons walks. However, we have already pointed out in the introduction the essential differences between our circuits and self-avoiding walks or polygons, which may render the use of the methods from \cite{MR2065628} ineffective here. Furthermore, the geometric trajectories of the circuits need to be determined in order to satisfy the constraint concerning the $N_j$. These circuits will also need to undergo eliminations to take into account the repetition of isometric circuits (see Section~\ref{isometrie}), as well as the satisfaction of a local constraint, which is not that of self-avoiding walks or polygons (see Section~\ref{contrainteslocales}). Finally, we consider it important, in addition to counting all the circuits, to present them all, for small values of $N$ at least.

The disadvantage of this method is that it goes first through the determination of all possible circuits, which represents a very large number of circuits, and then keeps only those which form a loop. Another method will be to determine all of the looping geometric paths, either by Cartesian product, or in a recursive and parallel way. From all of these paths, one retains only those which correspond to a circuit. This method seems to be as long as the one shown above, and equally presents the disadvantage of having to determine a large number of elements, only to keep but a portion.%

}

\iflanguage{french}{%

%%%%%%%%%%%%%%%%%%%%%%%%%%%%%%%%%%%%%%%%%%%%%%%%%%%%%%%%%%%%%
\subsection{Prise en compte des isométries}
\label{isometrie}

Commençons par étudier un exemple simple.
On choisira une petite valeur de 
$N$ sans se soucier des contraintes imposées par $N_j$.

\begin{example}
\label{examplesimulation500}
\input{simulations_circuit/simulation500}
\end{example}

De façon plus générale, on tracera tous les circuits obtenus, pour $N$ et $N_j$ donnés.

Une prise en compte des isométries directes se traduira par la comparaison 
de tous les circuits obtenus. Si deux d'entre eux possèdent les mêmes numéros signés de pièces,
à une permutation circulaire près,
l'un des deux sera éliminé.

Dans un second temps, une prise en compte des isométries indirectes sera effectuée.
De même, si deux circuits possèdent des mêmes numéros signés de pièces
qui sont opposés (pour les pièces non rectilignes),
à une permutation circulaire près,
l'un des deux sera éliminé. Pour prendre en compte toutes les isométries indirectes,
il sera aussi nécessaire d'éliminer les circuits en comparant aussi les indices avec des permutations
du type $N$, $N-1$, ..., 2, 1, ce qui revient à considérer que l'on parcourt le circuit dans le sens opposé. 
Dans ce cas, on remplacera  les numéros signés  de pièces $\pm \piecec$ par $\mp \piececb$ 
\ifcase \nopiecesix
(et $\pm \pieces$ par $\mp \piecesb$)
\or
\fi
et vice-versa.
Les deux pièces \piecec\ et  \piececb\ 
\ifcase \nopiecesix
(et $\pieces$ et $\piecesb$)
\or
\fi
sont en effet identiques ; seule l'orientation change.
Cette élimination sera légitime si les nombres de pièces disponibles de types
\piecec\ et  \piececb\ 
\ifcase \nopiecesix
(et $\pieces$ et $\piecesb$) sont identiques.
\or
sont identiques, ce qui sera toujours vrai par la suite (voir remarque \ref{autantpiececinqsix}).
\fi

\begin{remark}
\label{autantpiececinqsix}
Dans le cas des traditionnels polygones autoévitants, le nombre de carrés est nécessairement pair, ce qui n'est plus vrai ici.
Néanmoins, on peut affirmer que, dans tout circuit, 
\ifcase \nopiecesix
la somme des nombres de pièces de type \piecec\ et \pieces\ est égale à la somme 
des nombres de pièces de type \piececb\ et \piecesb.
\or
les nombres de pièces de types
\piecec\ et  \piececb\ 
sont égaux.
\fi
En effet, remarquons que, d'après le tableau \ref{typepiece},
\begin{itemize}
\item
les pièces de type
\pieceu, \pieced, \piecet\ et \pieceq\  relient entre eux deux points de même nature (milieux de côtés ou sommets de carrés) ;
\item
\ifcase \nopiecesix
les pièces de type
\piecec\ et \pieces\  relient un milieu à un sommet (dans cet ordre) ;
\or
les pièces de type
\piecec\ relient un milieu à un sommet  (dans cet ordre) ;
\fi
\item
\ifcase \nopiecesix
les pièces de type
\piececb\ et \piecesb\  relient un sommet  à un milieu  (dans cet ordre) ;
\or
les pièces de type
\piececb\ relient un sommet à un milieu   (dans cet ordre) ;
\fi
\end{itemize}
Un circuit part d'un point et revient au même point. Tous les points du circuits, correspondant aux extrémités des pièces
utilisées sont  soit des sommets, soit des milieux. Il en découle qu'il y a autant de pièces reliant un milieu à un sommet  (dans cet ordre) 
que de pièce reliant un sommet  à un milieu  (dans cet ordre). Sinon, le point de départ ne serait être de même nature que celui d'arrivée.
\end{remark}%

\begin{remark}
\label{remtranslation}
Notons que, dans le dénombrement des traditionnels polygones autoévitants, seules 
les translations sont prises en compte dans les isométries.
Notre problème de dénombrement est donc bien différent de celui de la littérature.
Dans les cas où les deux notions coïncident, cela implique 
que les configurations que l'on obtiendra seront \textit{a priori}
moins nombreuses que celles de la littérature.
\end{remark}%
}%
{%

%%%%%%%%%%%%%%%%%%%%%%%%%%%%%%%%%%%%%%%%%%%%%%%%%%%%%%%%%%%%%
\subsection{Consideration of the isometries}
\label{isometrie}

Let us begin by studying a simple example. We will choose a small value of $N$ without regard for the constraints imposed by $N_j$.

\begin{example}
\label{examplesimulation500}
%%%%%%%%%%%%%%%%%%%%%%%%%%%%%%%%%%%%%%%%%%%%%%%%%%%%%%
%\input{simulations_circuit/simulation500}
% fichier tex crée par MaTeXBuild02 le 22-Sep-2015 03:13:07
% à compiler avec 
% MaTeXBuild02('simulation500',0)
% après le fichier 'enumeration_construction_circuit.matex'

%%%%%%%%%%%%%%%%%%%%%%%%%%%%%%%%%%%%%%%%%%%%%%%%%%%%%%
%\input{./simulations_circuit/circuit_numerique/circuits_complets_5_piecmax_ev0}
% fichier crée par 'presentation_exhaustif_circuit_boucle.m' le 22-Sep-2015 03:13:18
\begin{figure}[h]
\centering
%%% sous figure 1
\subfigure[\label{circuits_complets_5_piecmax_ev01}]
{\epsfig{file=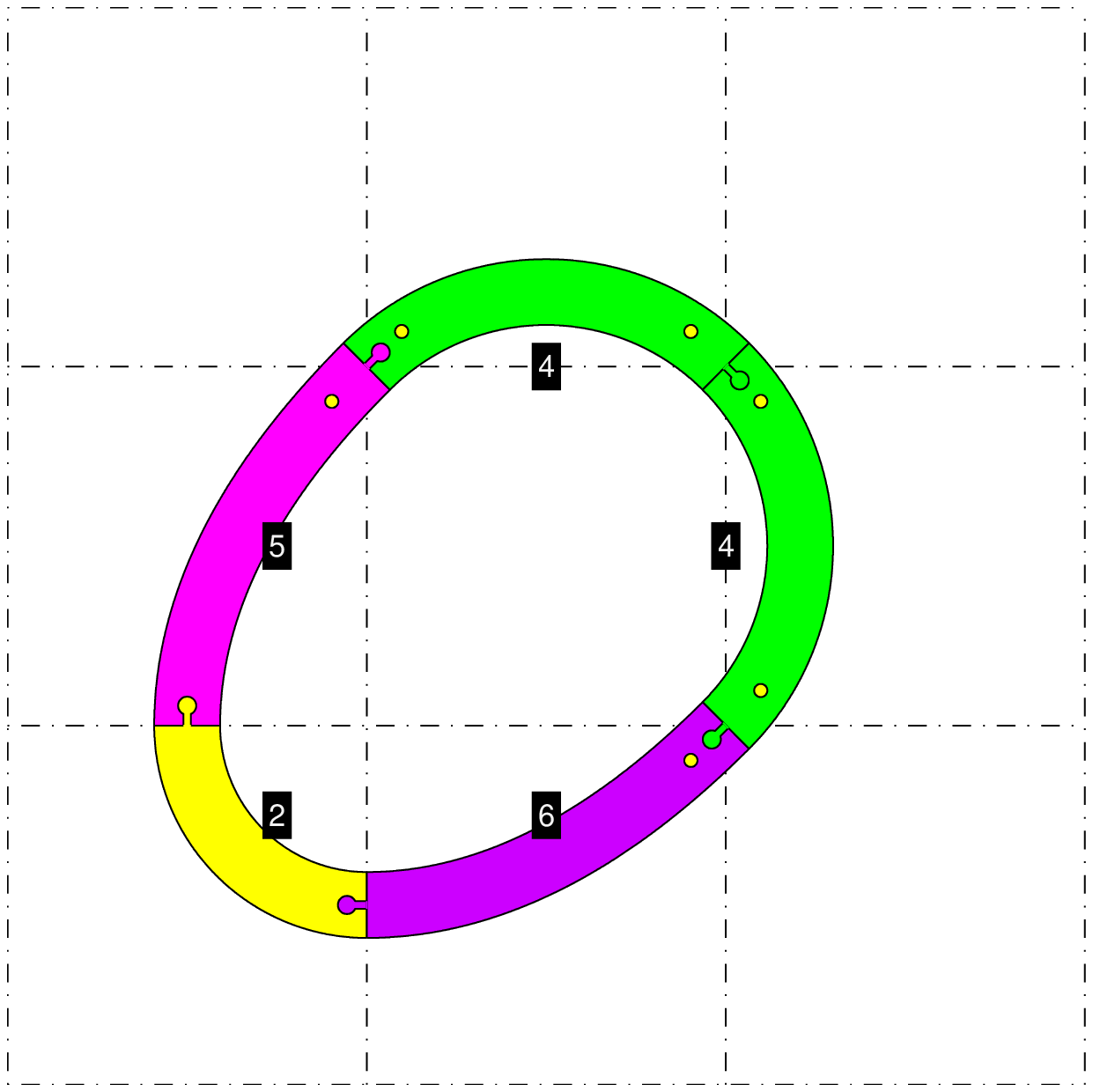, width=5 cm}}
\qquad
%%% sous figure 2
\subfigure[\label{circuits_complets_5_piecmax_ev02}]
{\epsfig{file=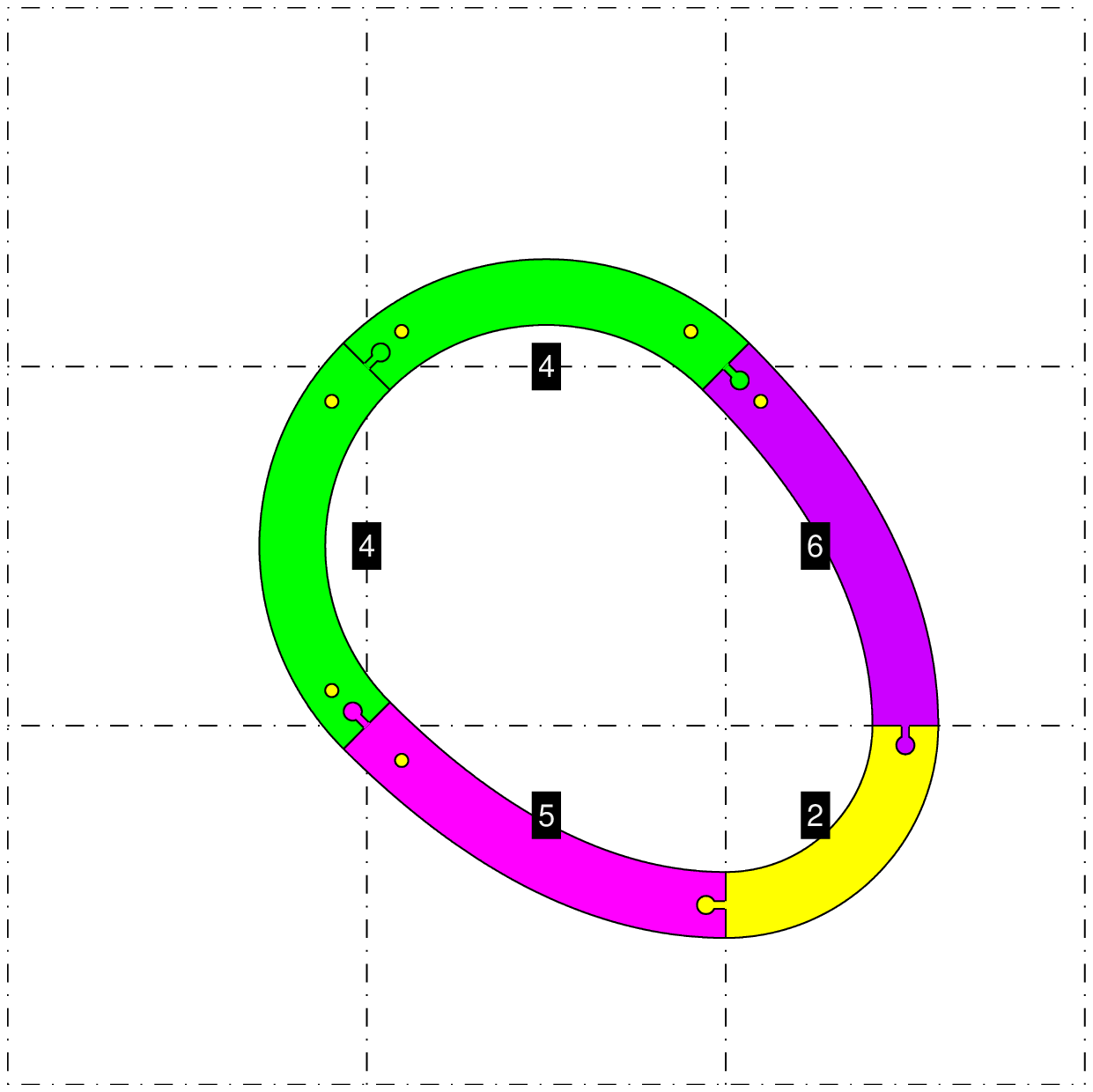, width=5 cm}}
\qquad
%%% sous figure 3
\subfigure[\label{circuits_complets_5_piecmax_ev03}]
{\epsfig{file=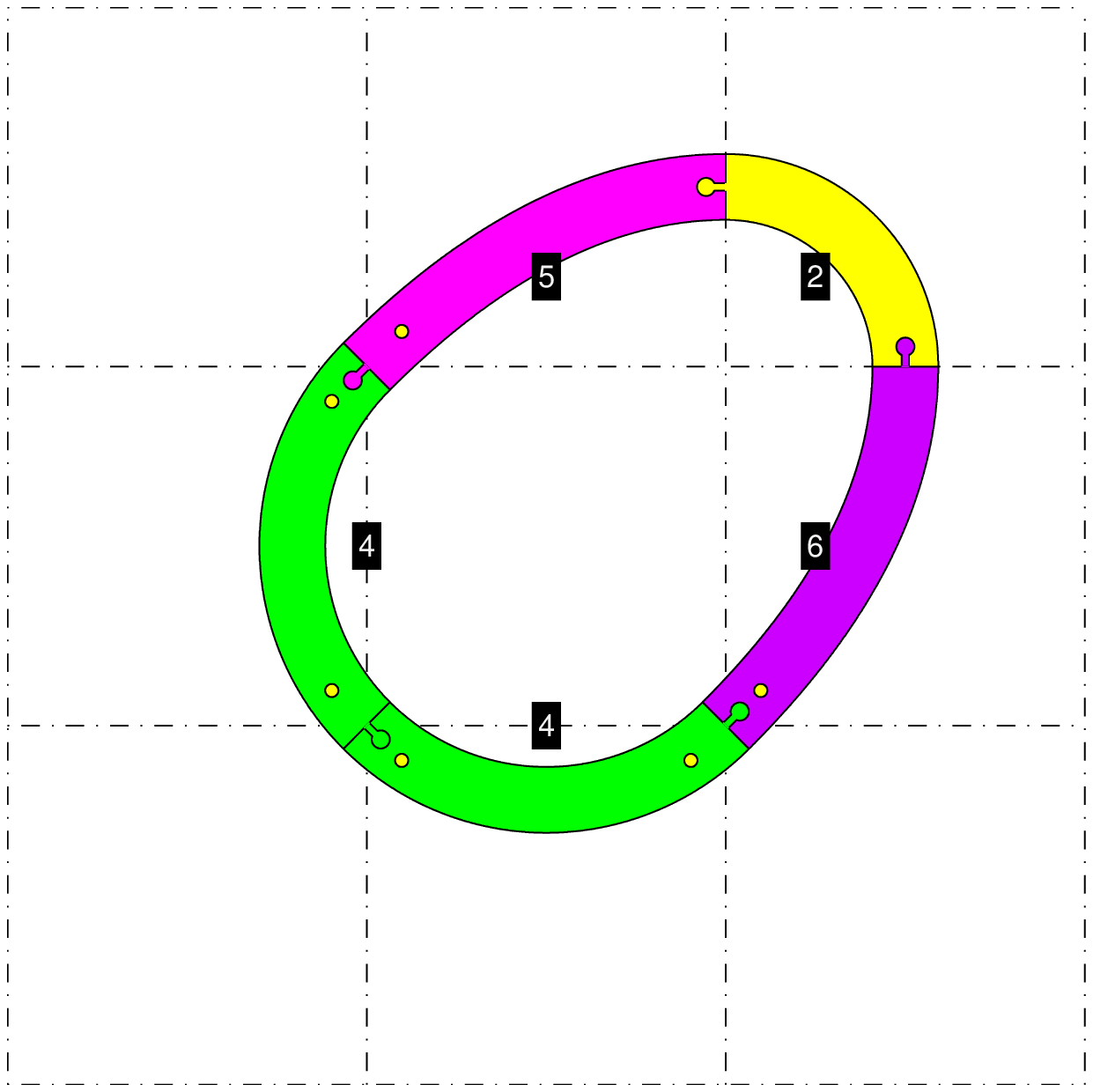, width=5 cm}}
\qquad
%%% sous figure 4
\subfigure[\label{circuits_complets_5_piecmax_ev04}]
{\epsfig{file=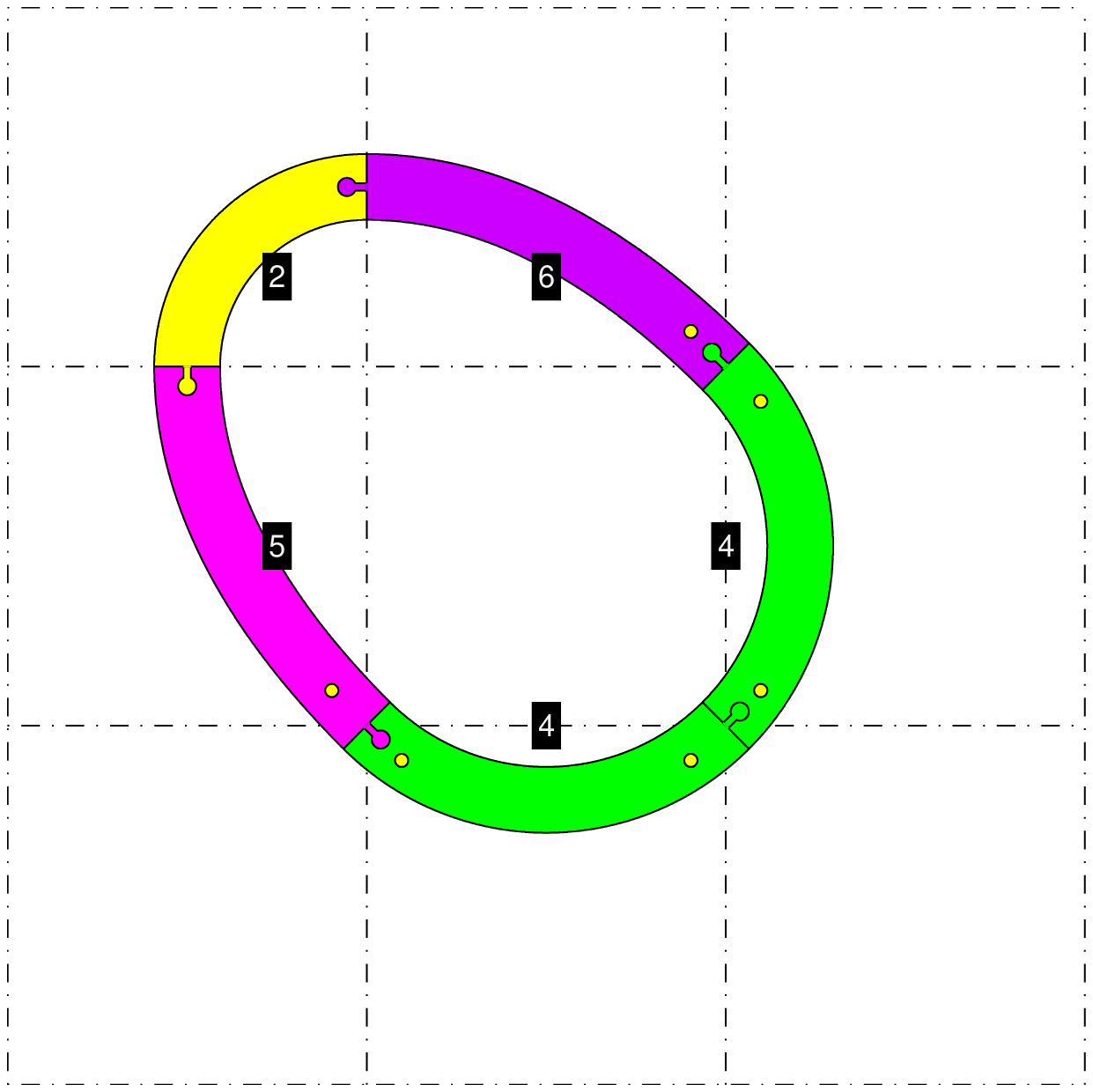, width=5 cm}}
\qquad
%%% sous figure 5
\subfigure[\label{circuits_complets_5_piecmax_ev05}]
{\epsfig{file=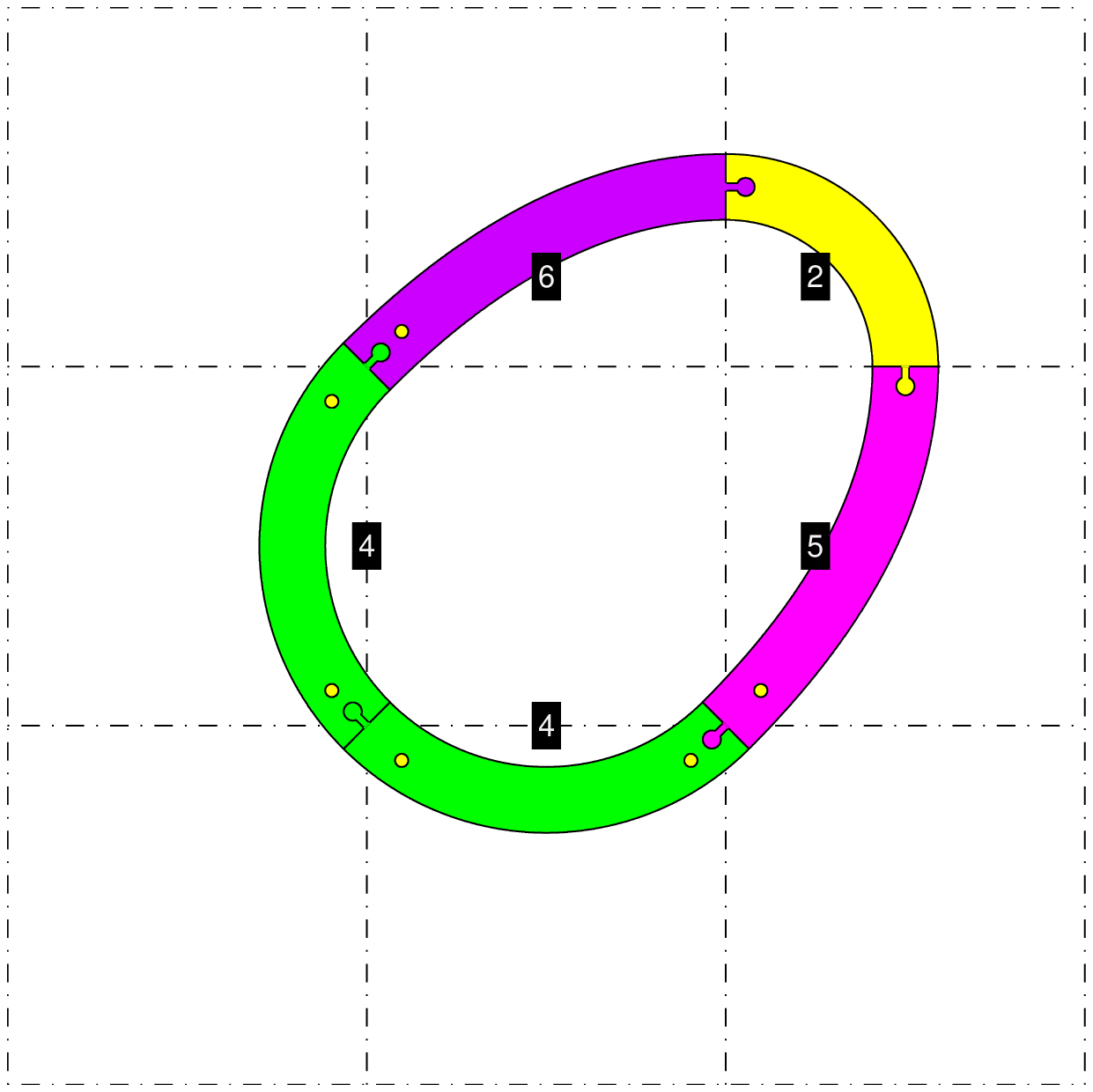, width=5 cm}}
\qquad
%%% sous figure 6
\subfigure[\label{circuits_complets_5_piecmax_ev06}]
{\epsfig{file=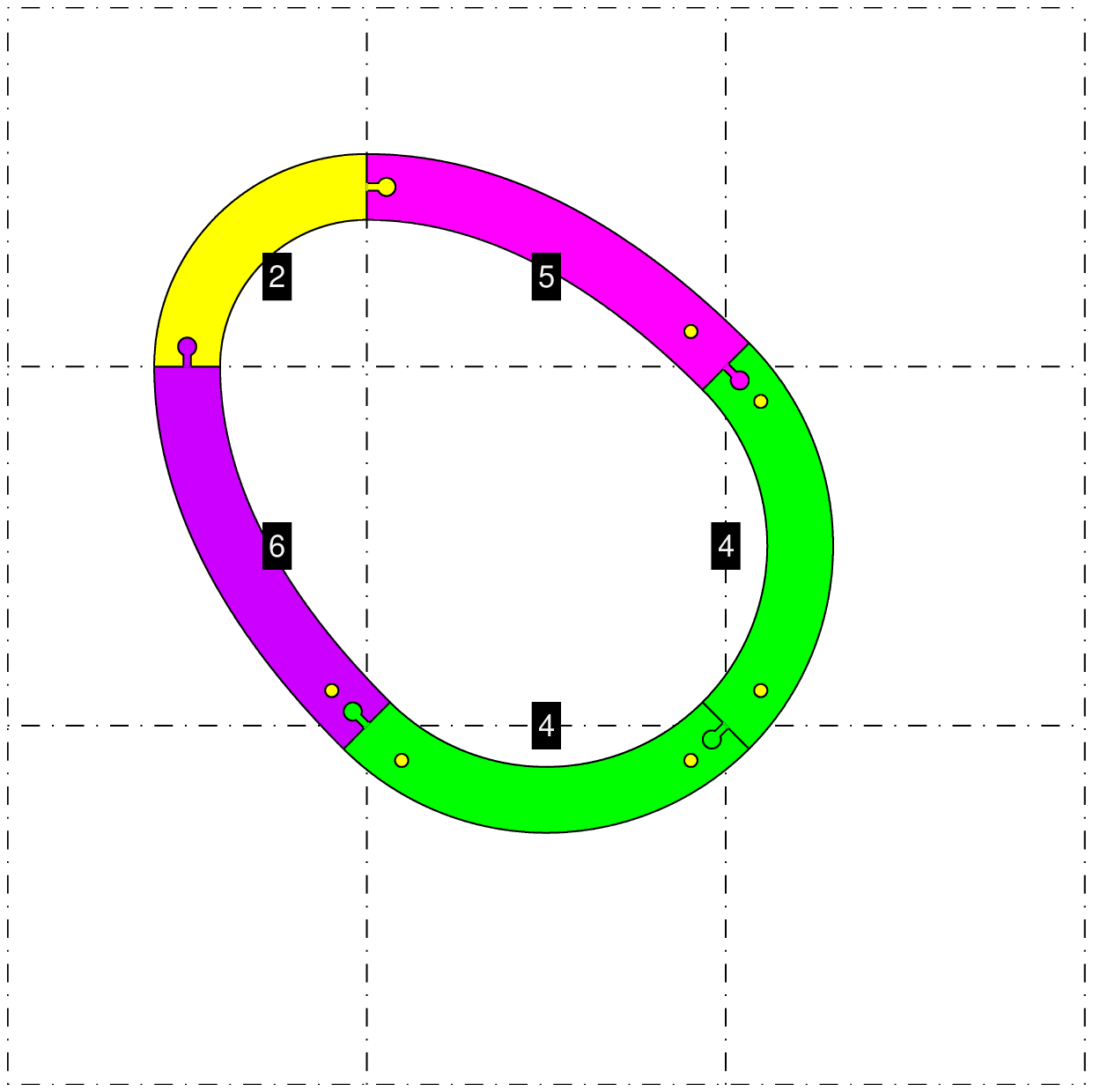, width=5 cm}}
\qquad
%%% sous figure 7
\subfigure[\label{circuits_complets_5_piecmax_ev07}]
{\epsfig{file=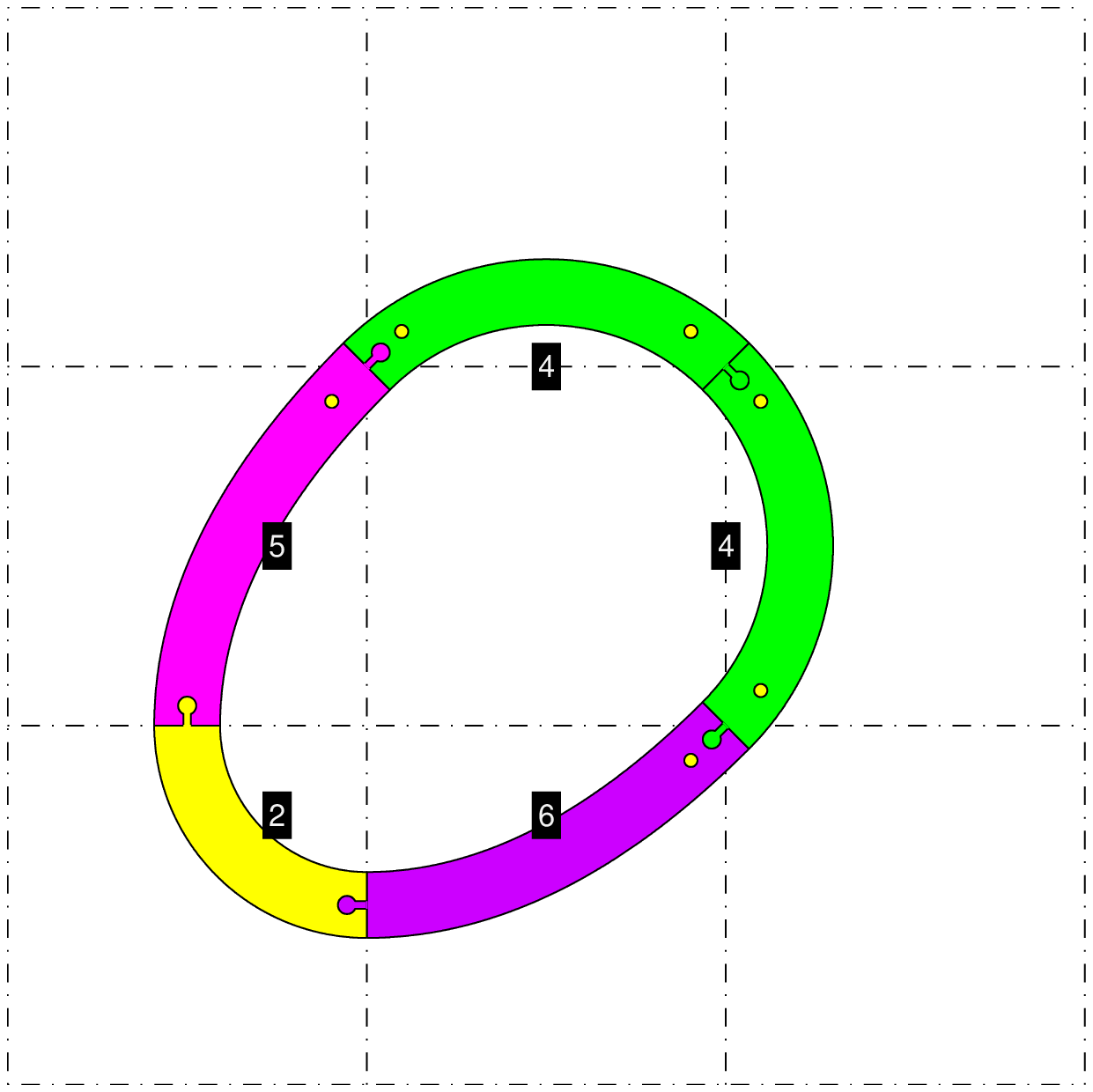, width=5 cm}}
\qquad
%%% sous figure 8
\subfigure[\label{circuits_complets_5_piecmax_ev08}]
{\epsfig{file=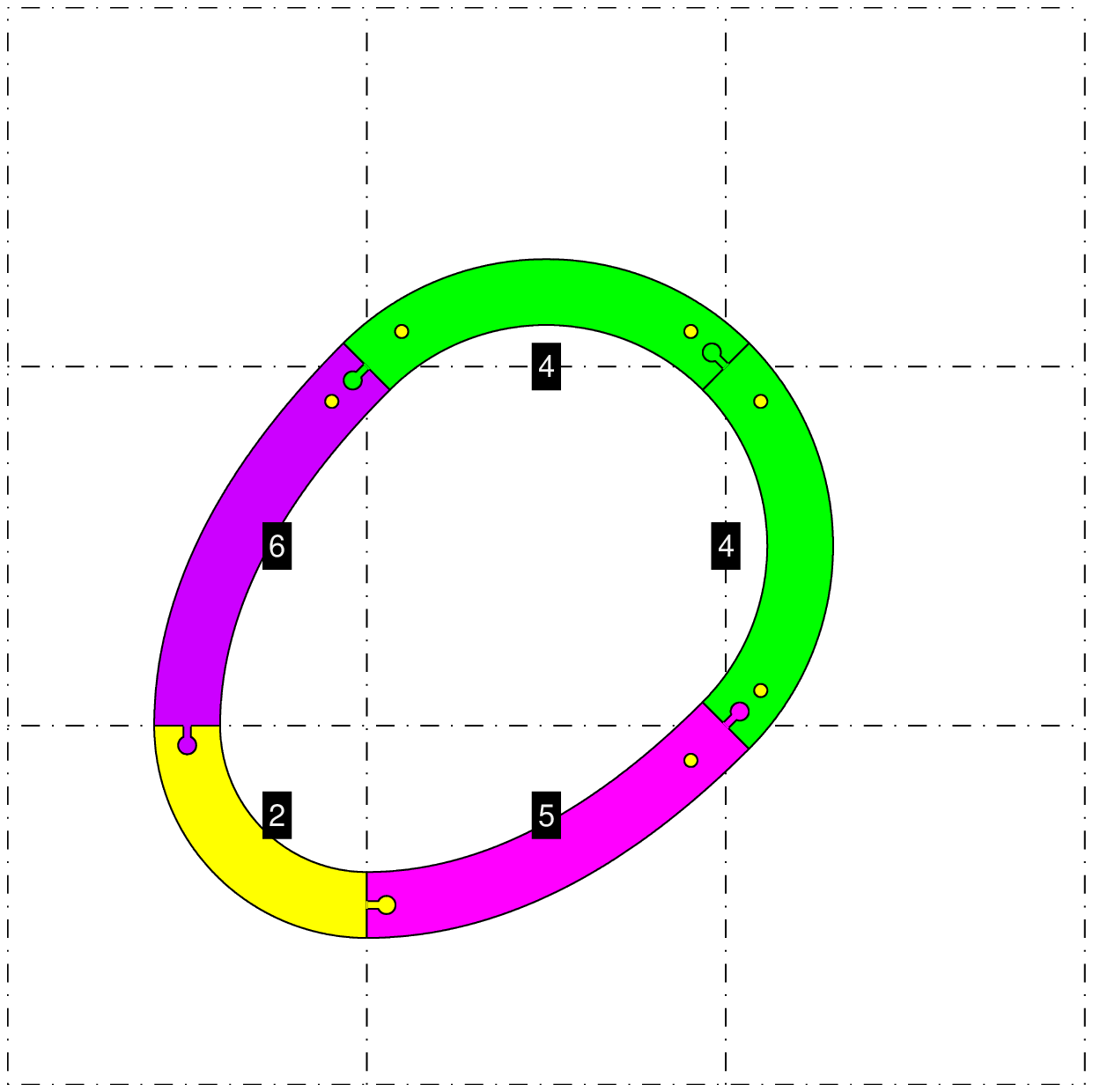, width=5 cm}}
\qquad
%%% sous figure 9
\subfigure[\label{circuits_complets_5_piecmax_ev09}]
{\epsfig{file=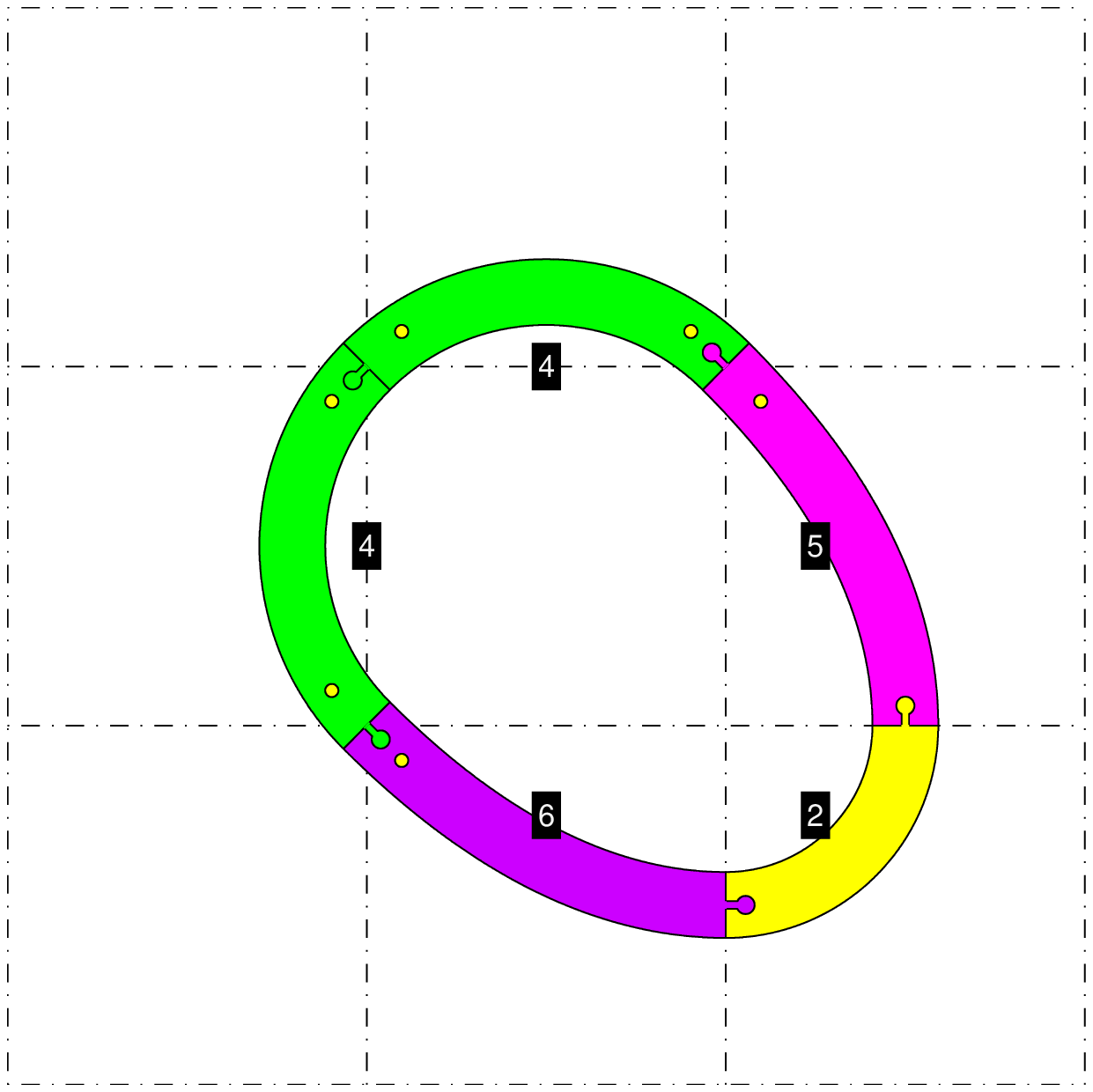, width=5 cm}}
\qquad
%%% sous figure 10
\subfigure[\label{circuits_complets_5_piecmax_ev010}]
{\epsfig{file=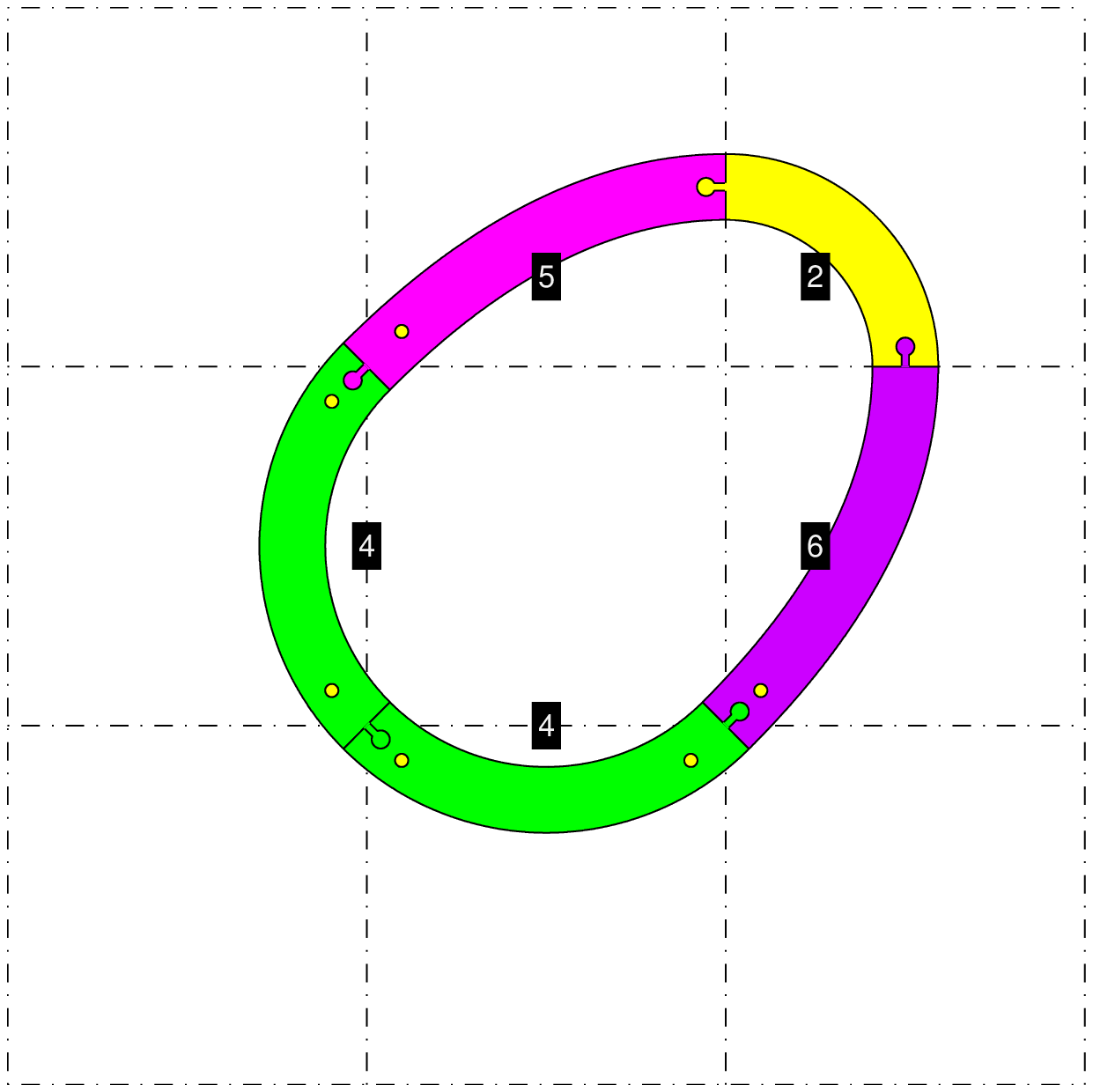, width=5 cm}}
\qquad
\caption{\label{circuits_complets_5_piecmax_ev0}\iflanguage{french}{Tous les 10 circuits retenus sur l'ensemble des 2000 circuits possibles}{All of the 10 circuits kept from the set of 2000 possible circuits}.}
\end{figure}
%%%%%%%%%%%%%%%%%%%%%%%%%%%%%%%%%%%%%%%%%%%%%%%%%%%%%%

\iflanguage{french}{%
Si on trace tous les  circuits réalisables %tous les circuits réalisables 
avec $N=5$ pièces
et $N_j=+\infty$, on obtient les  $10$ circuits de la figure \ref{circuits_complets_5_piecmax_ev0}.
Sur cette figure, on peut 
voir en fait se répéter deux circuits différents plusieurs fois. Dans chacun des deux ensembles de circuits,
on retrouve le même circuit à une isométrie directe près.
Deux circuits sont isométriques (cette isométrie étant directe) si et seulement si ils possèdent tous les deux 
les mêmes numéros signés (voir remarque \ref{remarkalpha}) de pièces, à une permutation circulaire près et au sens de parcours près.
Notons que la nombre total de circuits examinés est égal à $2000$.%
}{%
If we plot all the feasible circuits
 with $N=5$ pieces and $N_j=+\infty$, we obtain the $10$ circuits 
in Figure~\ref{circuits_complets_5_piecmax_ev0}. In this figure, one can in fact see two different circuits repeating several times. In each of the two sets of circuits, one finds the same circuit up to a direct isometry. Two circuits are isometric (this isometry being direct) if and only if they both possess the same signed number of pieces (see Remark~\ref{remarkalpha}), up to cyclic permutations 
%new 
and up to direction of travel. 
We note that the total number of circuits examined equals $2000$.%
}

%%%%%%%%%%%%%%%%%%%%%%%%%%%%%%%%%%%%%%%%%%%%%%%%%%%%%%
%\input{./simulations_circuit/circuit_numerique/circuits_complets_5_piecmax_ev1}
% fichier crée par 'presentation_exhaustif_circuit_boucle.m' le 22-Sep-2015 03:13:20
\begin{figure}[h]
\centering
%%% sous figure 1
\subfigure[\label{circuits_complets_5_piecmax_ev11}]
{\epsfig{file=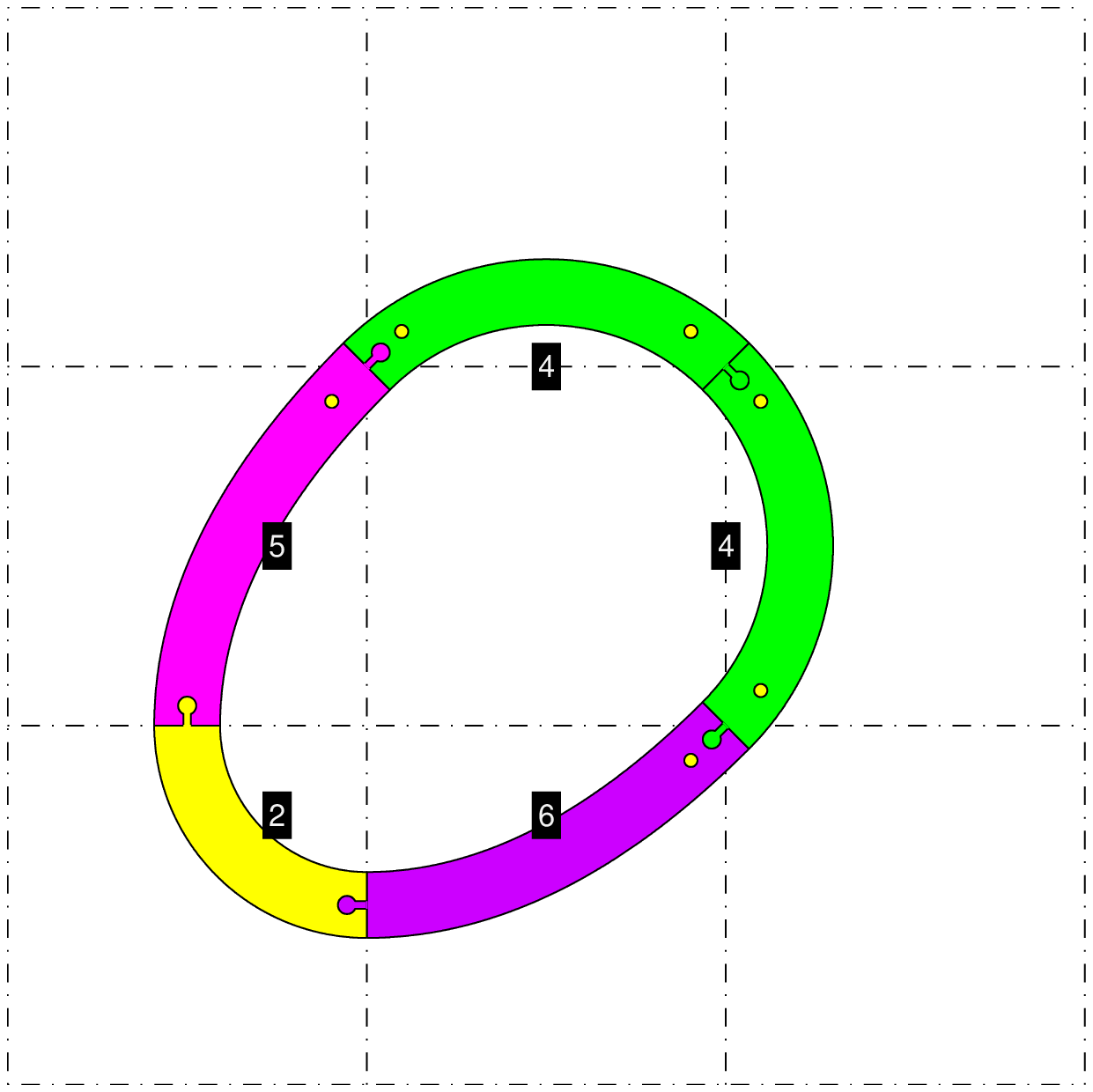, width=5 cm}}
\qquad
%%% sous figure 2
\subfigure[\label{circuits_complets_5_piecmax_ev12}]
{\epsfig{file=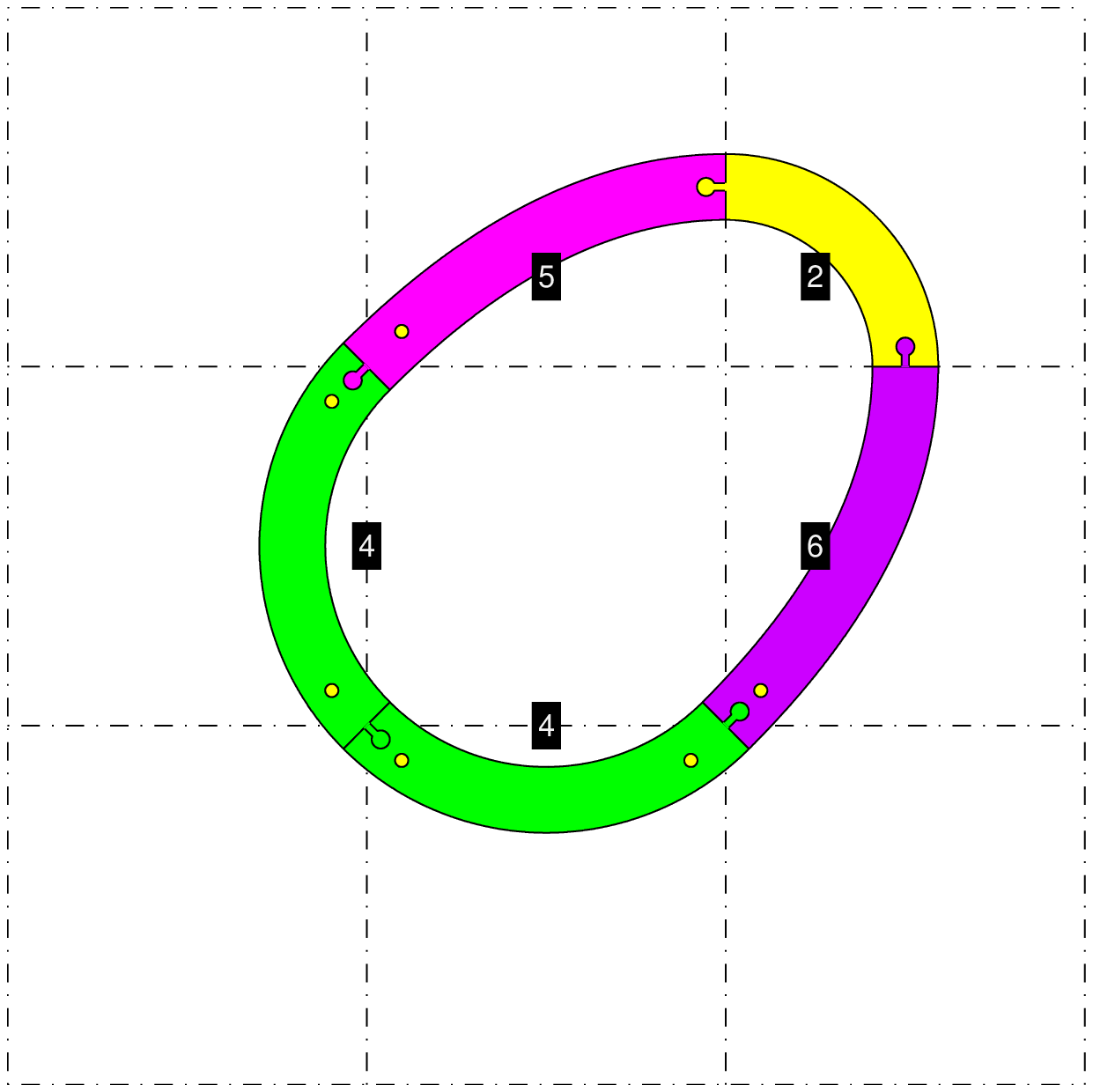, width=5 cm}}
\qquad
\caption{\label{circuits_complets_5_piecmax_ev1}\iflanguage{french}{Tous les 2 circuits retenus sur l'ensemble des 2000 circuits possibles}{All of the 2 circuits kept from the set of 2000 possible circuits}.}
\end{figure}
%%%%%%%%%%%%%%%%%%%%%%%%%%%%%%%%%%%%%%%%%%%%%%%%%%%%%%

\iflanguage{french}{%
Si on ne retient que ceux qui sont différents, à une isométrie directe près, on obtient
\ifcase 0
les  $2$ circuits 
\or
l'unique circuit 
\fi
de la figure \ref{circuits_complets_5_piecmax_ev1}.
Sur cette figure, on constate que le premier circuit est l'image du second par une isométrie indirecte.
Les numéros des pièces sont identiques, à une permutation circulaire près et au sens de parcours près ; 
de plus, pour prendre en compte cette isométrie indirecte,
il est nécessaire de remplacer les numéros  signés de pièces non rectilignes par leurs opposés.%
}{%
If we keep only those which are different, up to a direct isometry, we obtain 
\ifcase 0
the $2$ circuits 
\or
the unique circuit 
\fi
in Figure~\ref{circuits_complets_5_piecmax_ev1}. In this figure, one notices that the first circuit is the image of the second under an indirect isometry. The numbers of the pieces are identical, up to cyclic permutations
% new
and up to direction of travel.
Moreover, to take this indirect isometry into account, it is necessary to replace the signed numbers of curved pieces with their opposites.%
}

%%%%%%%%%%%%%%%%%%%%%%%%%%%%%%%%%%%%%%%%%%%%%%%%%%%%%%
%\input{./simulations_circuit/circuit_numerique/circuits_complets_5_piecmax_ev2}
% fichier crée par 'presentation_exhaustif_circuit_boucle.m' le 22-Sep-2015 03:13:21
\begin{figure}[h]
\begin{center}
\epsfig{file=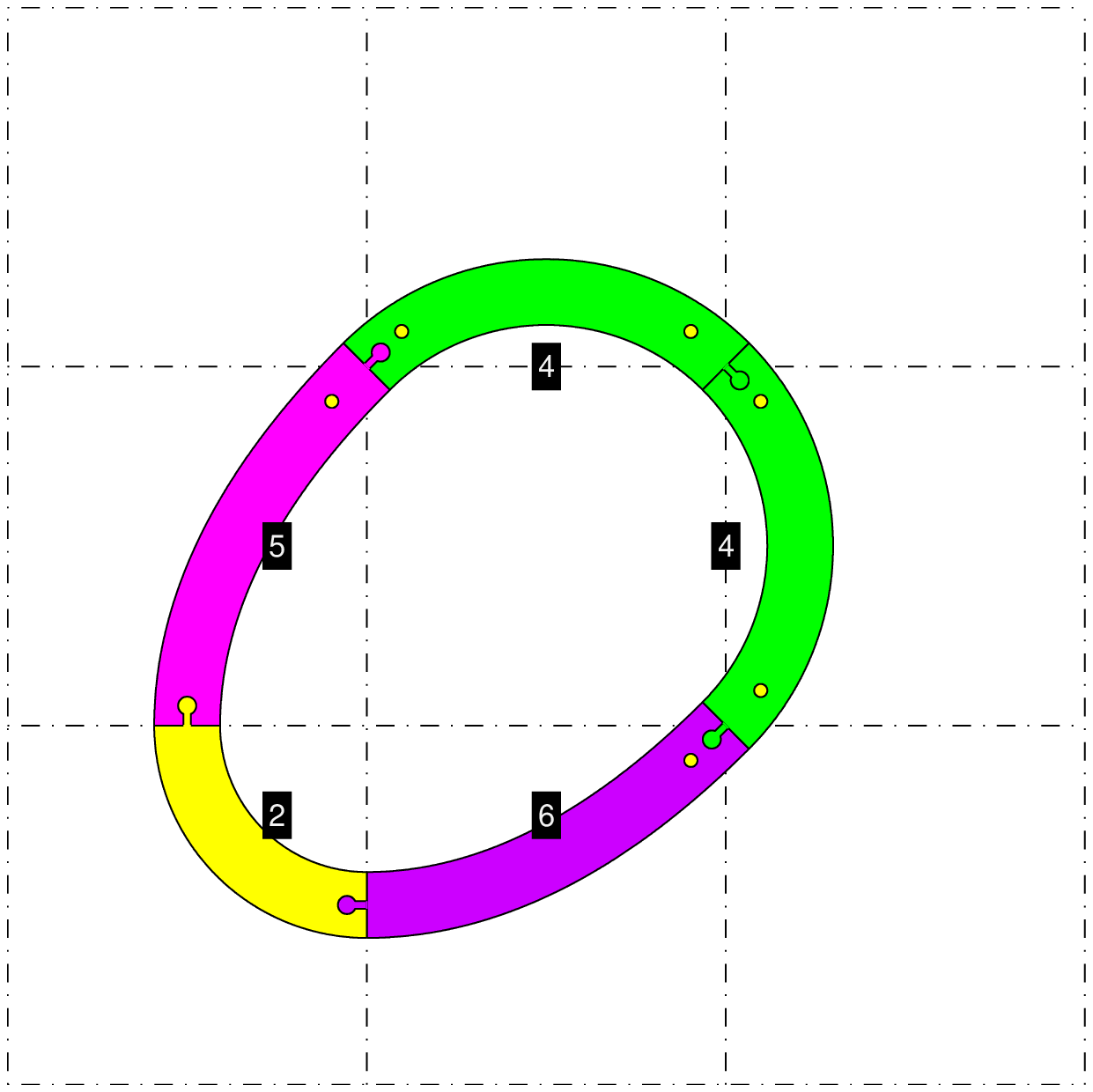, width=5 cm}
\end{center}
\caption{\label{circuits_complets_5_piecmax_ev2}\iflanguage{french}{Le seul circuit retenu sur l'ensemble des 2000 circuits possibles}{The sole circuit kept from the set of 2000 possible circuits}.}
\end{figure}
%%%%%%%%%%%%%%%%%%%%%%%%%%%%%%%%%%%%%%%%%%%%%%%%%%%%%%

\iflanguage{french}{%
Si  on ne retient  ceux qui sont
différents, à une isométrie près, on obtient
\ifcase 1
les  $1$ circuits 
\or
l'unique circuit 
\fi
de la figure \ref{circuits_complets_5_piecmax_ev2}.%
}{%
If we keep only those which are different up to an isometry, we obtain 
\ifcase 1
the $1$ circuits 
\or
the unique circuit 
\fi
in Figure~\ref{circuits_complets_5_piecmax_ev2}.%
}
%%%%%%%%%%%%%%%%%%%%%%%%%%%%%%%%%%%%%%%%%%%%%%%%%%%%%%
\end{example}

More generally, we draw all the circuits obtained for %a given $N$ and $N_j$.
$N$ and $N_j$ given .

A consideration of the direct isometries will result in the comparison of all of the obtained circuits. If two among them possess the same signed numbers of pieces, up to cyclic permutations, then one of the two will be eliminated.

Secondly, a consideration of the indirect isometries will be performed. Similarly, if two circuits possess the same signed numbers of pieces, but which are opposite (for the curved pieces), up to cyclic permutations, then one of the two will be eliminated. To take all of the indirect isometries into account, it will also be necessary to eliminate the circuits by also comparing the indices with permutations of the type $N$, $N-1$, ..., 2, 1, which amounts to considering the traversal of the circuit in the opposite direction. In this case, one will replace the signed number of pieces $\pm \piecec$ by $\mp \piececb$  
\ifcase \nopiecesix
(and $\pm \pieces$ by $\mp \piecesb$)
\or
\fi
and vice-versa. The two pieces \piecec\ and \piececb\ 
\ifcase \nopiecesix
(and $\pieces$ and $\piecesb$)
\or
\fi
are in effect identical; only the orientation changes. This elimination will be legitimate if the number of available pieces of types 
\piecec\ and  \piececb\  
\ifcase \nopiecesix
(and $\pieces$ and $\piecesb$) are identical
\or
are identical, , which will always be true in the following (see remark \ref{autantpiececinqsix}).
\fi

\begin{remark}
\label{autantpiececinqsix}
In the case of traditional self-avoiding polygons, the number of squares is necessarily even, which is no longer true here. Nonetheless, we can say that, in every circuit, 
\ifcase \nopiecesix
the sum of the numbers of pieces of type \piecec\ and \pieces\ equals the sum of the numbers of pieces of type \piececb\ and \piecesb.
\or
the numbers of pieces of types
\piecec\ and  \piececb\ 
are equal.
\fi
Indeed, note that, from Table \ref{typepiece},
\begin{itemize}
\item
the pieces of type
\pieceu, \pieced, \piecet\ and \pieceq\  connect together two points of the same nature (middles of sides or vertices of squares);
\item
\ifcase \nopiecesix
the pieces of type
\piecec\ and \pieces\  connect a middle to a vertex (in this order);
\or
the pieces of type
\piecec\ connect a middle to a vertex (in this order);
\fi
\item
\ifcase \nopiecesix
the pieces of type
\piececb\ and \piecesb\ connect a vertex to a middle (in this order);
\or
the pieces of type
\piececb\ connect a vertex to a middle (in this order);
\fi
\end{itemize}
A circuit departs from a point and returns to the same point. All of the points corresponding to the extremities of the pieces used in the circuit are either vertices or middles. It follows that that there are as many pieces connecting a middle to a vertex (in this order) as pieces connecting a vertex to a middle (in this order). Otherwise, the point of departure would not be of the same nature as the point of arrival.
\end{remark}%

\begin{remark}
\label{remtranslation}
Note that, in the enumeration of traditional self-avoiding polygons, only translations are taken into account in the isometries. Our enumeration problem is therefore quite different to the one in the literature. In the case where both notions coincide, this implies that the configurations that one will obtain will be \textit{a priori} less numerous than those in the literature.
\end{remark}%
}

\iflanguage{french}{%

%%%%%%%%%%%%%%%%%%%%%%%%%%%%%%%%%%%%%%%%%%%%%%%%%%%%%%%%%%%%%
\subsection{Prise en compte des contraintes locales de constructibilité}
\label{contrainteslocales}

Dans la recherche des chemins et des polygones autoévitants,  
une contrainte locale supplémentaire très importante est considérée :  
les carrés doivent être deux à deux distincts. Ici, cette contrainte n'est pas imposée, puisque seul compte le fait
de pouvoir réaliser des circuits  constructibles avec les pièces réelles.
Ces contraintes ne se manifestent pas sur l'exemple \ref{examplesimulation500}
puisque le petit nombre de pièces considéré ne fournit pas de circuits inconstructibles, mais ces contraintes interviendront dans de 
plus gros circuits, donnés en exemple plus bas.

Tout d'abord, il est nécessaire 
qu'hormis chacun des couples de pièces contiguës, qui ont donc une unique extrémité en commun, 
aucune des pièces n'ait d'extrémité en commun avec d'autres pièces qui ne soient pas continguës.
Si cette extrémité est un sommet, ce critère est facile à écrire et n'est pas détaillé ici.
Si cette extrémité  est un milieu et si deux pièces non contiguës ont cette extrémité en commun, elles appartiennent nécessairement au même carré, cas que l'on étudie maintenant.

Deux pièces peuvent appartenir au même carré à condition 
d'être disjointes (ou d'être à bords tangents).%

}%
{%

%%%%%%%%%%%%%%%%%%%%%%%%%%%%%%%%%%%%%%%%%%%%%%%%%%%%%%%%%%%%%
\subsection{Consideration of local constructibility constraints}
\label{contrainteslocales}

In the research into self-avoiding walks and polygons, 
%an additional local constraint 
a very important additional local constraint
is considered: the squares must be pairwise distinct. Here, this constraint isn't imposed, since only the fact of being able to produce circuits which are constructible with the real pieces counts. These constraints are not exhibited by Example~\ref{examplesimulation500}, since the small number of pieces considered does not provide unconstructible circuits, but these constraints come up in larger circuits, exemplified below.

First of all, it is necessary that, aside from each of the pairs of contiguous pieces, which therefore have a unique extremity in common, none of the pieces has an extremity in common with other pieces which aren't contiguous. If this extremity is a vertex, this criterion is easy to write, and is not detailed here. If this extremity is a middle, and if two non-contiguous pieces have this extremity in common, they necessarily belong to the same square, which is the case that we now study.

Two pieces may belong to the same square on the condition that they are disjoint (or have tangent edges).

}

\begin{figure}
\centering
%%% sous figure 1
\subfigure[\label{deux_pieces_meme_carrefigatheo}\iflanguage{french}{Cas où les deux lignes médianes sont nécessairement sécantes : $K_a=\{2\}$,
$K_b=\{0,4,5,6,7\}$}{Case where the two midlines are necessarily secant: $K_a=\{2\}$,
$K_b=\{0,4,5,6,7\}$}]
{\epsfig{file=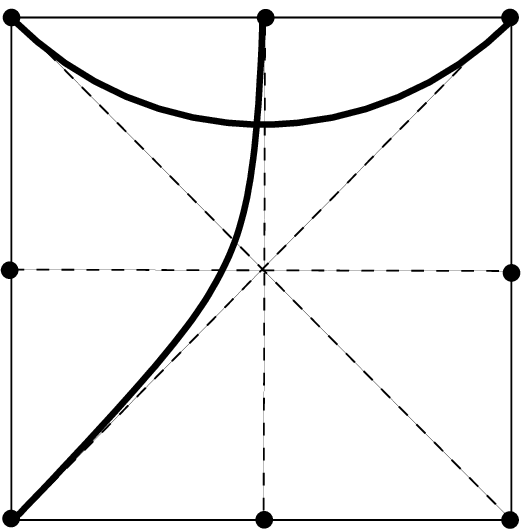, width=5cm}}
\qquad
%%% sous figure 2
\subfigure[\label{deux_pieces_meme_carrefigbtheo}\iflanguage{french}{Cas où les deux lignes médianes sont nécessairement disjointes}{Case where the two midlines are necessarily disjoint}]
{\epsfig{file=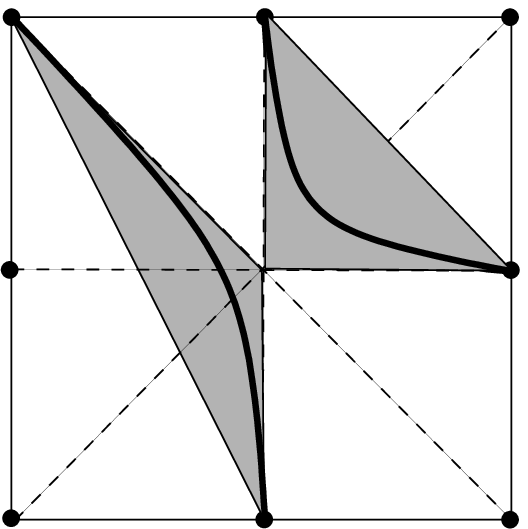, width=5cm}}
\caption{\label{deux_pieces_meme_carrefigtheo}\iflanguage{french}{Deux pièces dans un seul même carré.}{Two pieces within a single square.}}
\end{figure}

\ifcase \nopiecesix

\textbf{ATTENTION  !!!!!!!!!!!!!!!!!}
Début de raisonnement qui suit faux puisque \path|nopiece6|$=0$.

\or
\fi

\iflanguage{french}{%
On se donne donc deux pièces dans le même carré ${\mathcal{C}}_i$
dont on veut vérifier qu'elles sont disjointes.
Plusieurs cas se présentent : 
\begin{enumerate}
\item
\label{contloccasA}
Elles ont en commun au moins une extrémité : dans ce cas, elles ne sont pas disjointes.
\item
\label{contloccasB}
Leurs extrémités sont deux à deux distinctes.
\begin{enumerate}
\item
\label{contloccasBa}
Notons $P_1=A_i$ et $P_2=B_i$ (resp. $P'_1$ et $P'_2$)
les extrémités de la première (resp. seconde) pièce.
Les deux points  $P_1$ et $P_2$ appartiennent à la frontière $\partial  {\mathcal{C}}_i$ du carré ; ils définissent donc 
deux composantes connexes (pour la topologie induite) notée 
${\mathcal{P}}_a$ et ${\mathcal{P}}_b$. Si 
\begin{equation}
\label{ligncoupeante}
\text{($P'_1\in {\mathcal{P}}_a$ et $P'_2\in {\mathcal{P}}_b$)
ou 
($P'_2\in {\mathcal{P}}_a$ et $P'_1\in {\mathcal{P}}_b$)},
\end{equation}
alors les deux sommets de la seconde pièce sont de part et d'autre de la 
ligne médiane de la première pièce et par continuité des lignes médianes, elles  ont un point en commun et,
nécessairement, dans ce cas, les deux pièces  ne sont pas disjointes.
Réciproquement, si \eqref{ligncoupeante} n'a pas lieu, on peut montrer que les lignes médianes sont nécessairement disjointes
(voir cas \ref{contloccasBb}).

Précisons tout cela :
repérons chacune des extrémités $P_1$ et  $P_2$  de la première pièce par deux entiers $\kappa_1$ et $\kappa_2$  appartenant à $\{0,...,7\}$
de la façon suivante : si $\vec I$ désigne le premier vecteur de la base orthonormée choisie et $c_i$ est le centre du carré 
alors 
\begin{equation*}
\widehat{\left(\vec I,\overrightarrow{c_iP_1} \right)}=\frac{\kappa_1\pi}{4},\quad
\widehat{\left(\vec I,\overrightarrow{c_iP_2} \right)}=\frac{\kappa_2\pi}{4}.
\end{equation*}
On définit de même $\kappa'_1$ et $\kappa'_2$ pour la seconde pièce.
L'ensemble $\{0,...,7\}$
peut être partitionné de la façon suivante : $\{0,...,7\}=K_a\cup K_b\cup \{\kappa_1\}\cup \{\kappa_2\}$, de telle sorte que 
tous les sommets correspondant aux indices de $K_a$ (resp. $K_b$) soient successifs dans le carré.
La propriété \eqref{ligncoupeante}
est équivalente à 
\begin{equation}
\label{ligncoupe}
\text{($\kappa'_1\in K_a$ et $\kappa'_2\in K_b$)
ou ($\kappa'_2\in K_a$ et $\kappa'_1\in K_b$).}
\end{equation}
Si elle a lieu,  les deux pièces  ne sont donc pas disjointes.
Voir par exemple les figure \ref{deux_pieces_meme_carrefigatheo} et \ref{deux_pieces_meme_carrefiga} qui illustrent ce cas.
\item
\label{contloccasBb}
Supposons maintenant que \eqref{ligncoupe}
n'ait pas lieu. Dans ce cas, $\kappa'_1$ et $\kappa'_2$ appartiennent tous les deux soit à $K_a$ soit à $K_b$.
Si aucune des pièces n'est rectiligne, le cardinal de $K_a$ et $K_b$
est nécessairement dans $\{1,2\}$ ou dans $\{4,5\}$
Ainsi,
$\kappa'_1$ et $\kappa'_2$ ne peuvent appartenir 
au plus petit ensemble et on est nécessairement dans le cas où l'ensemble des sommets compris entre les deux extrémités
de la première pièce et l'ensemble des sommets compris entre les deux extrémités
de la seconde pièce sont disjoints.
Voir par exemple la figure \ref{deux_pieces_meme_carrefigbtheo}  qui illustre ce cas.
Une courbe de Bézier est nécessairement incluse dans le polygone de contrôle $A_ic_iB_i$.  Il en est de même pour les lignes
médianes des pièces circulaires. Ainsi, dans ce cas, chacune des courbes médiane est incluse dans deux triangles 
qui n'ont en commun que le centre $c_i$ du carré, par lequel ne passe aucune des lignes médiane (puisque les lignes médianes
rectilignes ne sont pas prises en compte).
Ce raisonnement est encore valable si l'une des lignes médianes est rectiligne.
Bref, dans ce cas, les lignes médianes des deux pièces 
sont disjointes.
\end{enumerate}
\end{enumerate}
\ifcase \nopiecesix
\textbf{ATTENTION  !!!!!!!!!!!!!!!!!}
Fin de raisonnement qui suit faux puisque \path|nopiece6|$=0$.
\or
\fi
Il ressort de tout cela que, si  les deux pièces étudiées sont confondues,  on est dans le cas \ref{contloccasA}.
Sinon, on est soit dans les cas \ref{contloccasA} ou \ref{contloccasBa}, auxquels cas les pièces ne sont pas disjointes.
Enfin, si on est dans le dernier cas \ref{contloccasBb}, les lignes médianes des deux pièces 
sont disjointes.
Dans ce cas, si on note  ${{\Gamma}}_i$ et ${{\Gamma}}'_i$
les deux courbes (qui sont dans le même carré), on peut considérer
\begin{equation}
\label{dmin}
\delta =\sqrt{\inf_{(M,M') \in {{\Gamma}}_i\times \in {{\Gamma}}'_i}
d^2(M,M')},
\end{equation}
où $d(M,M')$ est la distance euclidienne entre les points $M$ et $M'$.
Le couple $(M,M')$ décrit une partie compacte de $\Er^2$ et $d$ est continue ; cette borne inférieure
existe et est nécessairement atteinte en un couple de points $(M_0,M'_0)$ de ${{\Gamma}}_i\times {{\Gamma}}'_i$. 
Puisque les deux courbes ${{\Gamma}}_i$ et ${{\Gamma}}'_i$ sont disjointes, le nombre $\delta$ est nécessairement 
strictement positif. 
Par ailleurs, on peut montrer que pour tous les couples de courbes ${{\Gamma}}_i\times {{\Gamma}}'_i$
qui sont dans ce cas, ce réel $\delta$ est nécessairement atteint en un couple de points qui ne peut être au bord du carré. 
Dans ce cas, puisque $(M,M')\mapsto d^2(M,M')$ est une fonction différentiable, sa différentielle  y est nulle, ce qui se traduit
par la perpendicularité de la droite $(M_0M'_0)$ avec la tangente à la courbe ${{\Gamma}}_i$ (resp. ${{\Gamma}}'_i$)
au point $M_0$ (resp. $M'_0$).%
}%
{%
We therefore take two pieces in the same square ${\mathcal{C}}_i$, and we wish to verify that they are disjoint. Several cases present themselves:
\begin{enumerate}
\item
\label{contloccasA}
They have at least one extremity in common. In this case, they are not disjoint.
\item
\label{contloccasB}
Their extremities are pairwise distinct.
\begin{enumerate}
\item
\label{contloccasBa}
We denote by $P_1=A_i$ and $P_2=B_i$ (respectively $P'_1$ and $P'_2$) the extremities of the first (respectively second) piece. The two points $P_1$ and $P_2$ belong to the border $\partial  {\mathcal{C}}_i$ of the square. They therefore define two connected components (for the induced topology) denoted ${\mathcal{P}}_a$ and ${\mathcal{P}}_b$. If
\begin{equation}
\label{ligncoupeante}
\text{($P'_1\in {\mathcal{P}}_a$ and $P'_2\in {\mathcal{P}}_b$)
or 
($P'_2\in {\mathcal{P}}_a$ and $P'_1\in {\mathcal{P}}_b$)},
\end{equation}
then the two vertices of the second piece are on both sides of the midline of the first piece, and by continuity of the midlines, they have a point in common and, necessarily, in this case, the two pieces are not disjoint. Conversely, if \eqref{ligncoupeante} does not hold, one can show that the midlines are necessarily disjoint (see case~\ref{contloccasBb}).

Let us clarify this. We describe the location of the extremities $P_1$ and  $P_2$ of the first piece by two integers, $\kappa_1$ and $\kappa_2$, in $\{0,...,7\}$ in the following way: if $\vec I$ designates the first vector of the chosen orthonormal basis, and $c_i$ is the center of the square, then
\begin{equation*}
\widehat{\left(\vec I,\overrightarrow{c_iP_1} \right)}=\frac{\kappa_1\pi}{4},\quad
\widehat{\left(\vec I,\overrightarrow{c_iP_2} \right)}=\frac{\kappa_2\pi}{4}.
\end{equation*}
We similarly define $\kappa'_1$ and $\kappa'_2$ for the second piece. The set $\{0,...,7\}$ can be partitioned in the following way: $\{0,...,7\}=K_a\cup K_b\cup \{\kappa_1\}\cup \{\kappa_2\}$, such that all of the vertices corresponding to the indices of $K_a$ (respectively $K_b$) are consecutive in the square.
Property \eqref{ligncoupeante} is equivalent to
\begin{equation}
\label{ligncoupe}
\text{($\kappa'_1\in K_a$ and $\kappa'_2\in K_b$)
or ($\kappa'_2\in K_a$ and $\kappa'_1\in K_b$).}
\end{equation}
If this holds, the two pieces are not disjoint. See, for example, Figures~\ref{deux_pieces_meme_carrefigatheo} and~\ref{deux_pieces_meme_carrefiga}, which illustrate this case.
\item
\label{contloccasBb}
Let us now assume that \eqref{ligncoupe} does not hold. In this case, $\kappa'_1$ and $\kappa'_2$ both belong to either $K_a$ or to $K_b$. If none of the pieces is straight, the cardinality of $K_a$ and $K_b$ is necessarily in $\{1,2\}$ or in $\{4,5\}$.
Thus,
$\kappa'_1$ and $\kappa'_2$ cannot belong to the smaller set, and we are necessarily in the case where the set of vertices contained between the two extremities of the first piece and the set of vertices contained between the two extremities of the second piece are disjoint.
See for example Figure~\ref{deux_pieces_meme_carrefigbtheo}, which illustrates this case. A Bézier curve is necessarily included in the control polygon $A_ic_iB_i$.   It is the same for the midlines of the circular pieces. Thus, in this case, each of the middle curves is included in two triangles which only have the center $c_i$ of the square in common, through which none of the midlines pass (since the straight midlines are not considered). This reasoning is still valid if one of the midlines is straight. In short, in this case, the midlines of the two pieces are disjoint.
\ifcase \nopiecesix
\textbf{ATTENTION  !!!!!!!!!!!!!!!!!}
Fin de raisonnement qui suit faux puisque \path|nopiece6|$=0$.
\or
\fi
\end{enumerate}
\end{enumerate}
It thus follows that, if the two pieces under consideration are coincident, we are in case~\ref{contloccasA}. Otherwise we are either in case~\ref{contloccasA} or~\ref{contloccasBa}, in which the pieces are not disjoint. Finally, if we are in the last case, \ref{contloccasBb}, the midlines of the two pieces are disjoint.
In this case, if we write ${{\Gamma}}_i$ and ${{\Gamma}}'_i$ for the two curves (which are in the same square), one may consider
\begin{equation}
\label{dmin}
\delta =\sqrt{\inf_{(M,M') \in {{\Gamma}}_i\times \in {{\Gamma}}'_i}
d^2(M,M')},
\end{equation}
where $d(M,M')$ is the Euclidean distance between the points $M$ and $M'$. The pair $(M,M')$ describes a compact subset of $\Er^2$ and $d$ is continuous. This lower bound exists and is necessarily attained at a pair of points $(M_0,M'_0)$ of ${{\Gamma}}_i\times {{\Gamma}}'_i$. Since the two curves ${{\Gamma}}_i$ and ${{\Gamma}}'_i$ are disjoint, the number $\delta$ is necessarily strictly positive. In addition, one can show that for every pair of curves ${{\Gamma}}_i\times {{\Gamma}}'_i$ which fall into in this case, the real number $\delta$ is necessarily attained at a pair of points which cannot be at the edge of the square. In this case, since $(M,M')\mapsto d^2(M,M')$ is a differentiable function, its %derivative 
differential
is zero there, which results in the perpendicularity of the straight line $(M_0M'_0)$ to the tangent to the curve ${{\Gamma}}_i$ (respectively ${{\Gamma}}'_i$) at the point $M_0$ (respectively $M'_0$).% 
}
\ifcase \nopiecesix

\textbf{ATTENTION  !!!!!!!!!!!!!!!!!}
Début de raisonnement à reprendre numériquement puisque \path|nopiece6|$=0$.

\or
\fi
%%%%%%%%%%%%%%%%%%%%%%%%%%%%%%%%%%%%%%%%%%%%%%
%\input{deux_pieces_meme_carre}
% fichier tex crée par MaTeXBuild02 (manuellement) le 21-Sep-2015 8:35:11
% à compilier avec 
% MaTeXBuild02('deux_pieces_meme_carre',0)
% après enumeration_circuit
% Compilation longue

\begin{figure}
\psfrag{se coupent (test indice)}{}
\psfrag{ne se coupe pas}{}
\centering
%%% sous figure 1
\subfigure[\label{deux_pieces_meme_carrefiga}\iflanguage{french}{Cas où les deux lignes médianes sont nécessairement sécantes}{Case where the two midlines are necessarily secant}]
{\epsfig{file=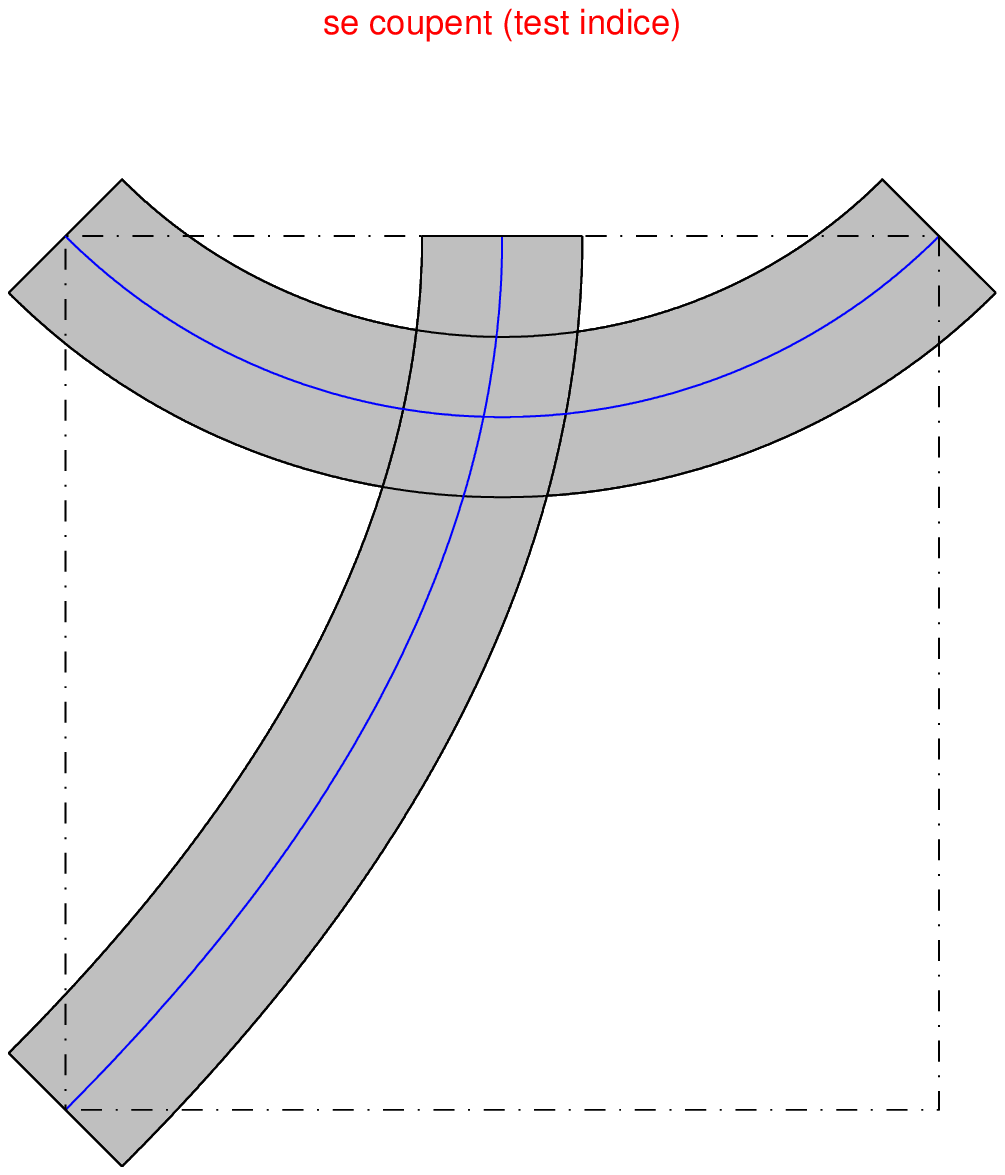, width=5cm}}
\qquad
%%% sous figure 2
\subfigure[\label{deux_pieces_meme_carrefigb}\iflanguage{french}{Configuration où la distance entre les deux lignes médianes est la plus petite}{Configuration where the distance between the two midlines is the smallest}]
{\epsfig{file=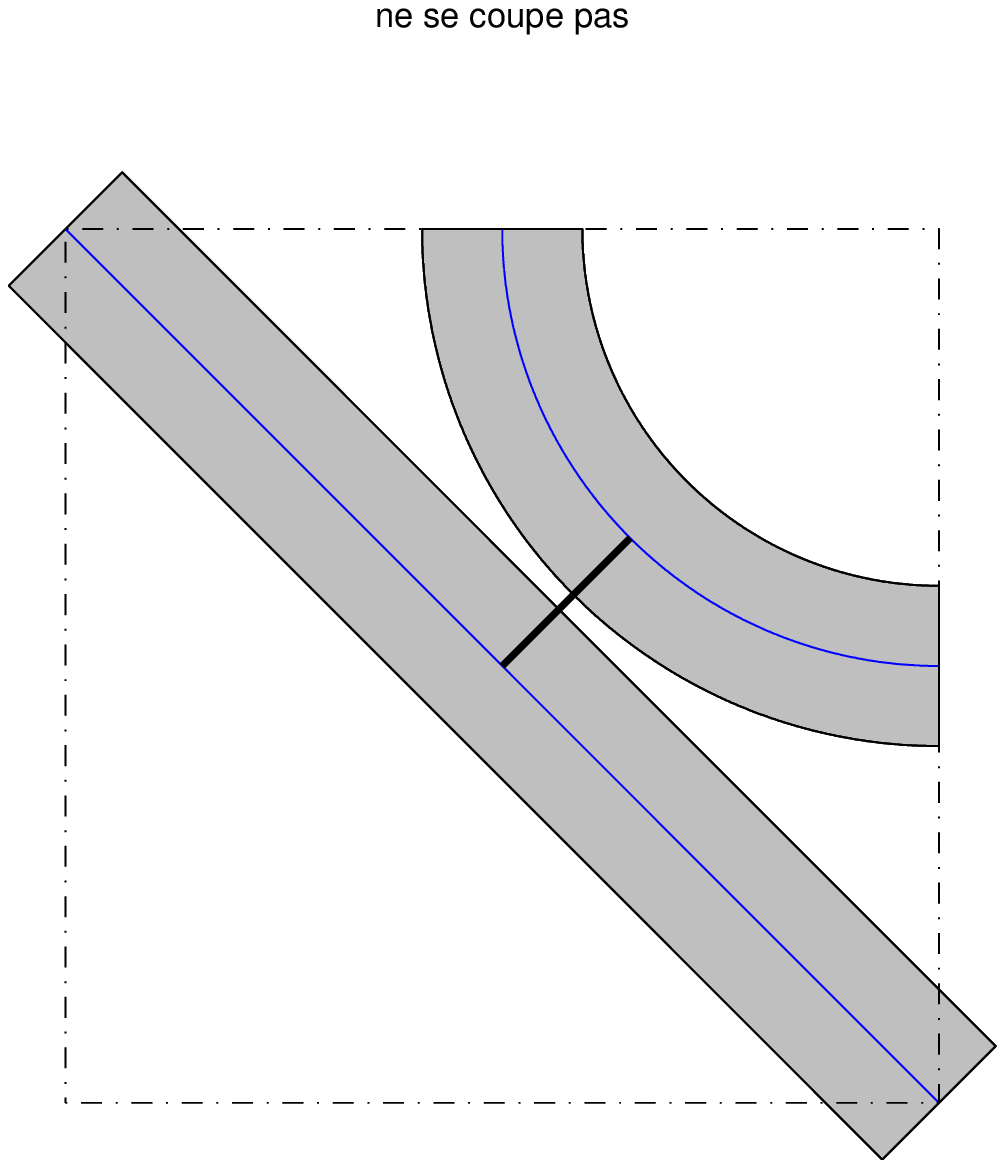, width=5cm}}
\caption{\label{deux_pieces_meme_carrefig}\iflanguage{french}{Deux pièces dans un seul même carré.}{Two pieces within the same square.}}
\end{figure}

% Commandes locales utilisées pour le commentaire
\newcommand{\dpmcchomatot}{1}
\newcommand{\dpmcdmin}{0.20711}
\newcommand{\dpmceccart}{0.51493}
\newcommand{\dpmcminint}{1}
\newcommand{\dpmcntot}{1600}
\newcommand{\dpmcnoerreur}{1}
\newcommand{\dpmcdminexact}{\frac{1}{2}\left(\sqrt 2 -1\right)}
%%%%%%%%%%%%%%%%%%%%%%%%%%%%%%%%%%%%%%%%%%%%%%
\iflanguage{french}{%
Informatiquement, toutes ces propriété ont été vérifiées
en faisant un balayage sur tous les couples possibles  de courbes  ${{\Gamma}}_i\times {{\Gamma}}'_i$,
ce qui représente $\dpmcntot$  cas à étudier.
Nous avons déterminé le couple de courbes qui correspond à la plus petite distance $\delta$ possible, donnée
par (dans le cas d'un carré de longueur unité) : 
\begin{equation}
\label{dminmin}
\delta_{\min} =\dpmcdmin,
\end{equation}
qui correspond à la configuration de la figure \ref{deux_pieces_meme_carrefigb}.
Cette expression est égale à 
\begin{equation}
\label{dminminexa}
\delta_{\min} =\dpmcdminexact.
\end{equation}
Par construction des bords et par les propriétés de perpendicularité vues, à distance constante égale à la demi-largeur $e/2$ du 
rail, deux rails seront disjoints si et seulement si 
\begin{equation}
\label{dminminexateste}
\delta_{\min}\geq e.
\end{equation}
Le cas $\delta_{\min}=e$ correspond au cas où les bords des deux rails sont tangents et est encore acceptable\footnote{Si la conception des pièces est parfaite,
car aucune imperfection n'est dans ce cas autorisée. En pratique, on préfèrera donc choisir par sécurité $\delta_{\min}> e$.}.
Le choix d'une section standard, compatible avec les véhicules miniatures de type \text{Brio} \textregistered, 
correspond à $e$ donné par \eqref{dminminexaeval}.
%\begin{equation}
%\label{dminminexaeval}
%e=ateval{affiche_nombre02(2*(e+lambda+mu),3,6,5)}.
%0.18349.
%\end{equation}
On vérifie donc que \eqref{dminminexateste} a lieu, ce qui signifie, qu'avec le choix de section fait,
tous les couples de courbes qui ne sont pas dans le cas où les lignes médianes se coupent nécessairement
(cas \ref{contloccasA} ou \ref{contloccasBa}) donnent des situations où les deux pièces sont disjointes
(comme dans le cas indiqué par la figure \ref{deux_pieces_meme_carrefigb}).
Notons que dans le cas où la plus petite distance est atteinte, la distance la plus petite entre les deux bords est donnée par 
$\delta_{\min}-e\geq 0$, ce qui numériquement donne (en multipliant par la longueur de référence donnée par \eqref{longueurl}) :
\begin{equation*}
\xi=
(\delta_{\min}-e)L=\dpmceccart \text{ cm},
\end{equation*}
ce qui est très peu, finalement, eu égard à la valeur donnée par \eqref{longueurl} !%
}%
{%
Computationally, all of these properties have been verified by conducting a sweep of all the possible pairs of curves ${{\Gamma}}_i\times {{\Gamma}}'_i$, which represents $\dpmcntot$ cases to study. We have determined the pair of curves which correspond to the smallest possible distance $\delta$, given by (in the case of a square of unit length): 
\begin{equation}
\label{dminmin}
\delta_{\min} =\dpmcdmin,
\end{equation}
which corresponds to the configuration in Figure~\ref{deux_pieces_meme_carrefigb}. This expression equals
\begin{equation}
\label{dminminexa}
\delta_{\min} =\dpmcdminexact.
\end{equation}
By construction of the edges, and by the perpendicularity properties seen above, at a constant distance equal to the half-width $e/2$ of the piece, two pieces will be disjoint if and only if
\begin{equation}
\label{dminminexateste}
\delta_{\min}\geq e.
\end{equation}
The case $\delta_{\min}=e$ corresponds to the case where the edges of both rails are tangent, and is still acceptable\footnote{If the design of the pieces is perfect, as no imperfection is permitted. In practice therefore, to be safe, one will prefer to choose $\delta_{\min}> e$.}.
The choice of a standard cross-section, compatible with the \text{Brio}\textregistered-type miniature vehicles, corresponds to 
$e$ given by \eqref{dminminexaeval}.
%\begin{equation}
%\label{dminminexaeval}
%e=ateval{affiche_nombre02(2*(e+lambda+mu),3,6,5)}.
%0.18349.
%\end{equation}
We therefore verify that~\eqref{dminminexateste} holds, which means that with the chosen cross-section, every pair of curves which are not in the case where the midlines necessarily cut across each other (case \ref{contloccasA} or \ref{contloccasBa}), gives rise to a situation where the two pieces are disjoint, as in the case indicated by Figure~\ref{deux_pieces_meme_carrefigb}). We note that in the case where the smallest distance is attained, the smallest distance between the two edges is given by
$\delta_{\min}-e\geq 0$, which (multiplying by the reference length given by~\eqref{longueurl}) numerically gives:
\begin{equation*}
\xi=
(\delta_{\min}-e)L=\dpmceccart \text{ cm},
\end{equation*}
which is very small in the end, with respect to the value given by~\eqref{longueurl}! %
}

\ifcase \dpmcchomatot

\textbf{ATTENTION  !!!!!!!!!!!!!!!!!}
Le programme a detecté de cas foireux, à reprendre.

\or
\fi

\ifcase \dpmcminint

\textbf{ATTENTION  !!!!!!!!!!!!!!!!!}
Le programme a detecté de cas où le minimum est atteint au bord du carré, à reprendre.

\or
\fi

\ifcase\dpmcnoerreur

\textbf{ATTENTION  !!!!!!!!!!!!!!!!!}
La configuration minimale ne correspond pas à $\delta_{\min} =\dpmcdminexact$. 
\`A reprendre. 

\or
\fi

\iflanguage{french}{%
Notons aussi que l'on pourrait procéder, par récursivité par exemple, à l'étude complète où plusieurs pièces 
appartiennent au même carré, ce qui a été évité ici, puisque ce cas est très rare.
\ifcase \nopiecesix

\textbf{ATTENTION  !!!!!!!!!!!!!!!!!}
Fin de raisonnement à reprendre numériquement puisque \path|nopiece6|$=0$.

\or
\fi

Dans le brevet 
\ifcase \choano
\cite{brevetJB,pctJB},
\or
\cite{brevetJBano,pctJBano},
\fi
la prise en compte d'aiguillages, de ponts, de croisements a été prévue ; il suffit que ces éléments
soient aussi inclus dans des carrés et vérifient les principes de construction. Dans un premier temps,
seules des pièces planes et simples (sans aiguillages ni croisements) ont été réalisées. 
Pour des circuits comprenant d'autres éléments que ceux-là, cette étude de contraintes locales serait donc à reconsidérer.
Voir section \ref{geneaiguillage}.

Concluons par deux exemples montrant des circuits inconstructibles ou constructibles, mis en évidence informatiquement.%
}{%
We also note that one could proceed, by recursion for example, to the complete study where several pieces belong to the same square, which was avoided here, since this case is very rare.
\ifcase \nopiecesix

\textbf{ATTENTION  !!!!!!!!!!!!!!!!!}
Fin de raisonnement à reprendre numériquement puisque \path|nopiece6|$=0$.

\or
\fi

In the patent 
\ifcase \choano
\cite{brevetJB,pctJB},
\or
\cite{brevetJBano,pctJBano},
\fi
the inclusion of switches, bridges, and crossings was considered; it suffices that these elements are also included in squares and satisfy the principles of construction. In the first instance, only planar and simple (without switches or crossings) pieces have been realised. For circuits including other element than these, this study of local constraints will therefore need to be reconsidered. See Section~\ref{geneaiguillage}.

We conclude by two examples showing unconstructable and constructible circuits, demonstrated computationally.%
}

\begin{example}
\label{examplesimulation410}
%%%%%%%%%%%%%%%%%%%%%%%%%%%%%%%%%%%%%%%%%%%%%%%%%%%%%%%
%\input{simulations_circuit/simulation410}
% fichier tex crée par MaTeXBuild02 le 03-Sep-2015 06:04:09
% A compiler avec MaTeXBuild02('simulation410',0)
% après le fichier 'enumeration_construction_circuit.matex'
% Sortie statique, pas très pertinent avec matex !!

\begin{figure}
\centering
%%% sous figure 1
\subfigure[\label{figexamplesimulation410a}\iflanguage{french}{Un circuit non constructible}{A non-constructible circuit}]
{\epsfig{file=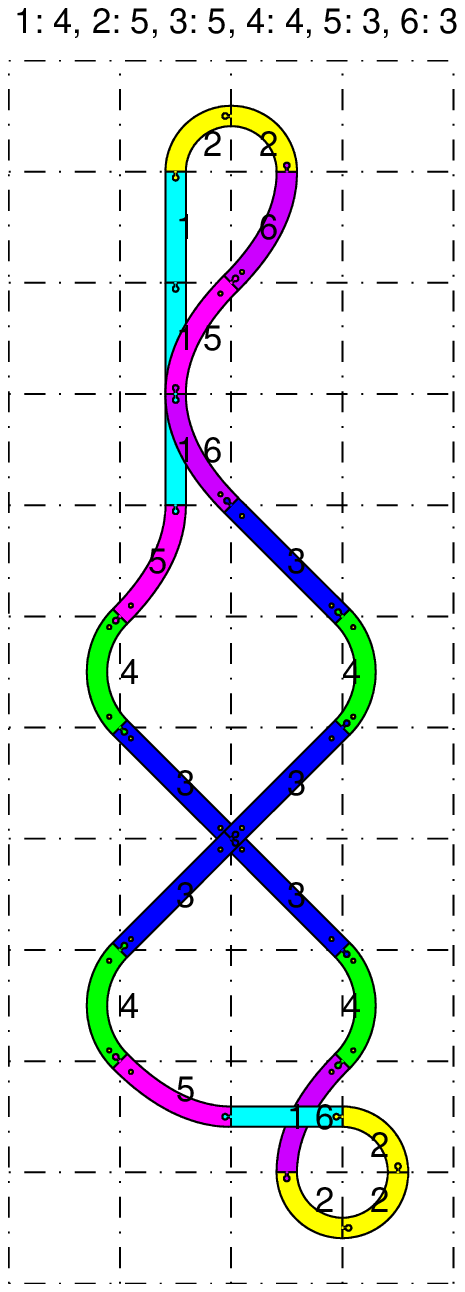, width=4 cm}}
\qquad
%%% sous figure 2
\subfigure[\label{figexamplesimulation410b}\iflanguage{french}{Un circuit constructible}{A constructible circuit}]
{\epsfig{file=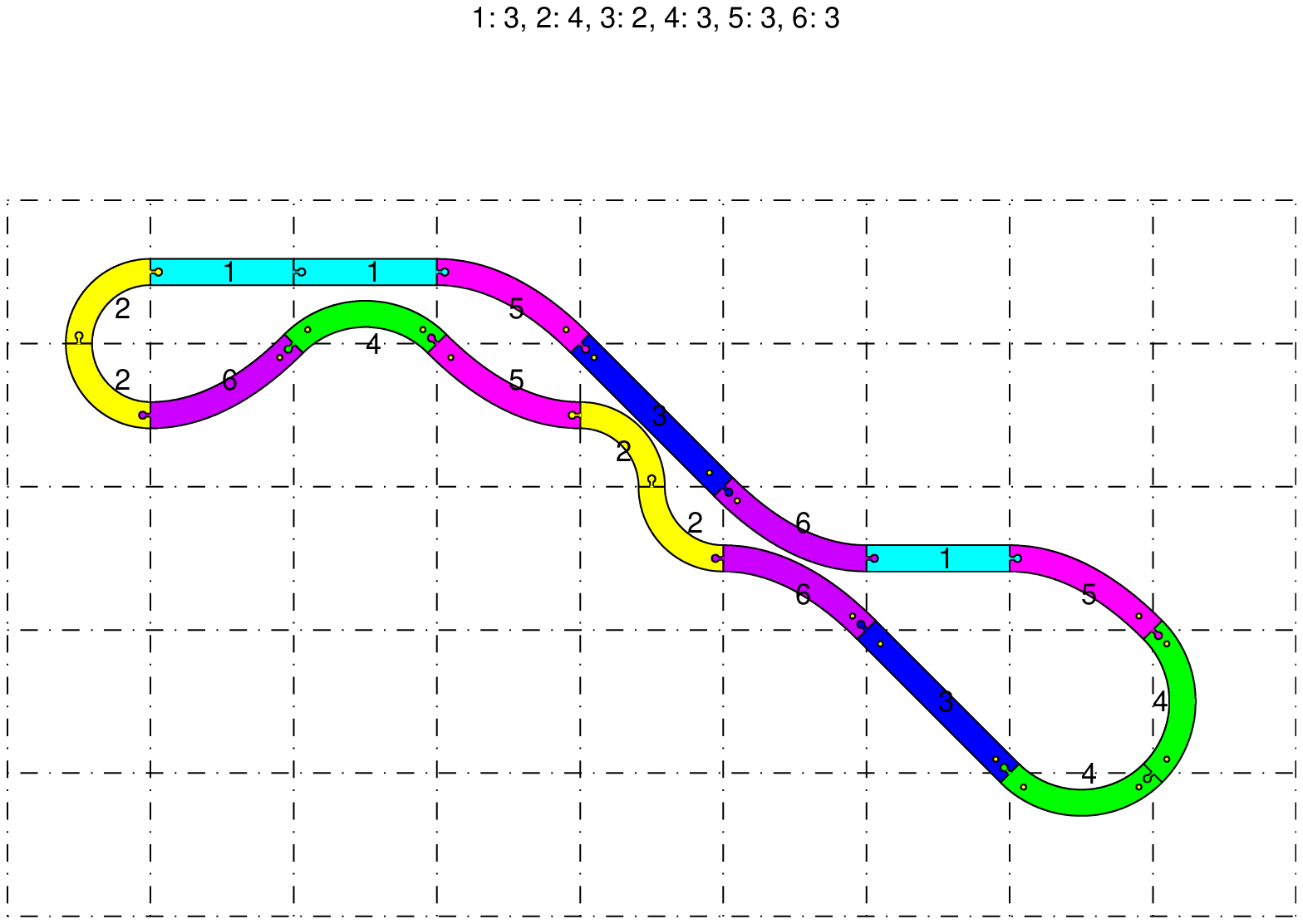, width=12 cm}}
\caption{\label{figexamplesimulation410}\iflanguage{french}{Deux circuits.}{Two circuits.}}
\end{figure}

\iflanguage{french}{Sur la figure \ref{figexamplesimulation410a}, on a choisi trois exemples de pièces non disjointes, couvrant les trois cas vus. 
Au contraire, en figure \ref{figexamplesimulation410b}, le circuit est constructible avec trois carrés dans lesquels apparaissent 
à chaque fois deux pièces disjointes.}{In Figure~\ref{figexamplesimulation410a}, we chose three examples of non-disjoint pieces, covering the three cases seen above. In contrast, in Figure \ref{figexamplesimulation410b}, the circuit is constructible with three squares in which two disjoint pieces appear each time.}
%%%%%%%%%%%%%%%%%%%%%%%%%%%%%%%%%%%%%%%%%%%%%%%%%%%%%%%
\end{example}

\iflanguage{french}{%

%%%%%%%%%%%%%%%%%%%%%%%%%%%%%%%%%%%%%%%%%%%%%%%%%%%%%%%%%%%%%
\subsection{Définition adoptée}
\label{defslw}

\begin{definition}[\slw] 
On appellera un \slw, un chemin $\Gamma$ de classe $\mathcal{C}^1$ de $\Er^2$, 
défini comme en section \ref{principecourbe}.
Ce chemin est défini 
 à une isométrie près, à une permutation circulaire près et 
à un sens de parcours  près.
Par abus de langage, nous appellerons aussi \slw, ce qui permet de définir le chemin $\Gamma$, c'est-à-dire, respectivement : 
\begin{itemize}
\item
la suite des centres des carrés  ${(c_i)}_{1\leq i \leq N}$ occupés par le chemin $\Gamma$,
vérifiant toutes les contraintes données dans 
les  sections 
\ref{enumerationtouscircuit},
\ref{isometrie} et 
\ref{contrainteslocales}, cette suite étant définie à une isométrie près, à une permutation circulaire près et 
à un sens de parcours  près ; 
\item
la suite des angles ${(\alpha_i)}_{1\leq i \leq N}$ définis en sections
\ref{enumerationpositionprobleme}
 et 
\ref{enumerationtouscircuit}, cette suite  
étant définie, à une multiplication par $-1$ près, à une permutation circulaire près et 
à un sens de parcours  près ; 
\item
la suite des numéros signés de pièces ${(p_i)}_{1\leq i \leq N}$ définis en section 
\ref{enumerationpositionprobleme}, cette suite  
étant définie 
au changement  
$\pm \piecec$ par $\mp \piececb$ 
\ifcase \nopiecesix
(et $\pm \pieces$ par $\mp \piecesb$)
\or
\fi
et vice-versa près,
à une permutation circulaire près et 
à un sens de parcours  près. 
\end{itemize}
\end{definition}%

Notons aussi qu'un \slw\ dépend  des $N_j$ et de la largeur $e$ qui, dans tout cet article est choisie  
inférieure  à la largeur critique 
égale au maximum de $e_0$ définie par 
\ifcase \nopiecesix
\eqref{maxlar6}
\or
\eqref{maxlar5}
\fi
et de 
$\delta_{\min}$
défini par 
\eqref{dminmin} et 
\eqref{dminminexa}.
Dans le cas légèrement différent où $e>e_0$, $e$ doit rester inférieur à $1/2$ ; dans ce dernier cas, 
les règles de connexion, donnée en section \ref{contrainteslocales}, sont à modifier, ce qui est aussi pris en compte dans les algorithmes utilisés.
\ifcase \nopiecesix
\or
Si on décide d'adjoindre les pièces de type \pieces\ et \piecesb, $e_0$ est alors donné par \eqref{maxlar6}.
\fi.%

\begin{definition}[circuit]
Dans toute la suite de cet article, un circuit est donc l'ensemble des pièces géométriques 
qui s'appuie sur la construction géométrique d'un \slw.
\end{definition}
}{%

%%%%%%%%%%%%%%%%%%%%%%%%%%%%%%%%%%%%%%%%%%%%%%%%%%%%%%%%%%%%%
\subsection{Adopted definitions}
\label{defslw}

\begin{definition}[\slw] 
We will call a \slw, a path $\Gamma$ of class $\mathcal{C}^1$ in $\Er^2$, defined in Section \ref{principecourbe}. 
This path is defined up to isometries, up to a cyclic permutations and up to direction of travel.
In an abuse of terminology, we will also call a \slw\ that which allows us to define the path $\Gamma$, i.e, any of: 
\begin{itemize}
\item
the sequence of the centres of the squares ${(c_i)}_{1\leq i \leq N}$ occupied by the path $\Gamma$, satisfying all of the constraints given in sections
\ref{enumerationtouscircuit},
\ref{isometrie} and 
\ref{contrainteslocales}, this sequence being defined up to isometries, up to cyclic permutations and up to direction of travel;
\item
the sequence of angles ${(\alpha_i)}_{1\leq i \leq N}$ defined in sections
\ref{enumerationpositionprobleme}
 and 
\ref{enumerationtouscircuit}, this sequence being defined up to multiplication by $-1$, up to cyclic permutations and up to direction of travel; 
\item
The sequence of signed piece numbers ${(p_i)}_{1\leq i \leq N}$ defined in Section \ref{enumerationpositionprobleme}, this sequence being defined up to exchange of $\pm \piecec$ with $\mp \piececb$,
\ifcase \nopiecesix
(and $\pm \pieces$ with $\mp \piecesb$),
\or
\fi
up to cyclic permutations and up to direction of travel.
\end{itemize}
\end{definition}%

Note as well that a \slw\ depends on $N_j$ and on the width $e$ which, in this article, is chosen to be less than the critical width equal to the maximum of $e_0$ defined by 
\ifcase \nopiecesix
\eqref{maxlar6}
\or
\eqref{maxlar5}
\fi
and
$\delta_{\min}$
defined by 
\eqref{dminmin} and 
\eqref{dminminexa}.
In the slightly different case where $e>e_0$, $e$ must remain less than $1/2$; in the latter case, the connection rules, given in Section \ref{contrainteslocales}, are to be modified, which is also taken into account in the algorithms used.
\ifcase \nopiecesix
\or
If we decide to include the pieces of types \pieces\ and \piecesb, $e_0$ is then given by \eqref{maxlar6}.
\fi.%

\begin{definition}[circuit]
% Pour max, en anglais, laisser ici "circuit", puisque "J'ai traduit "circuit" par "track" lorsqu'il parle du jeu, et "circuit" lorsqu'il parle de la mathématique"
In the remainder of this article, a circuit is therefore the set of geometric pieces which is based on the geometric construction of a \slw.
\end{definition}
}

\iflanguage{french}{%

%%%%%%%%%%%%%%%%%%%%%%%%%%%%%%%%%%%%%%%%%%%%%%%%%%%%%%%%%%%%%
\subsection{Limitation informatique et complexité des algorithmes}
\label{limitationinfo}

Dans \cite{MR2065628}, le nombre de chemins autoévitants
a pu être déterminé de façon exact jusqu'à $N=71$ ; l'auteur  obtient 
$4$ $190$ $893$ $020$ $903$ $935$ $054$ $619$ $120$ $005$ $916$  chemins ! 
Dans \cite{MR2902304}, le nombre de polygones autoévitants
a pu être déterminé de façon exacte jusqu'à $N=130$ ; l'auteur  obtient 
$17$ $ 076$ $ 613$ $ 429$ $ 289$ $ 025$ $ 223$ $ 970 $ $687$ $ 974 $ $244 $ $417 $ $384$ $ 681$ $ 143$ $ 572 $ $320$
polygones !
Malheureusement, comme signalé plus haut, ces méthodes parallélisables n'ont pas pu être implémentées ici. 
Tous les algorithmes présentés ont été programmés sous Matlab \textregistered.
Deux versions ont été prévues :  la première est vectorielle (donc parallélisable) et évite l'usage des boucles, 
ce qui est relativement rapide. Cependant, les tableaux utilisés sont vite de taille très importante, ainsi que 
le nombre total de circuits à étudier.
Jusqu'à $N=9$, ces calculs sont possibles. Au-delà, la taille mémoire est trop importante.
Il faut passer alors à des calculs par boucles, qui sont beaucoup plus longs, mais qui évitent 
de stocker de grands tableaux correspondant aux circuits possibles.%
}{%

%%%%%%%%%%%%%%%%%%%%%%%%%%%%%%%%%%%%%%%%%%%%%%%%%%%%%%%%%%%%%
\subsection{Computational limitations and algorithmic complexity}
\label{limitationinfo}

In \cite{MR2065628}, the number of self-avoiding walks was able to be exactly determined up to $N=71$; the author obtained 
$4$ $190$ $893$ $020$ $903$ $935$ $054$ $619$ $120$ $005$ $916$
paths!
In \cite{MR2902304}, the number of self-avoiding polygons could be exactly determined up to $N=130$; the author obtains 
$17$ $ 076$ $ 613$ $ 429$ $ 289$ $ 025$ $ 223$ $ 970 $ $687$ $ 974 $ $244 $ $417 $ $384$ $ 681$ $ 143$ $ 572 $ $320$
polygons!
Unfortunately, as noted above, these parallelizable methods could not be implemented here. All of the algorithms presented have been programmed in Matlab \textregistered. Two versions were planned: the first is vectorial (and therefore parallelizable), and avoids the use of loops, which is relatively fast. However, the tables used are quickly of a significant size, as well as the total number of circuits to be studied. Up to $N=9$, these calculations are possible. Beyond that, the memory size is too great. It is necessary to move on to calculations with loops, which are much longer, but which avoid the storage of large tables corresponding to possible circuits.%
}
\ifcase 0
\ifcase 0
\iflanguage{french}{%
Jusqu'à  $N=11$, les calculs sont raisonnables.
Au-delà de $N=11$,
il aurait fallu plus de $4$ jours de calcul
(voir section \ref{enumerationpetit}) et les calculs n'ont pas été menés.%
}{%
Up to $N=11$, the calculations are reasonable. Beyond $N=11$, more than 
$4$
days of calculation were required (see Section~\ref{enumerationpetit}), and the calculations were not carried out.%
}
\ifcase 0
\or
\textbf{ATTENTION  !!!!!!!!!!!!!!!!!}
Durée inférieure à une heure !!
\fi
\or
\textbf{ATTENTION  !!!!!!!!!!!!!!!!!}
On ne peut déterminer la complexité : Règle non valable
\fi
\or
\textbf{ATTENTION  !!!!!!!!!!!!!!!!!}
On ne peut déterminer la complexité : Nombre de pièce max trop petit.
\fi

\iflanguage{french}{%
Les énumérations exhaustives des circuits possibles  font appel à des produits cartésiens d'ensembles finis et sont 
donc de complexité en $\mathcal{O}(A^N)$, ce qui, de toute façon, limite en théorie les calculs informatiques.

Dans 
\cite{%
MR2883859,%
MR1985492,%
MR1718791,%
Guttmann2012,%
Guttmann2012b,%
MR2902304},
une estimation du nombre de polygones autoévitants est donnée 
pour $N\to+\infty$ :
\begin{equation}
\label{estime_nombre_circuit}
q(N)  \thicksim A\mu^N N^{\gamma-1},
\end{equation}
où
$\mu$ est appelée la constante de connectivité, $\gamma$ l'exposant  critique et $A$ l'amplitude critique.
%where $\mu$ is called the connective constant, $\gamma$ is a critical exponent and $A$ a critical amplitude.
La valeur de $\mu$ proposée est la même que celle correspondant aux chemins autoévitants (voir \cite{MR1985492}) :
\begin{subequations}
\label{valeurmugamma}
\begin{equation}
\label{valeurmu}
\mu \approx 2.638.
\end{equation}
La valeur de $\gamma$ correspondant aux réseaux carrés est donné par
\begin{equation}
\label{valeurgamma}
\gamma-1 \approx -\frac{5}{2}, 
\end{equation}
\end{subequations}
et enfin, on a (voir \cite{MR1985492})
\begin{equation}
\label{valeurA}
A \approx 
0.0795 774 715.
\end{equation}%
L'estimation \eqref{estime_nombre_circuit}
n'est pas utilisable en l'état ici, puisque nous avons vu que la recherche des circuits et des chemins et des polygones  autoévitants
n'était pas exactement identique ; néanmoins, nous utiliserons abusivement cette approximation pour 
évaluer le nombre de circuits, pour des valeurs plus importantes de $N$ en section  \ref{estimation}.%
}{%
The exhaustive enumeration of possible circuits  makes use of Cartesian products of finite sets, and is therefore of complexity $\mathcal{O}(A^N)$, which, in any case, limits computer calculations in theory. 

In 
\cite{%
MR2883859,%
MR1985492,%
MR1718791,%
Guttmann2012,%
Guttmann2012b,%
MR2902304},
an estimation of the number of self-avoiding polygons is given for $N\to+\infty$:
\begin{equation}
\label{estime_nombre_circuit}
q(N)  \thicksim A\mu^N N^{\gamma-1},
\end{equation}
where
$\mu$ is called the connective constant, $\gamma$ a critical exponent and $A$ a critical amplitude.
The proposed value of $\mu$ is the same as the one corresponding to self-avoiding walks (see \cite{MR1985492}):
\begin{subequations}
\label{valeurmugamma}
\begin{equation}
\label{valeurmu}
\mu \approx 2.638.
\end{equation}
The value of $\gamma$ corresponding to square lattices is given by
\begin{equation}
\label{valeurgamma}
\gamma-1 \approx -\frac{5}{2}, 
\end{equation}
\end{subequations}
and finally, we have (see \cite{MR1985492})
\begin{equation}
\label{valeurA}
A \approx 
0.0795 774 715.
\end{equation}%
The estimation \eqref{estime_nombre_circuit} cannot be used here as is, since we have seen that the search for self-avoiding walks and polygons is not identical. Nevertheless, we will improperly use this approximation to evaluate the number of circuits for larger values of $N$ in Section~\ref{estimation}.%
}

%\begin{remark}
%M\label{newrem02}
%\iflanguage{french}{%
%Comme dans la généralisation de la remarque \ref{newrem01}, notons que, dans 
%\cite[tableau 3]{MR2104301}, Jensen propose la valeur suivante de $\mu$ dans le cas de triangles :%
%}{%
%As in the generalization of Remark~\ref{newrem01}, we note that, in \cite[tableau 3]{MR2104301}, Jensen proposes the following value of 
%$\mu$ in the case of triangles:%
%}
%\begin{equation}
%\label{valeurmubis}
%\mu\approx
%4.1
%\end{equation}
%\end{remark}

\begin{remark}
\label{newrem02}
\iflanguage{french}{%
Comme dans la généralisation de la remarque \ref{newrem01}, notons que, dans 
\cite[tableau 3]{MR2104301}, Jensen propose la valeur suivante de $\mu$ dans le cas de réseau de Kagomé  :%
}{%
As in the generalization of Remark~\ref{newrem01}, we note that, in \cite[tableau 3]{MR2104301}, Jensen proposes the following value of $\mu$ in the case of Kagomé lattice:%
}
\begin{equation}
\label{valeurmubis}
\mu\approx
2.560.
\end{equation}
\end{remark}

%%%%%%%%%%%%%%%%%%%%%%
% Attention : merdebeta
% peut-on virer cela !!!!
% TESTER ET VIRER 
%%%%%%%%%%%%%%%%%%%%%%

%\ifcase \choarticle
%\clearpage
%\or
%\fi

\iflanguage{french}{%

%%%%%%%%%%%%%%%%%%%%%%%%%%%%%%%%%%%%%%%%%%%%%%%%%%%%%%%%%%%%%
\subsection{Quelques exemples d'énumération de circuits}
\label{exempleenumeration}

Reprenons quelques exemples similaires à l'exemple \ref{examplesimulation500}
en utilisant cette fois tous les algorithmes de sélection de circuits proposés.%
}{%

%%%%%%%%%%%%%%%%%%%%%%%%%%%%%%%%%%%%%%%%%%%%%%%%%%%%%%%%%%%%%
\subsection{Some examples of the enumeration of circuits}
\label{exempleenumeration}

We reconsider some examples similar to Example~\ref{examplesimulation500}, this time using all of the proposed circuit selection algorithms.%
}

\begin{example}
\label{examplesimulation510}
%%%%%%%%%%%%%%%%%%%%%%%%%%%%%%%%%%%%%%%%%%%%%%%%%%%%%%%%%%%%%
%\input{simulations_circuit/simulation510}
% fichier tex crée par MaTeXBuild02 le 04-Sep-2015 09:12:15
% à compiler avec 
% MaTeXBuild02('simulation510',0)
% après le fichier 'enumeration_construction_circuit.matex'

%%%%%%%%%%%%%%%%%%%%%%%%%%%%%%%%%%%%%%%%%%%%%%%%%%%%%%%%%%%%%
%\input{./simulations_circuit/circuit_numerique/circuits_complets_7_piecmax_ev3}
% fichier crée par 'presentation_exhaustif_circuit_boucle.m' le 04-Sep-2015 09:12:17
\begin{figure}[h]
\centering
%%% sous figure 1
\subfigure[\label{circuits_complets_7_piecmax_ev31}]
{\epsfig{file=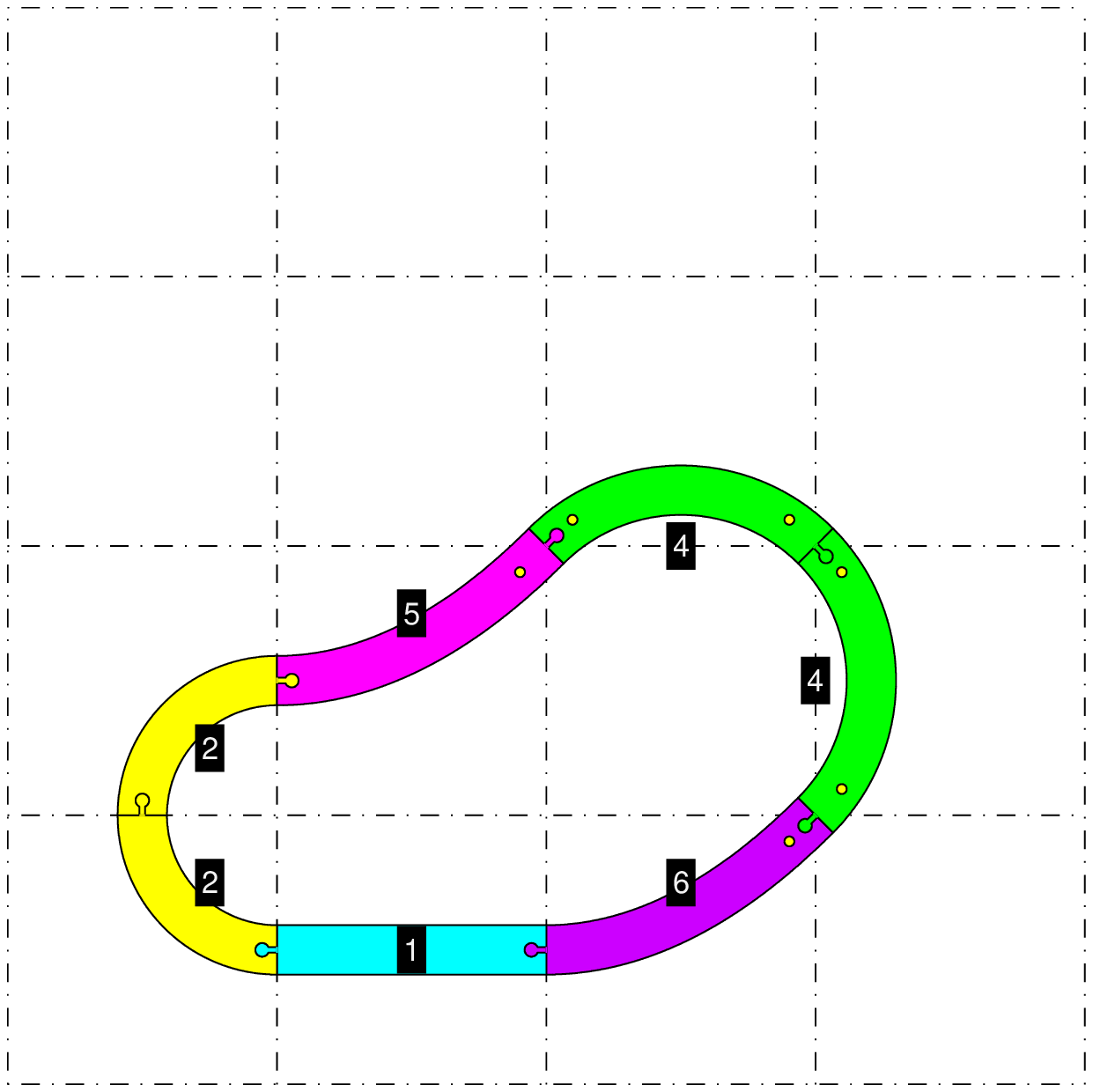, width=5 cm}}
\qquad
%%% sous figure 2
\subfigure[\label{circuits_complets_7_piecmax_ev32}]
{\epsfig{file=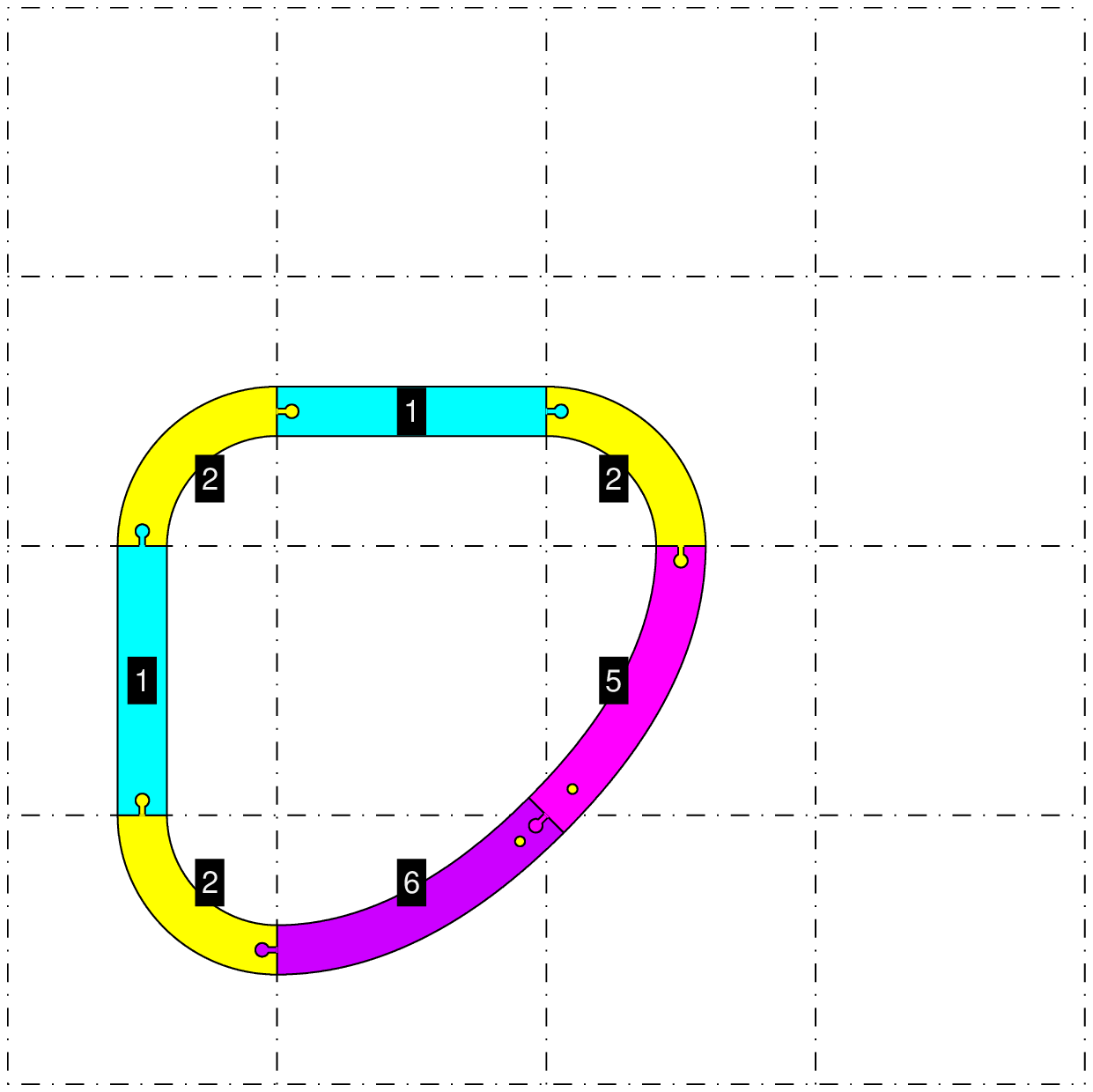, width=5 cm}}
\qquad
%%% sous figure 3
\subfigure[\label{circuits_complets_7_piecmax_ev33}]
{\epsfig{file=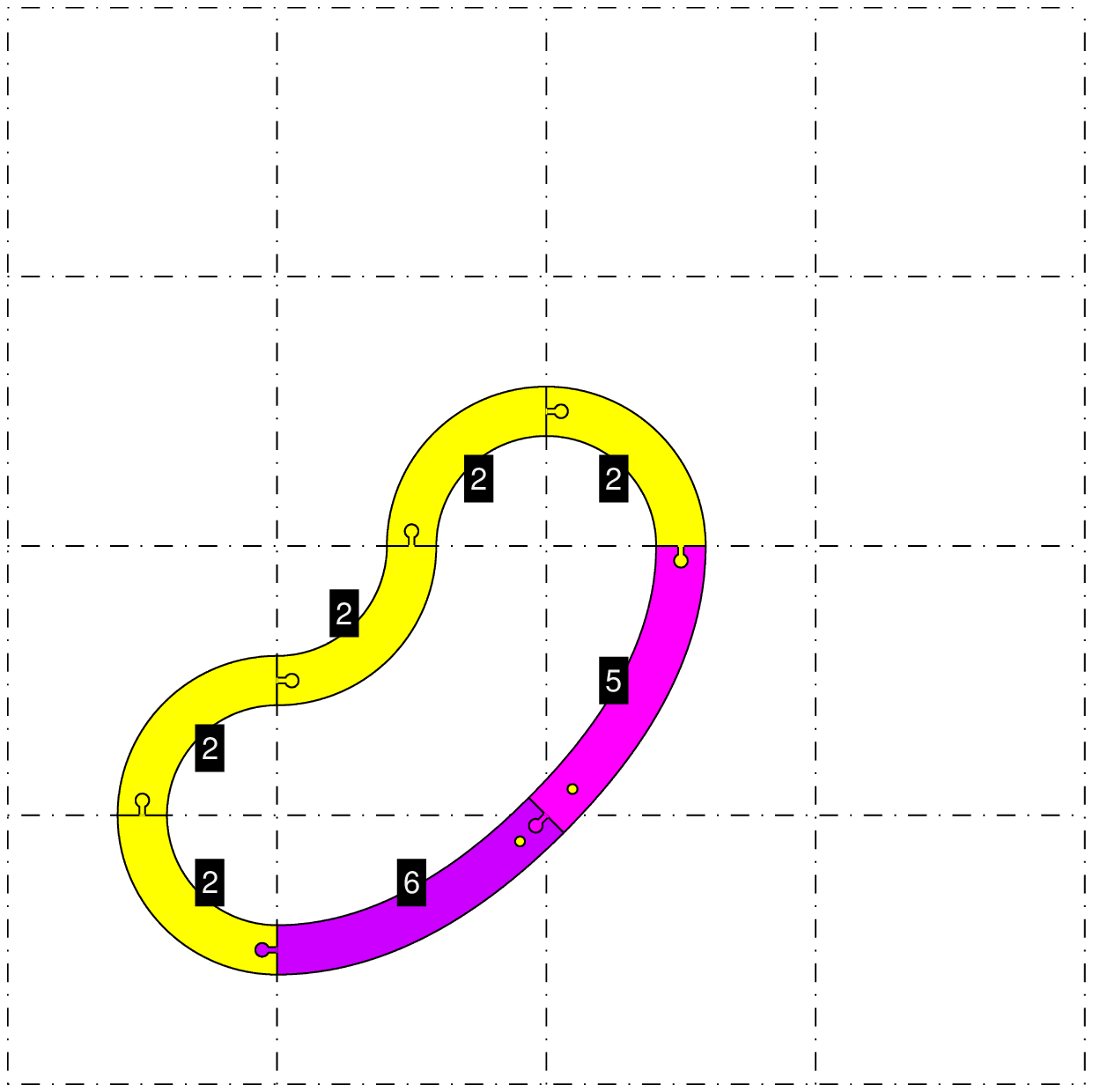, width=5 cm}}
\qquad
%%% sous figure 4
\subfigure[\label{circuits_complets_7_piecmax_ev34}]
{\epsfig{file=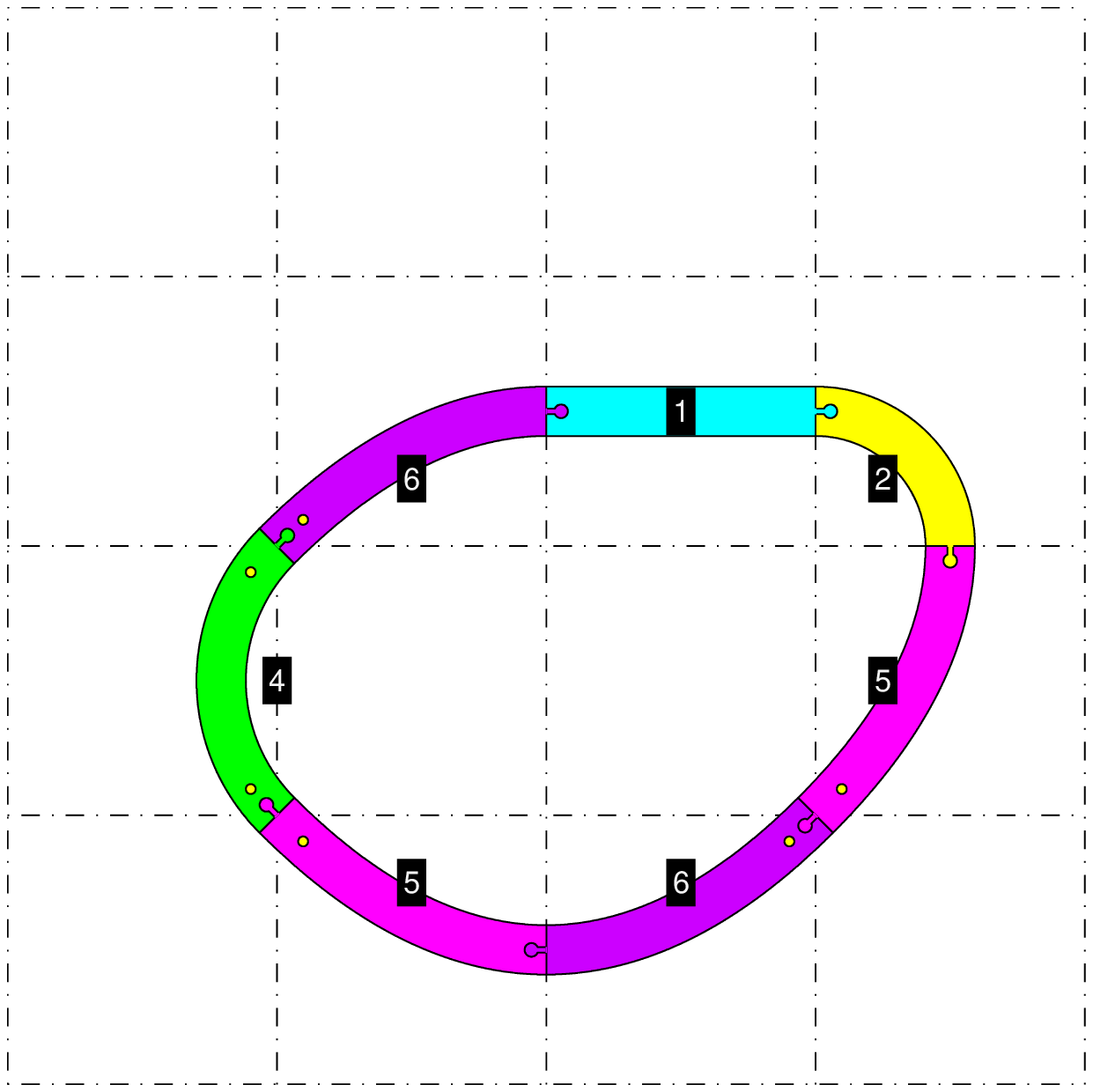, width=5 cm}}
\qquad
%%% sous figure 5
\subfigure[\label{circuits_complets_7_piecmax_ev35}]
{\epsfig{file=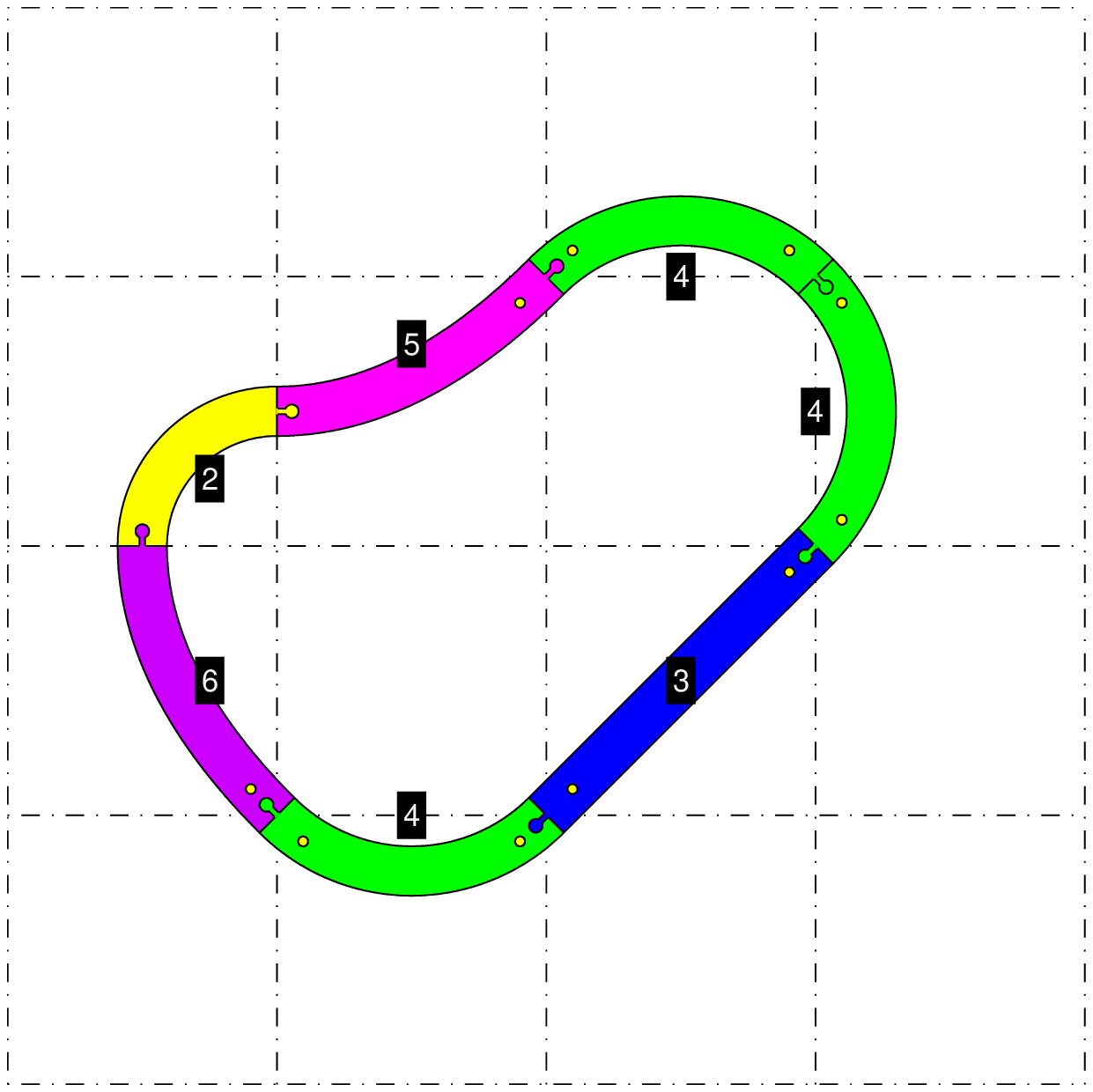, width=5 cm}}
\qquad
%%% sous figure 6
\subfigure[\label{circuits_complets_7_piecmax_ev36}]
{\epsfig{file=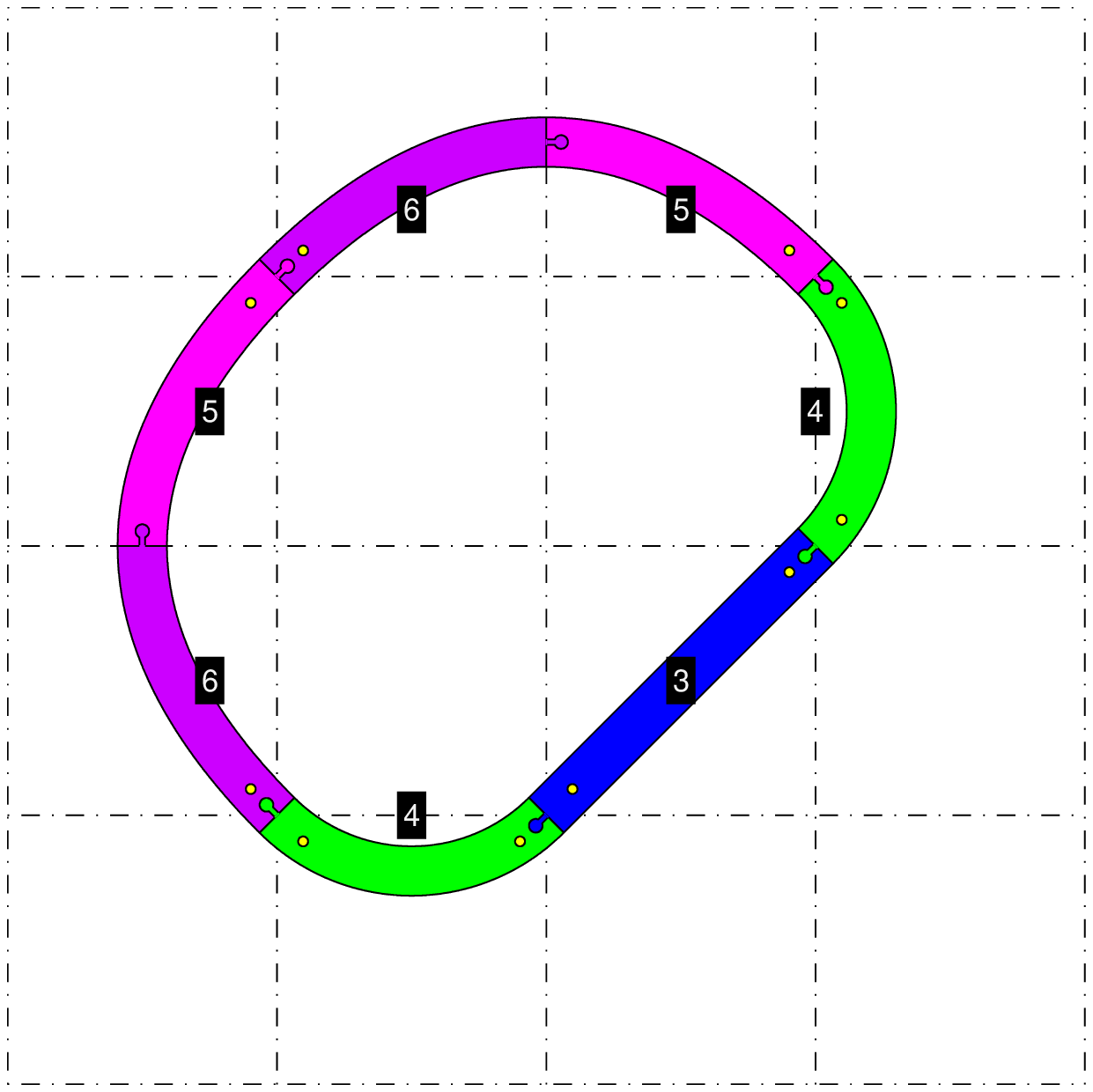, width=5 cm}}
\qquad
%%% sous figure 7
\subfigure[\label{circuits_complets_7_piecmax_ev37}]
{\epsfig{file=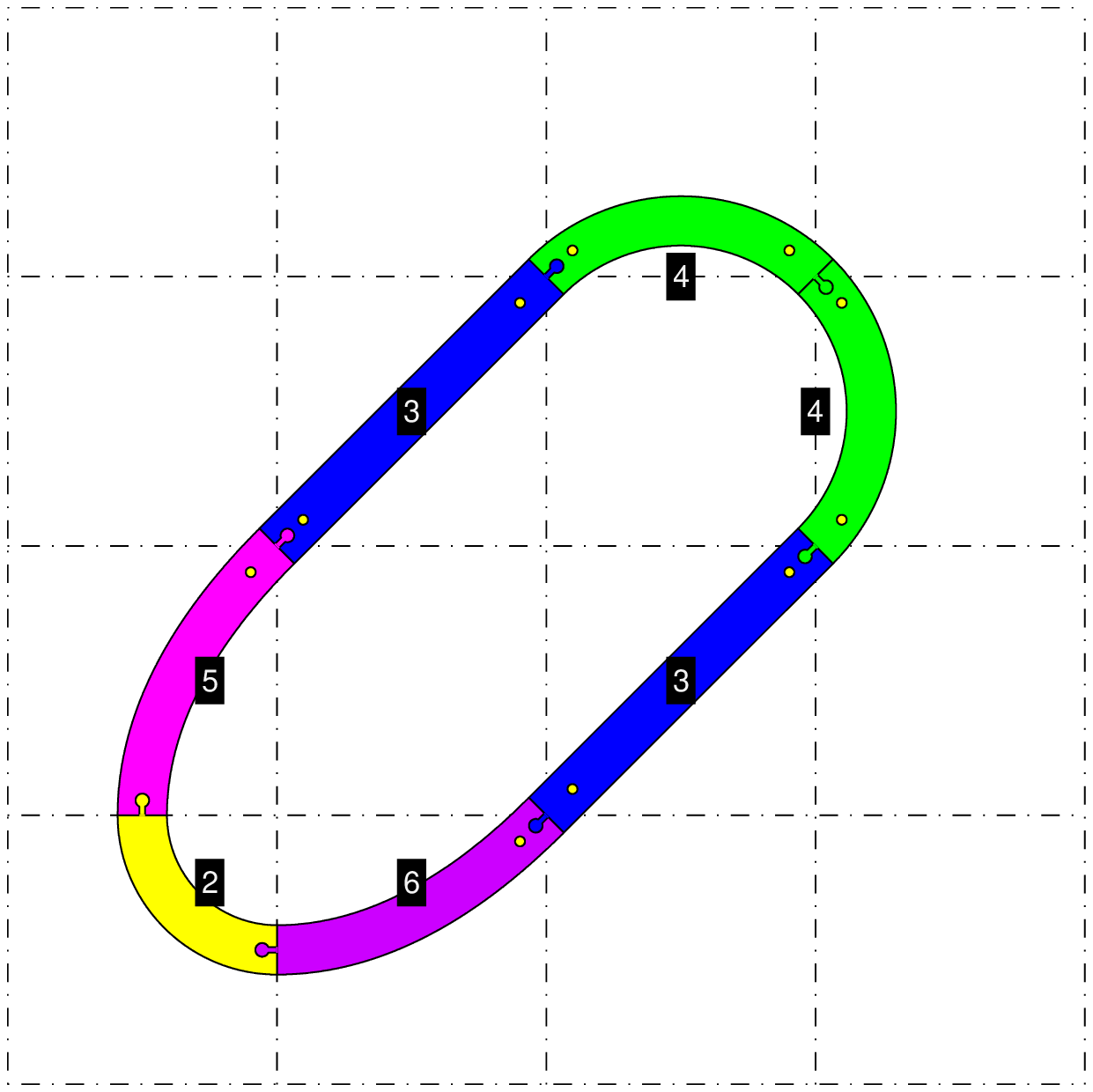, width=5 cm}}
\qquad
\caption{\label{circuits_complets_7_piecmax_ev3}\iflanguage{french}{Tous les 7 circuits retenus sur l'ensemble des 50000 circuits possibles}{All of the 7 circuits kept from the set of 50000 possible circuits}.}
\end{figure}
%%%%%%%%%%%%%%%%%%%%%%%%%%%%%%%%%%%%%%%%%%%%%%%%%%%%%%%%%%%%%

\iflanguage{french}{%
Si on trace tous les  circuits réalisables %tous les circuits réalisables 
avec $N=7$ pièces
et $N_j=+\infty$, on obtient les  $7$ circuits de la figure \ref{circuits_complets_7_piecmax_ev3}. 
Sur cette figure, seuls les circuits, tous différents à une isométrie près et constructibles ont été tracés.%
}{%
If we draw 
all the feasible circuits
%all of the feasible circuits 
with $N=7$ pieces and $N_j=+\infty$, we obtain the $7$ circuits in Figure~\ref{circuits_complets_7_piecmax_ev3}. In this figure, only the circuits which are all different up to an isometry have been drawn.%
}

%%%%%%%%%%%%%%%%%%%%%%%%%%%%%%%%%%%%%%%%%%%%%%%%%%%%%%%%%%%%%
\end{example}

\begin{example}
\label{examplesimulation515}
%%%%%%%%%%%%%%%%%%%%%%%%%%%%%%%%%%%%%%%%%%%%%%%%%%%%%%%%%%%%%
%\input{simulations_circuit/simulation515}
% fichier tex crée par MaTeXBuild02 le 19-Feb-2016 16:33:48
% à compiler avec 
% MaTeXBuild02('simulation515',0)
% après le fichier 'enumeration_construction_circuit.matex'

%%%%%%%%%%%%%%%%%%%%%%%%%%%%%%%%%%%%%%%%%%%%%%%%%%%%%%%%%%%%%
%\input{./simulations_circuit/circuit_numerique/circuits_complets_8_piecmax_ev3}
% fichier crée par 'presentation_exhaustif_circuit_boucle.m' le 19-Feb-2016 16:33:52
\begin{figure}[h]
\centering
%%% sous figure 1
\subfigure[\label{circuits_complets_8_piecmax_ev31}]
{\epsfig{file=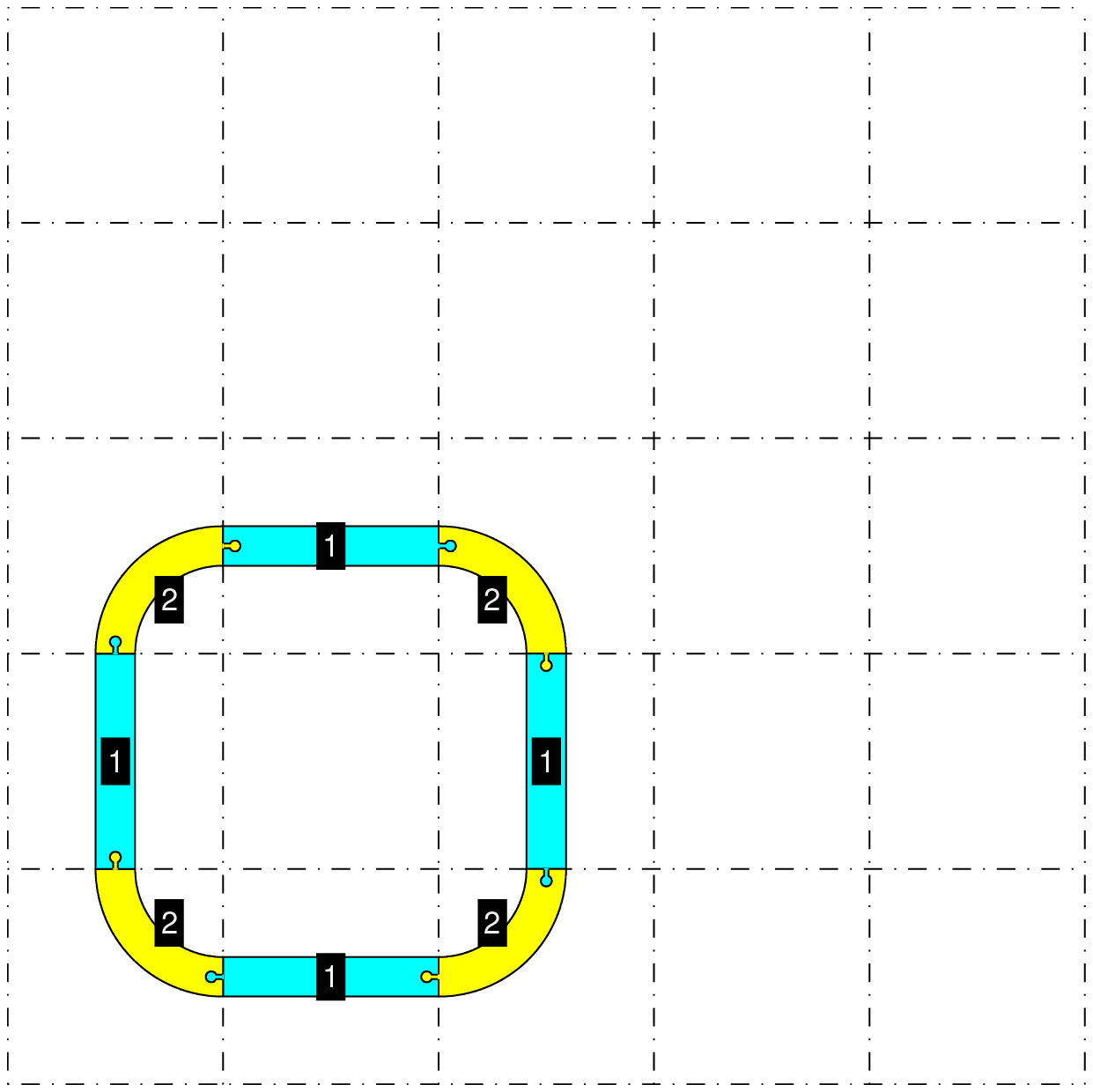, width=5 cm}}
\qquad
%%% sous figure 2
\subfigure[\label{circuits_complets_8_piecmax_ev32}]
{\epsfig{file=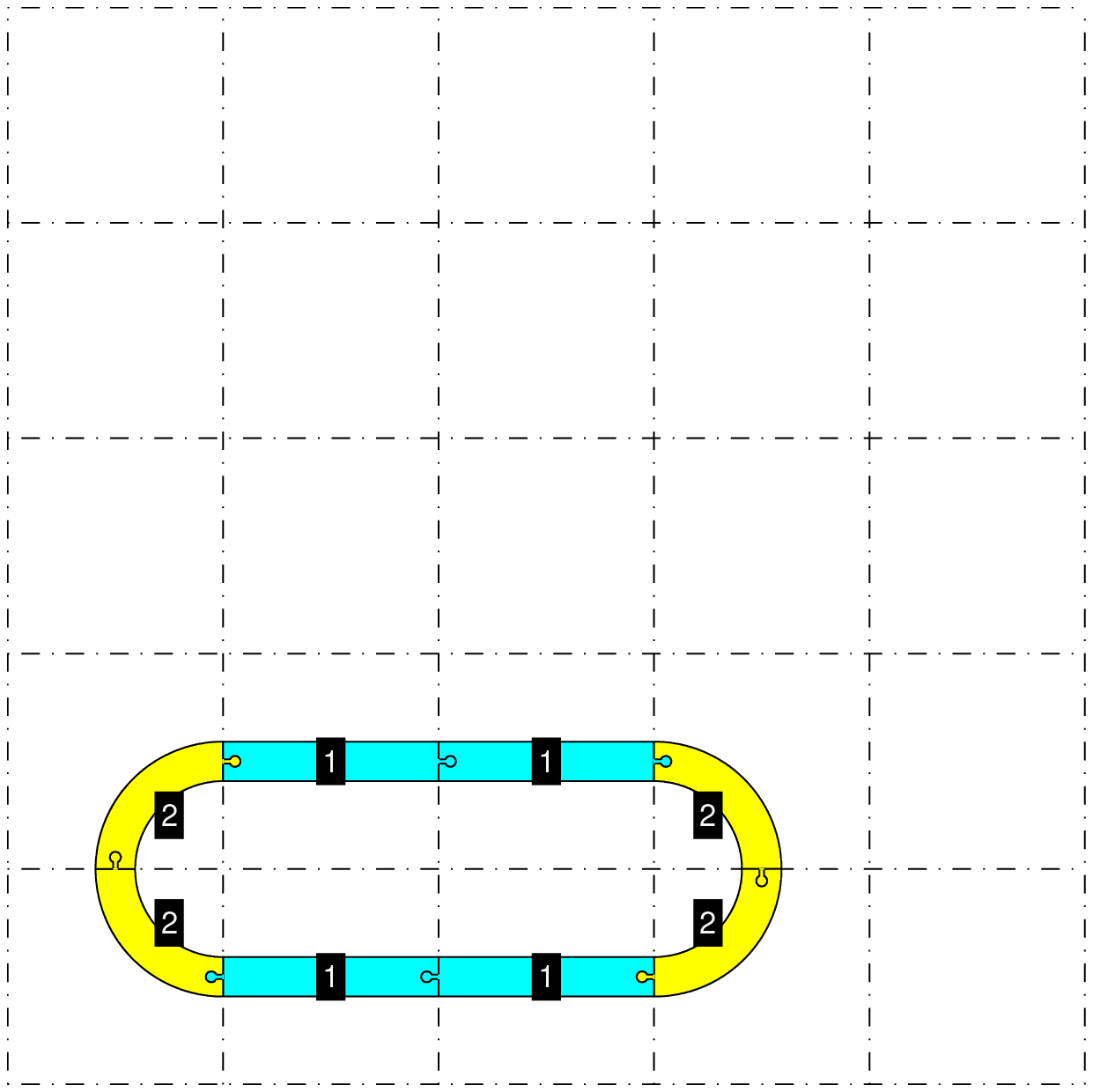, width=5 cm}}
\qquad
%%% sous figure 3
\subfigure[\label{circuits_complets_8_piecmax_ev33}]
{\epsfig{file=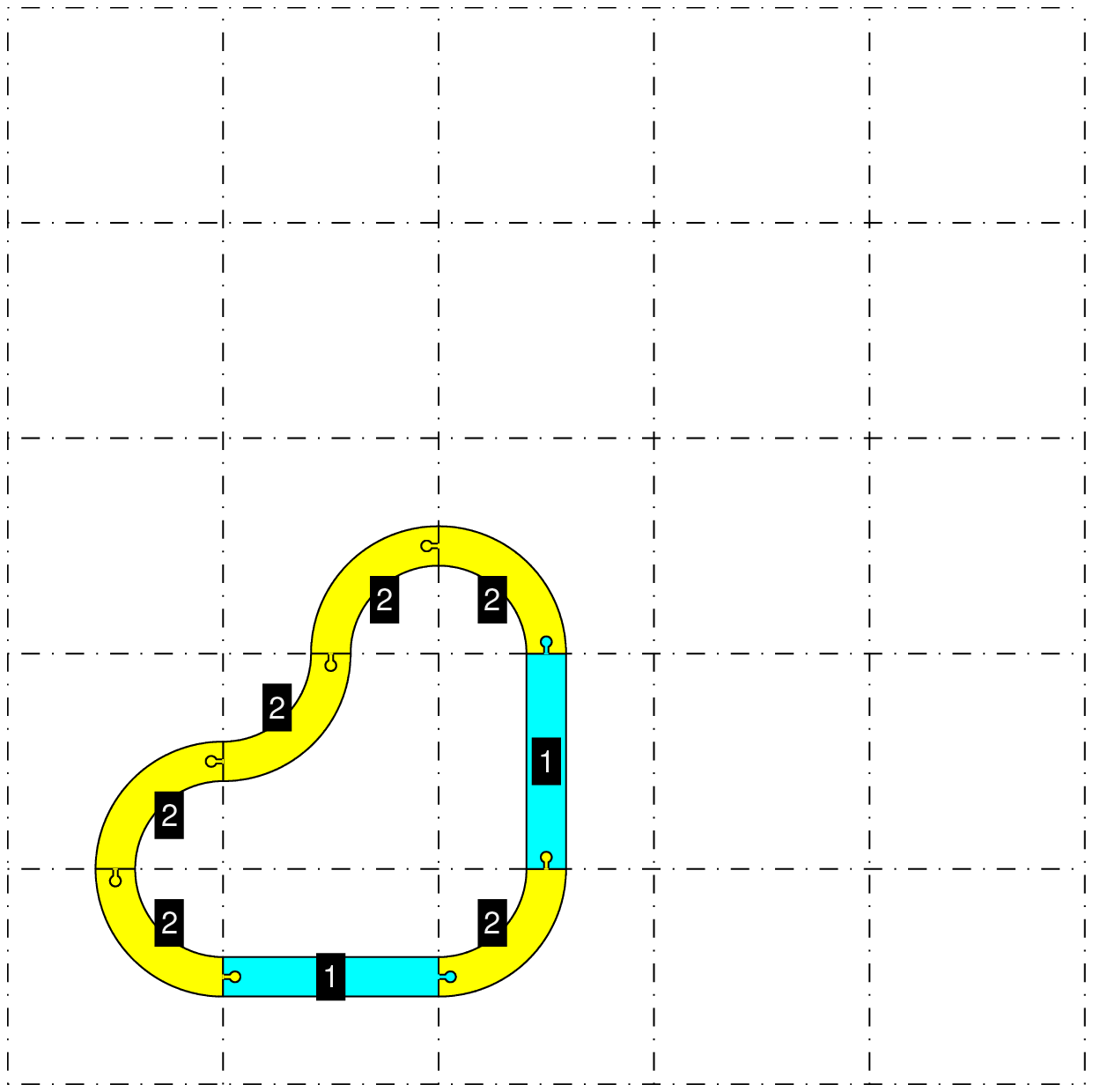, width=5 cm}}
\qquad
%%% sous figure 4
\subfigure[\label{circuits_complets_8_piecmax_ev34}]
{\epsfig{file=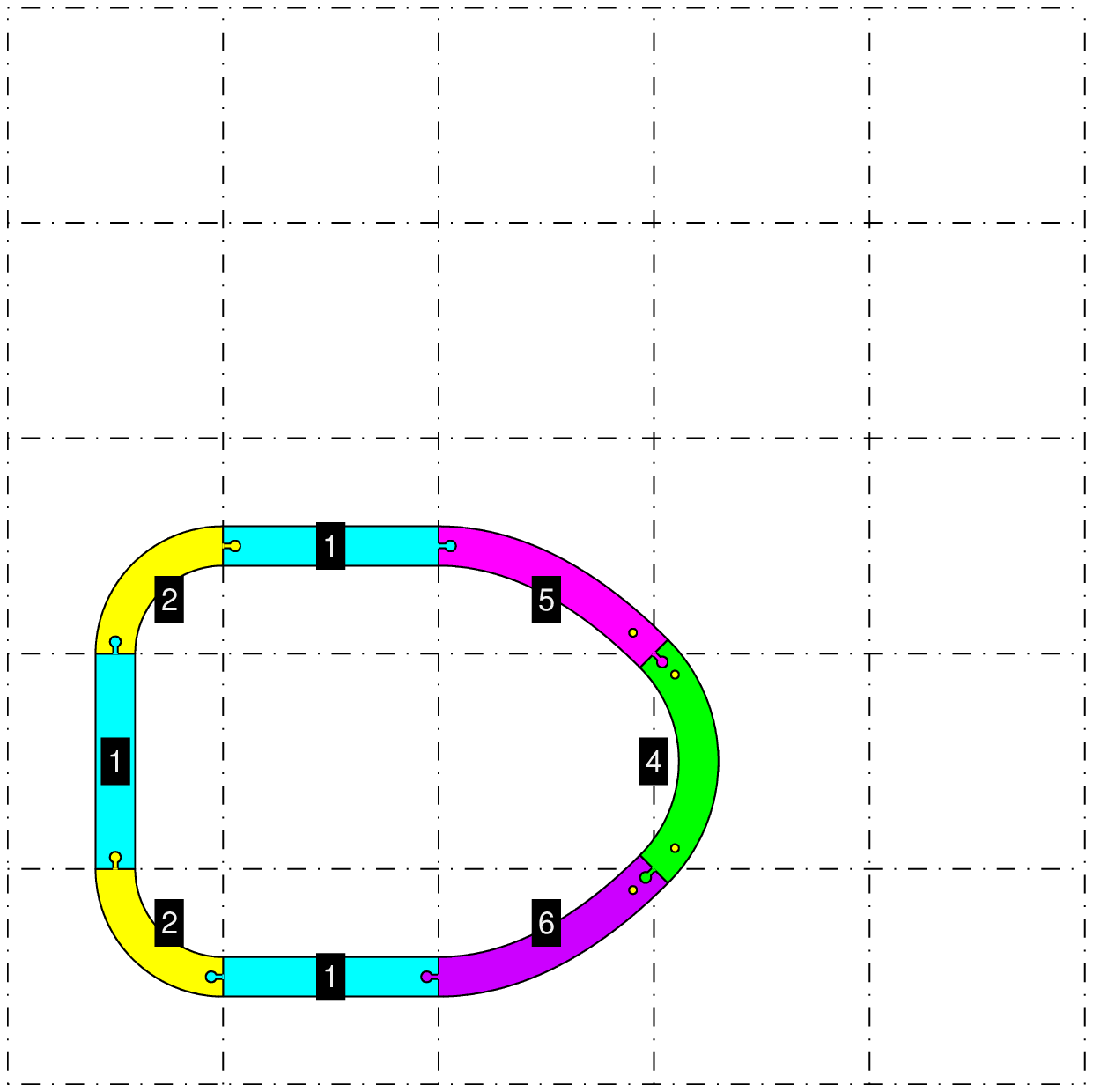, width=5 cm}}
\qquad
%%% sous figure 5
\subfigure[\label{circuits_complets_8_piecmax_ev35}]
{\epsfig{file=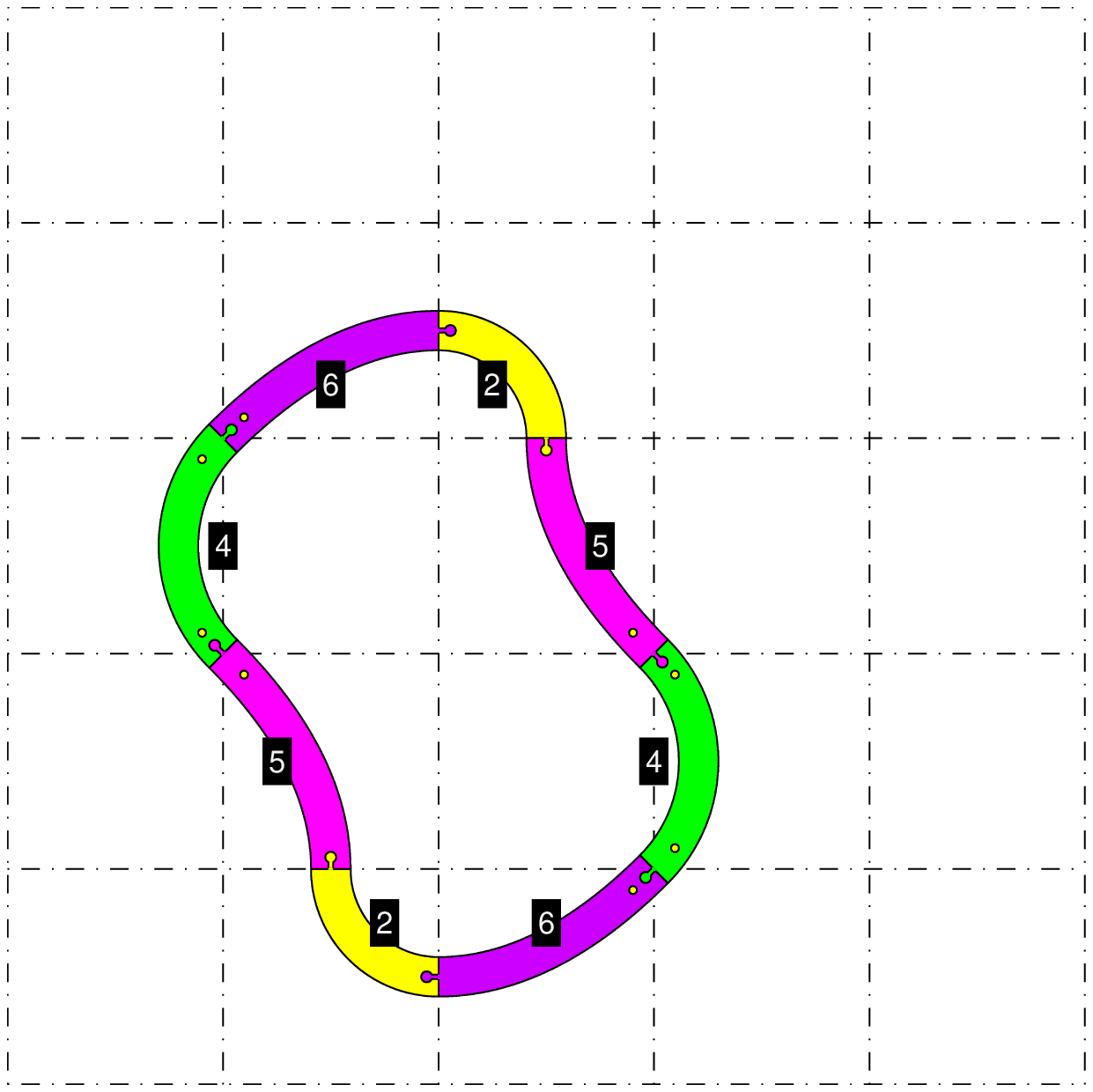, width=5 cm}}
\qquad
%%% sous figure 6
\subfigure[\label{circuits_complets_8_piecmax_ev36}]
{\epsfig{file=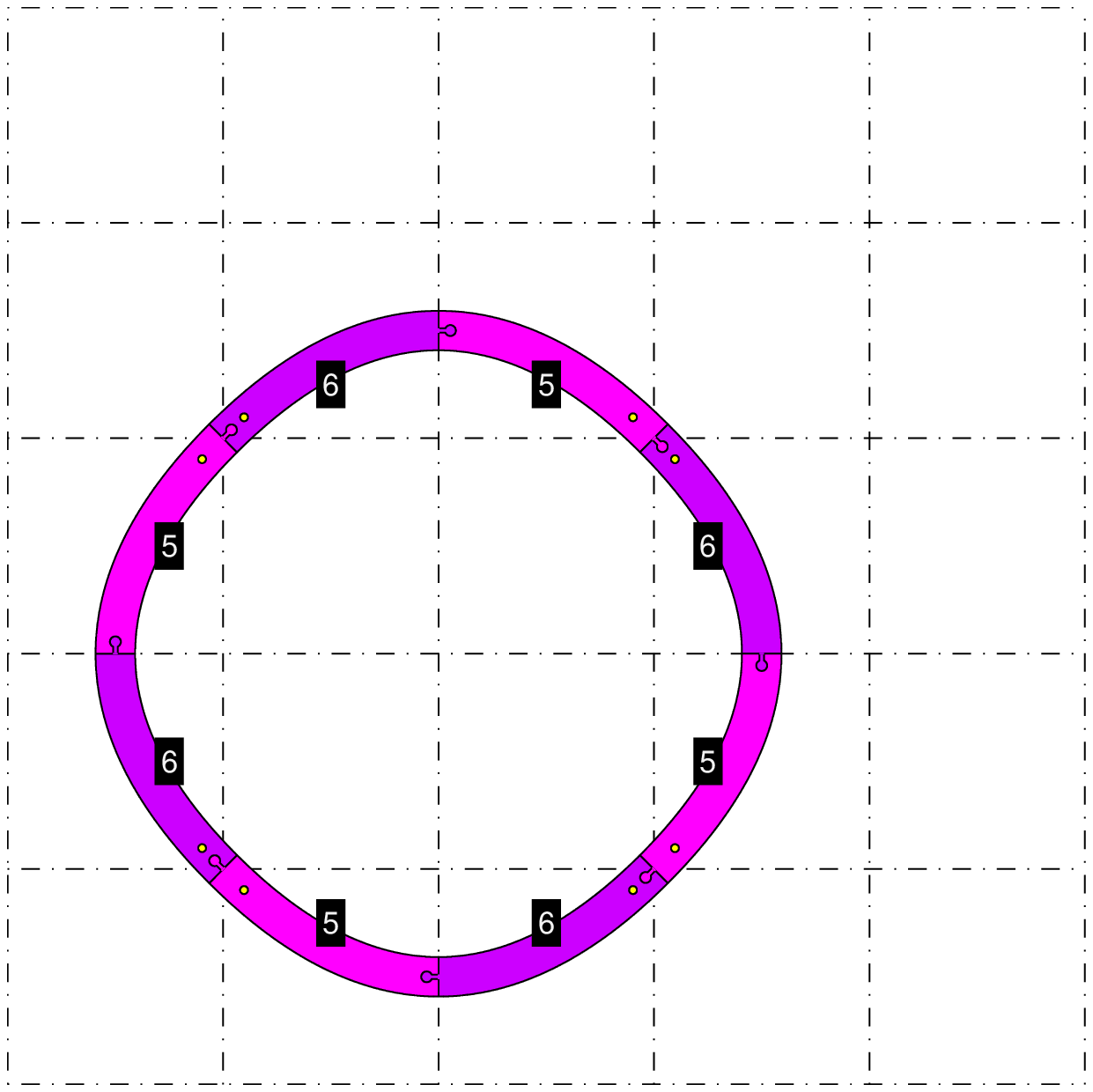, width=5 cm}}
\qquad
%%% sous figure 7
\subfigure[\label{circuits_complets_8_piecmax_ev37}]
{\epsfig{file=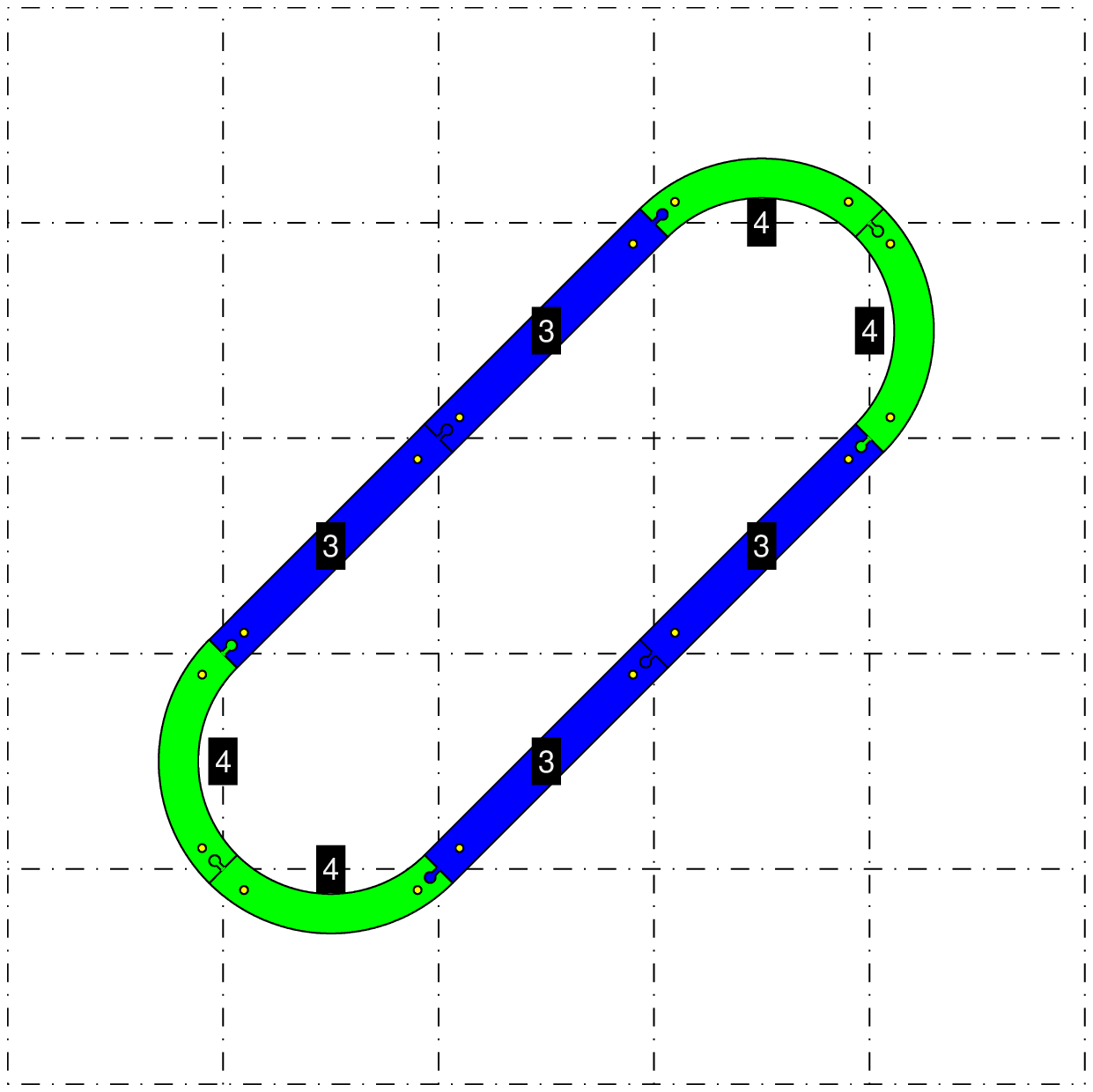, width=5 cm}}
\qquad
%%% sous figure 8
\subfigure[\label{circuits_complets_8_piecmax_ev38}]
{\epsfig{file=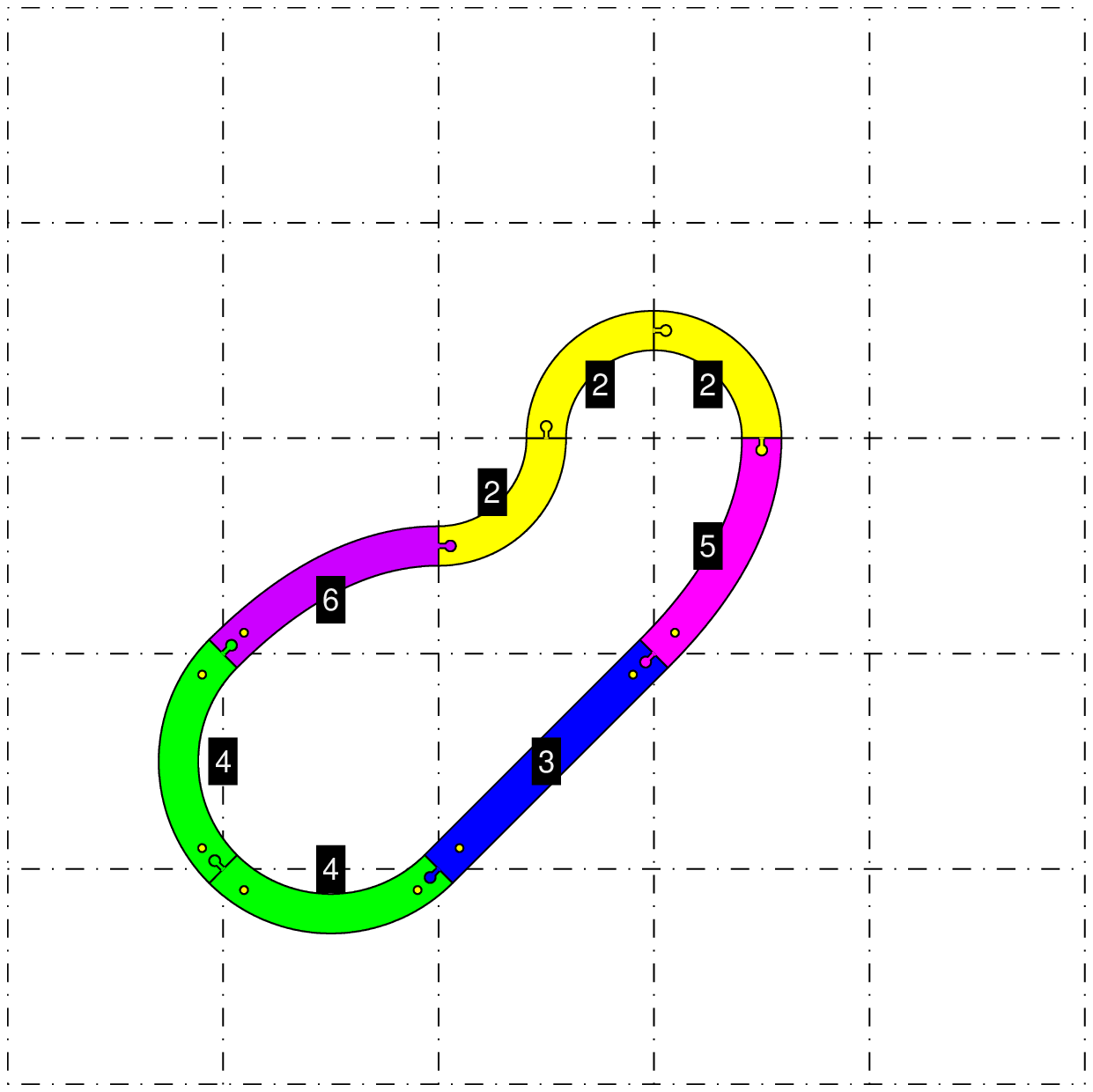, width=5 cm}}
\qquad
%%% sous figure 9
\subfigure[\label{circuits_complets_8_piecmax_ev39}]
{\epsfig{file=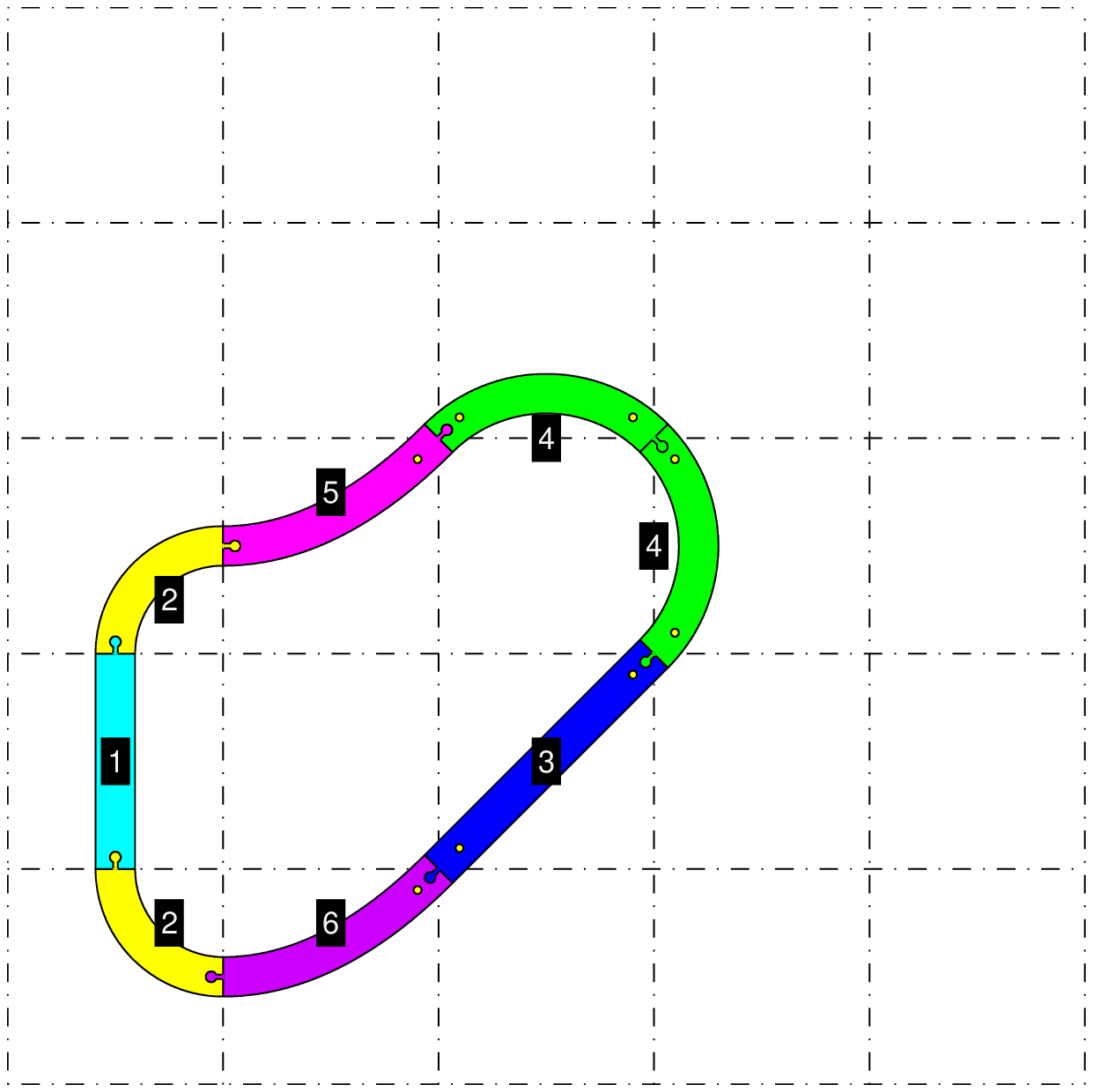, width=5 cm}}
\qquad
%%% sous figure 10
\subfigure[\label{circuits_complets_8_piecmax_ev310}]
{\epsfig{file=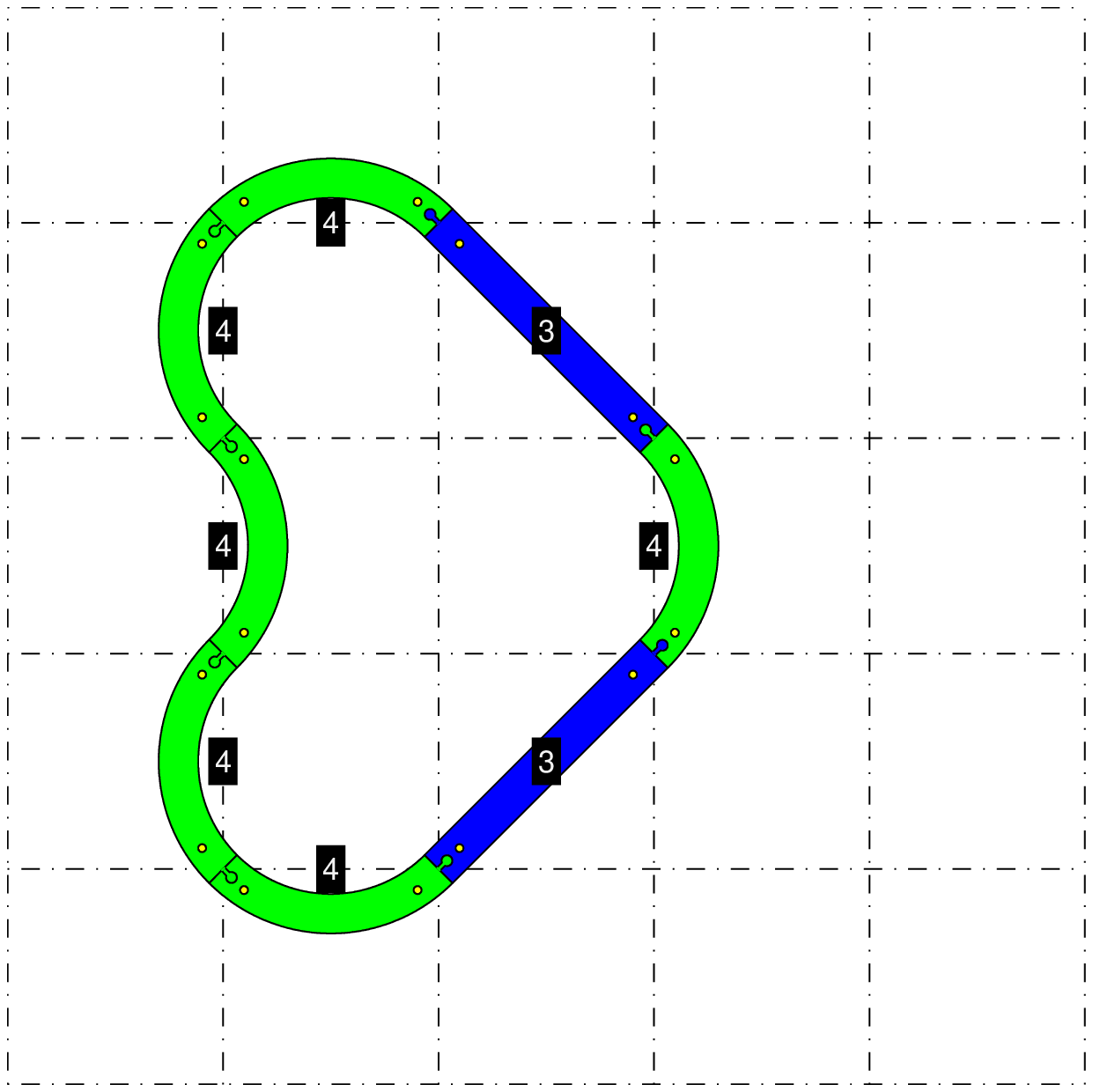, width=5 cm}}
\qquad
\caption{\label{circuits_complets_8_piecmax_ev3}\iflanguage{french}{10 des 33 circuits retenus sur l'ensemble des 250000 circuits possibles}{10 of the 33 circuits retained from the set of 250000 possible circuits}.}
\end{figure}
%%%%%%%%%%%%%%%%%%%%%%%%%%%%%%%%%%%%%%%%%%%%%%%%%%%%%%%%%%%%%

\iflanguage{french}{%
Si on trace quelques uns des circuits réalisables %tous les circuits réalisables 
avec $N=8$ pièces
et $N_j=+\infty$, on obtient les  $10$ circuits de la figure \ref{circuits_complets_8_piecmax_ev3}. 
Sur cette figure, seuls les circuits, tous différents à une isométrie près et constructibles ont été tracés.
Notons que les figures 
\ref{circuits_complets_8_piecmax_ev31},
\ref{circuits_complets_8_piecmax_ev32} et 
\ref{circuits_complets_8_piecmax_ev33}
correspondent à des circuits où ne sont utilisées que les pièces \pieceu\ et \pieced.
Nous reviendrons plus loin sur ces circuits particuliers.%
}{%
If we draw 
some of the feasible circuits
%all of the feasible circuits 
with $N=8$ pieces and $N_j=+\infty$, we obtain the $10$ circuits in Figure~\ref{circuits_complets_8_piecmax_ev3}. In this figure, only the circuits which are all different up to an isometry have been drawn.
%  Attention, anglais non traduit par ml !!
Note that Figures 
\ref{circuits_complets_8_piecmax_ev31},
\ref{circuits_complets_8_piecmax_ev32} and
\ref{circuits_complets_8_piecmax_ev33}
correspond to circuits which are used only the parts \pieceu\ et \pieced. We shall return to these particular circuits.%
}

%%%%%%%%%%%%%%%%%%%%%%%%%%%%%%%%%%%%%%%%%%%%%%%%%%%%%%%%%%%%%
\end{example}

\begin{example}
\label{examplesimulation520}
%%%%%%%%%%%%%%%%%%%%%%%%%%%%%%%%%%%%%%%%%%%%%%%%%%%%%%%%%%%%%
%\input{simulations_circuit/simulation520}
% fichier tex crée par MaTeXBuild02 le 04-Sep-2015 09:12:17
% à compiler avec 
% MaTeXBuild02('simulation520',0)
% après le fichier 'enumeration_construction_circuit.matex'

%%%%%%%%%%%%%%%%%%%%%%%%%%%%%%%%%%%%%%%%%%%%%%%%%%%%%%%%%%%%%
%\input{./simulations_circuit/circuit_numerique/circuits_complets_9_piecmax_ev3}
% fichier crée par 'presentation_exhaustif_circuit_boucle.m' le 04-Sep-2015 09:12:40
\begin{figure}[h]
\centering
%%% sous figure 1
\subfigure[\label{circuits_complets_9_piecmax_ev31}]
{\epsfig{file=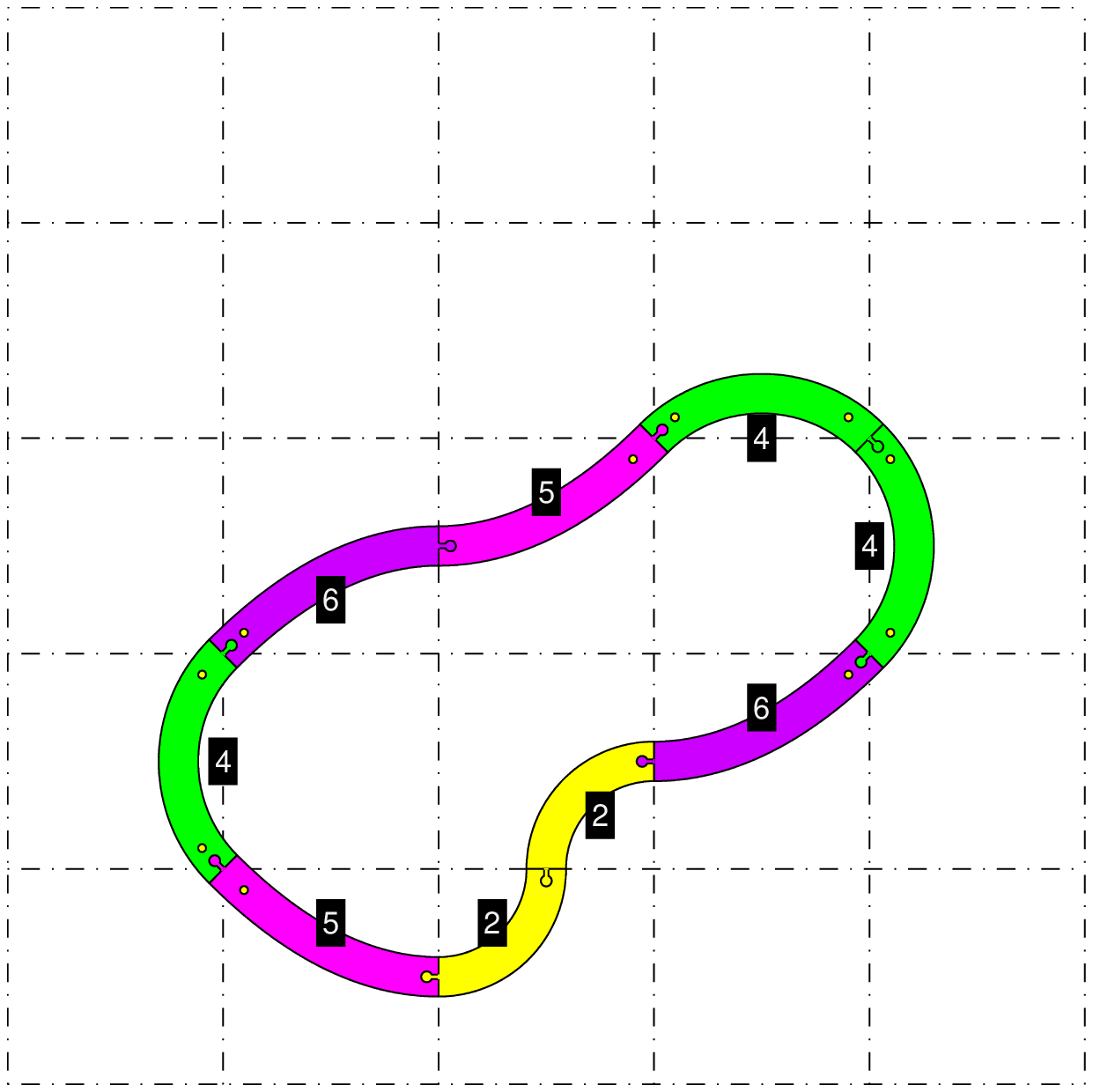, width=5 cm}}
\qquad
%%% sous figure 2
\subfigure[\label{circuits_complets_9_piecmax_ev32}]
{\epsfig{file=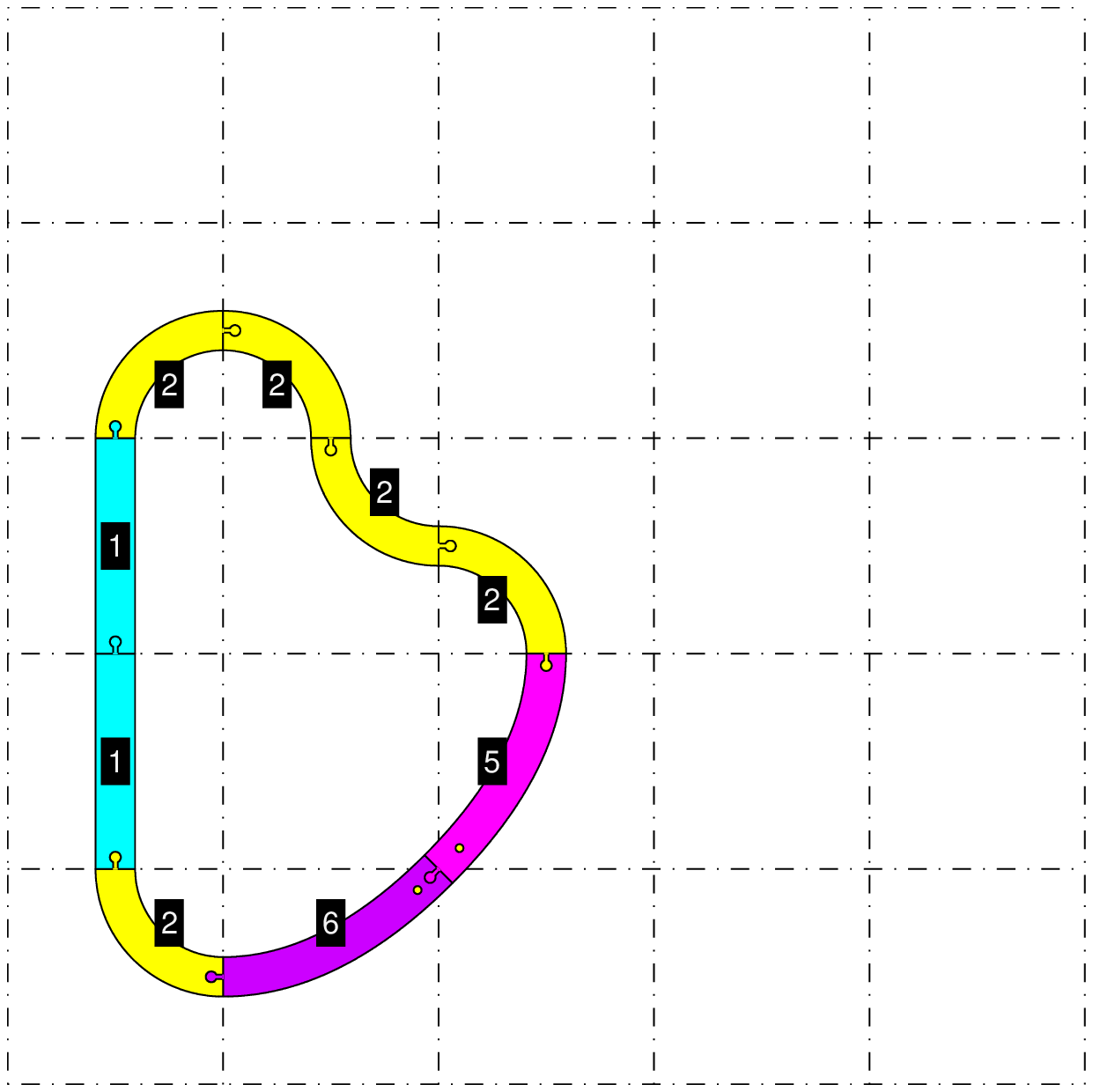, width=5 cm}}
\qquad
%%% sous figure 3
\subfigure[\label{circuits_complets_9_piecmax_ev33}]
{\epsfig{file=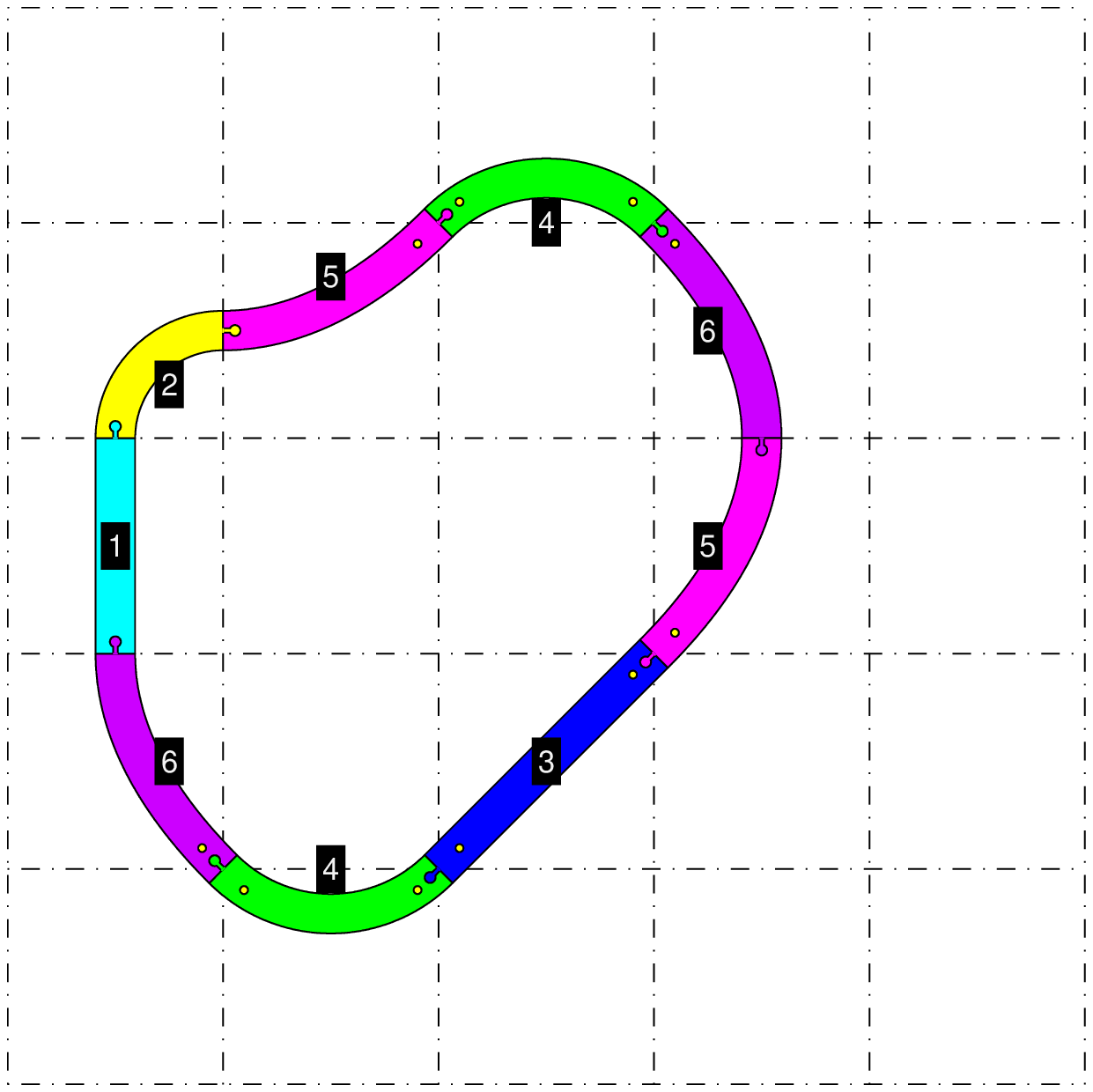, width=5 cm}}
\qquad
%%% sous figure 4
\subfigure[\label{circuits_complets_9_piecmax_ev34}]
{\epsfig{file=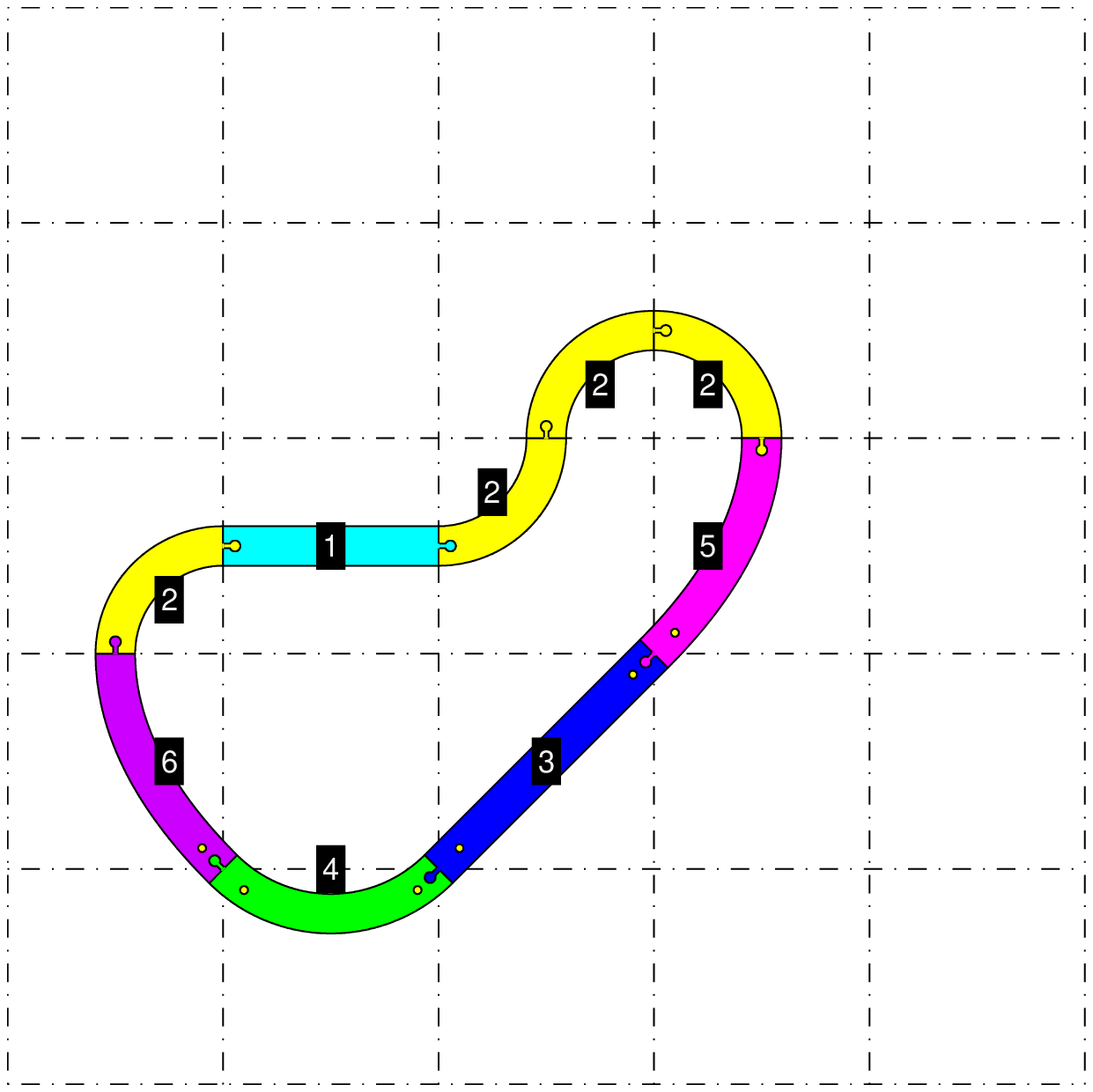, width=5 cm}}
\qquad
%%% sous figure 5
\subfigure[\label{circuits_complets_9_piecmax_ev35}]
{\epsfig{file=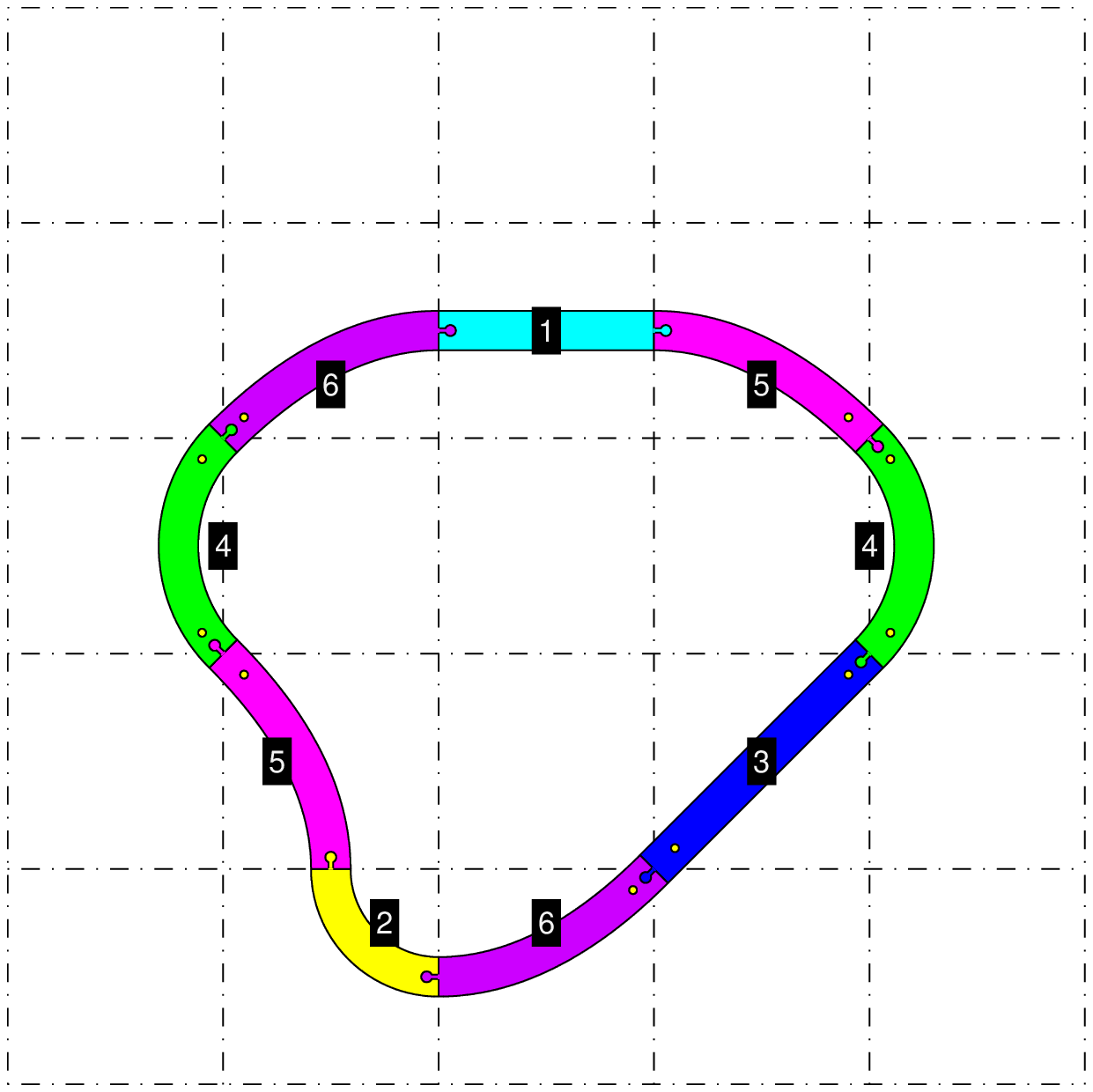, width=5 cm}}
\qquad
%%% sous figure 6
\subfigure[\label{circuits_complets_9_piecmax_ev36}]
{\epsfig{file=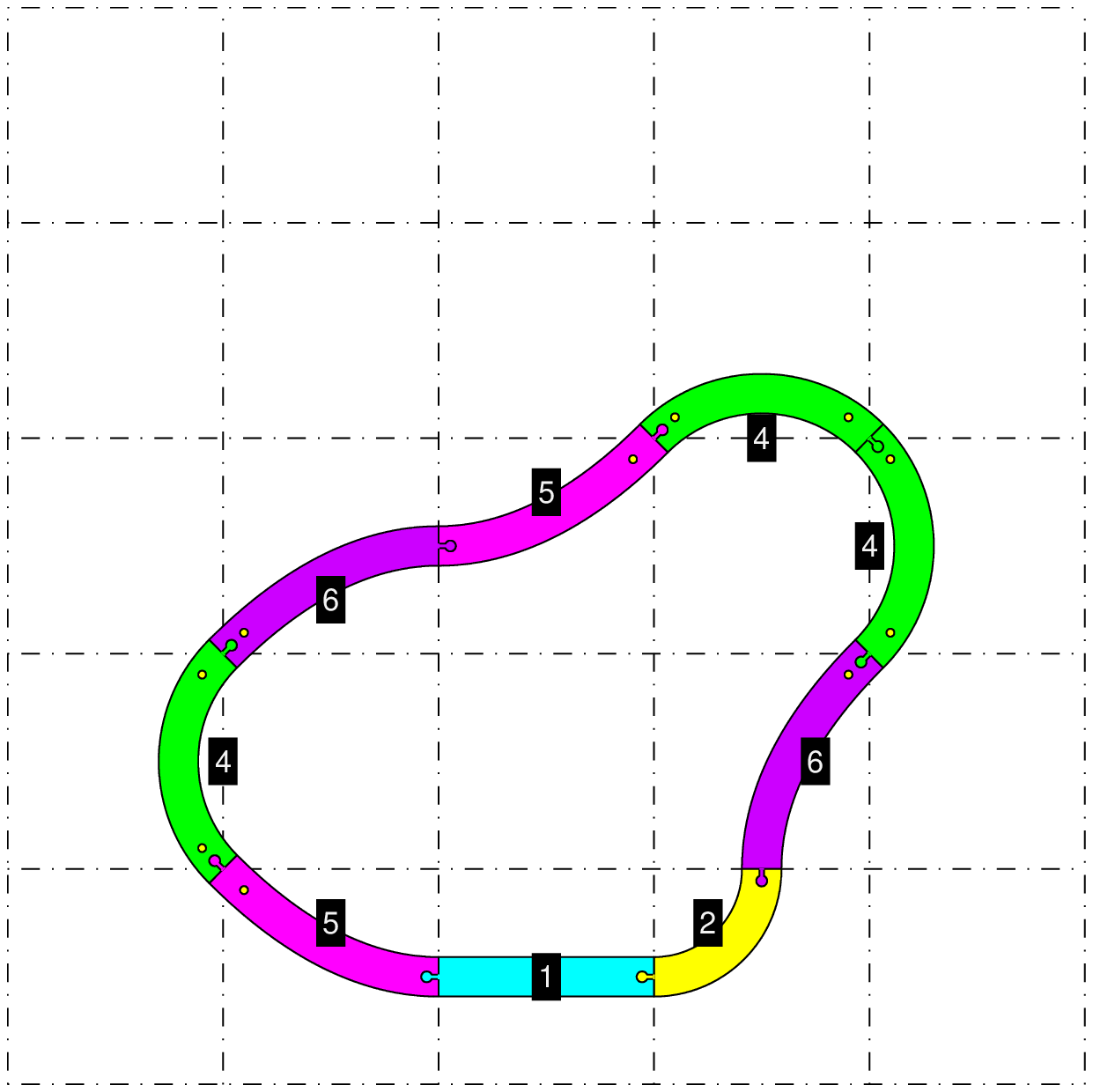, width=5 cm}}
\qquad
%%% sous figure 7
\subfigure[\label{circuits_complets_9_piecmax_ev37}]
{\epsfig{file=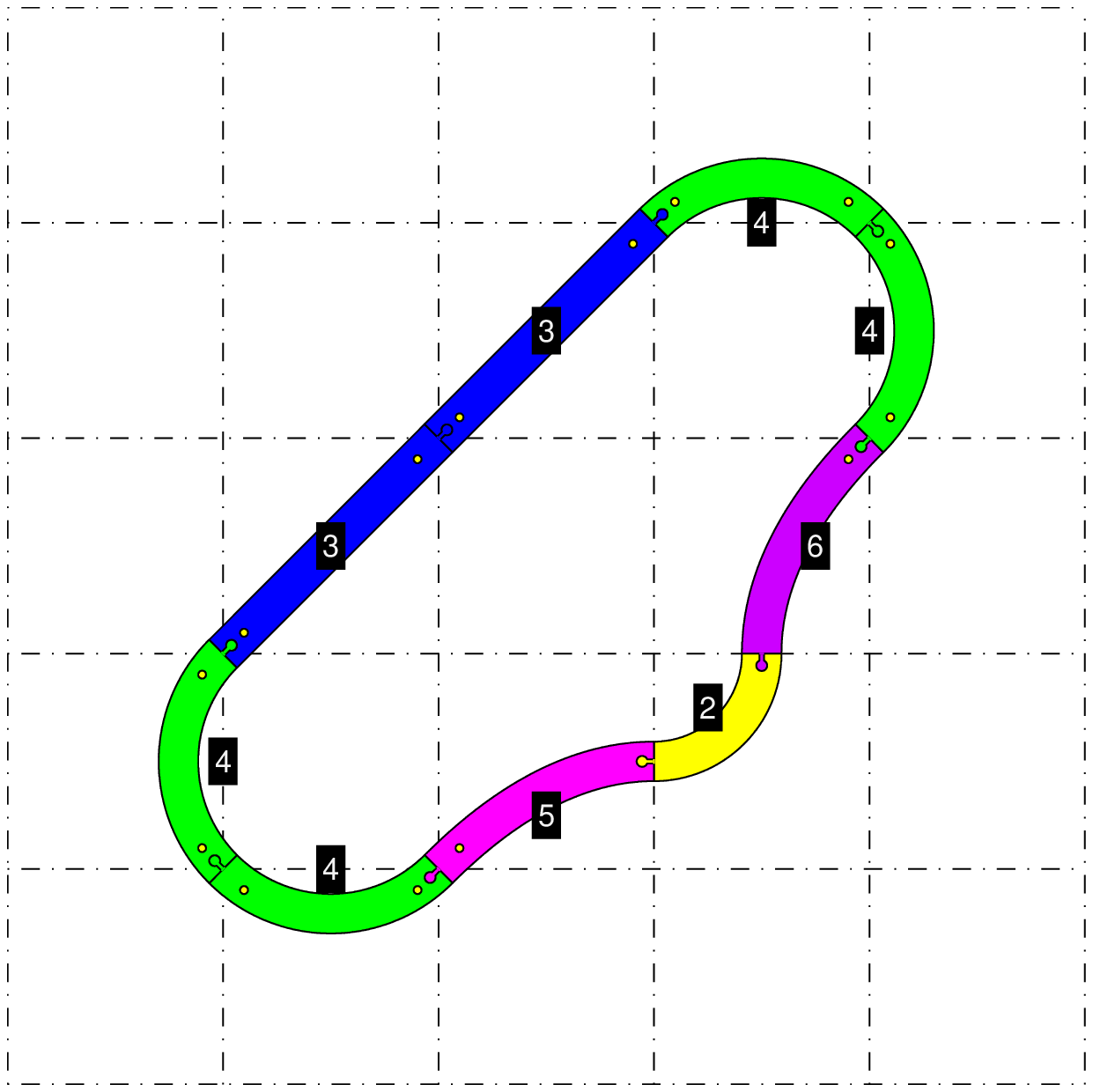, width=5 cm}}
\qquad
%%% sous figure 8
\subfigure[\label{circuits_complets_9_piecmax_ev38}]
{\epsfig{file=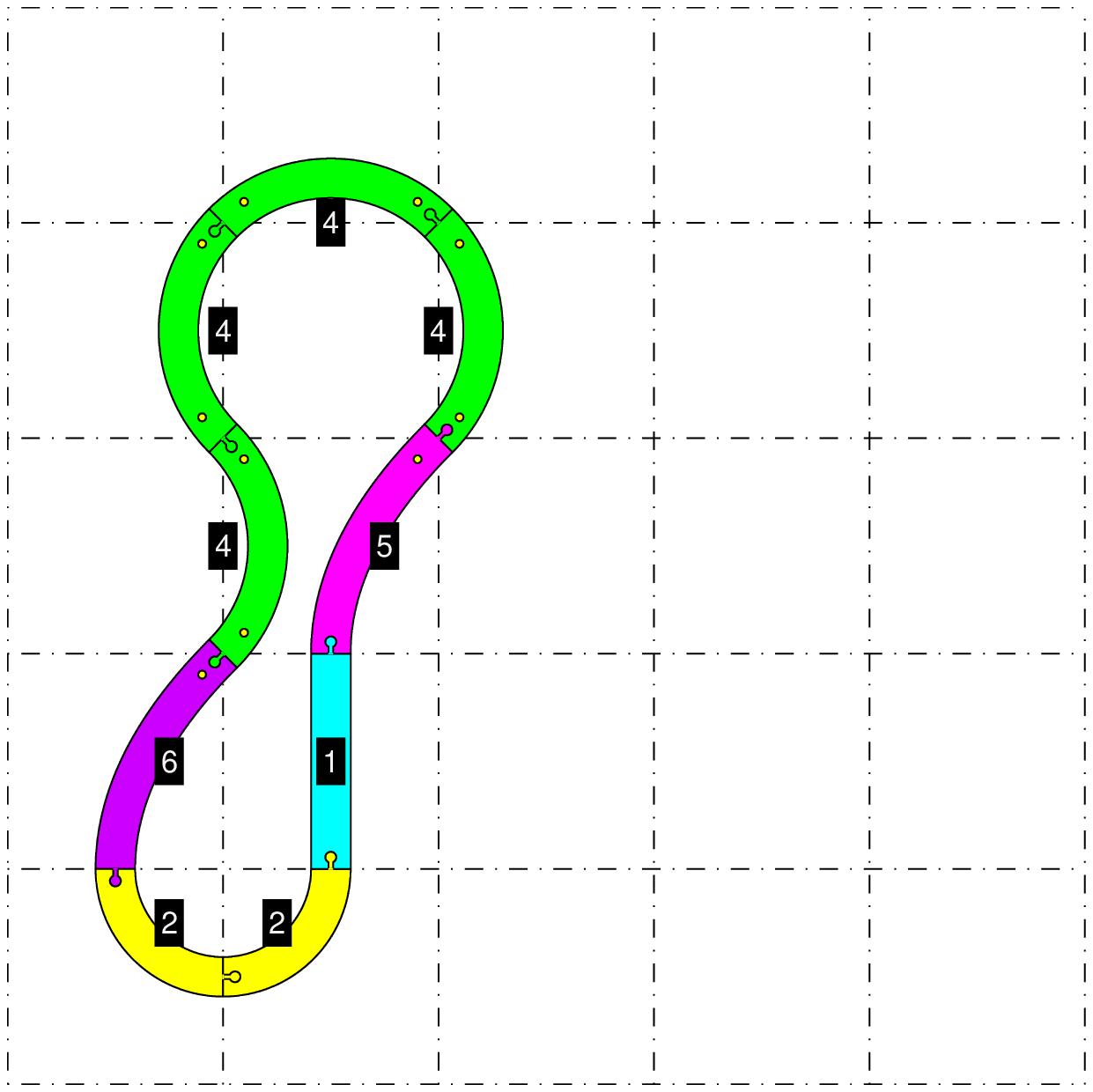, width=5 cm}}
\qquad
%%% sous figure 9
\subfigure[\label{circuits_complets_9_piecmax_ev39}]
{\epsfig{file=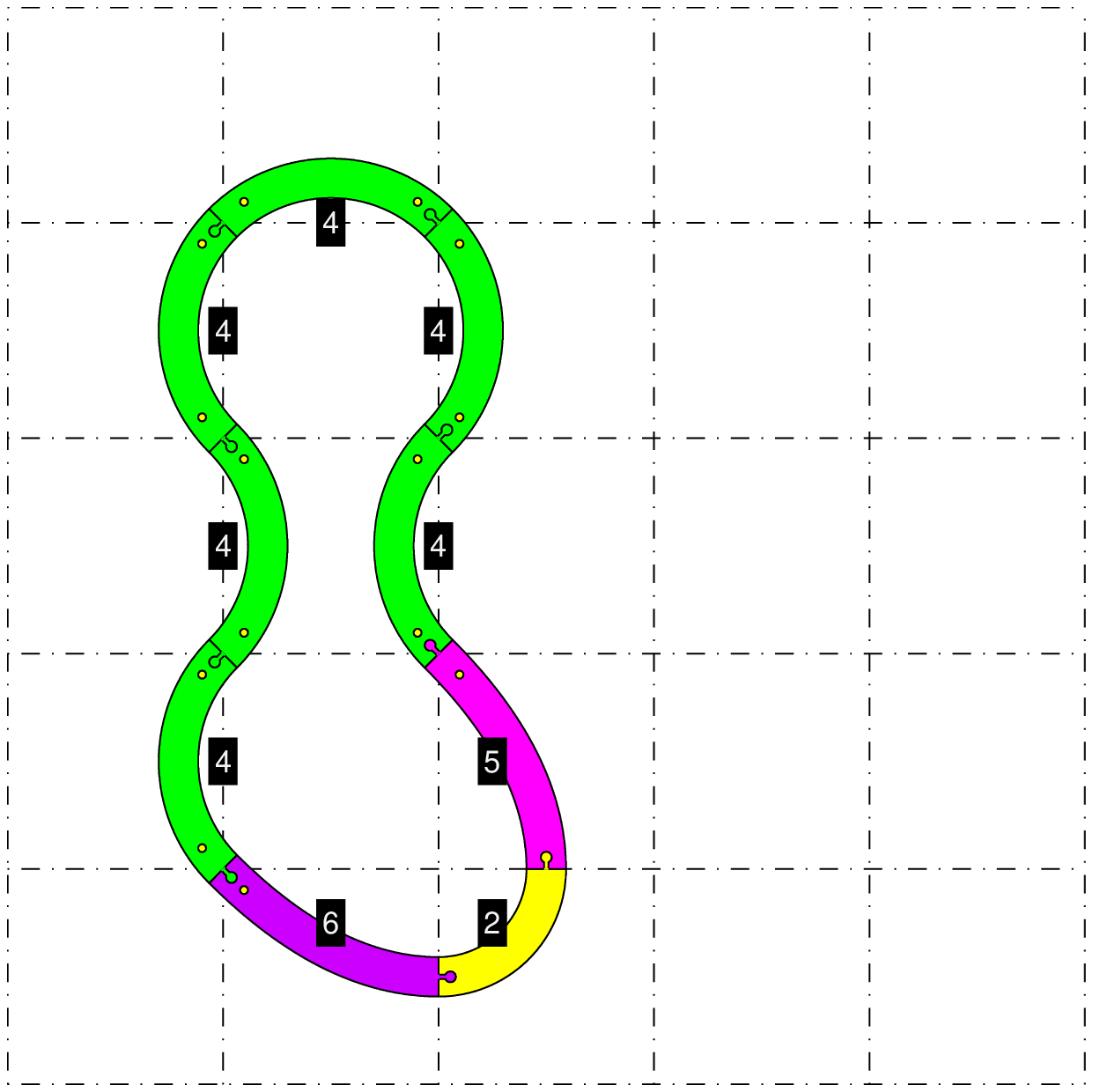, width=5 cm}}
\qquad
%%% sous figure 10
\subfigure[\label{circuits_complets_9_piecmax_ev310}]
{\epsfig{file=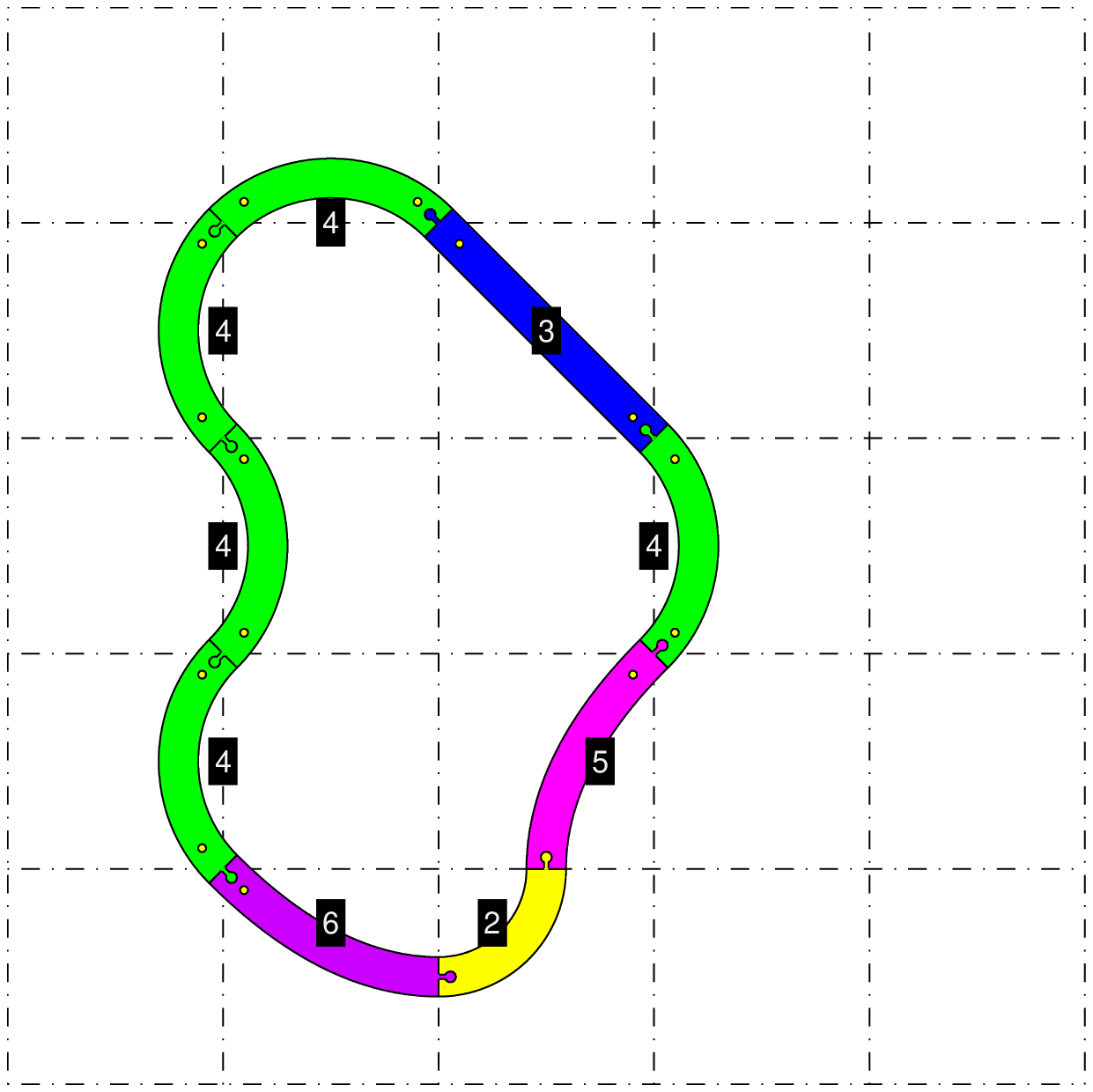, width=5 cm}}
\qquad
\caption{\label{circuits_complets_9_piecmax_ev3}\iflanguage{french}{10 des 74 circuits retenus sur l'ensemble des 1250000 circuits possibles}{10 of the 74 circuits retained from the set of 1250000 possible circuits}.}
\end{figure}
%%%%%%%%%%%%%%%%%%%%%%%%%%%%%%%%%%%%%%%%%%%%%%%%%%%%%%%%%%%%%

\iflanguage{french}{%
Si on trace quelques uns des circuits réalisables %tous les circuits réalisables 
avec $N=9$ pièces
et $N_j=+\infty$, on obtient les  $10$ circuits de la figure \ref{circuits_complets_9_piecmax_ev3}. 
Sur cette figure, seuls les circuits, tous différents à une isométrie près et constructibles ont été tracés.
Pour plusieurs d'entre eux, on peut observer des pièces qui appartiennent au même carré en étant disjointes.%
}{%
If we draw 
some of the feasible circuits
%some of the feasible circuits 
with $N=9$ pieces and $N_j=+\infty$, we obtain the  $10$ circuits in Figure~\ref{circuits_complets_9_piecmax_ev3}. In this figure, only the circuits which are all different up to an isometry have been drawn. For several among them, one may observe some pieces which belong to the same square being disjoint.%
}

%%%%%%%%%%%%%%%%%%%%%%%%%%%%%%%%%%%%%%%%%%%%%%%%%%%%%%%%%%%%%
\end{example}

% Gardé tout de même pour l'article
%\ifcase \choarticle
\begin{example}
\label{examplesimulation525}
%%%%%%%%%%%%%%%%%%%%%%%%%%%%%%%%%%%%%%%%%%%%%%%%%%%%%%%%%%%%%
%\input{simulations_circuit/simulation525}
% fichier tex crée par MaTeXBuild02 le 04-Sep-2015 09:11:52
% à compiler avec 
% MaTeXBuild02('simulation525',0)
% après le fichier 'enumeration_construction_circuit.matex'

%%%%%%%%%%%%%%%%%%%%%%%%%%%%%%%%%%%%%%%%%%%%%%%%%%%%%%%%%%%%%
%\input{./simulations_circuit/circuit_numerique/circuits_complets_9_piecmax_ev4}
% fichier crée par 'presentation_exhaustif_circuit_boucle.m' le 04-Sep-2015 09:12:15
\begin{figure}[h]
\centering
%%% sous figure 1
\subfigure[\label{circuits_complets_9_piecmax_ev41}]
{\epsfig{file=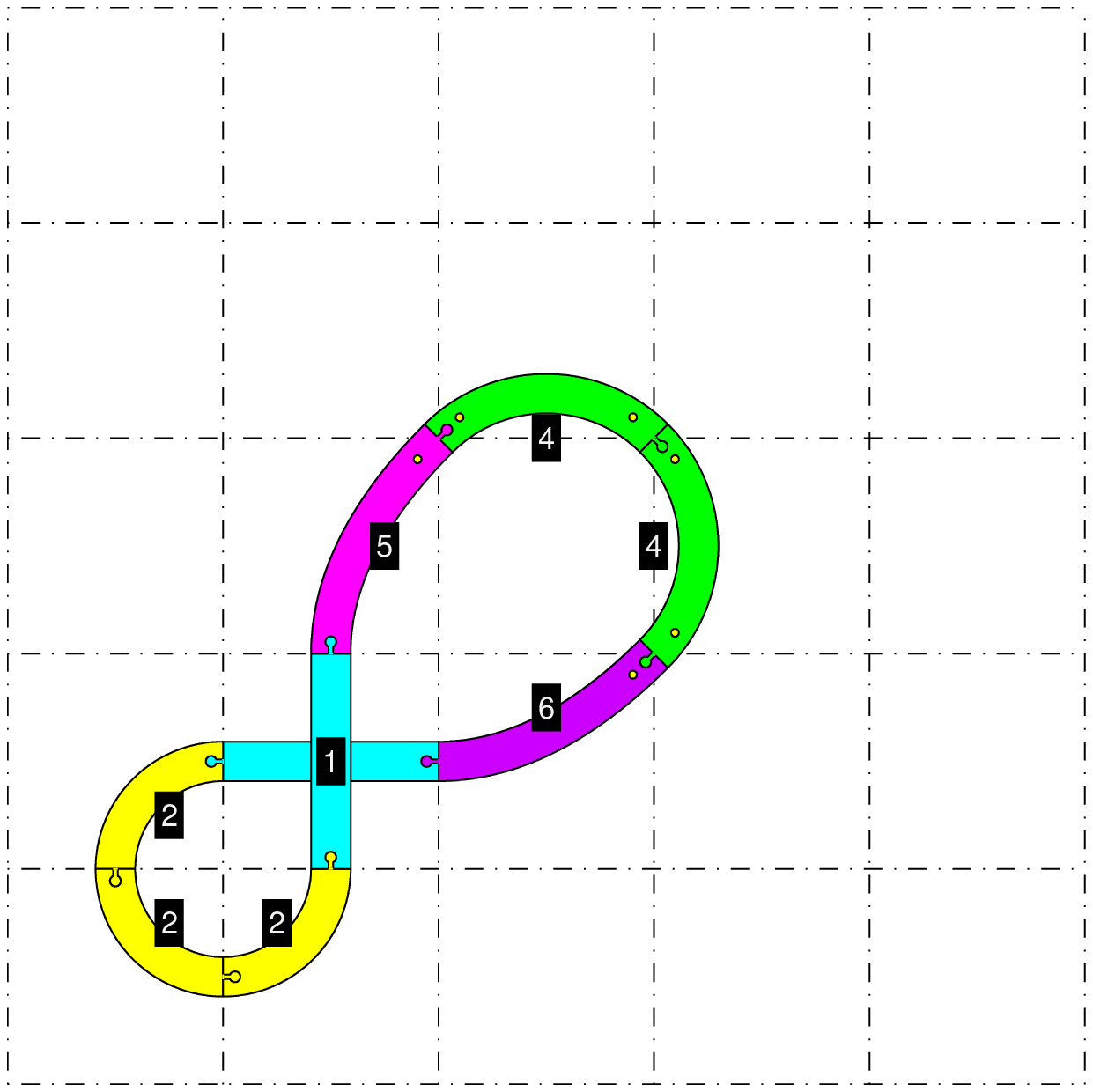, width=5 cm}}
\qquad
%%% sous figure 2
\subfigure[\label{circuits_complets_9_piecmax_ev42}]
{\epsfig{file=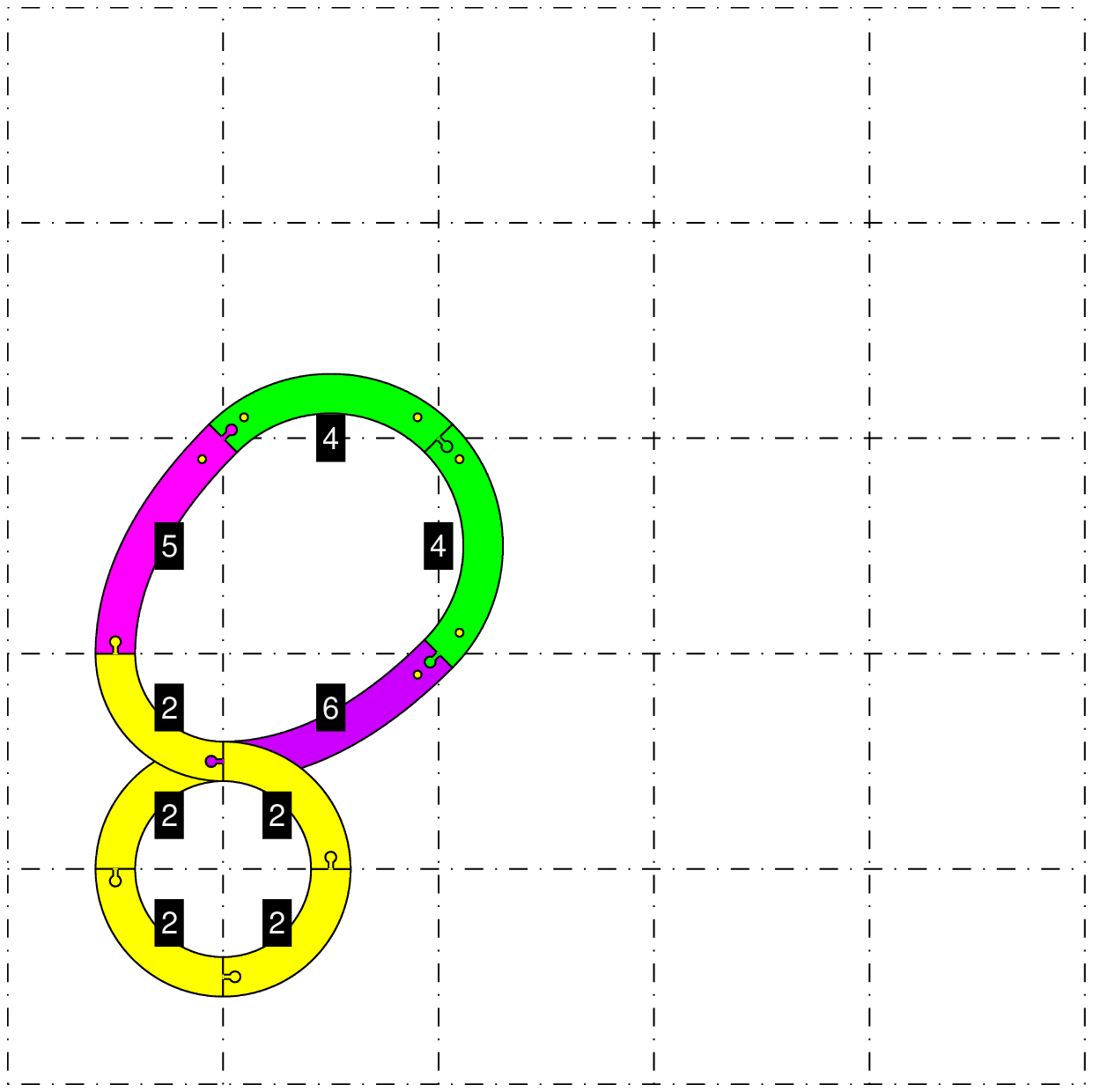, width=5 cm}}
\qquad
%%% sous figure 3
\subfigure[\label{circuits_complets_9_piecmax_ev43}]
{\epsfig{file=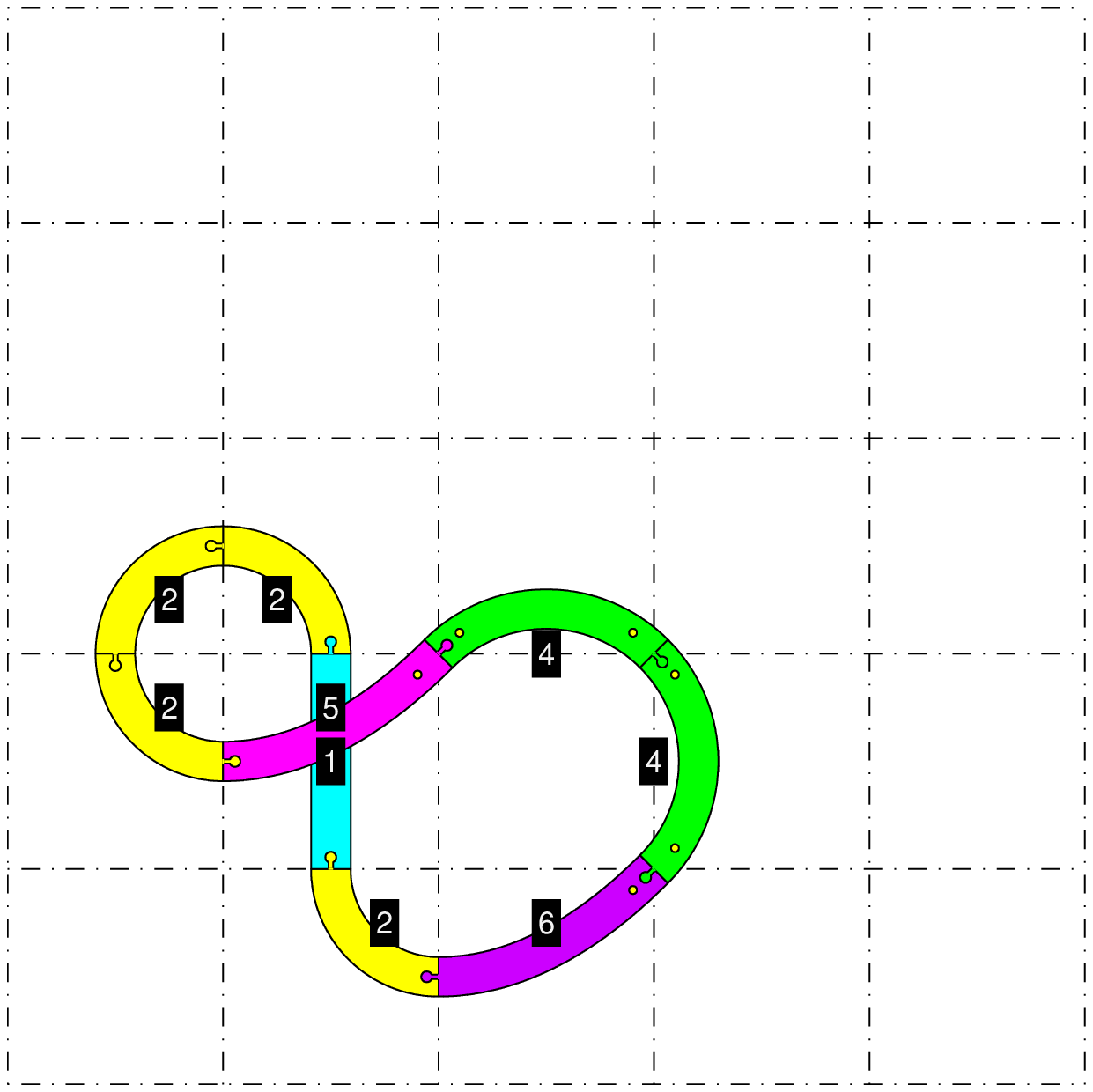, width=5 cm}}
\qquad
%%% sous figure 4
\subfigure[\label{circuits_complets_9_piecmax_ev44}]
{\epsfig{file=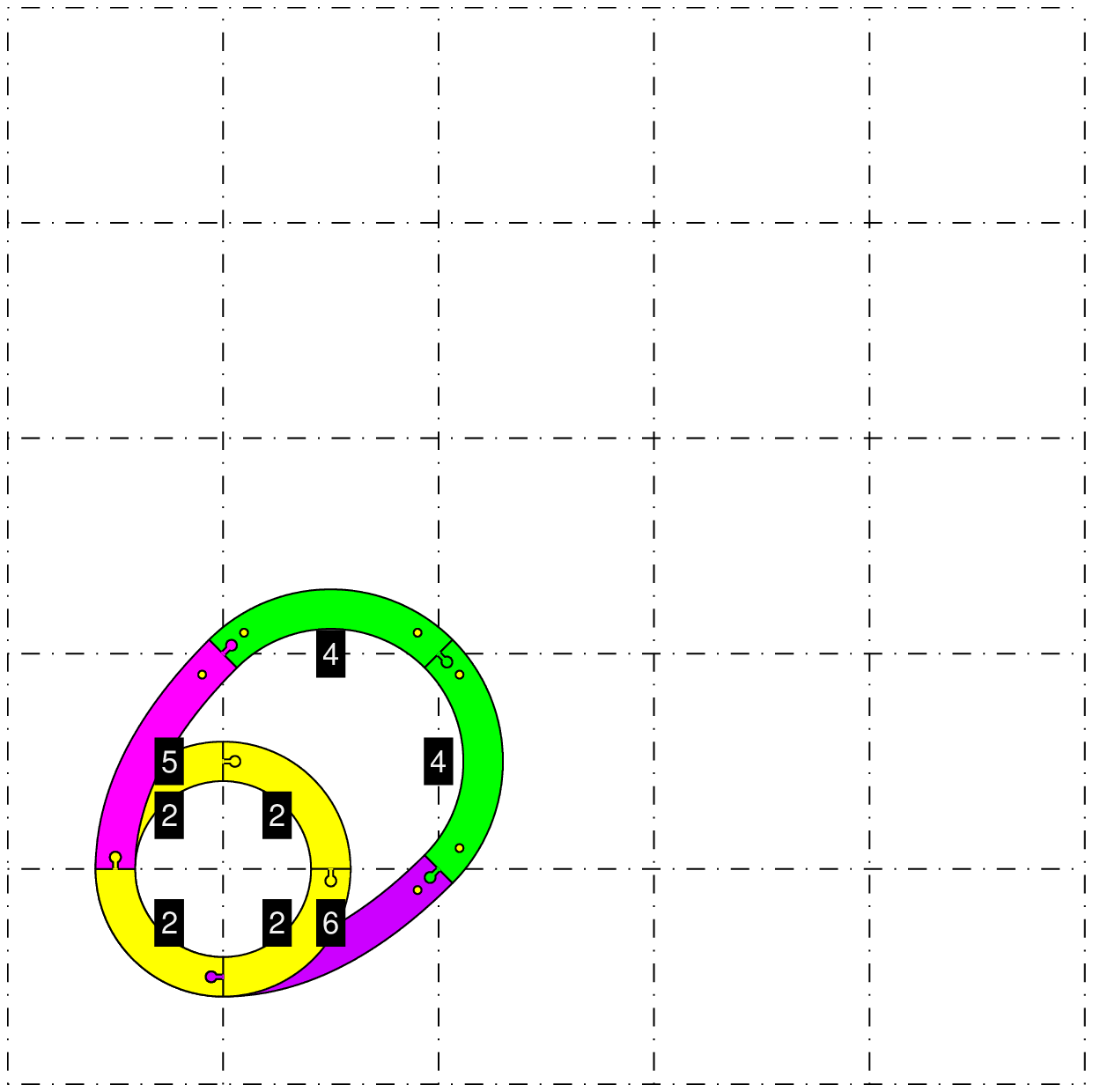, width=5 cm}}
\qquad
%%% sous figure 5
\subfigure[\label{circuits_complets_9_piecmax_ev45}]
{\epsfig{file=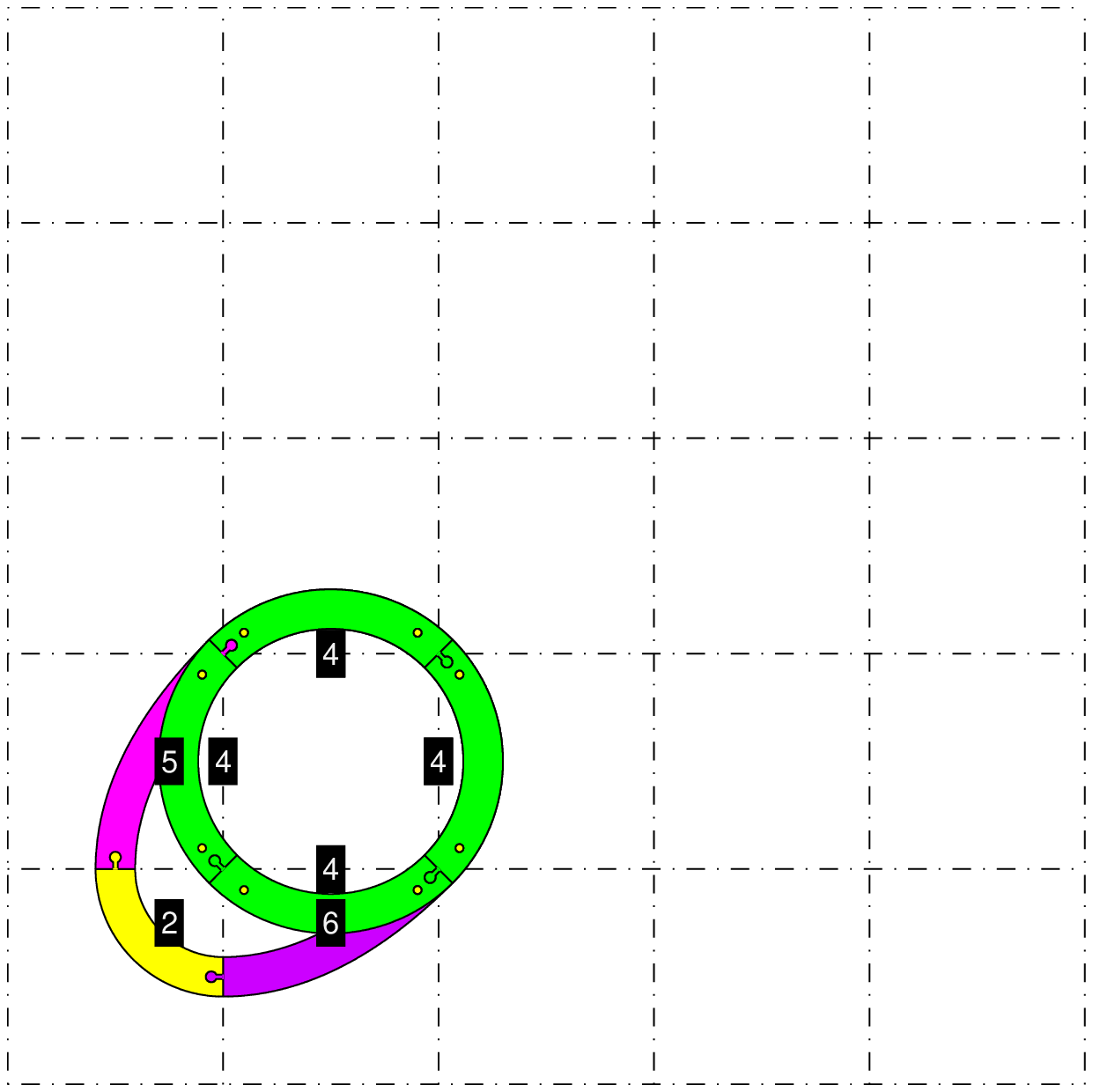, width=5 cm}}
\qquad
%%% sous figure 6
\subfigure[\label{circuits_complets_9_piecmax_ev46}]
{\epsfig{file=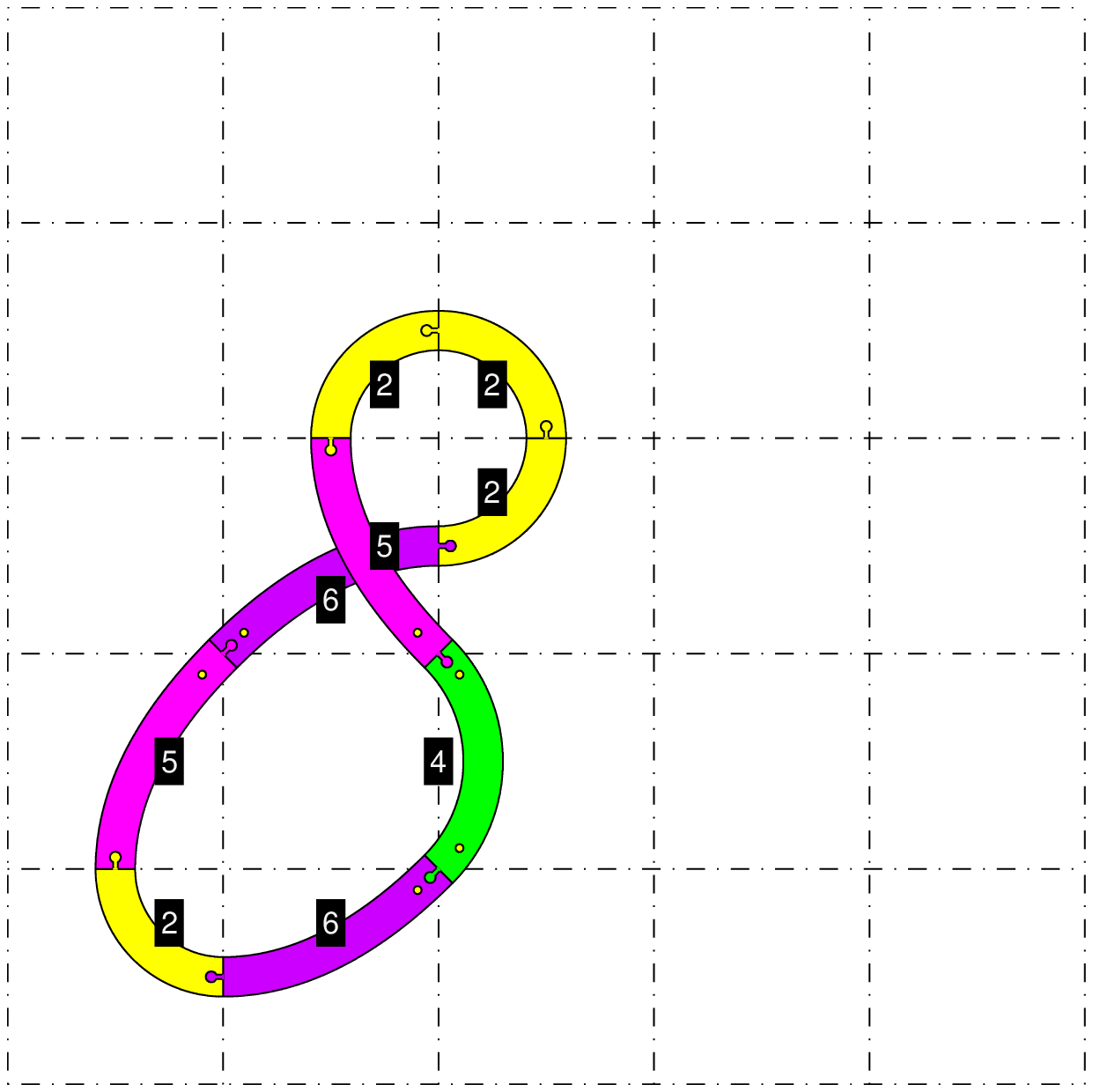, width=5 cm}}
\qquad
%%% sous figure 7
\subfigure[\label{circuits_complets_9_piecmax_ev47}]
{\epsfig{file=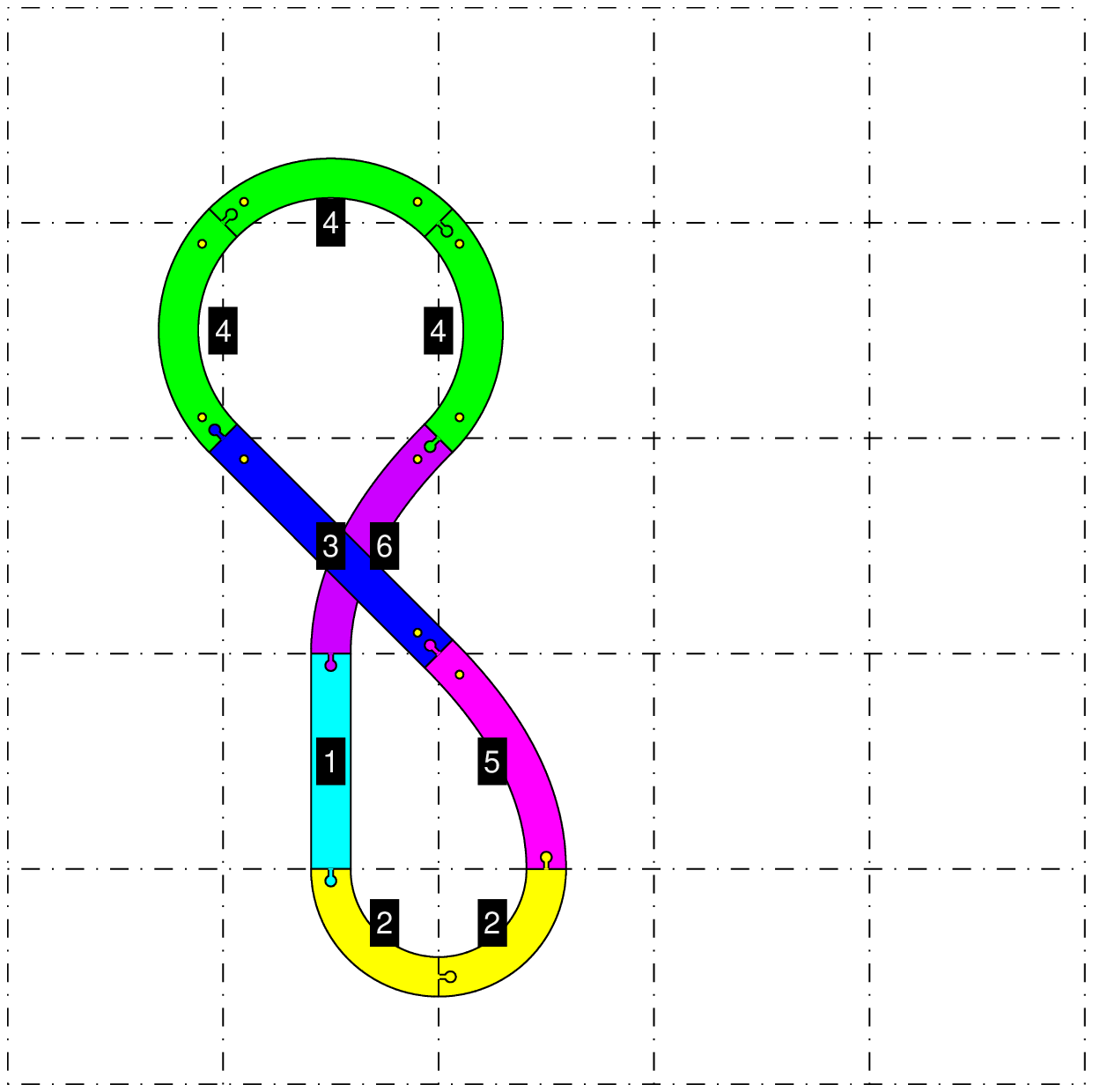, width=5 cm}}
\qquad
%%% sous figure 8
\subfigure[\label{circuits_complets_9_piecmax_ev48}]
{\epsfig{file=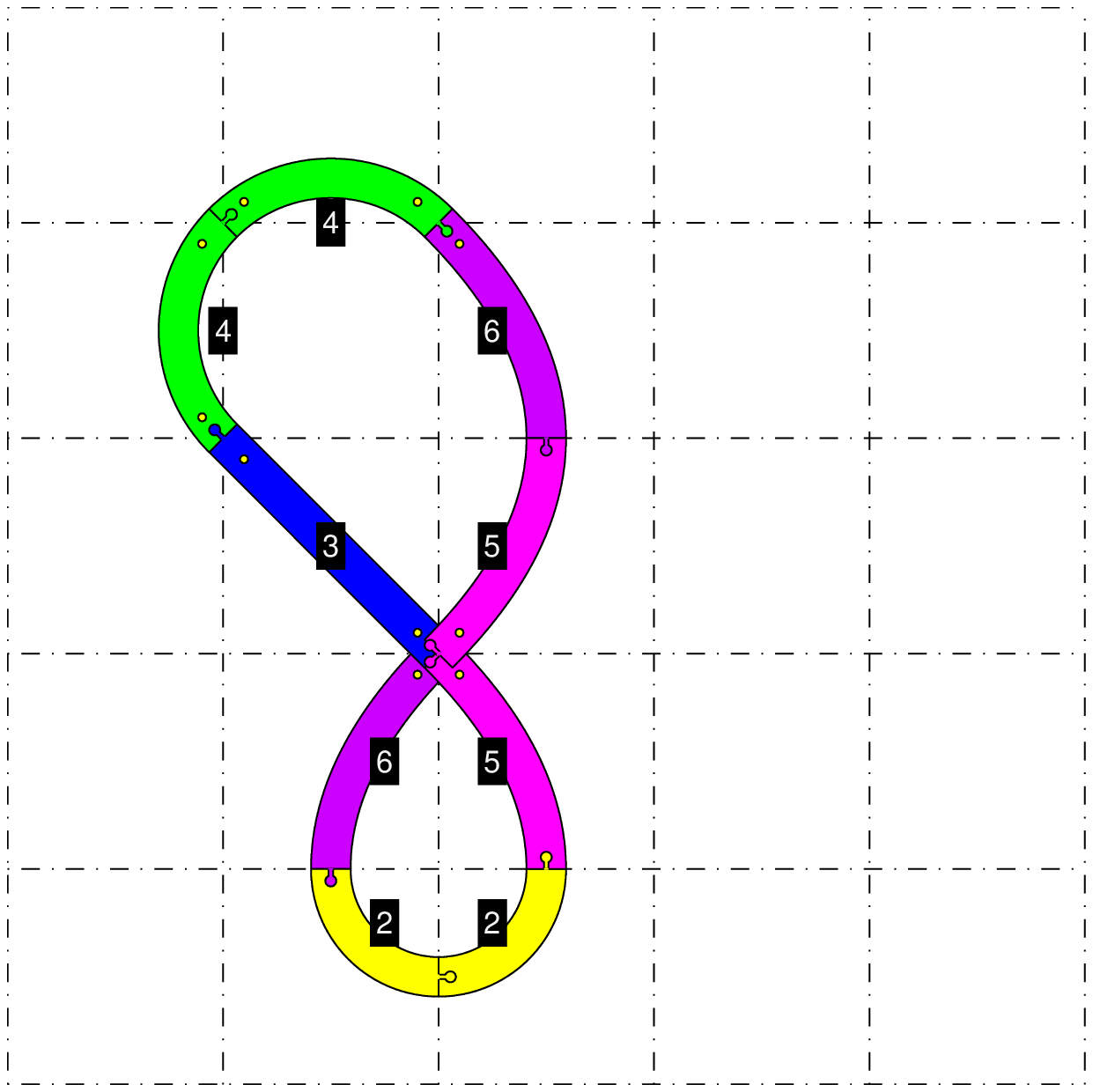, width=5 cm}}
\qquad
%%% sous figure 9
\subfigure[\label{circuits_complets_9_piecmax_ev49}]
{\epsfig{file=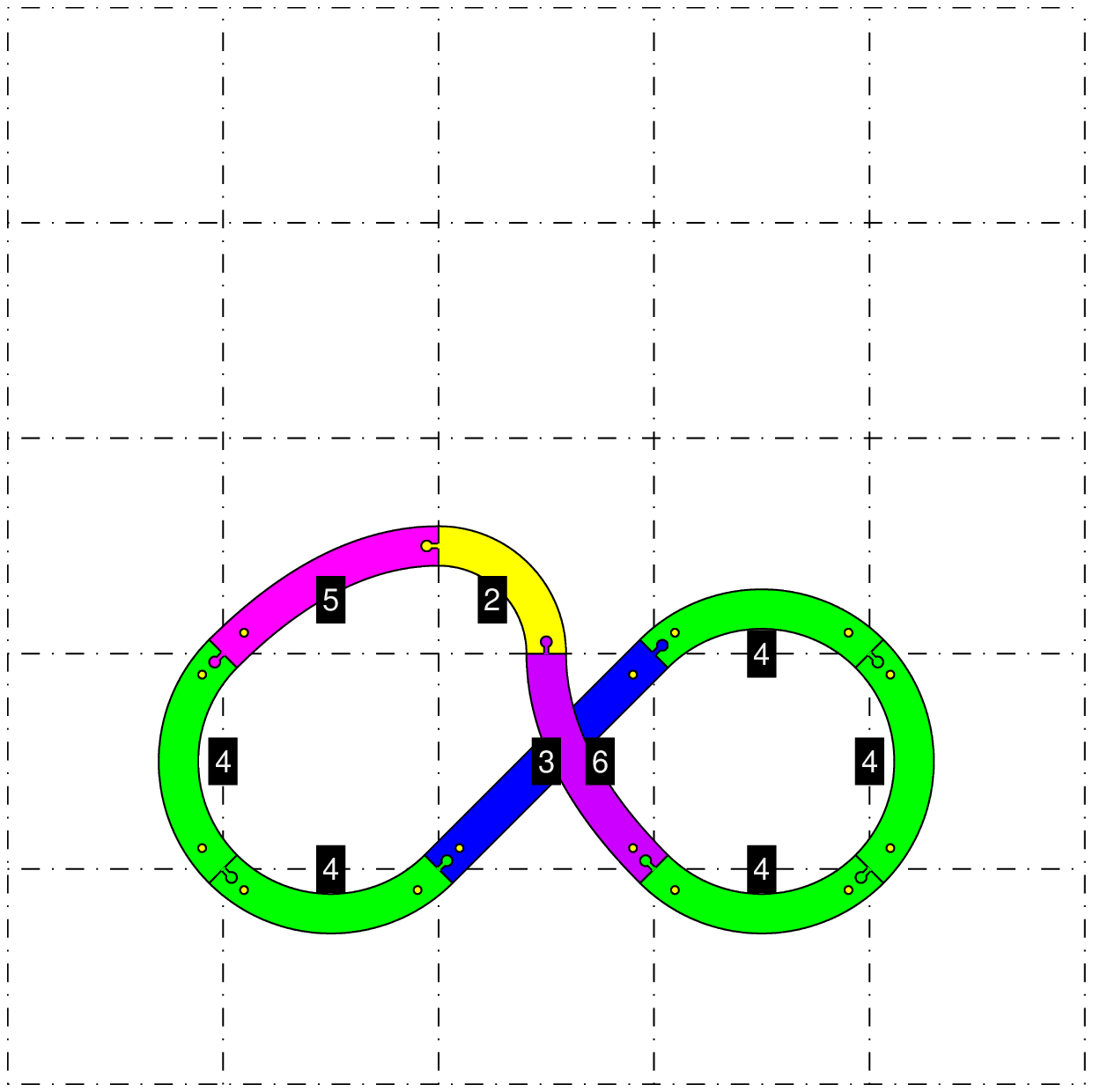, width=5 cm}}
\qquad
%%% sous figure 10
\subfigure[\label{circuits_complets_9_piecmax_ev410}]
{\epsfig{file=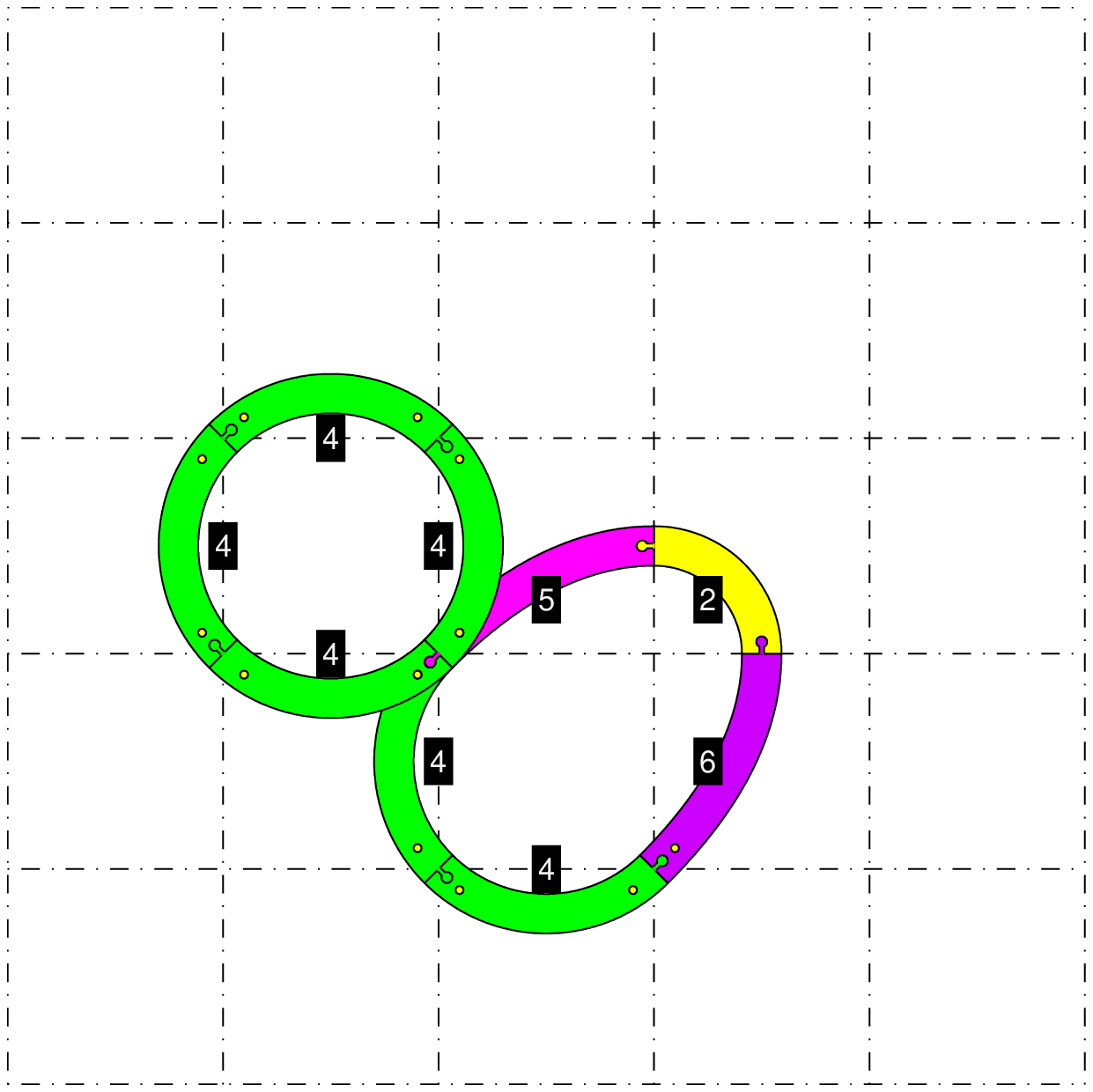, width=5 cm}}
\qquad
\caption{\label{circuits_complets_9_piecmax_ev4}\iflanguage{french}{10 des 18 circuits retenus sur l'ensemble des 1250000 circuits possibles}{10 of the 18 circuits retained from the set of 1250000 possible circuits}.}
\end{figure}

%%%%%%%%%%%%%%%%%%%%%%%%%%%%%%%%%%%%%%%%%%%%%%%%%%%%%%%%%%%%%

\iflanguage{french}{%
Si on trace quelques uns des circuits réalisables %tous les circuits réalisables 
avec $N=9$ pièces
et $N_j=+\infty$, on obtient les  $10$ circuits de la figure \ref{circuits_complets_9_piecmax_ev4}. 
Nous nous sommes amusés  à ne retenir dans cet exemple  que les circuits non constructibles, 
 ce qui reste loisible puisque tous les algorithmes ont été implémentés informatiquement\footnote{Le lecteur averti
constera que certaines figures ont des prises mâles/femelles mal dessinées, puisque le programme de dessin
n'est prévu en fait que pour traiter des circuits constructibles.}.%
}{%
If we draw 
some of the feasible circuits
%some of the feasible circuits 
with  $N=9$ pieces and $N_j=+\infty$, we obtain the $10$ circuits in 
Figure~\ref{circuits_complets_9_piecmax_ev4}. In this example we content ourselves to keep only the non-constructible circuits, which remain permissible since all of the algorithms have been implemented computationally\footnote{The alert reader will notice that certain figures have male/female connectors drawn on them, since the drawing program was in fact conceived to treat only constructible circuits.}.%
}

%%%%%%%%%%%%%%%%%%%%%%%%%%%%%%%%%%%%%%%%%%%%%%%%%%%%%%%%%%%%%
\end{example}
%\or
%\fi

\begin{example}
\label{examplesimulation530}
%%%%%%%%%%%%%%%%%%%%%%%%%%%%%%%%%%%%%%%%%%%%%%%%%%%%%%%%%%%%%
%\input{simulations_circuit/simulation530}
% fichier tex crée par MaTeXBuild02 le 04-Sep-2015 09:59:45
% à compiler avec 
% MaTeXBuild02('simulation530',0)
% après le fichier 'enumeration_construction_circuit.matex'
% Attention, longue compilation

%%%%%%%%%%%%%%%%%%%%%%%%%%%%%%%%%%%%%%%%%%%%%%%%%%%%%%%%%%%%%
%\input{./simulations_circuit/circuit_numerique/circuits_complets_11_piecmax_ev3}
% fichier crée par 'presentation_exhaustif_circuit_boucle.m' le 04-Sep-2015 14:37:10
\begin{figure}[h]
\centering
%%% sous figure 1
\subfigure[\label{circuits_complets_11_piecmax_ev31}]
{\epsfig{file=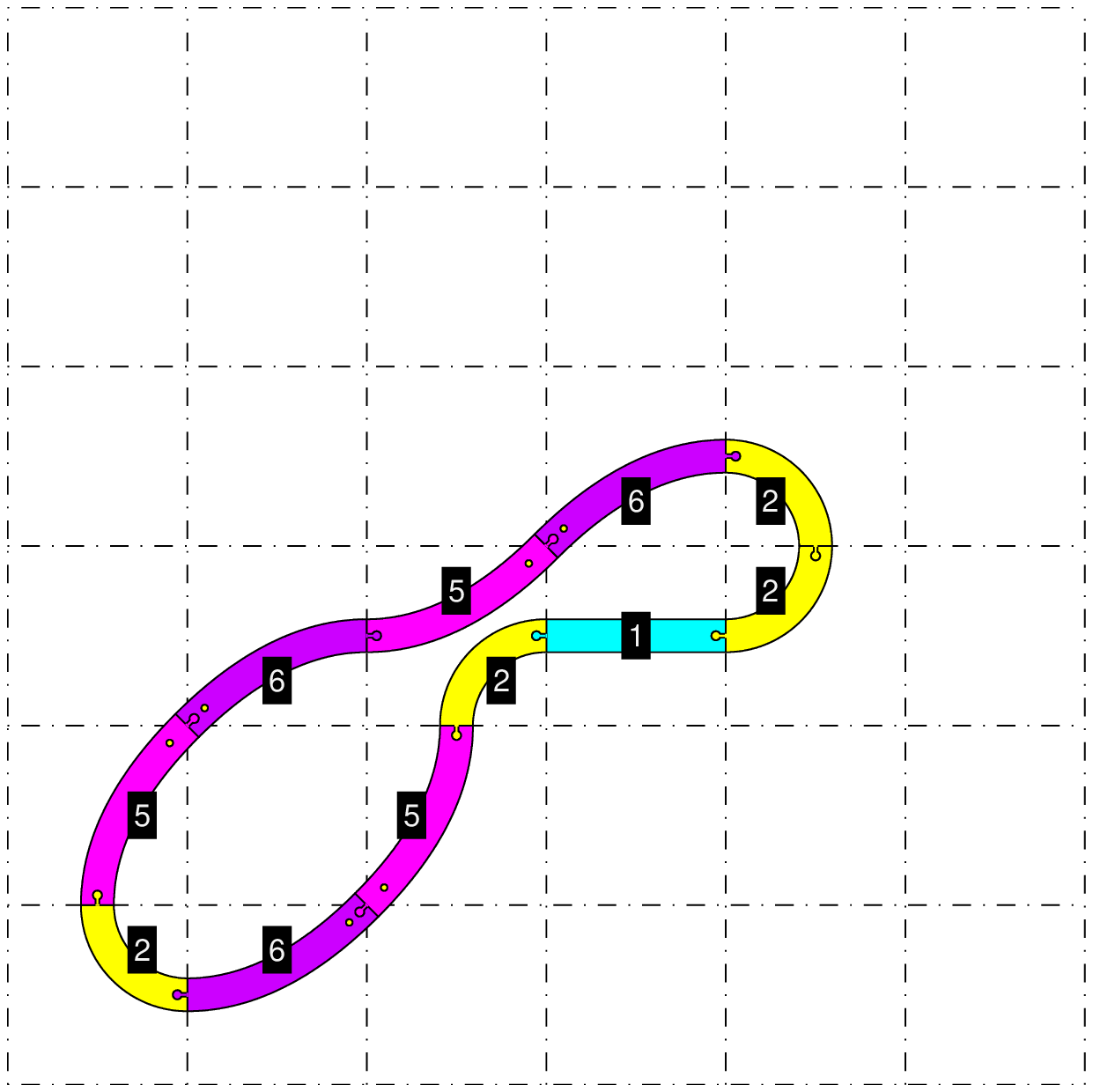, width=5 cm}}
\qquad
%%% sous figure 2
\subfigure[\label{circuits_complets_11_piecmax_ev32}]
{\epsfig{file=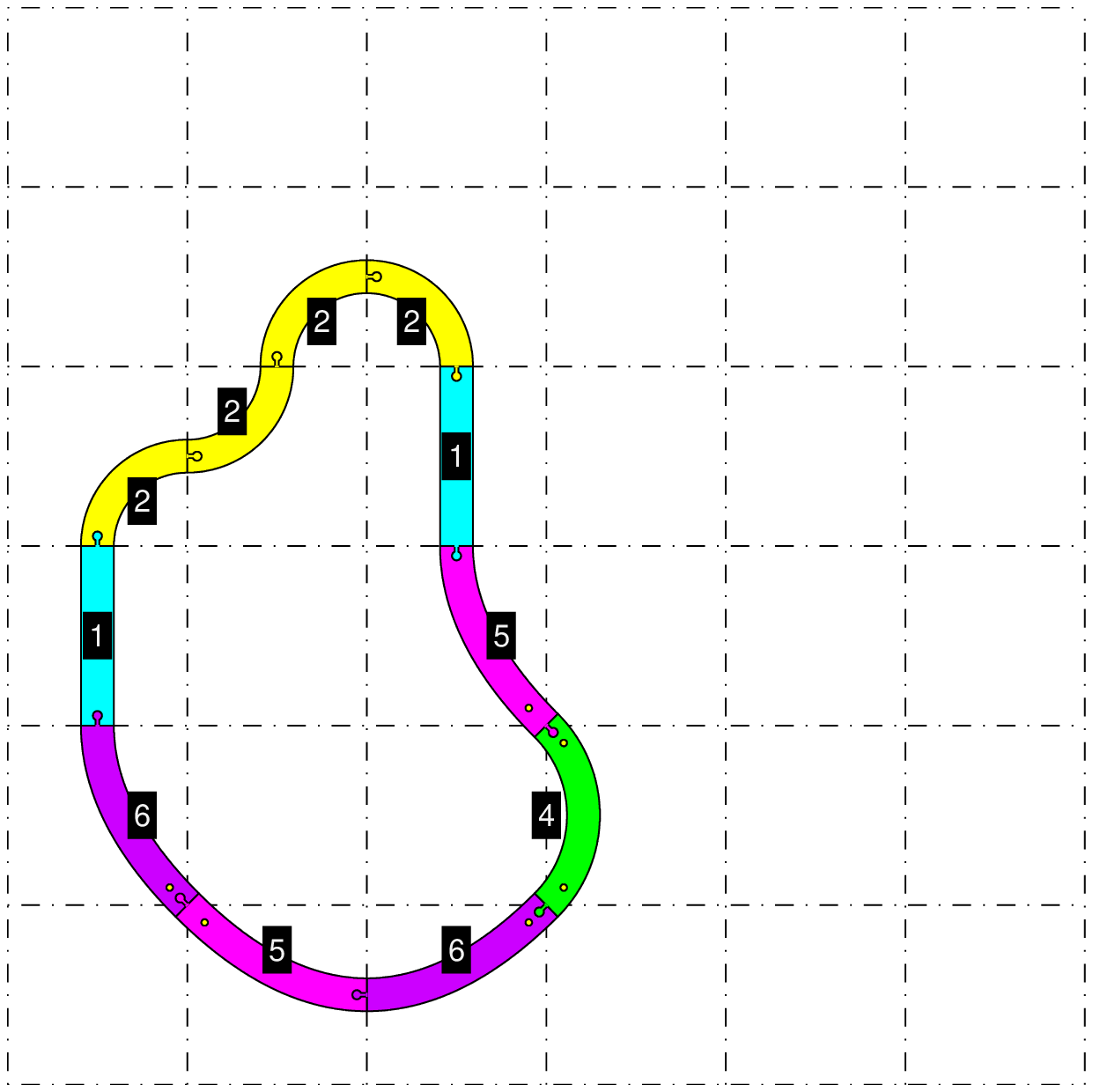, width=5 cm}}
\qquad
%%% sous figure 3
\subfigure[\label{circuits_complets_11_piecmax_ev33}]
{\epsfig{file=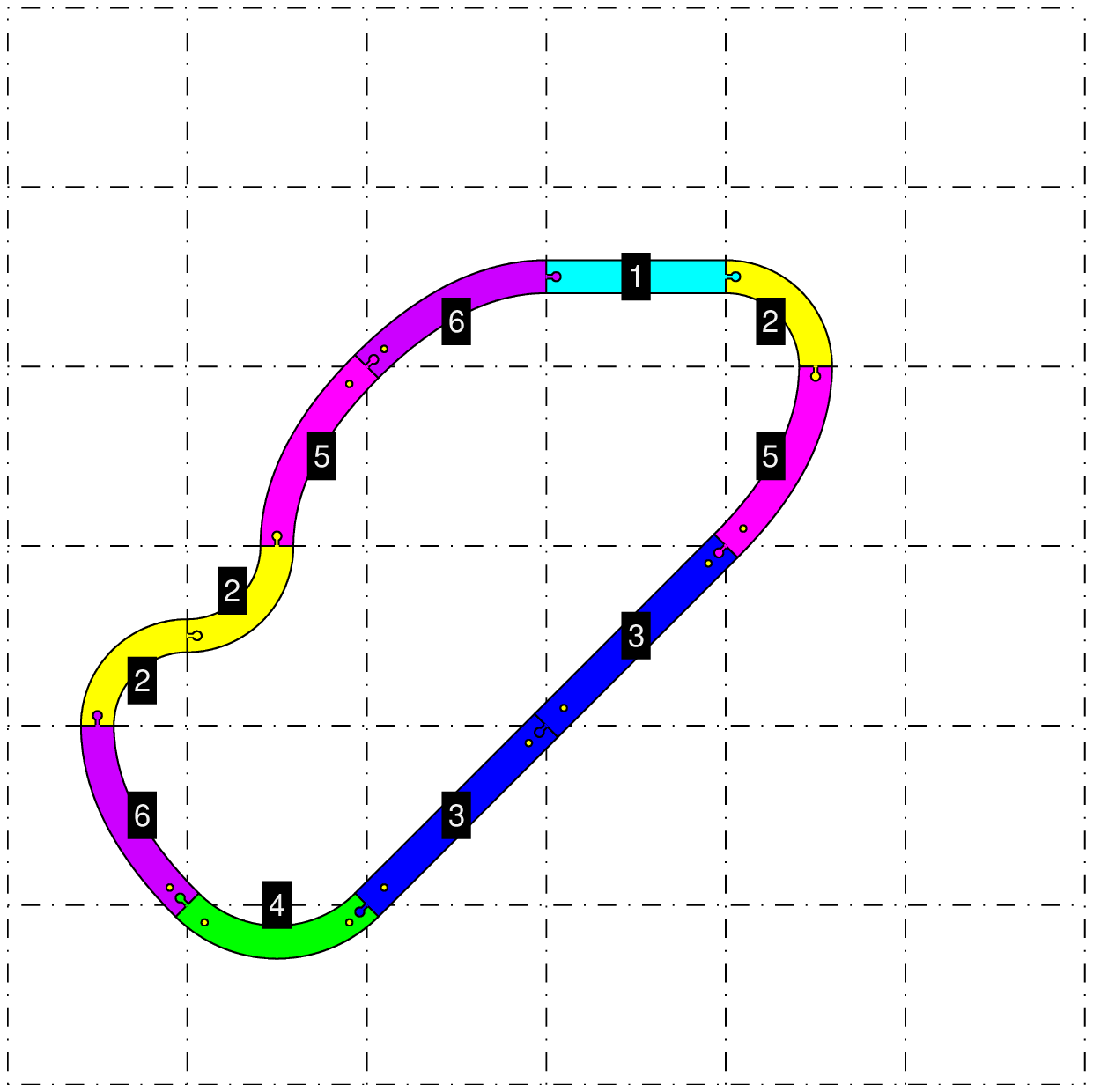, width=5 cm}}
\qquad
%%% sous figure 4
\subfigure[\label{circuits_complets_11_piecmax_ev34}]
{\epsfig{file=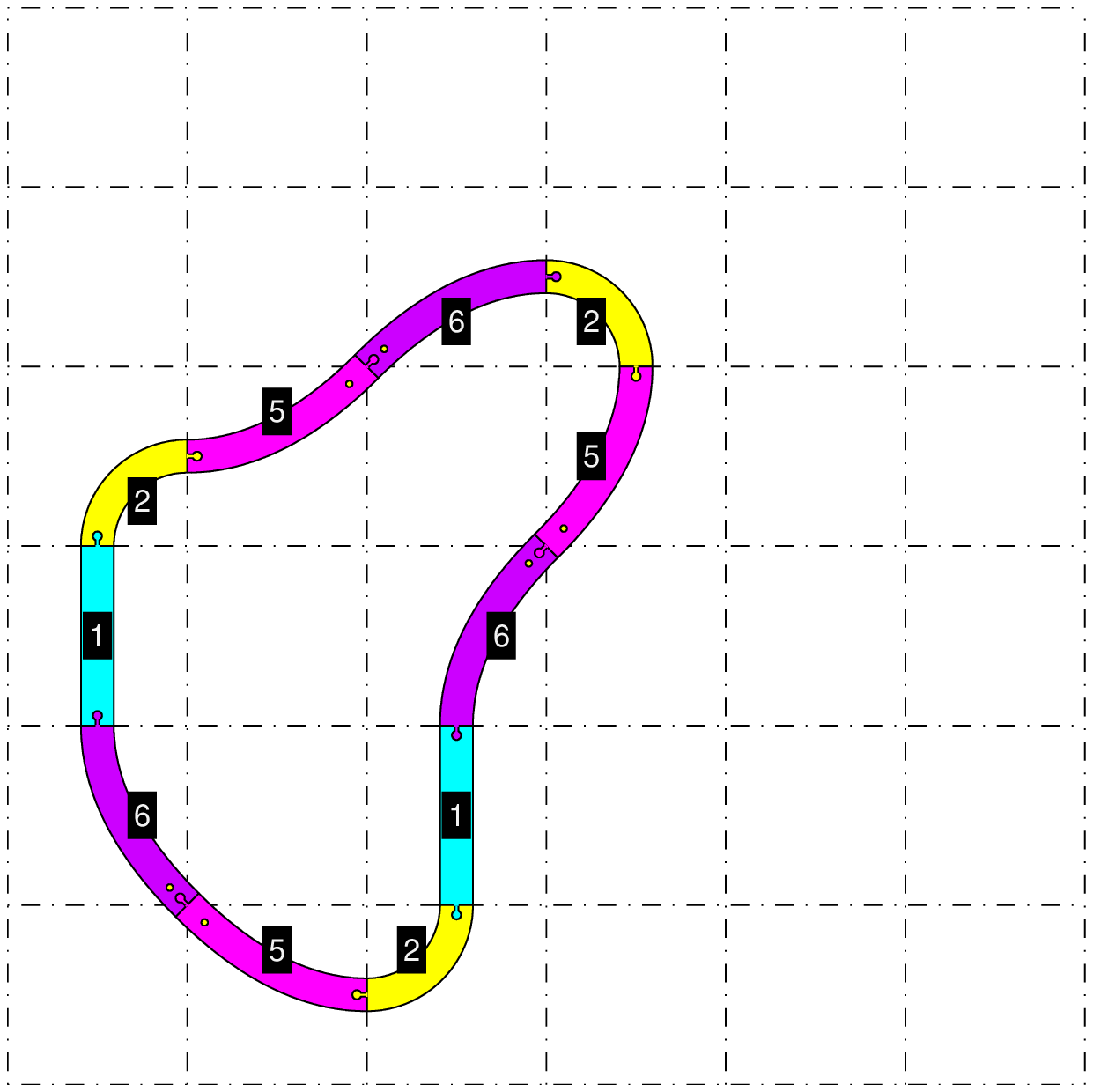, width=5 cm}}
\qquad
%%% sous figure 5
\subfigure[\label{circuits_complets_11_piecmax_ev35}]
{\epsfig{file=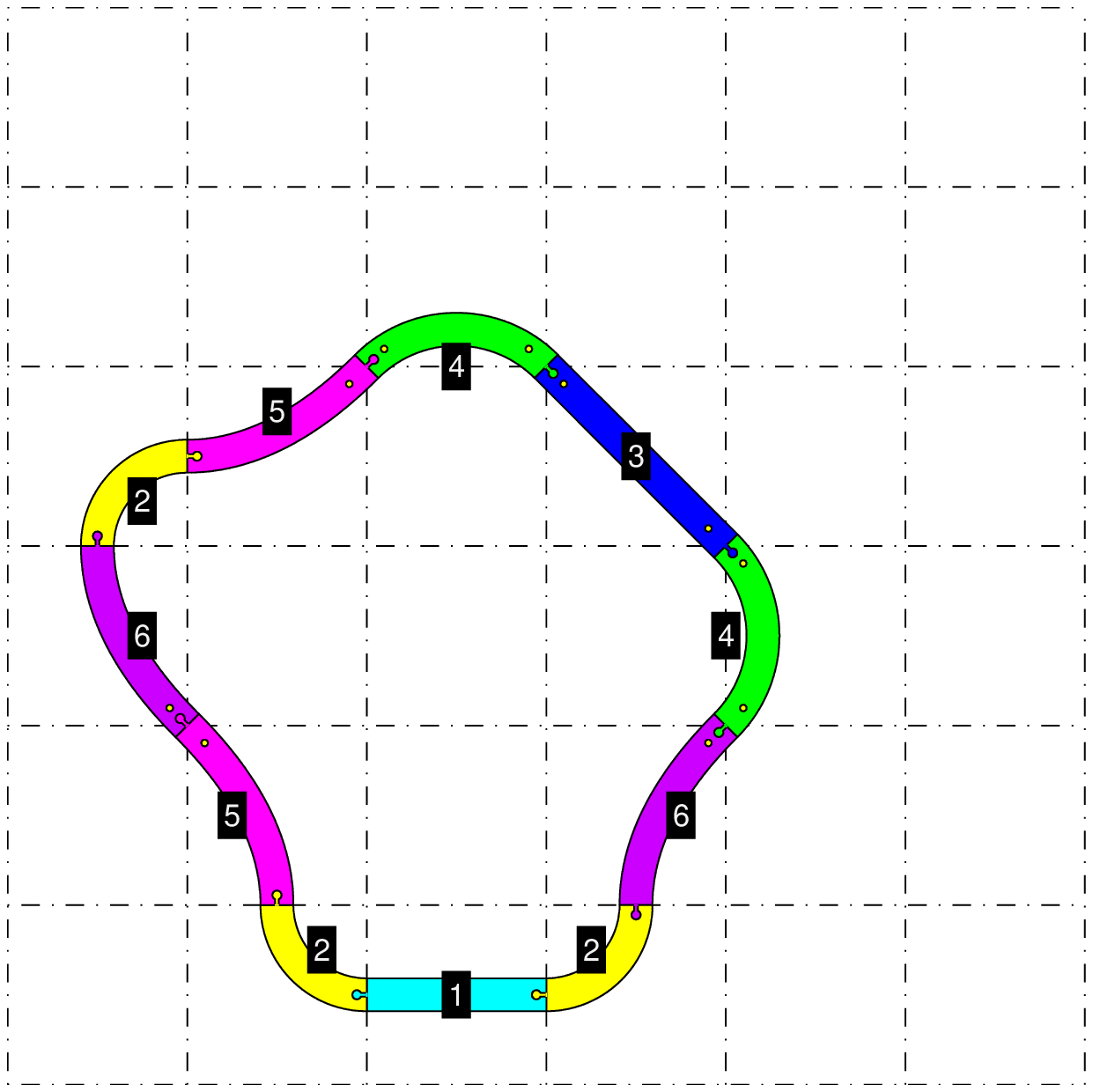, width=5 cm}}
\qquad
%%% sous figure 6
\subfigure[\label{circuits_complets_11_piecmax_ev36}]
{\epsfig{file=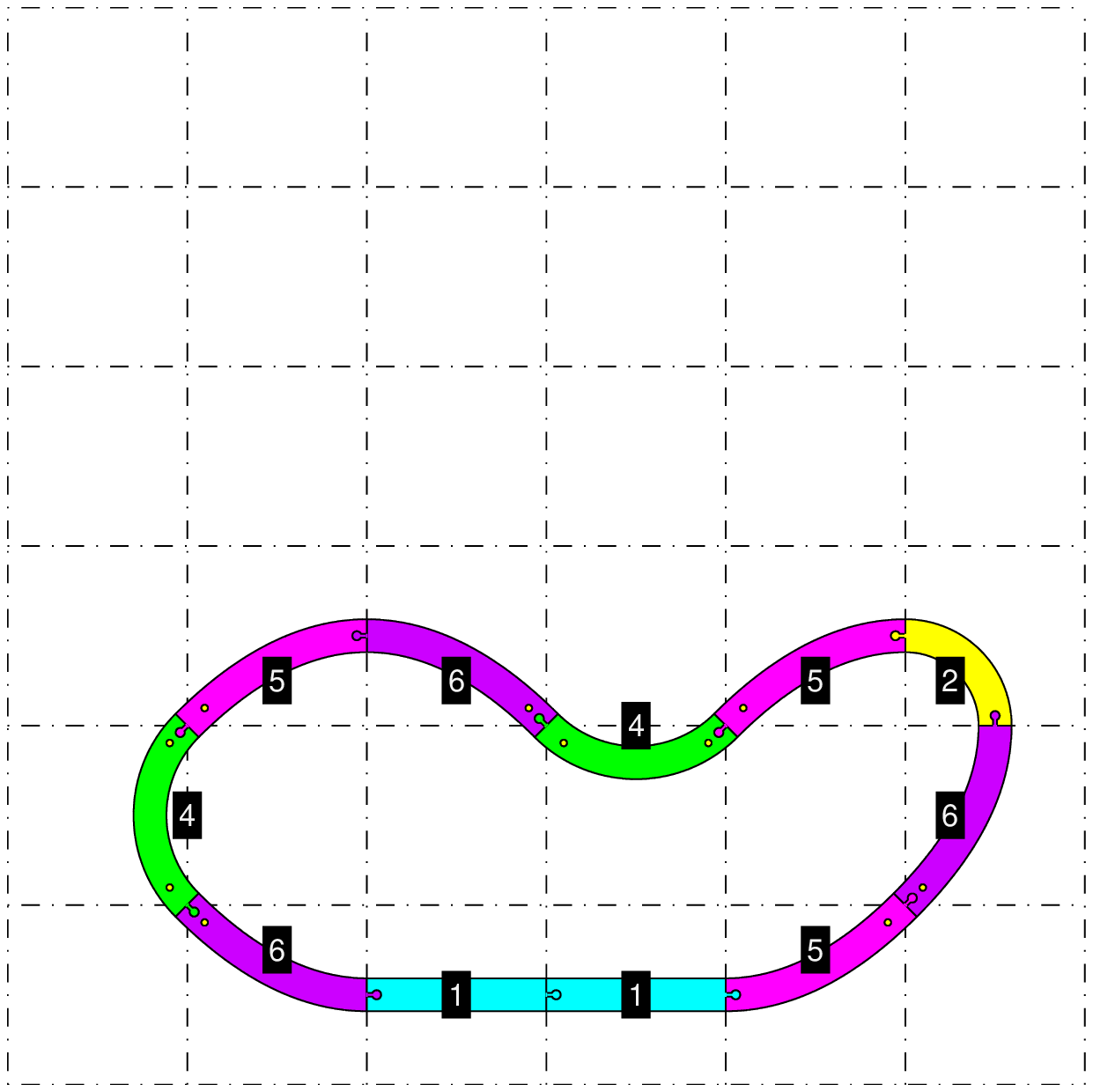, width=5 cm}}
\qquad
%%% sous figure 7
\subfigure[\label{circuits_complets_11_piecmax_ev37}]
{\epsfig{file=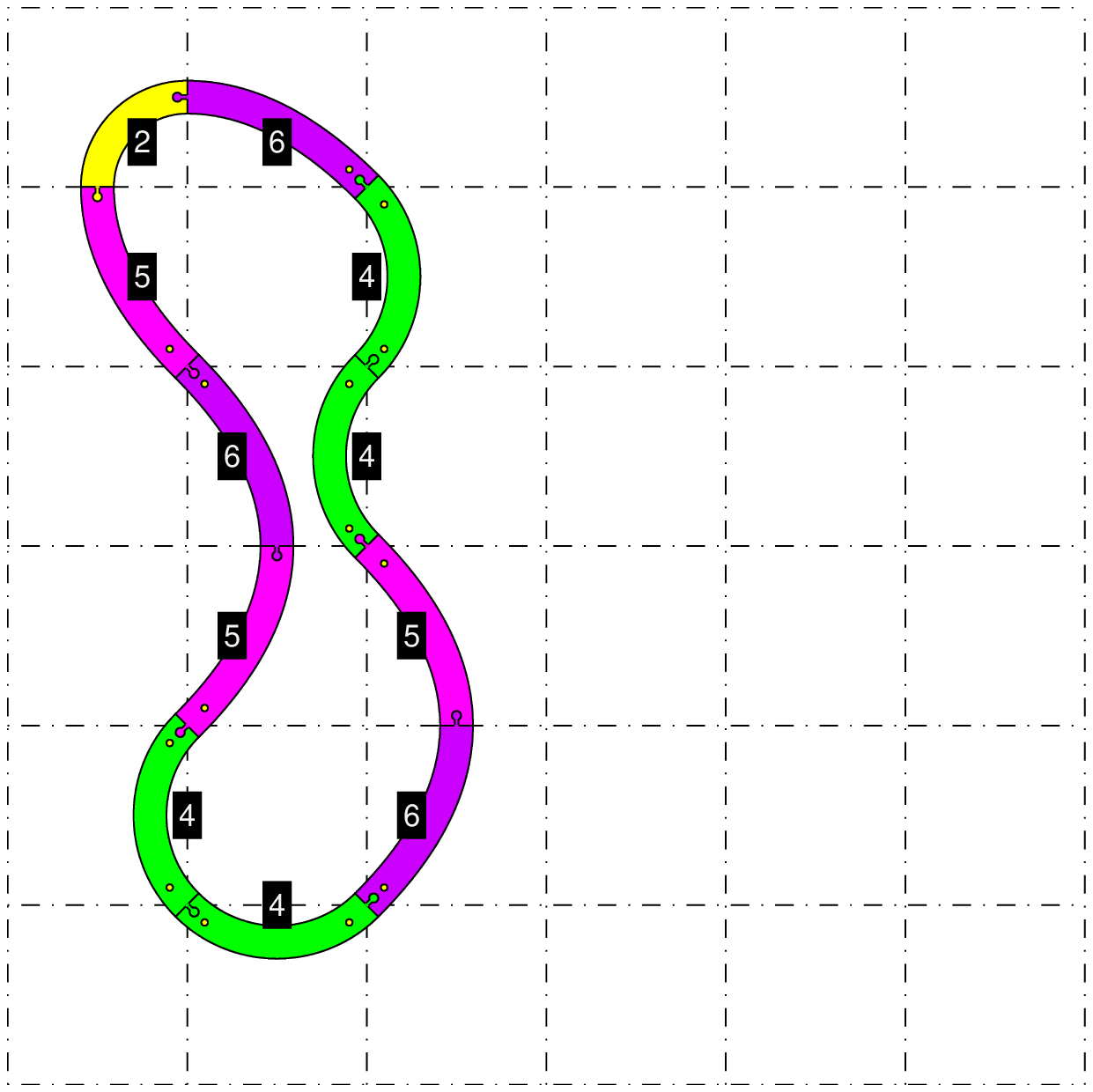, width=5 cm}}
\qquad
%%% sous figure 8
\subfigure[\label{circuits_complets_11_piecmax_ev38}]
{\epsfig{file=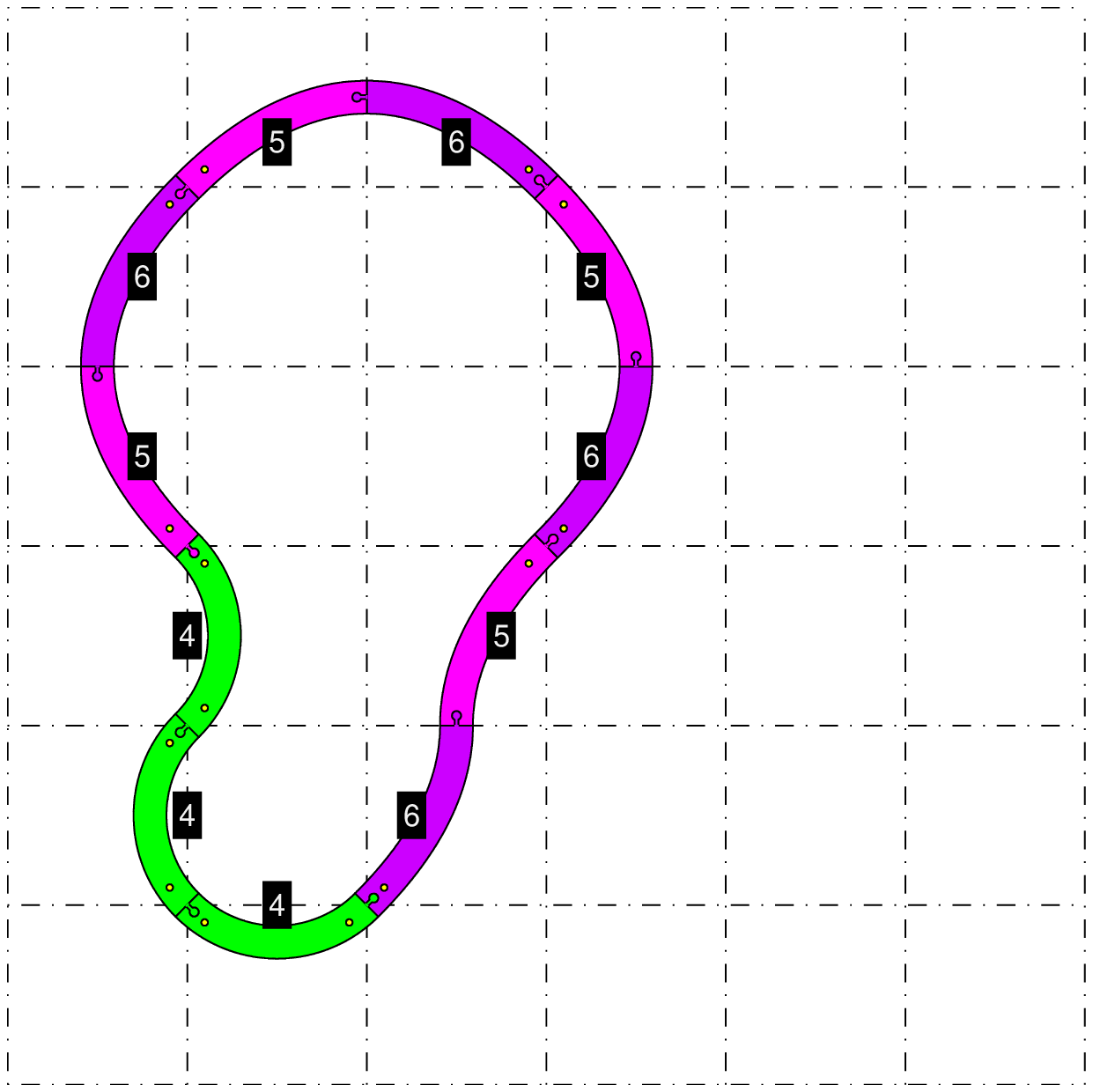, width=5 cm}}
\qquad
%%% sous figure 9
\subfigure[\label{circuits_complets_11_piecmax_ev39}]
{\epsfig{file=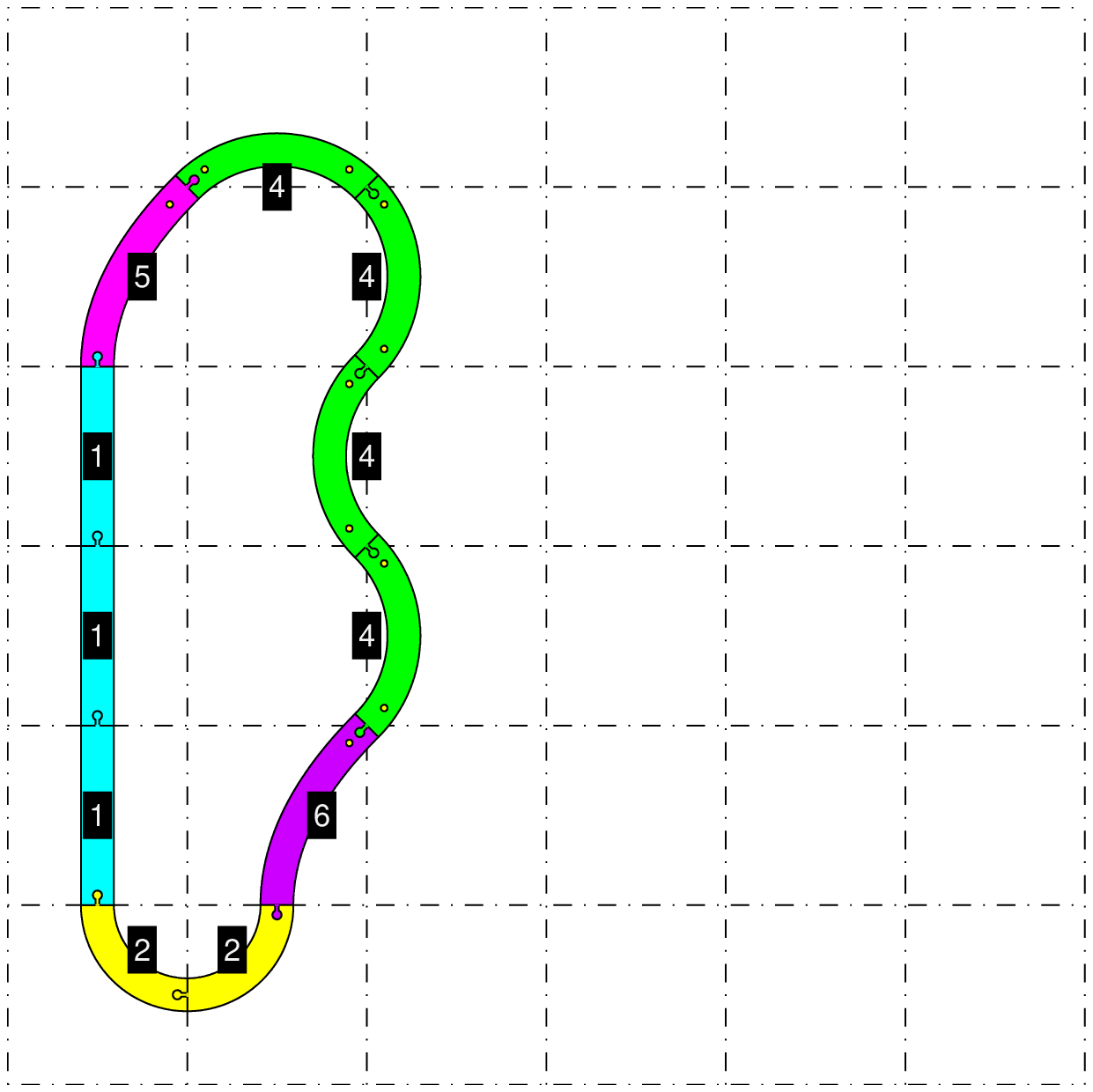, width=5 cm}}
\qquad
%%% sous figure 10
\subfigure[\label{circuits_complets_11_piecmax_ev310}]
{\epsfig{file=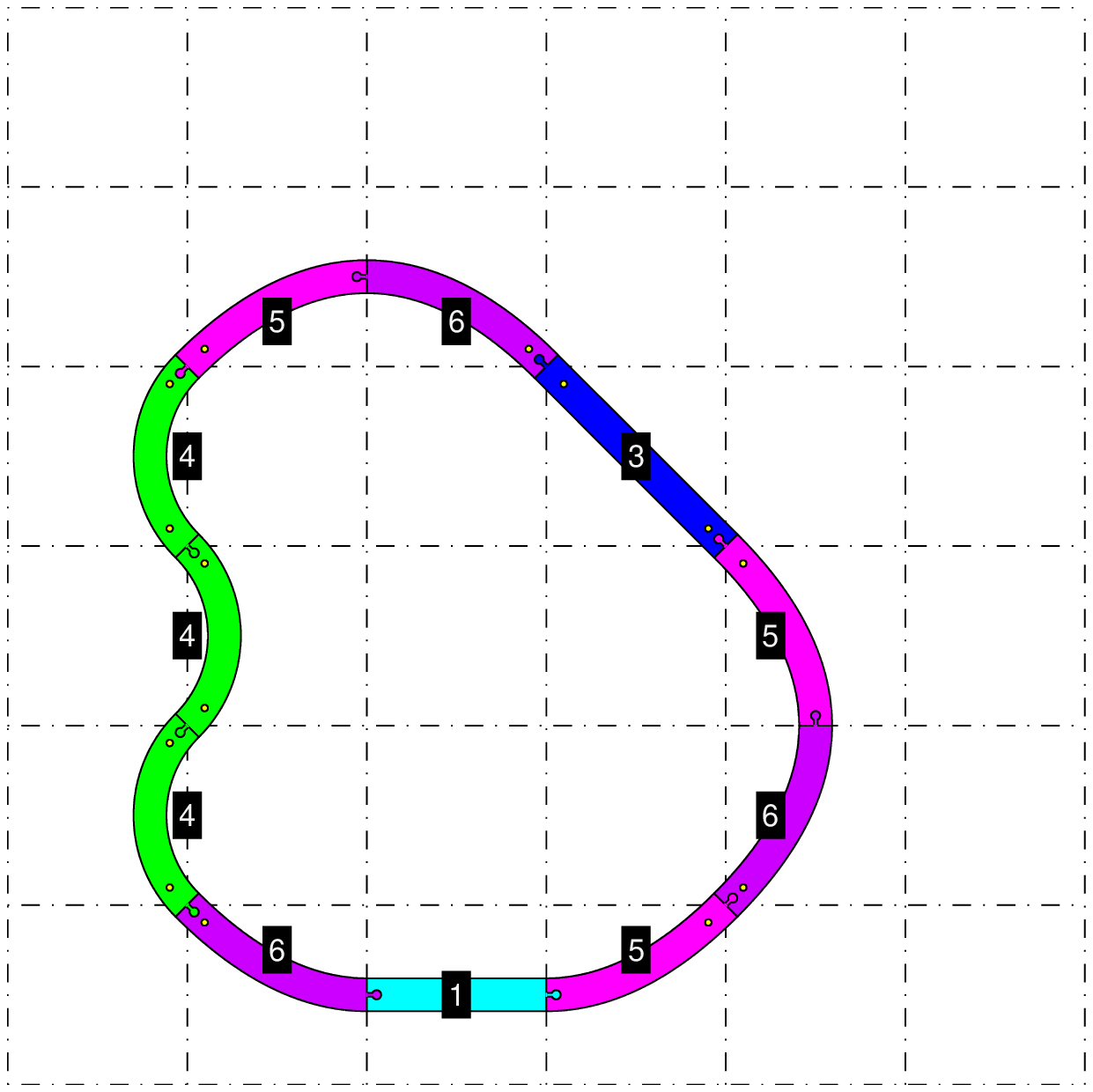, width=5 cm}}
\qquad
\caption{\label{circuits_complets_11_piecmax_ev3}\iflanguage{french}{10 des 753 circuits retenus sur l'ensemble des 31250000 circuits possibles}{10 of the 753 circuits retained from the set of 31250000 possible circuits}.}
\end{figure}
%%%%%%%%%%%%%%%%%%%%%%%%%%%%%%%%%%%%%%%%%%%%%%%%%%%%%%%%%%%%%

\iflanguage{french}{%
Concluons par le calcul correspondant à la valeur maximale de $N$, calculable raison\-na\-ble\-ment sur le plan informatique.

Si on trace quelques uns des circuits réalisables %tous les circuits réalisables 
avec $N=11$ pièces
et $N_j=4$, ce qui correspond aux boîtes de jeux, distribuées par \textit{Easyloop},
on obtient les  $10$ circuits de la figure \ref{circuits_complets_11_piecmax_ev3}.%
}{%
We conclude with the calculation corresponding to the maximal value of $N$ which can be reasonably calculated computationally.

If we draw 
 some of the feasible circuits
%some of the feasible circuits 
with $N=11$ pieces and $N_j=4$, which corresponds to the toy-box distributed by \textit{Easyloop}, we obtain the 
$10$ circuits in Figure~\ref{circuits_complets_11_piecmax_ev3}.% 
}

%Sur cette figure, nous n'avons conservé que des circuits contenant au moins une pièce de chaque type pour en 
%montrer les aspects variés.
%ne marche pas à cause du message suivant :
%100% du calcul en cours fait
%??? Error using ==> print at 315
%Need a handle to a Figure object.
%
%Error in ==> dessine_circuit02_1p1 at 324
%        print('-depsc',[nomfig,'.eps']);
%
%Error in ==> dessine_multi_circuitbis at 59
%            dessine_circuit02_1p1(ctot(:,1:nombpiec(i),i),ptot(:,1:nombpiec(i),i),(itot(1:nombpiec(i),i)).',(rottot(1:nombpiec(i),i)).',(symetot(1:nombpiec(i),i)).',...%
%
%Error in ==> presentation_exhaustif_circuit_boucle at 129
%        dessine_multi_circuitbis...
% Par la suite, prévoir le cas où il n'y a pas de sortie !!!!!!
% Bizarre tout de même, marche pour N=9, 10 ; Histoire de parité ?

%%%%%%%%%%%%%%%%%%%%%%%%%%%%%%%%%%%%%%%%%%%%%%%%%%%%%%%%%%%%%
\end{example}

\iflanguage{french}{%

%%%%%%%%%%%%%%%%%%%%%%%%%%%%%%%%%%%%%%%%%%%%%%%%%%%%%%%%%%%%%
\subsection{Comparaison avec les systèmes traditionnels et la théorie classique des polygones autoévitants (pavage carré)}
\label{compartradi}

Les systèmes traditionnels, comme \textit{Brio} \textregistered\
proposent une multitude de forme de rails, qui sont tous 
circulaires ou droits, mais jamais paraboliques. Voir
{\footnotesize{\url{http://www.woodenrailway.info/track/trackmath.html}}}.
Naturellement, on ne peut comparer les circuits étudiés % système \textit{Easyloop} \textregistered\
avec ces types de systèmes qui ne présentent pas de plans variés extensibles et modulables à volonté.
Cependant, dans des tels systèmes existent en particulier des huitièmes de cercles, qui assemblés deux par deux, 
donnent un quart de cercle, dont le rayon est égal à  l'une des longueurs des rails droits. Autrement dit, les pièces 
 des circuits étudiés %système \textit{Easyloop} \textregistered\, 
numéros \pieceu\ et \pieced\, utilisées seules\footnote{On peut aussi considérer les pièces
\piecet\ et \pieceq\ qui sont homothétiques des pièces  \pieceu\ et \pieced\ avec un rapport $\sqrt{2}$.}
permettent de créer des circuits simples qui pourraient être créés avec les systèmes traditionnels. Ces circuits
ont des pièces droites qui ne peuvent qu'être perpendiculaires entre elles : les formes sont moins variées et surtout, le nombre
de possibilités de circuits offerts est beaucoup moins grand que ceux de notre système (voir sections \ref{enumerationpetit} et \ref{estimation}).

Choisissons donc maintenant de montrer un circuit, uniquement formé des pièces \pieceu\ et \pieced .
Notons que, dans ce cas, les carrés adjacents ne peuvent avoir en commun qu'un côté et que l'on est très proche du 
cas des polygones autoévitants, mais des pièces dans le même carré peuvent aussi coexister. 
Notons aussi que le nombre de pièces utilisées est nécessairement pair exactement comme dans le cas des polygone autoévitants.%
}{%

%%%%%%%%%%%%%%%%%%%%%%%%%%%%%%%%%%%%%%%%%%%%%%%%%%%%%%%%%%%%%
%\subsection{Comparison with traditional systems}
\subsection{Comparison with traditional systems and the classical theory of self-avoiding polygons (square tiling)}
\label{compartradi}

Traditional systems, such as  \textit{Brio} \textregistered, offer a multitude of shapes of track pieces, which are all circular or straight, but never parabolic. See
{\footnotesize{\url{http://www.woodenrailway.info/track/trackmath.html}}}.
Naturally, one cannot compare the studied circuits %\textit{Easyloop} \textregistered\ system 
with these types of system, which do not offer designs which can be varied, scaled and modulated at will. However, in some such systems there exist in particular eighths of a circle, which assembled in pairs gives a quarter-circle, whose radius %diameter 
is equal to one of the lengths of the straight track pieces. In other words, the pieces of the  
studied systems%\textit{Easyloop} \textregistered\ system 
numbered \pieceu\ and \pieced, used alone\footnote{One may also consider pieces \piecet\ and \pieceq\ which are homothetic to pieces \pieceu\ and \pieced\ with the ratio $\sqrt{2}$.}, allow the creation of simple circuits which could be created with the traditional systems. These circuits have straight pieces which can only be perpendicular to each other. The forms are less varied and above all, the number of possible circuits offered is much smaller than those of our system (see Sections \ref{enumerationpetit} and \ref{estimation}).

We now choose to show a circuit formed solely from pieces \pieceu\ and \pieced. We note that, in this case, adjacent squares may only have one common side, and that this is very close to the case of self-avoiding polygons, but pieces in the same square may also coexist. We also note that the number of pieces used is necessarily even, exactly in the case of self-avoiding polygons.%

}

\begin{example}\
\label{examplesimulation600}
%%%%%%%%%%%%%%%%%%%%%%%%%%%%%%%%%%%%%%%%%%%%%%%%%
%\input{simulations_circuit/simulation_new_600}
% fichier tex crée par MaTeXBuild02 le 11-Feb-2016 10:09:53
% à compiler avec 
% MaTeXBuild02('simulation_new_600',0)
% après le fichier 'enumeration_construction_circuit_new.matex'

%%%%%%%%%%%%%%%%%%%%%%%%%%%%%%%%%%%%%%%%%%%%%%%%%
%\input{./simulations_circuit/circuit_numerique/circuits_complets_4_piecmax_ev3brio_new}
% fichier crée par 'presentation_exhaustif_circuit_boucle.m' le 11-Feb-2016 10:09:53
\begin{figure}[h]
\begin{center}
\epsfig{file=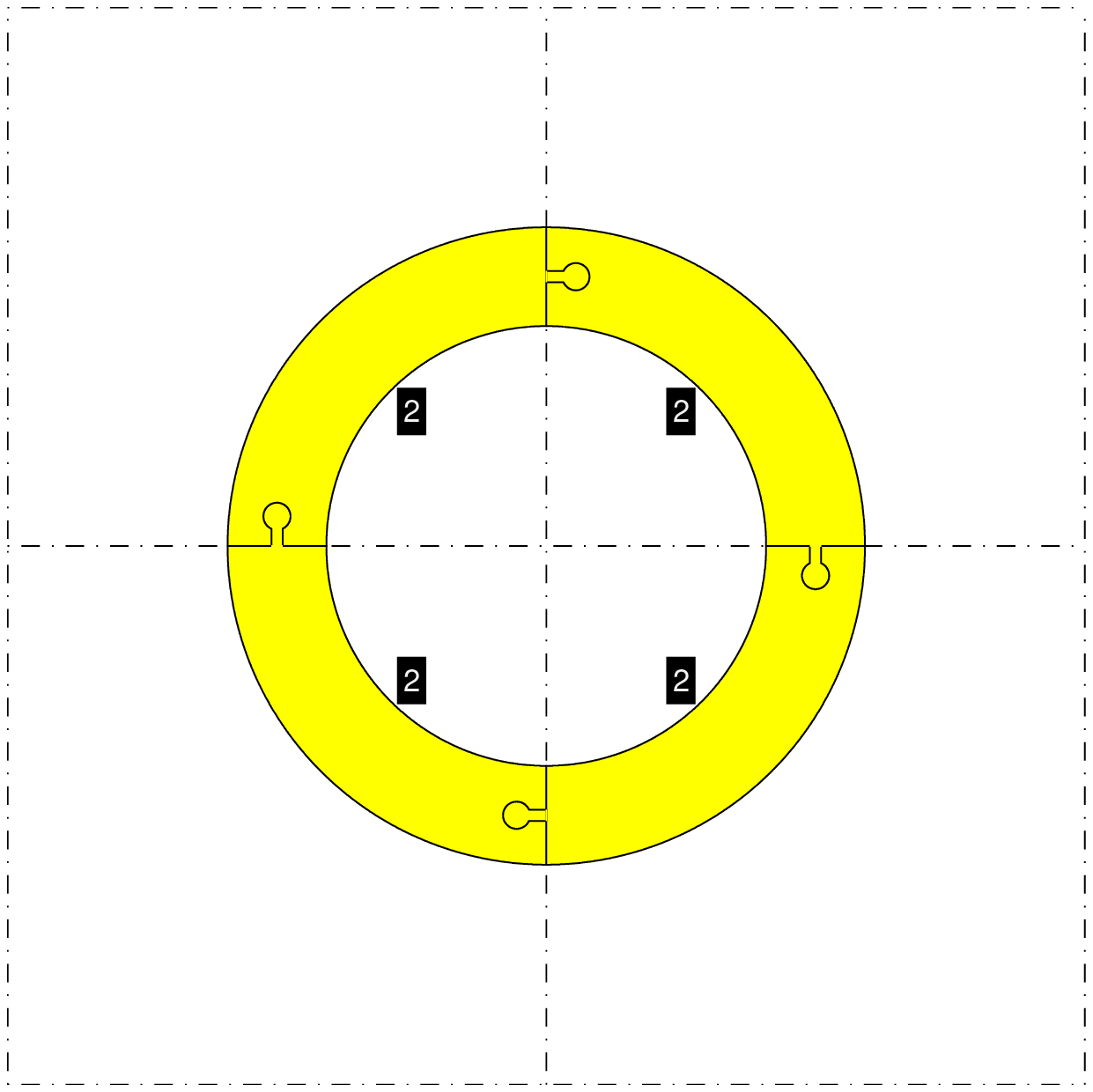, width=5 cm}
\end{center}
\caption{\label{circuits_complets_4_piecmax_ev3brio_new}\iflanguage{french}{Le seul circuit retenu sur l'ensemble des 400 circuits possibles}{The sole circuit kept from the set of 400 possible circuits}.}
\end{figure}

%%%%%%%%%%%%%%%%%%%%%%%%%%%%%%%%%%%%%%%%%%%%%%%%%

\iflanguage{french}{%
On trace l'unique circuit réalisable %tous les circuits réalisables 
avec $N=4$ pièces. On choisit 
$N_j=+\infty$, si $j\in \{1,2\}$ et nul sinon.
On obtient les  $1$ circuits de la figure \ref{circuits_complets_4_piecmax_ev3brio_new}.%
}{%
We draw 
the sole feasible circuit
%all of the feasible circuits 
with $N=4$ pieces. We choose $N_j=+\infty$ if $j\in \{1,2\}$, and zero otherwise. We obtain the $1$ circuits in Figure 
\ref{circuits_complets_4_piecmax_ev3brio_new}.%
}

%%%%%%%%%%%%%%%%%%%%%%%%%%%%%%%%%%%%%%%%%%%%%%%%%
\end{example}

\begin{example}\
\label{examplesimulation601}
%%%%%%%%%%%%%%%%%%%%%%%%%%%%%%%%%%%%%%%%%%%%%%%%%
%\input{simulations_circuit/simulation_new_601}
% fichier tex crée par MaTeXBuild02 le 11-Feb-2016 10:13:46
% à compiler avec 
% MaTeXBuild02('simulation_new_601',0)
% après le fichier 'enumeration_construction_circuit_new.matex'

%%%%%%%%%%%%%%%%%%%%%%%%%%%%%%%%%%%%%%%%%%%%%%%%%
%\input{./simulations_circuit/circuit_numerique/circuits_complets_6_piecmax_ev3brio_new}
% fichier crée par 'presentation_exhaustif_circuit_boucle.m' le 11-Feb-2016 10:13:47
\begin{figure}[h]
\begin{center}
\epsfig{file=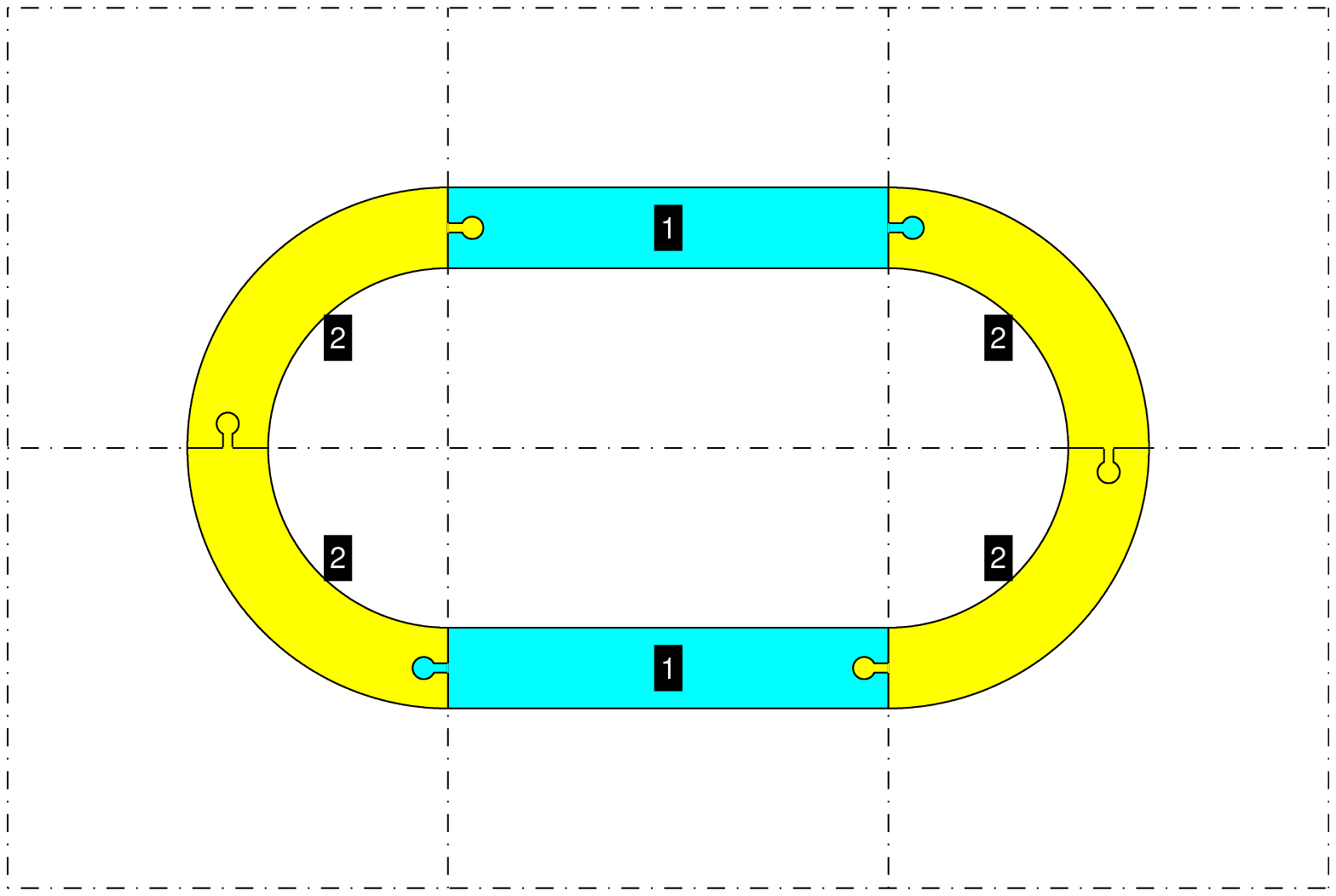, width=5 cm}
\end{center}
\caption{\label{circuits_complets_6_piecmax_ev3brio_new}\iflanguage{french}{Le seul circuit retenu sur l'ensemble des 10000 circuits possibles}{The sole circuit kept from the set of 10000 possible circuits}.}
\end{figure}
%%%%%%%%%%%%%%%%%%%%%%%%%%%%%%%%%%%%%%%%%%%%%%%%%

\iflanguage{french}{%
% Attention, label non automatique apparemment !!!
Comme dans l'exemple \ref{examplesimulation600},
on trace l'unique circuit réalisable %tous les circuits réalisables 
avec $N=6$ pièces. On choisit aussi 
$N_j=+\infty$, si $j\in \{1,2\}$ et nul sinon.
On obtient les  $1$ circuits de la figure \ref{circuits_complets_6_piecmax_ev3brio_new}.%
}{%
As in the example \ref{examplesimulation600},
we draw 
the sole feasible circuit
%all of the feasible circuits 
with $N=6$ pieces. We choose also $N_j=+\infty$ if $j\in \{1,2\}$, and zero otherwise. We obtain the $1$ circuits in Figure 
\ref{circuits_complets_6_piecmax_ev3brio_new}.%
}
%%%%%%%%%%%%%%%%%%%%%%%%%%%%%%%%%%%%%%%%%%%%%%%%%
\end{example}

\begin{example}\
\label{examplesimulation602}
%%%%%%%%%%%%%%%%%%%%%%%%%%%%%%%%%%%%%%%%%%%%%%%%%
%\input{simulations_circuit/simulation_new_602}
% fichier tex crée par MaTeXBuild02 le 19-Feb-2016 16:41:44
% à compiler avec 
% MaTeXBuild02('simulation_new_602',0)
% après le fichier 'enumeration_construction_circuit_new.matex'

%%%%%%%%%%%%%%%%%%%%%%%%%%%%%%%%%%%%%%%%%%%%%%%%%
%\input{./simulations_circuit/circuit_numerique/circuits_complets_8_piecmax_ev3brio_new}
% fichier crée par 'presentation_exhaustif_circuit_boucle.m' le 19-Feb-2016 16:41:46
\begin{figure}[h]
\centering
%%% sous figure 1
\subfigure[\label{circuits_complets_8_piecmax_ev3brio_new1}]
{\epsfig{file=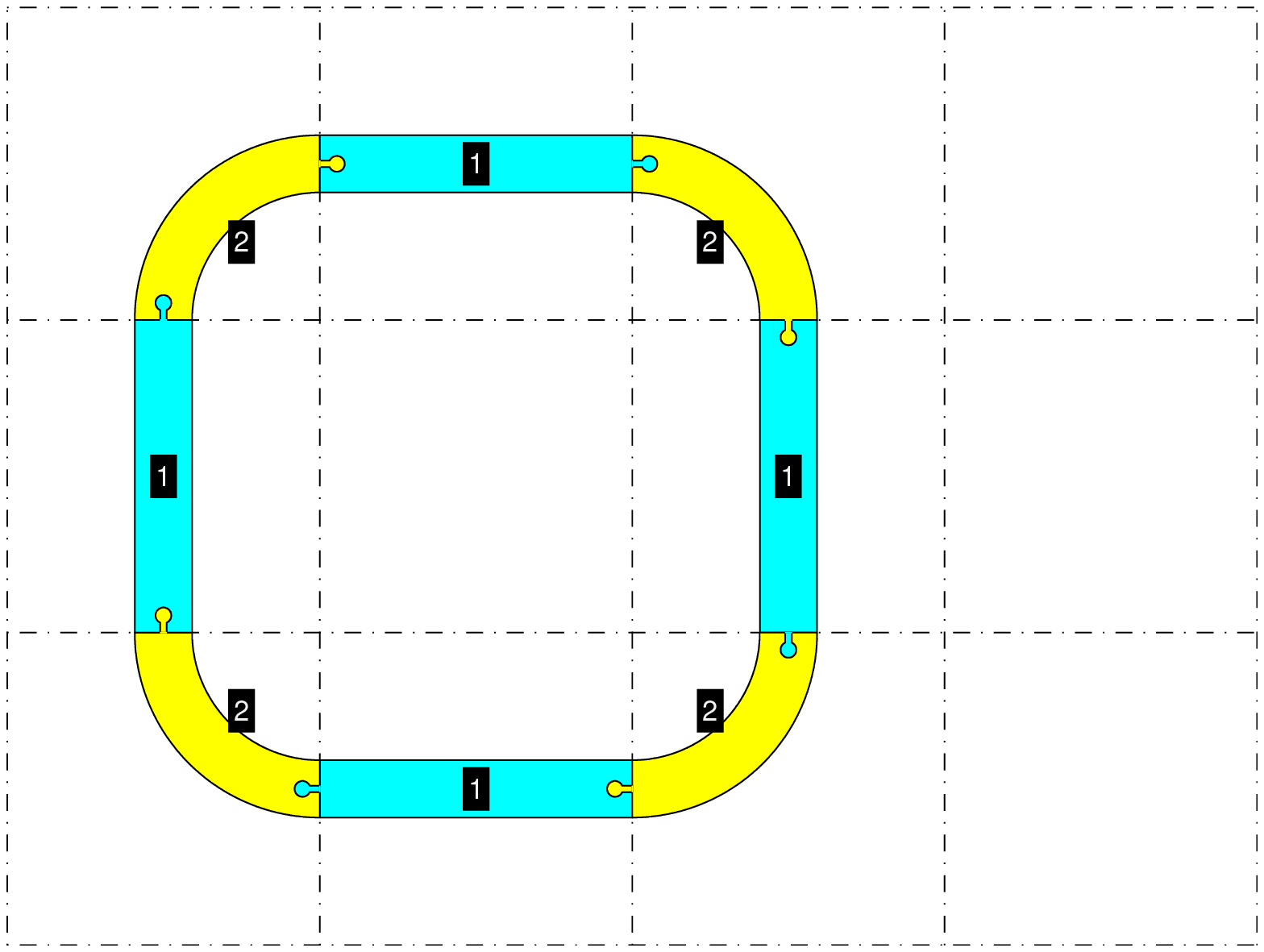, width=5 cm}}
\qquad
%%% sous figure 2
\subfigure[\label{circuits_complets_8_piecmax_ev3brio_new2}]
{\epsfig{file=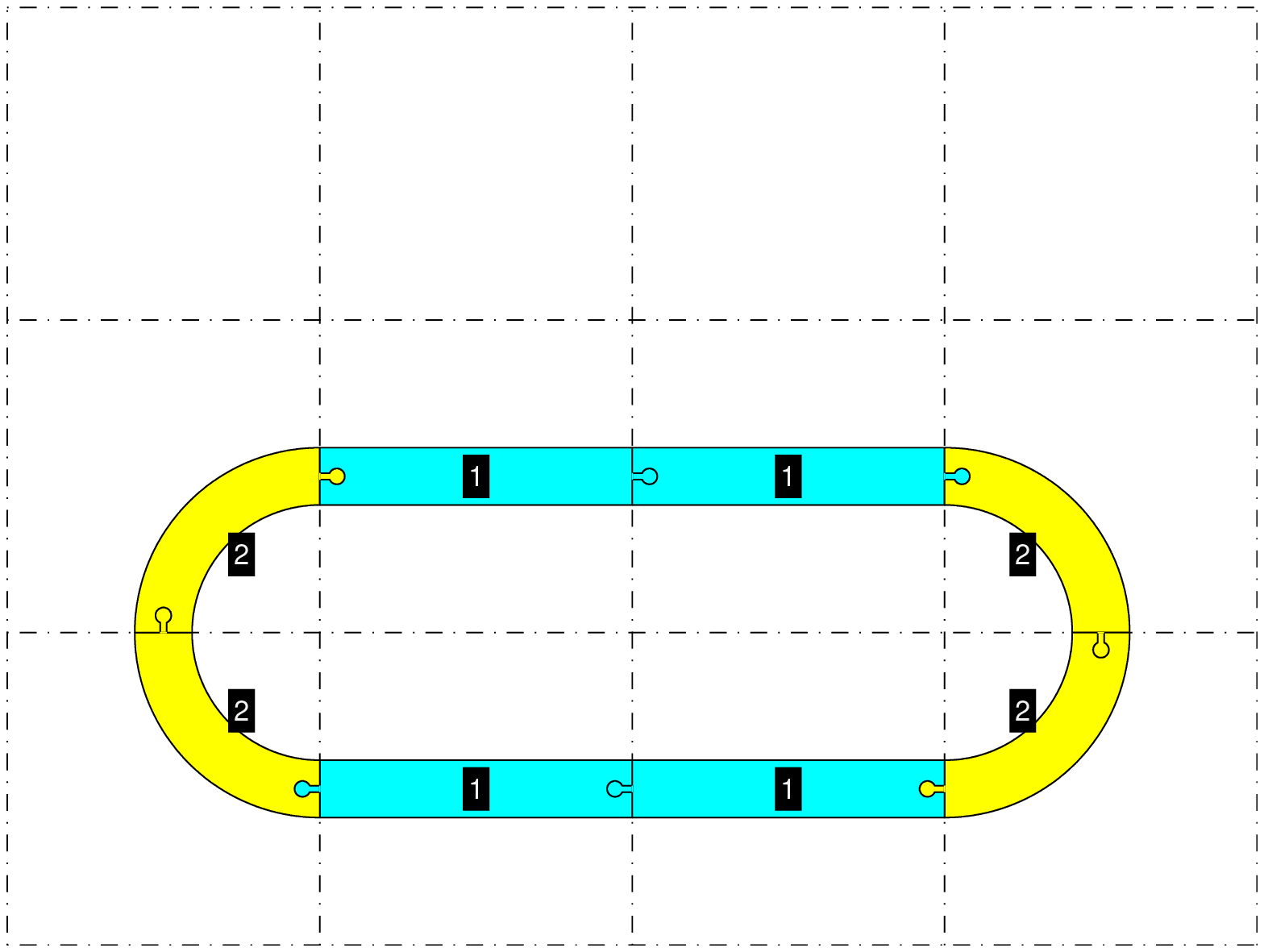, width=5 cm}}
\qquad
%%% sous figure 3
\subfigure[\label{circuits_complets_8_piecmax_ev3brio_new3}]
{\epsfig{file=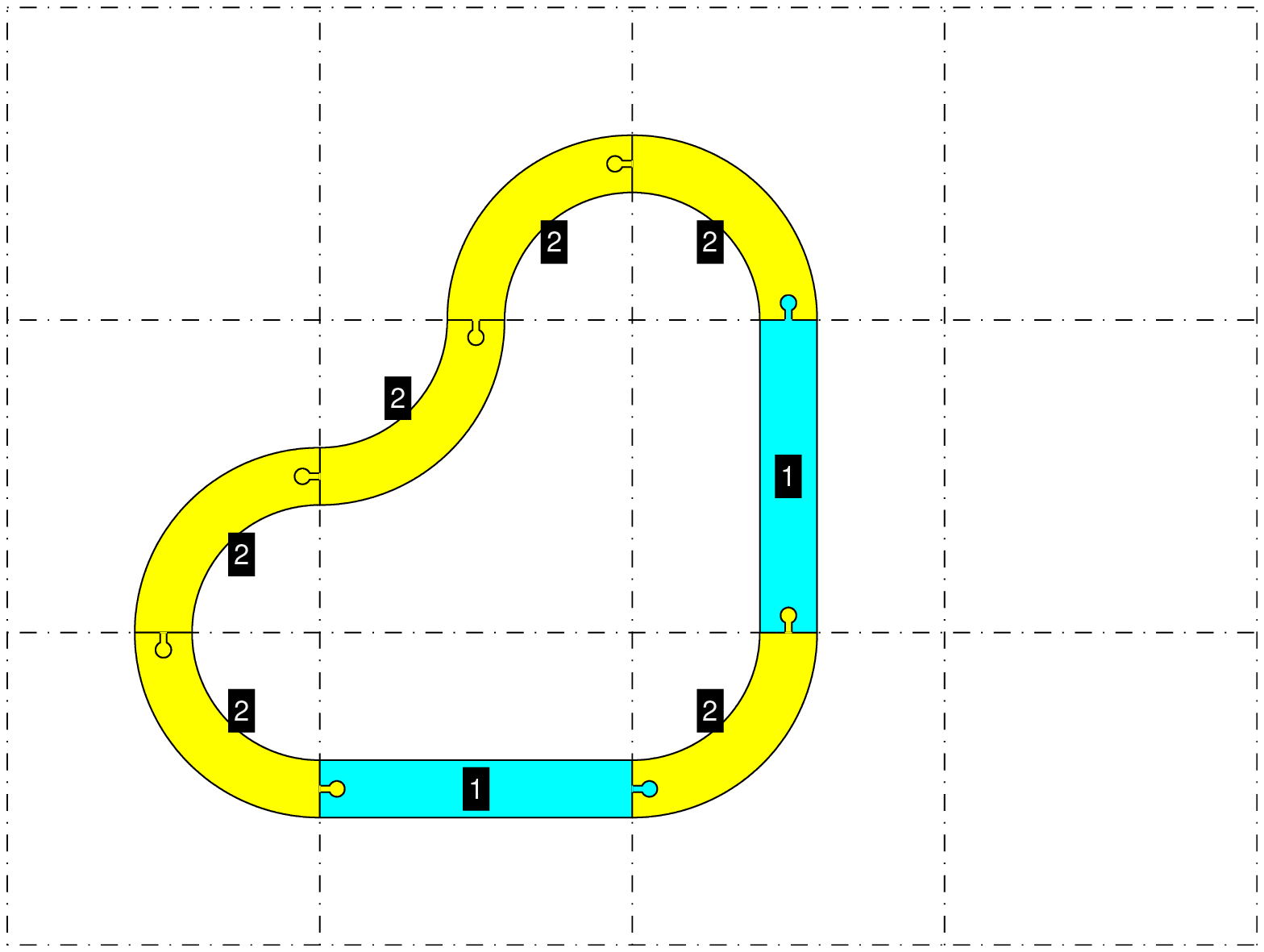, width=5 cm}}
\qquad
%%% sous figure 4
\subfigure[\label{circuits_complets_8_piecmax_ev3brio_new4}]
{\epsfig{file=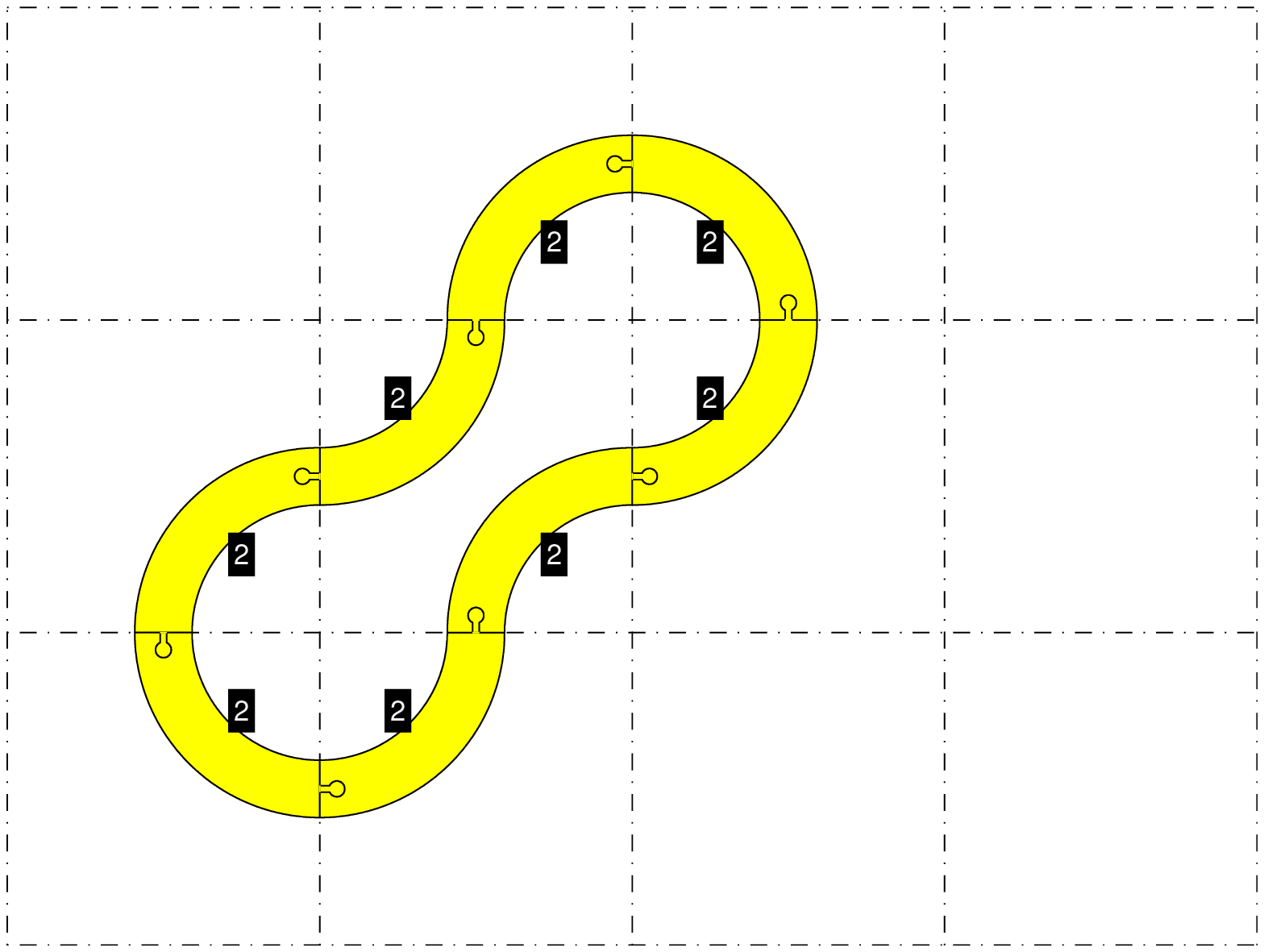, width=5 cm}}
\qquad
\caption{\label{circuits_complets_8_piecmax_ev3brio_new}\iflanguage{french}{Tous les 4 circuits retenus sur l'ensemble des 250000 circuits possibles}{All of the 4 circuits kept from the set of 250000 possible circuits}.}
\end{figure}
%%%%%%%%%%%%%%%%%%%%%%%%%%%%%%%%%%%%%%%%%%%%%%%%%

\iflanguage{french}{%
% Attention, label non automatique apparemment !!!
Comme dans l'exemple \ref{examplesimulation600},
on trace tous les  circuits réalisables %tous les circuits réalisables 
avec $N=8$ pièces. 
On obtient les  $4$ circuits de la figure \ref{circuits_complets_8_piecmax_ev3brio_new}.
Notons que l'on retrouve naturellement quelques-uns des circuits 
à $8$ pièces de l'exemple \ref{examplesimulation515}
(voir figures
\ref{circuits_complets_8_piecmax_ev31},
\ref{circuits_complets_8_piecmax_ev32} et 
\ref{circuits_complets_8_piecmax_ev33}).%
}{%
As in the example \ref{examplesimulation600},
we draw 
all the feasible circuits
%all of the feasible circuits 
with $N=8$ pieces. 
We obtain the $4$ circuits in Figure 
\ref{circuits_complets_8_piecmax_ev3brio_new}.
% attention anglais non traduit par ml
Note that we find naturally some of circuits with $8$ pieces of Example \ref{examplesimulation515}
(see Figures
\ref{circuits_complets_8_piecmax_ev31},
\ref{circuits_complets_8_piecmax_ev32} and 
\ref{circuits_complets_8_piecmax_ev33}).%
}
%%%%%%%%%%%%%%%%%%%%%%%%%%%%%%%%%%%%%%%%%%%%%%%%%
\end{example}

\begin{example}\
\label{examplesimulation603}
%%%%%%%%%%%%%%%%%%%%%%%%%%%%%%%%%%%%%%%%%%%%%%%%%
%\input{simulations_circuit/simulation_new_603}
% fichier tex crée par MaTeXBuild02 le 11-Feb-2016 16:03:54
% à compiler avec 
% MaTeXBuild02('simulation_new_603',0)
% après le fichier 'enumeration_construction_circuit_new.matex'

%%%%%%%%%%%%%%%%%%%%%%%%%%%%%%%%%%%%%%%%%%%%%%%%%
%\input{./simulations_circuit/circuit_numerique/circuits_complets_10_piecmax_ev3brio_new}
% fichier crée par 'presentation_exhaustif_circuit_boucle.m' le 11-Feb-2016 16:12:58
\begin{figure}[h]
\centering
%%% sous figure 1
\subfigure[\label{circuits_complets_10_piecmax_ev3brio_new1}]
{\epsfig{file=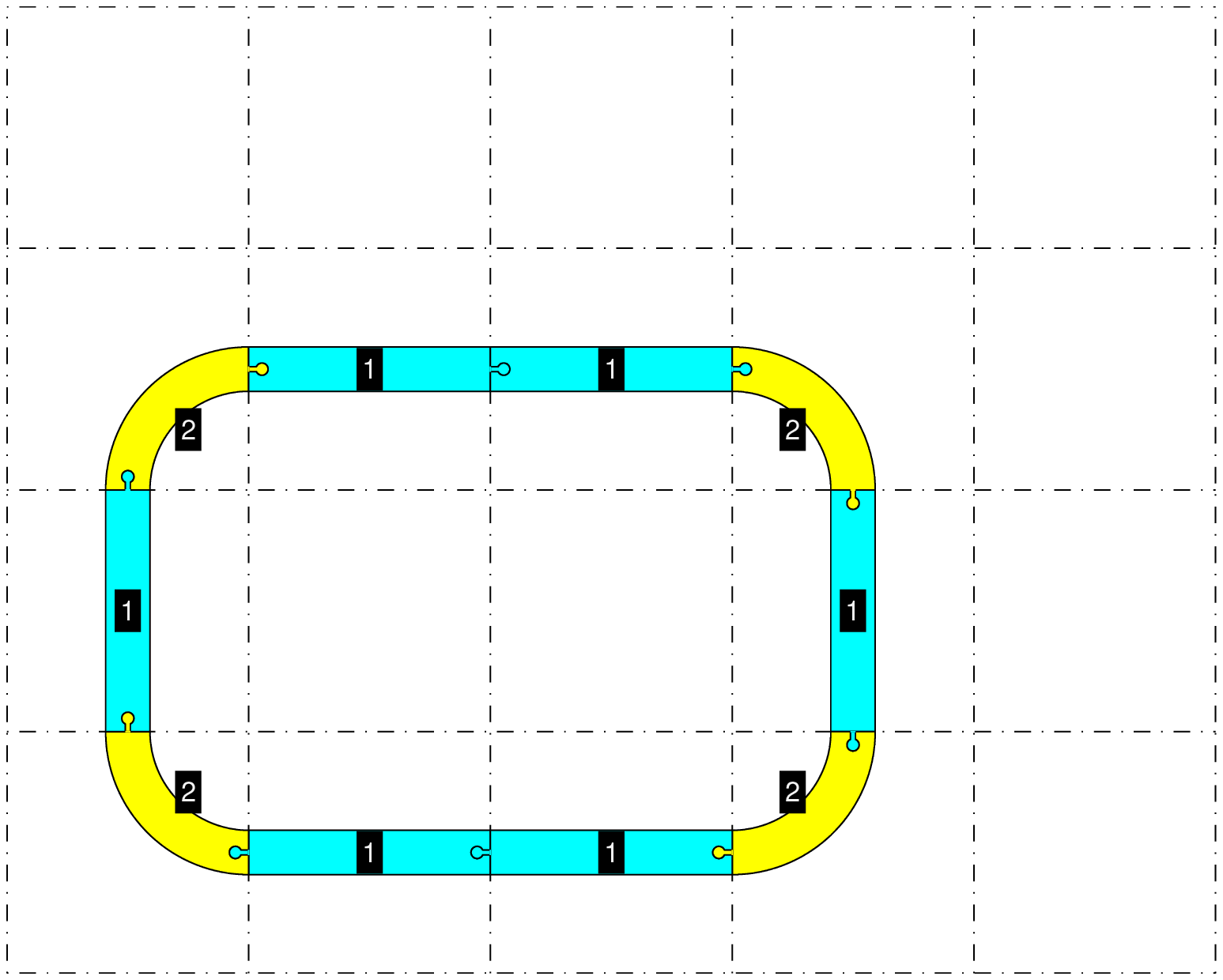, width=5 cm}}
\qquad
%%% sous figure 2
\subfigure[\label{circuits_complets_10_piecmax_ev3brio_new2}]
{\epsfig{file=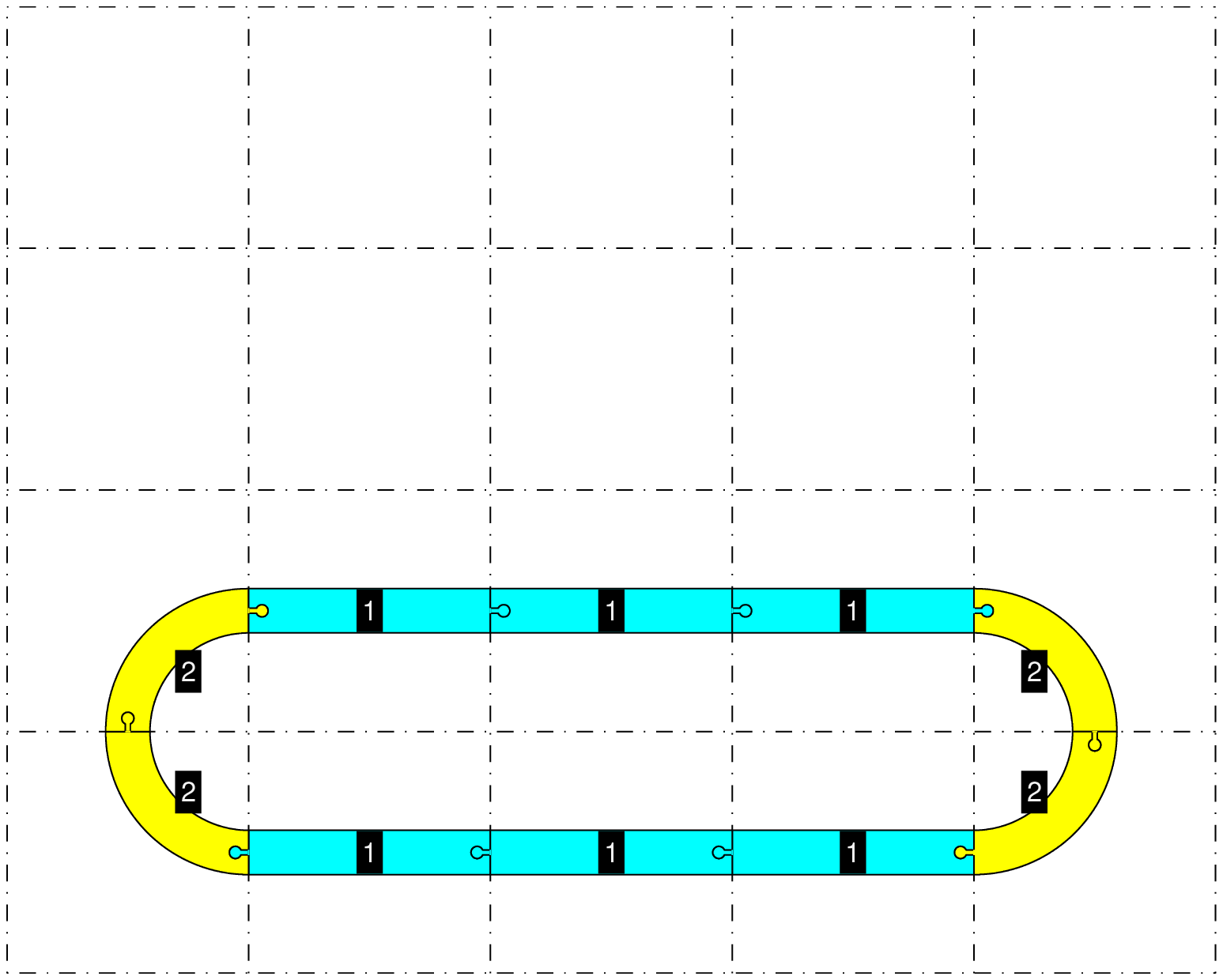, width=5 cm}}
\qquad
%%% sous figure 3
\subfigure[\label{circuits_complets_10_piecmax_ev3brio_new3}]
{\epsfig{file=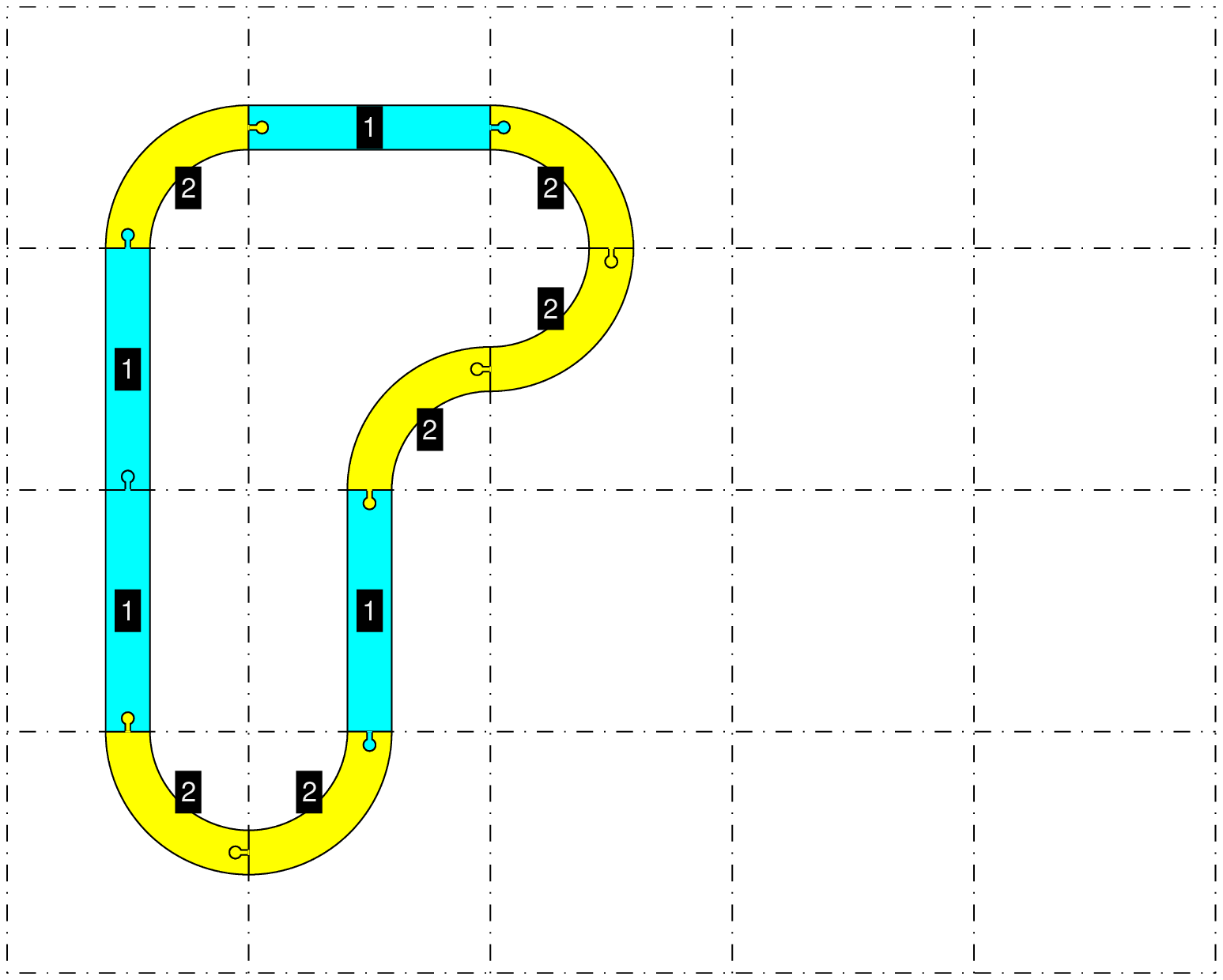, width=5 cm}}
\qquad
%%% sous figure 4
\subfigure[\label{circuits_complets_10_piecmax_ev3brio_new4}]
{\epsfig{file=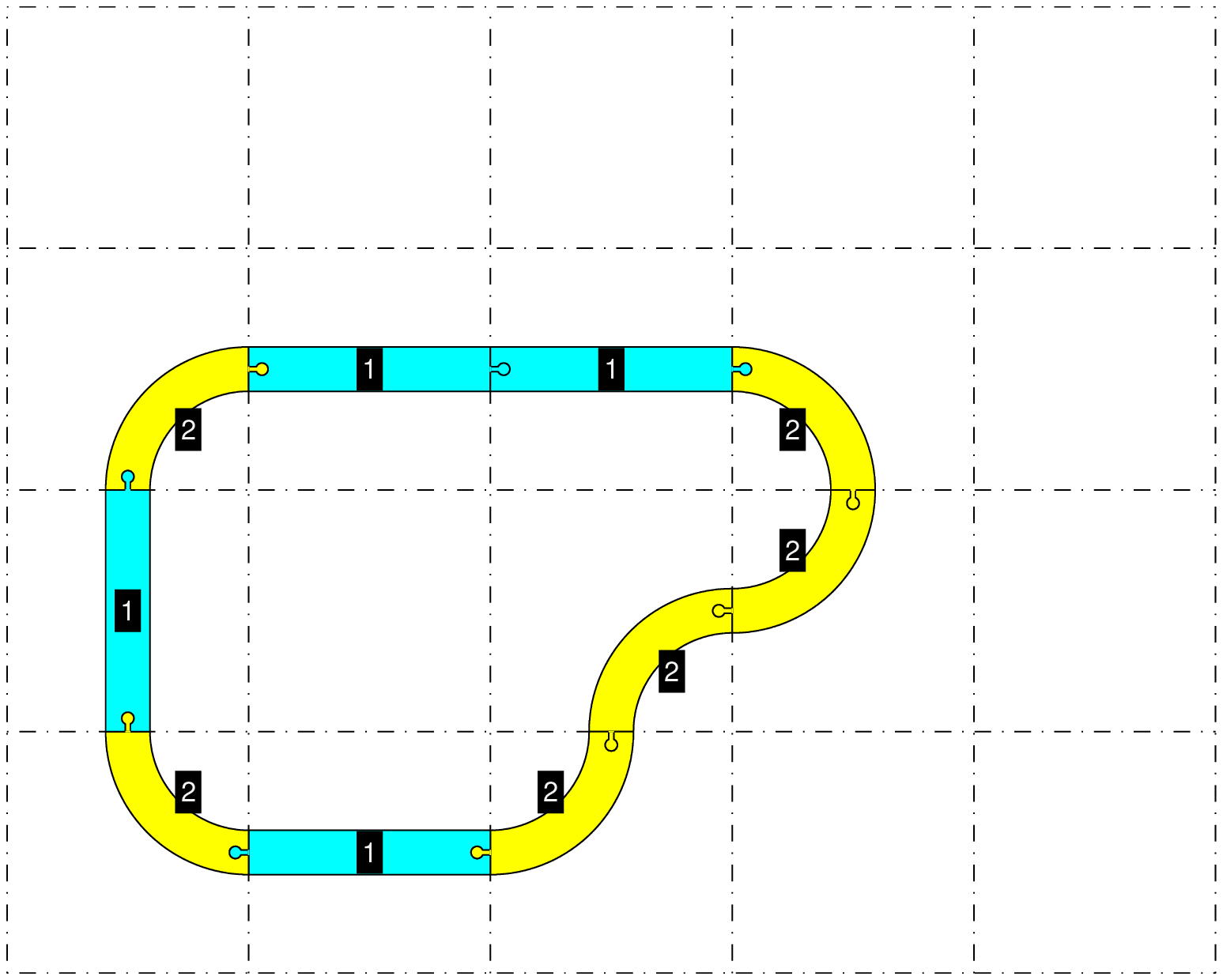, width=5 cm}}
\qquad
%%% sous figure 5
\subfigure[\label{circuits_complets_10_piecmax_ev3brio_new5}]
{\epsfig{file=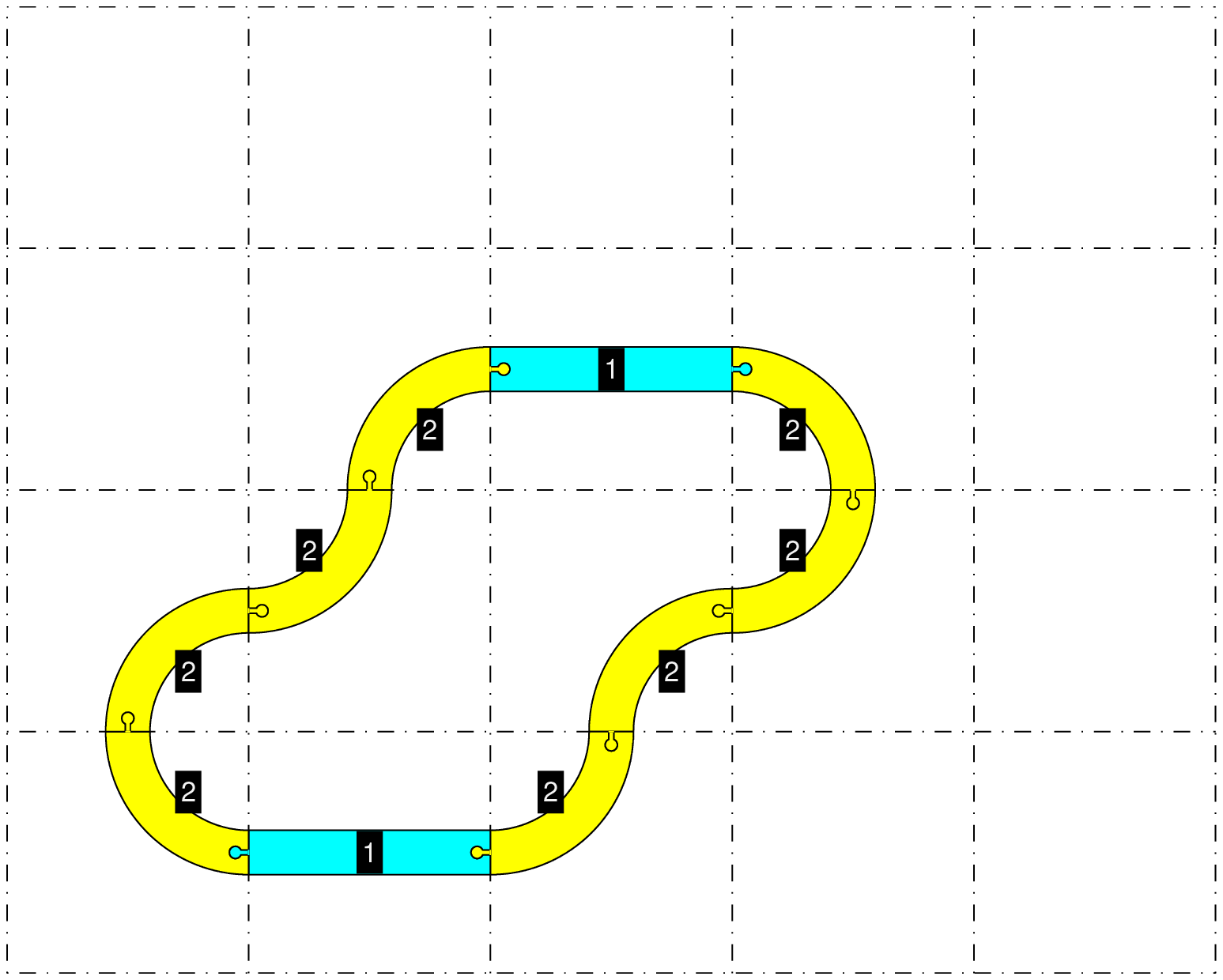, width=5 cm}}
\qquad
%%% sous figure 6
\subfigure[\label{circuits_complets_10_piecmax_ev3brio_new6}]
{\epsfig{file=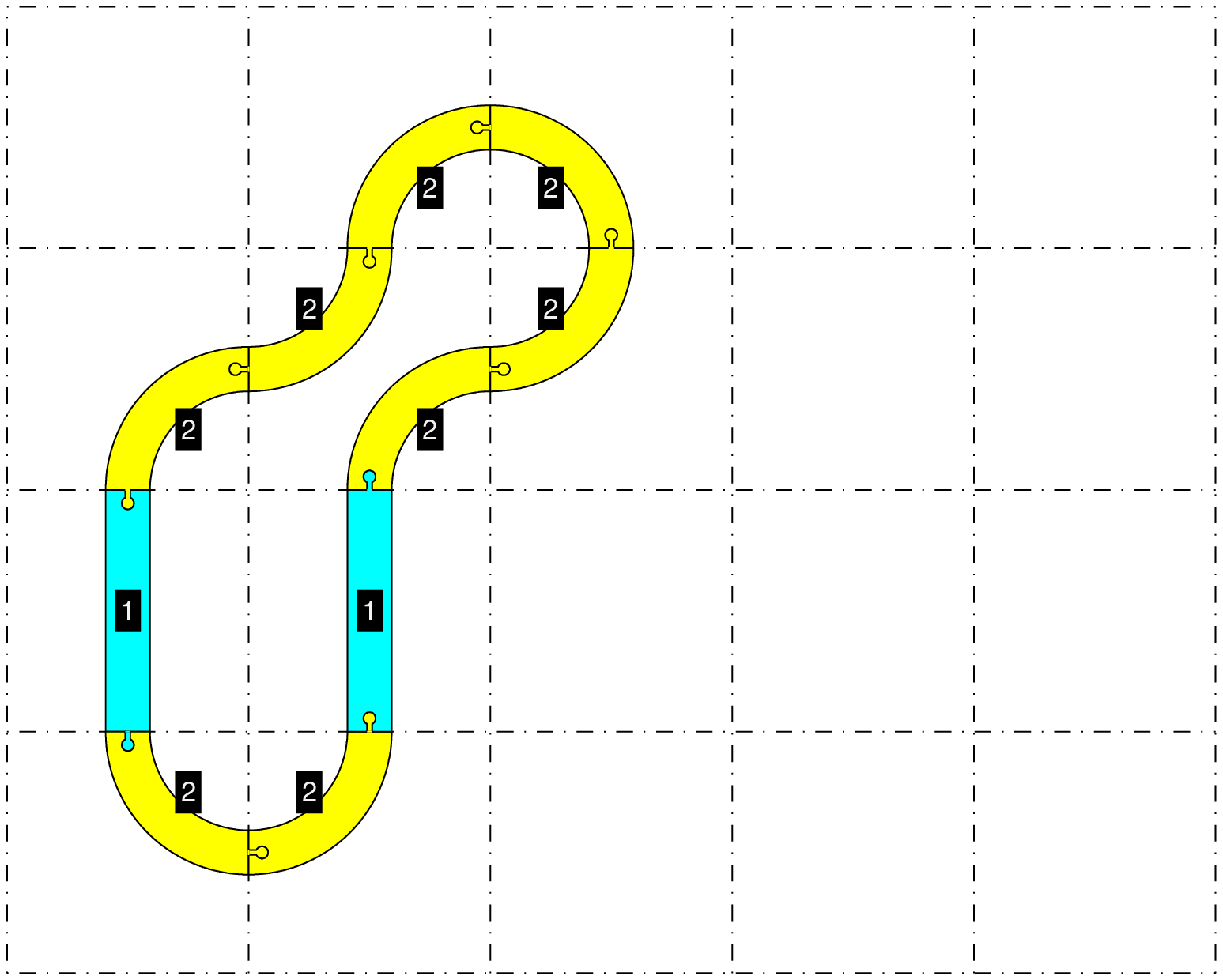, width=5 cm}}
\qquad
%%% sous figure 7
\subfigure[\label{circuits_complets_10_piecmax_ev3brio_new7}]
{\epsfig{file=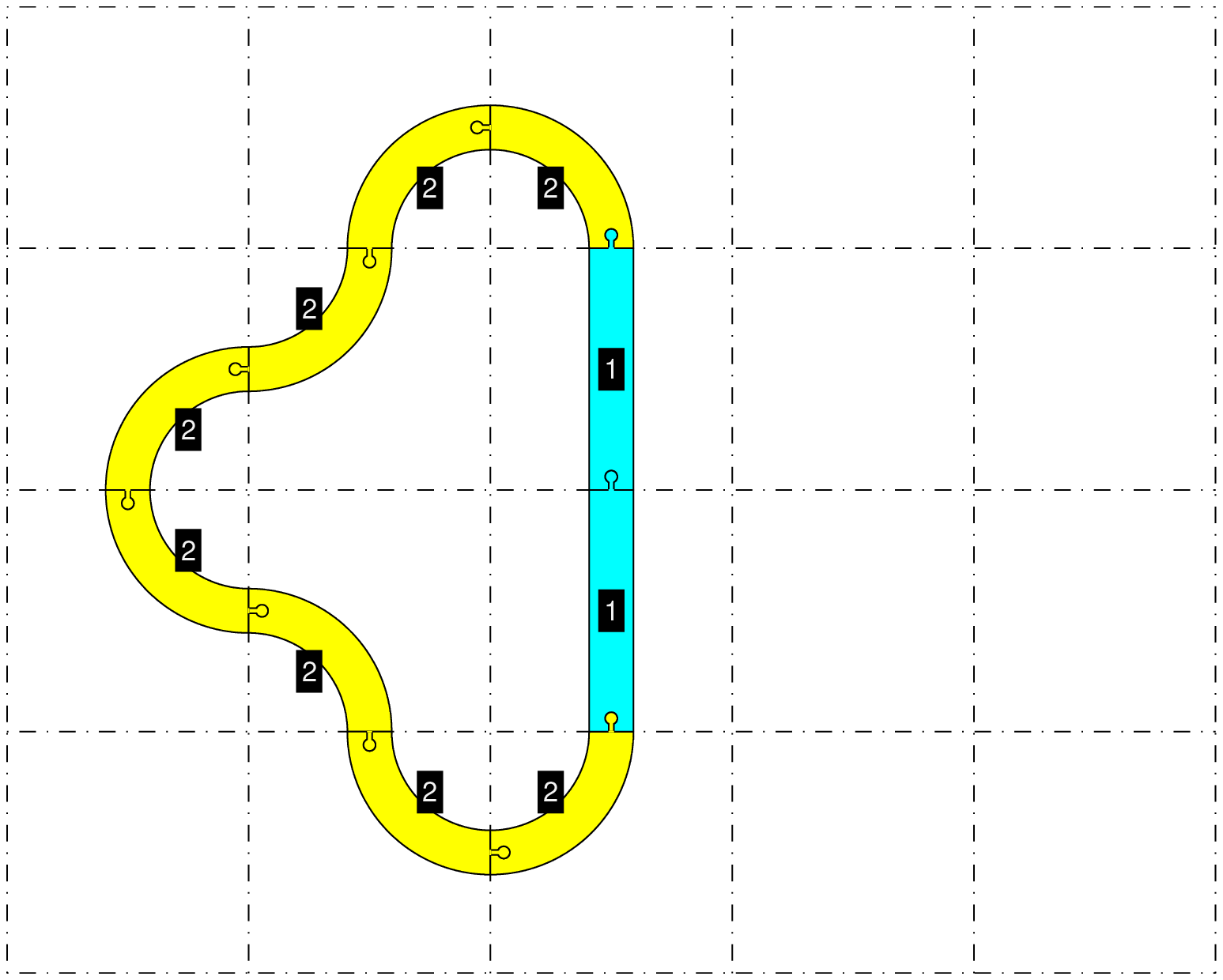, width=5 cm}}
\qquad
\caption{\label{circuits_complets_10_piecmax_ev3brio_new}\iflanguage{french}{Tous les 7 circuits retenus sur l'ensemble des 6250000 circuits possibles}{All of the 7 circuits kept from the set of 6250000 possible circuits}.}
\end{figure}
%%%%%%%%%%%%%%%%%%%%%%%%%%%%%%%%%%%%%%%%%%%%%%%%%

\iflanguage{french}{%
% Attention, label non automatique apparemment !!!
Comme dans l'exemple \ref{examplesimulation600},
on trace tous les  circuits réalisables %tous les circuits réalisables 
avec $N=10$ pièces. 
On obtient les  $7$ circuits de la figure \ref{circuits_complets_10_piecmax_ev3brio_new}.%
}{%
As in the example \ref{examplesimulation600},
we draw 
all the feasible circuits
%all of the feasible circuits 
with $N=10$ pieces. 
We obtain the $7$ circuits in Figure 
\ref{circuits_complets_10_piecmax_ev3brio_new}.%
}
%%%%%%%%%%%%%%%%%%%%%%%%%%%%%%%%%%%%%%%%%%%%%%%%%
\end{example}

%%%%%%%%%%%%%%%%%%%%%%%%%%%%%%%%%%%%%%%%%%%%%%%%%%%%%%%%%%%%%%%%%%%%%%%%%%%%%%%%%%%%%%%%%%%%%%%%%%%%%%%%%%%%%%%%%%%%%%%%%%%
%%%%%%%%%%%%%%%%%%%%%%%%%%%%%%%%%%%%%%%%%%%%%%%%%%%%%%%%%%%%%%%%%%%%%%%%%%%%%%%%%%%%%%%%%%%%%%%%%%%%%%%%%%%%%%%%%%%%%%%%%%%
\clearpage
%%%%%%%%%%%%%%%%%%%%%%%%%%%%%%%%%%%%%%%%%%%%%%%%%%%%%%%%%%%%%%%%%%%%%%%%%%%%%%%%%%%%%%%%%%%%%%%%%%%%%%%%%%%%%%%%%%%%%%%%%%%
%%%%%%%%%%%%%%%%%%%%%%%%%%%%%%%%%%%%%%%%%%%%%%%%%%%%%%%%%%%%%%%%%%%%%%%%%%%%%%%%%%%%%%%%%%%%%%%%%%%%%%%%%%%%%%%%%%%%%%%%%%%

\iflanguage{french}{%
%%%%%%%%%%%%%%%%%%%%%%%%%%%%%%%%%%%%%%%%%%%%%%%%%%%%%%%%%%%%%
\subsection{Détermination du nombre de circuits}
\label{enumerationpetit}%
}{%
%%%%%%%%%%%%%%%%%%%%%%%%%%%%%%%%%%%%%%%%%%%%%%%%%%%%%%%%%%%%%
\subsection{Determination of the number of circuits}
\label{enumerationpetit}
}

%%% Vérifier passage de 11 à 12

%\ifcase \nopiecesix
%\textbf{ERREUR !!!!!}
%Calculs prévus pour \path|nopiece6|=1.
%\or
%\fi

% table crée par stocke_tableau_tex le 29-Mar-2016 14:34:39
\begin{table}[h]
\begin{center}
\begin{tabular}{|r|r|r|r|r|r|r|}
\hline
$N$&possible circuits&looping circuits&direct isom.&isometries&constructible
\\
\hline
\hline
$1$ & $2$ & $0$ & $0$ & $0$ & $0$
\\
\hline
$2$ & $16$ & $0$ & $0$ & $0$ & $0$
\\
\hline
$3$ & $80$ & $0$ & $0$ & $0$ & $0$
\\
\hline
$4$ & $400$ & $4$ & $4$ & $2$ & $2$
\\
\hline
$5$ & $2000$ & $10$ & $2$ & $1$ & $1$
\\
\hline
$6$ & $10000$ & $36$ & $10$ & $5$ & $5$
\\
\hline
$7$ & $50000$ & $140$ & $20$ & $7$ & $7$
\\
\hline
$8$ & $250000$ & $664$ & $107$ & $41$ & $33$
\\
\hline
$9$ & $1250000$ & $2988$ & $332$ & $92$ & $74$
\\
\hline
$10$ & $6250000$ & $13910$ & $1466$ & $428$ & $304$
\\
\hline
$11$ & $31250000$ & $64592$ & $5872$ & $1512$ & $986$
\\
\hline
\end{tabular}
\vspace{1 cm}
\caption{\label{enumeration_construction_circuit_tab_new01}
Numbers of circuits corresponding to $N_j=+\infty$}
\end{center}
\end{table}

%\iflanguage{french}{\input{enumeration_construction_circuit_tab_new01}}{\input{enumeration_construction_circuit_eng_tab_new01}}

\iflanguage{french}{%
Dans le tableau \ref{enumeration_construction_circuit_tab_new01},
ont été donnés successivement (pour $N_j=+\infty$).
\begin{itemize}
\item
le nombre de circuits possibles examinés, parmi lesquels ont été retenus uniquement ceux qui se rebouclent ;
\item
le nombre de circuits qui se rebouclent ;
\item
le nombre de circuits qui se rebouclent, à une isométrie directe près ;
\item
le nombre de circuits qui se rebouclent, à une isométrie quelconque près ;
\item
et enfin, le nombre de circuits qui se rebouclent, à une isométrie quelconque  près en ne retenant que les circuits constructibles.
\end{itemize}%
}{%
In Table~\ref{enumeration_construction_circuit_tab_new01}, 
we give in succession (for $N_j=+\infty$).
\begin{itemize}
\item
the number of possible circuits examined, from which only those which form a loop were retained;
\item
the number of looping circuits;
\item
the number of looping circuits, up to a direct isometry;
\item
the number of looping circuits, up to any isometry;
\item
and finally, the number of looping circuits, up to any isometry, and only keeping constructible circuits.
\end{itemize}%
}

\iflanguage{french}{%
Notons que pour les valeurs de $N$ inférieures à $3$, aucun circuit n'existe.
Pour les valeurs de $N$ inférieures à 
$7$, tous les circuits sont constructibles.%
}{%

We note that for values of $N$ smaller than $3$, no circuits exist. For values of $N$ smaller than 
$7$, all of the circuits are constructible.%
}

% table crée par stocke_tableau_tex le 29-Mar-2016 14:34:39
\begin{table}[h]
\begin{center}
\begin{tabular}{|r|r|r|r|r|r|}
\hline
$N$&looping circuits&direct isom.&isometries&constructible
\\
\hline
\hline
$4$ & $4$ & $4$ & $2$ & $2$
\\
\hline
$5$ & $10$ & $2$ & $1$ & $1$
\\
\hline
$6$ & $36$ & $10$ & $5$ & $5$
\\
\hline
$7$ & $126$ & $18$ & $6$ & $6$
\\
\hline
$8$ & $564$ & $89$ & $30$ & $28$
\\
\hline
$9$ & $2358$ & $262$ & $72$ & $63$
\\
\hline
$10$ & $10710$ & $1126$ & $324$ & $244$
\\
\hline
$11$ & $45034$ & $4094$ & $1046$ & $753$
\\
\hline
\end{tabular}
\vspace{1 cm}
\caption{\label{enumeration_construction_circuit_tab_new02}
The (non-zero) numbers of circuits corresponding to $N_j=4$}
\end{center}
\end{table}

%\iflanguage{french}{\input{enumeration_construction_circuit_tab_new02}}{\input{enumeration_construction_circuit_eng_tab_new02}}

\iflanguage{french}{%
Dans le tableau \ref{enumeration_construction_circuit_tab_new02} ont été donnés les nombres (non nuls) de  circuits, correspondant
aux boîtes distribuées \textit{Easyloop} où  $N_j=4$.
Notons que jusqu'à $6$, les valeurs pour $N_j=4$ et $N_j=+\infty$ sont égales.%
}{%
In Table~\ref{enumeration_construction_circuit_tab_new02} we give the (non-zero) numbers of circuits corresponding to the distributed \textit{Easyloop} boxes, where $N_j=4$. We note that up to $6$, the values for $N_j=4$ and $N_j=+\infty$ are equal.%
}

% table crée par stocke_tableau_tex le 29-Mar-2016 14:34:39
\begin{table}[h]
\begin{center}
\begin{tabular}{|r|r|r|r|r|}
\hline
$N$&Self-avoiding polygons&traditional \textit{Brio} system&\textit{Easyloop} system
\\
\hline
\hline
$4$ & $1$ & $1$ & $2$
\\
\hline
$5$ & $0$ & $0$ & $1$
\\
\hline
$6$ & $2$ & $1$ & $5$
\\
\hline
$7$ & $0$ & $0$ & $6$
\\
\hline
$8$ & $7$ & $4$ & $28$
\\
\hline
$9$ & $0$ & $0$ & $63$
\\
\hline
$10$ & $28$ & $7$ & $244$
\\
\hline
$11$ & $0$ & $0$ & $753$
\\
\hline
\end{tabular}
\vspace{1 cm}
\caption{\label{enumeration_construction_circuit_tab_new03}
The numbers of self-avoiding polygons, the (non-zero) numbers of circuits for traditional  \textit{Brio} system and the \textit{Easyloop} system}
\end{center}
\end{table}

%\iflanguage{french}{\input{enumeration_construction_circuit_tab_new03}}{\input{enumeration_construction_circuit_eng_tab_new03}}

\iflanguage{french}{%
Enfin, dans le tableau \ref{enumeration_construction_circuit_tab_new03}
les nombres de polygones autoévitants, correspondant à un pavage carré, issu de \cite[tableau p. 396]{MR1197356} 
ou 
\url{http://oeis.org/A002931/b002931.txt}
ainsi qu'une 
comparaison entre les systèmes traditionnels (voir section \ref{compartradi})
et le système \textit{Easyloop} sont proposés ; pour ces derniers, seul le nombre de circuits constructibles à une isométrie près, est affiché.
Le nombre $N_j$ est égal à 
\ifcase 1
$12$ si $j=1$,
$12$ si $j=2$ et nul sinon.
\or
$12$ si $j=1,2$ et nul sinon.
\fi

Les circuits traditionnels sont très proches des polygones autoévitants, hormis les deux différences suivantes,
déjà évoquée plus haut : les isométries autorisées sont toutes les isométries et un carré  peut être emprunté plusieurs 
fois par le circuit. Le point commun est que pour $N$ impair, le nombre obtenu est nul.

%%%%%%%%%%%%%%%%%%%%%%%%%%%%%%%%%%%%%%%%%%%%%%%%
%\input{simulations_circuit/simulation_new_610}
%%%%%%%%%%%%%%%%%%%%%%%%%%%%%%%%%%%%%%%%%%%%%%%%

Si maintenant, on compare la dernière colonne du tableau \ref{enumeration_construction_circuit_tab_new03} avec les deux premières, 
on constate que 
la comparaison est donc en faveur du système \textit{Easyloop}.%

Le nombre d'opération est  en 
\begin{equation}
\label{nombreoperation}
\mathcal{O}\left(5^N\right)
\end{equation}
Ainsi, le passage de $N$ à $N+1$ multiplie le nombre d'opération par $5$.
Tous les résultats de cette section ont été réalisés
\ifcase 0
en un temps supérieur à $23$ heures.
De remplacer $11$ par 
$12$
aurait donc exigé  plus de $4$ jours de calcul  !
\or
\textbf{ATTENTION  !!!!!!!!!!!!!!!!!}
Durée inférieure à une heure !!
\fi
}{%
Finally, in Table~\ref{enumeration_construction_circuit_tab_new03}, 
the numbers of self-avoiding polygons, corresponding to a square lattice (see \cite[table p. 396]{MR1197356}) 
or
\url{http://oeis.org/A002931/b002931.txt}
and 
a comparison between traditional systems (see Section~\ref{compartradi}) and the \textit{Easyloop} system are proposed. 
For the \textit{Easyloop} system,
only the number of constructible circuits up to an isometry is displayed. The number $N_j$ equals 
\ifcase 1
$12$ if  $j=1$,
$12$ if  $j=2$, and zero otherwise .
\or
$12$ if  $j=1,2$, and zero otherwise.
\fi

Traditional circuits are very close to self-avoiding polygons, except for the following two differences already mentioned above: the permitted isometries are all the isometries of a square, and a square may be used multiple times by the circuit. The common point is that for $N$ odd, the number obtained is zero.

%%%%%%%%%%%%%%%%%%%%%%%%%%%%%%%%%%%%%%%%%%%%%%%%%%%
%\input{simulations_circuit/simulation_new_610}
% fichier tex crée par MaTeXBuild02 le 29-Feb-2016 15:57:03
% à compiler avec 
% MaTeXBuild02('simulation_new_610',0)
% après le fichier 'enumeration_construction_circuit_new.matex'

% attention, suite non automatique

\begin{figure}[h]
\begin{center}
\epsfig{file=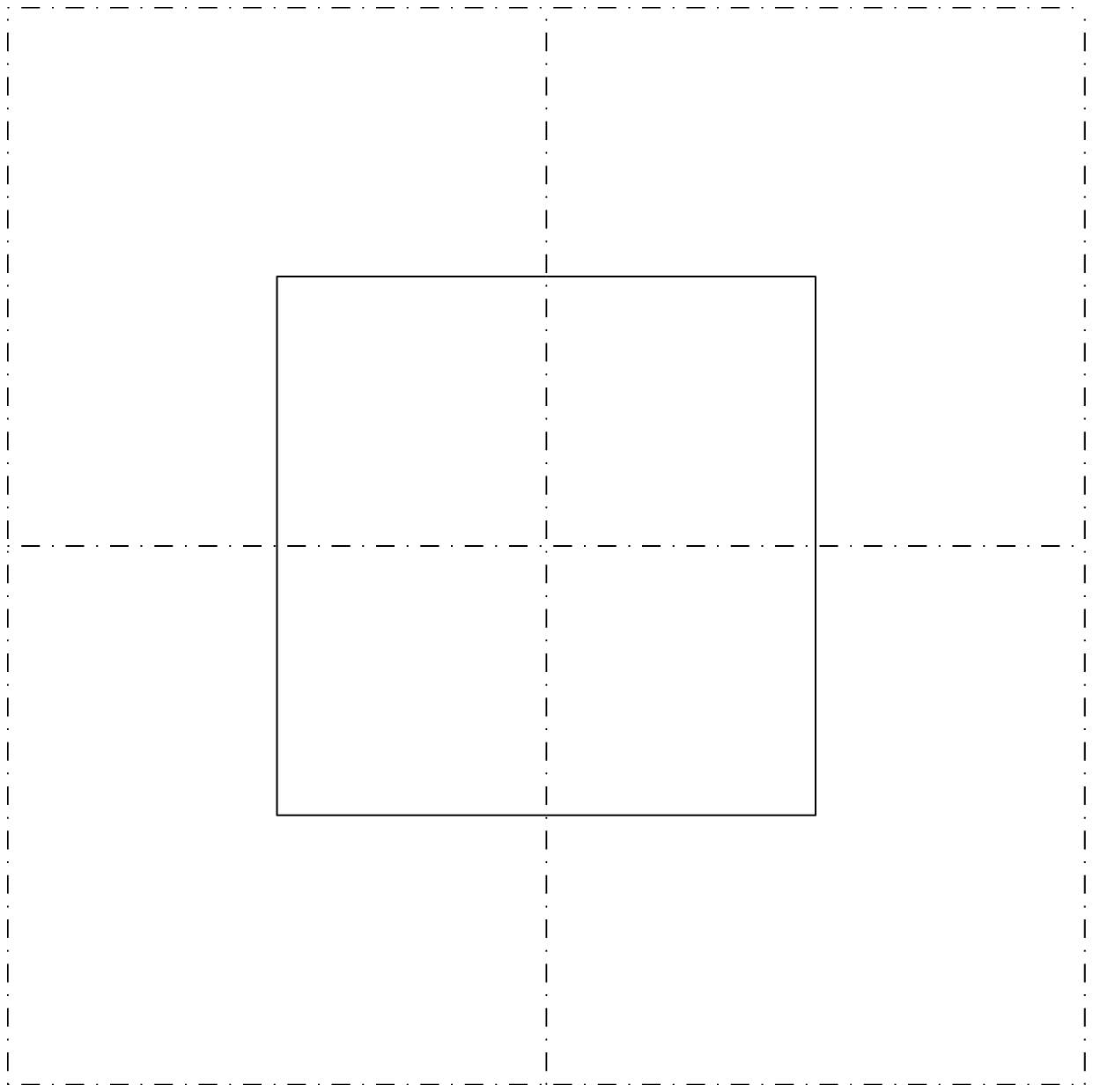, width=3.5 cm}
\end{center}
\caption{\label{circuits_complets_4_piecmax_ev3brio_new_saw}\iflanguage{french}{Le seul circuit à 4 pièces  tracé sous forme de polygone}{The sole circuit 
with 4 pieces,  
drawn under polygons form}.}
\end{figure}

\begin{figure}[h]
\begin{center}
\epsfig{file=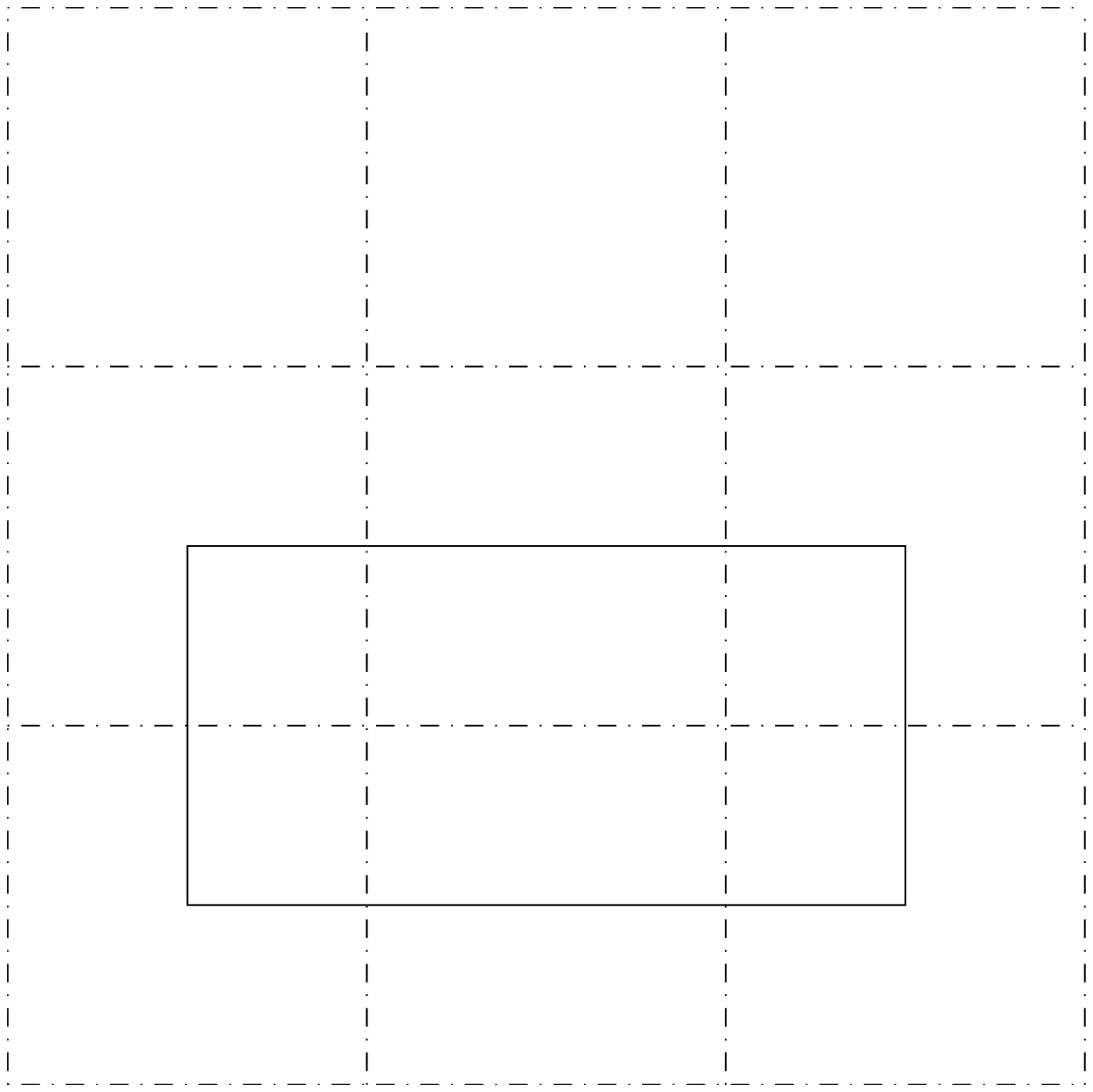, width=3.5 cm}
\end{center}
\caption{\label{circuits_complets_6_piecmax_ev3brio_new_saw}\iflanguage{french}{Le seul circuit à 6 pièces  tracé sous forme de polygone}{The sole circuit 
with 6 pieces,  
drawn under polygons form}.}
\end{figure}

\begin{figure}[h]
\centering
%%% sous figure 1
\subfigure[\label{circuits_complets_8_piecmax_ev3brio_new_saw1}]
{\epsfig{file=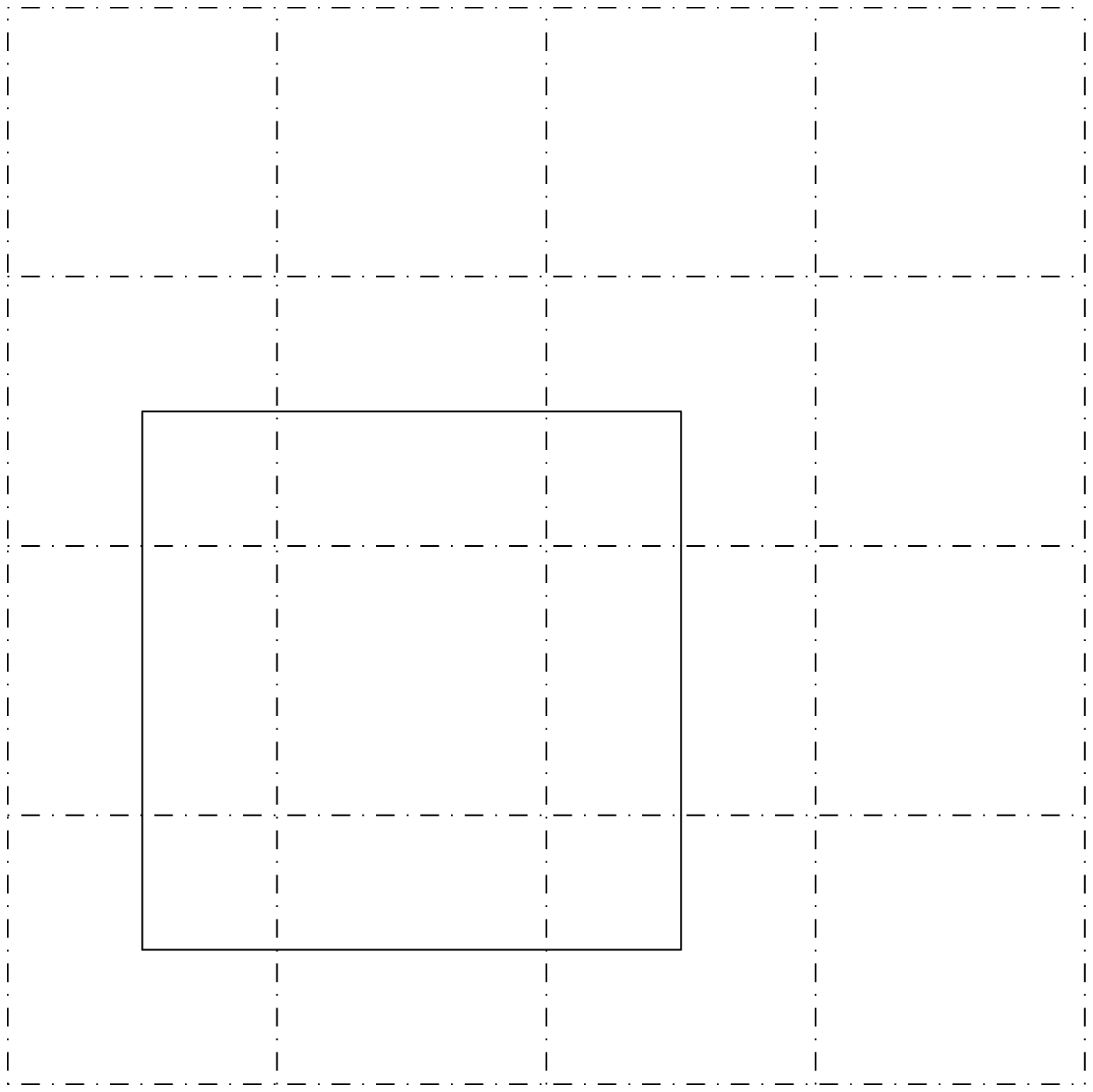, width=3.5 cm}}
\qquad
%%% sous figure 2
\subfigure[\label{circuits_complets_8_piecmax_ev3brio_new_saw2}]
{\epsfig{file=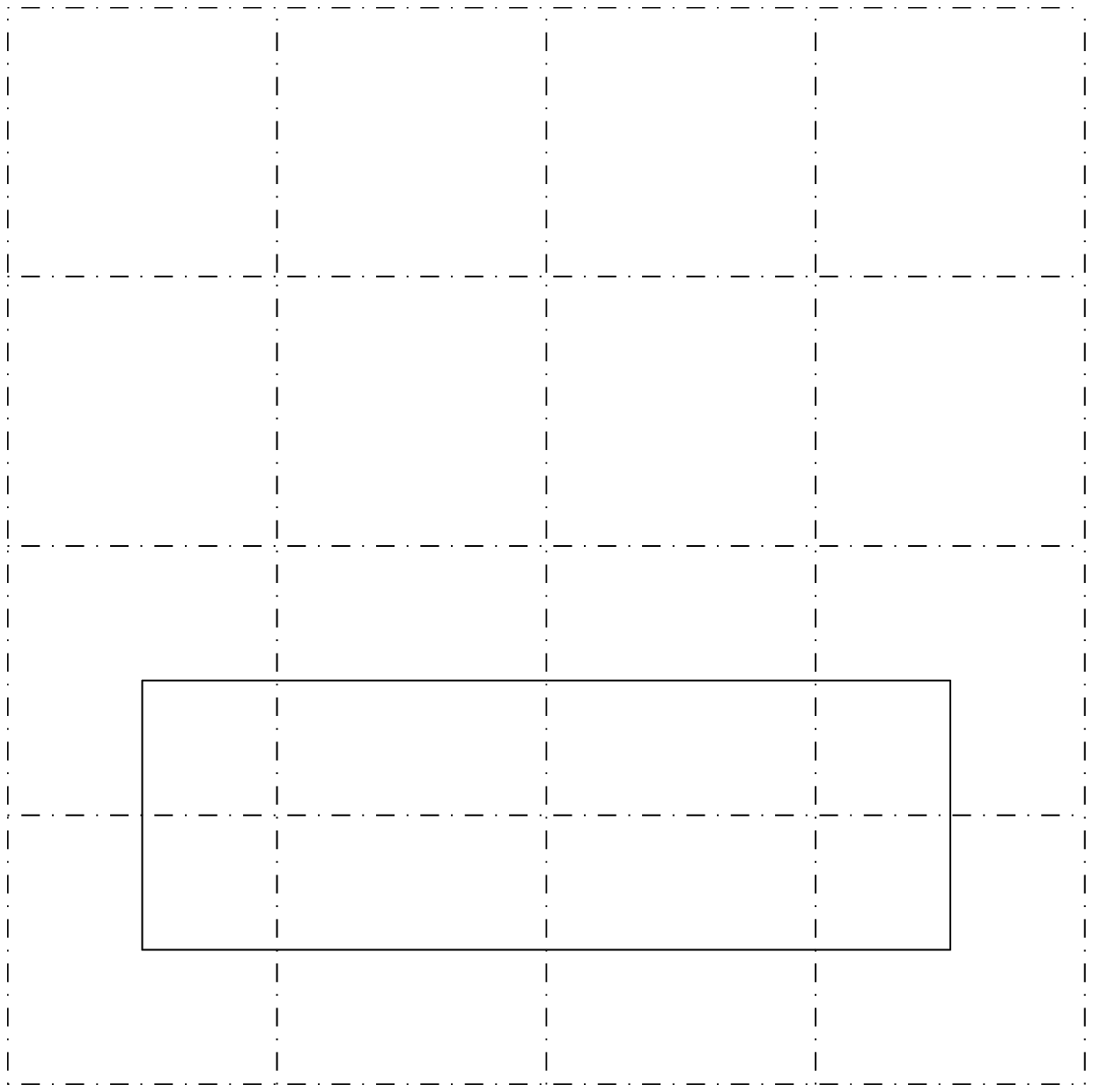, width=3.5 cm}}
\qquad
%%% sous figure 3
\subfigure[\label{circuits_complets_8_piecmax_ev3brio_new_saw3}]
{\epsfig{file=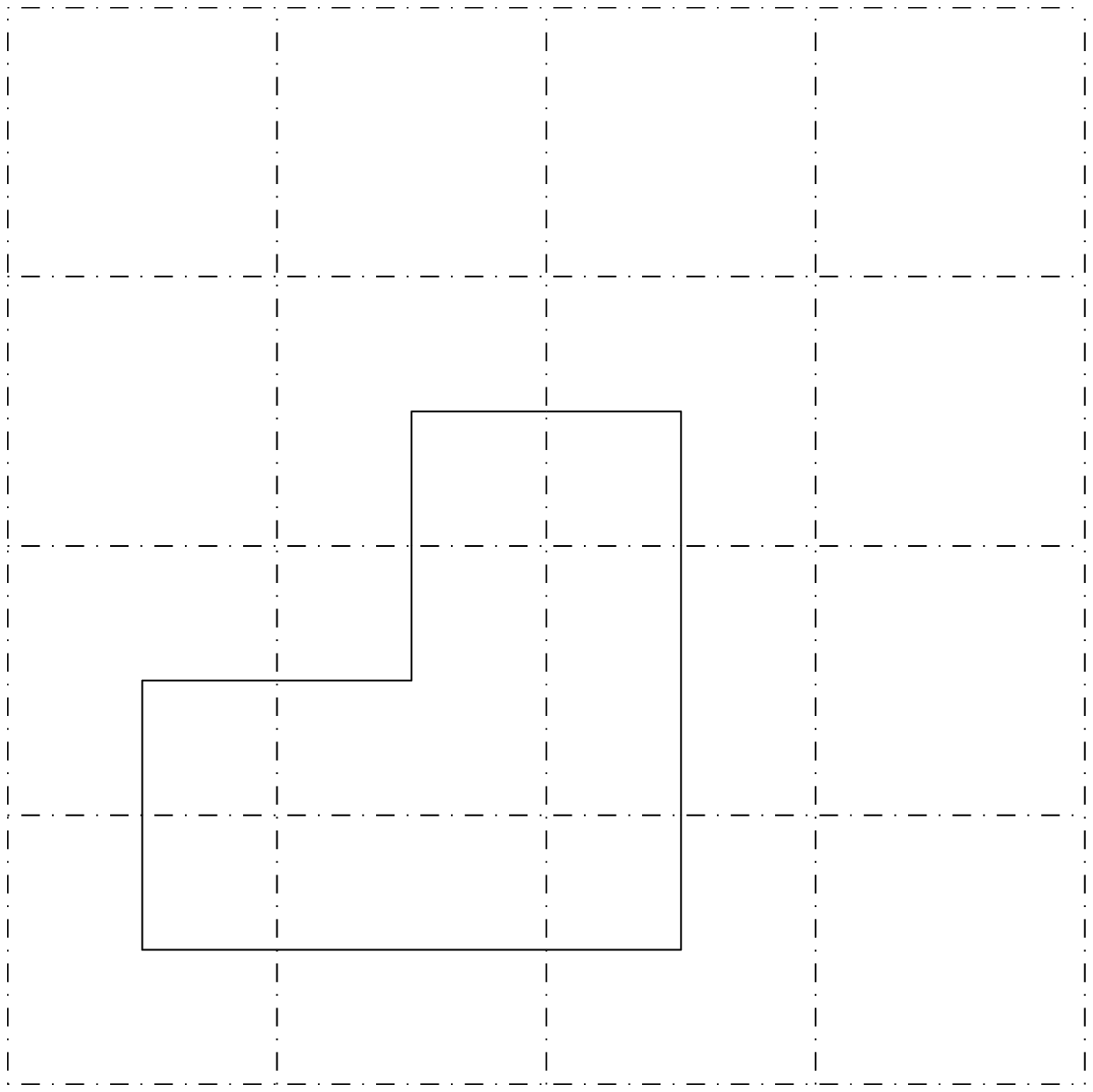, width=3.5 cm}}
\qquad
%%% sous figure 4
\subfigure[\label{circuits_complets_8_piecmax_ev3brio_new_saw4}]
{\epsfig{file=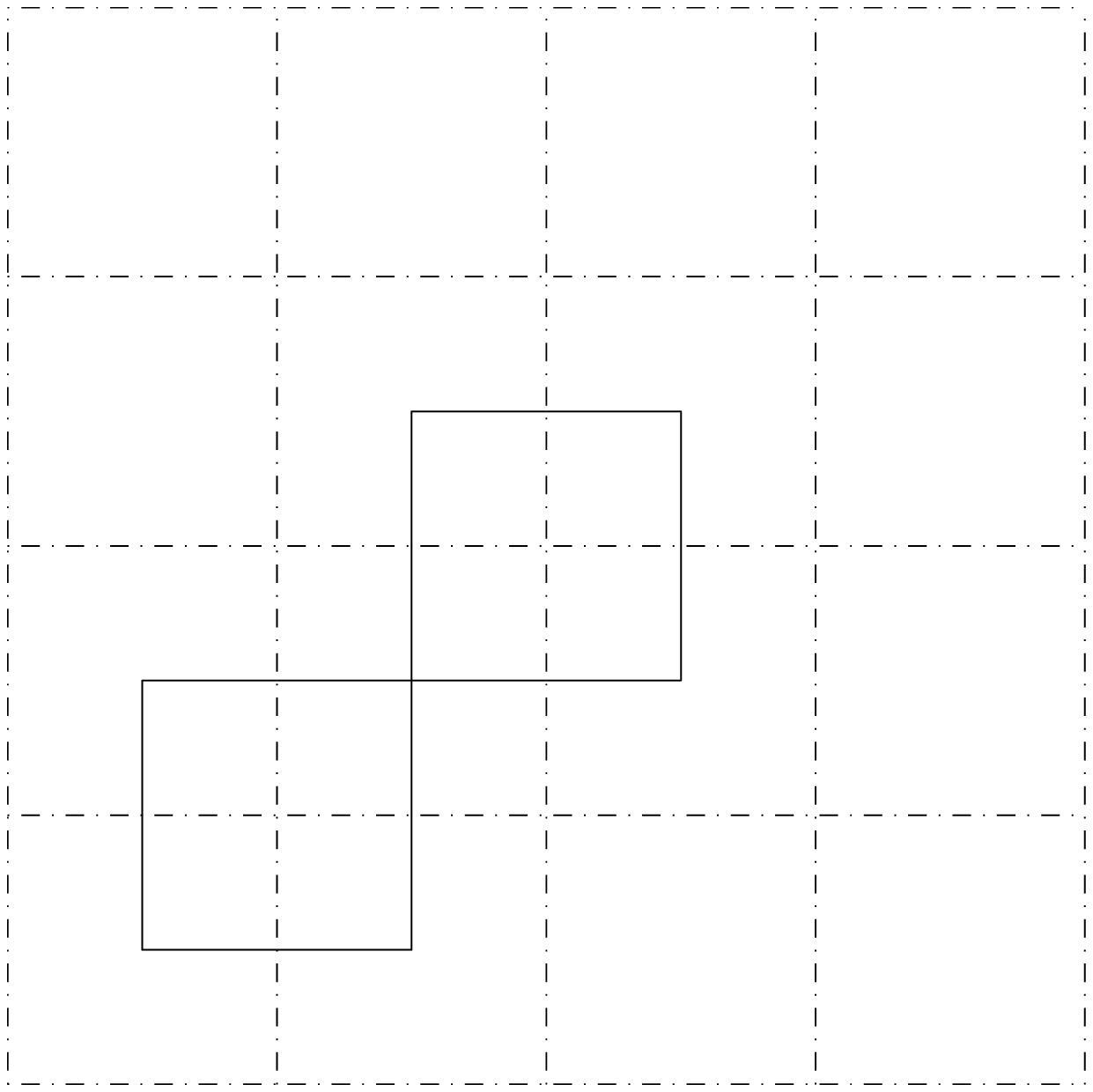, width=3.5 cm}}
\qquad
\caption{\label{circuits_complets_8_piecmax_ev3brio_new_saw}\iflanguage{french}{Les 4 circuits à 8 pièces  tracés sous forme de polygone}{The 4 circuits 
with 8 pieces,  
drawn under polygons form}.}
\end{figure}

\begin{figure}[h]
\centering
%%% sous figure 1
\subfigure[\label{circuits_complets_10_piecmax_ev3brio_new_saw1}]
{\epsfig{file=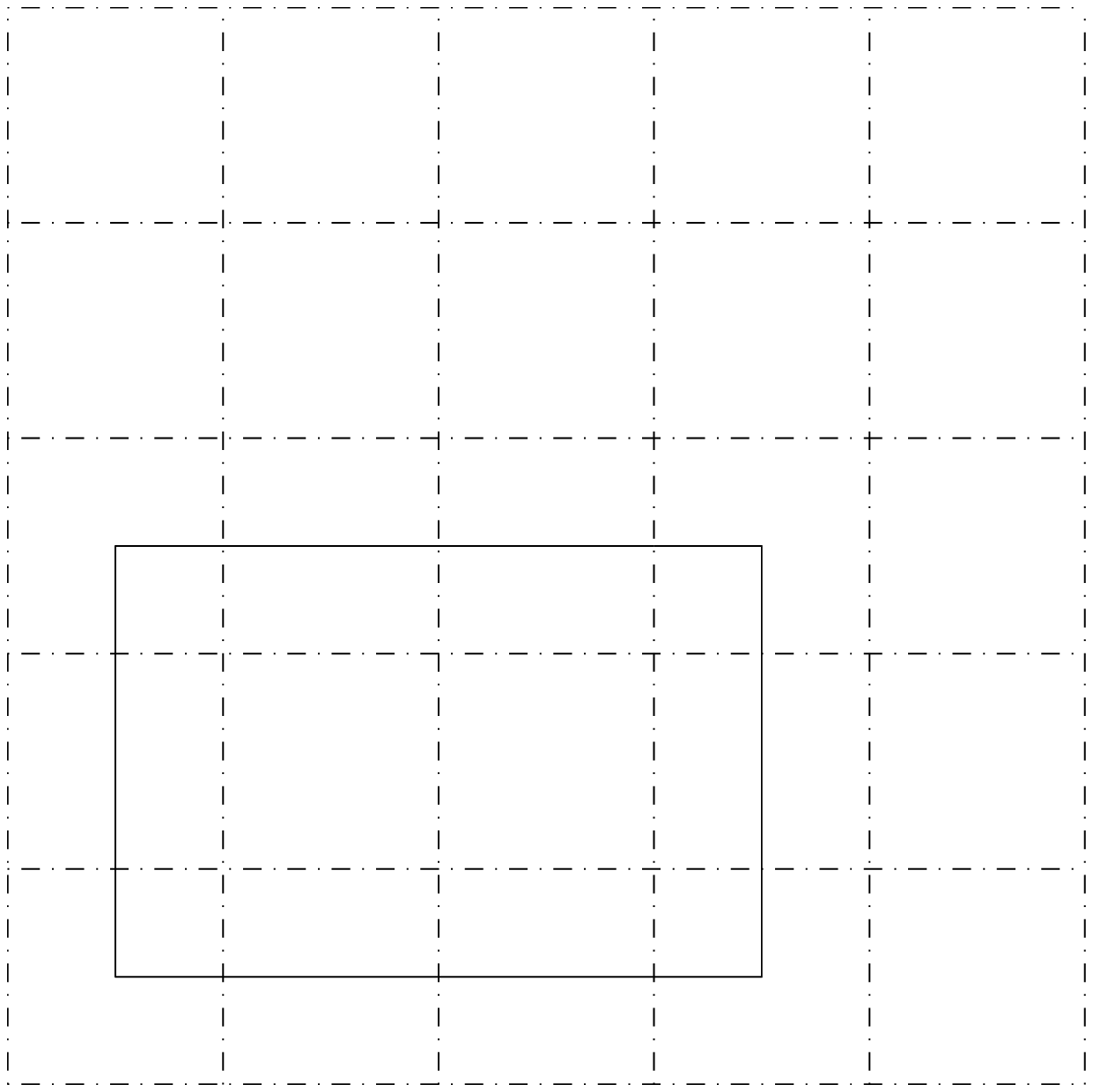, width=3.5 cm}}
\qquad
%%% sous figure 2
\subfigure[\label{circuits_complets_10_piecmax_ev3brio_new_saw2}]
{\epsfig{file=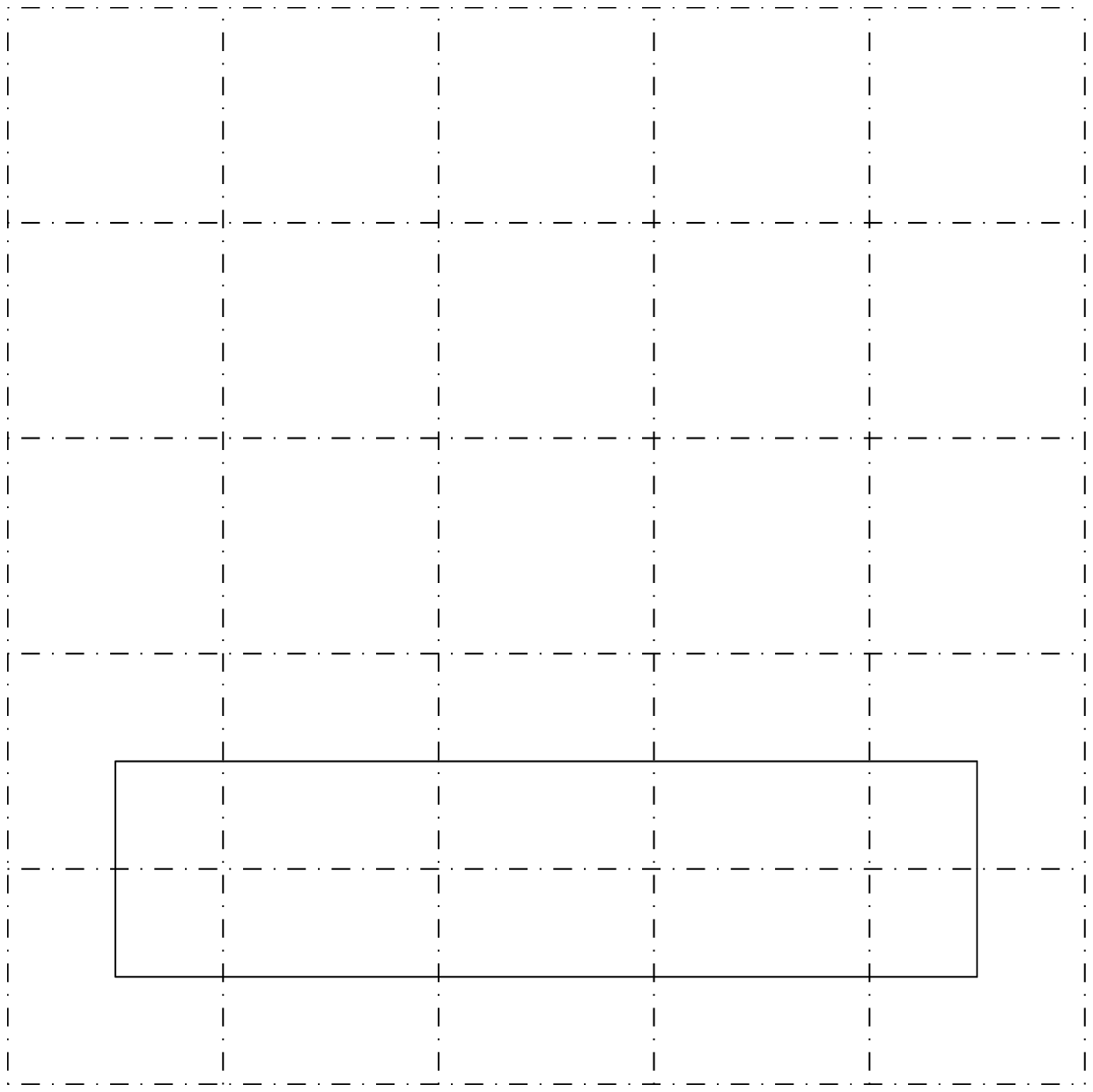, width=3.5 cm}}
\qquad
%%% sous figure 3
\subfigure[\label{circuits_complets_10_piecmax_ev3brio_new_saw3}]
{\epsfig{file=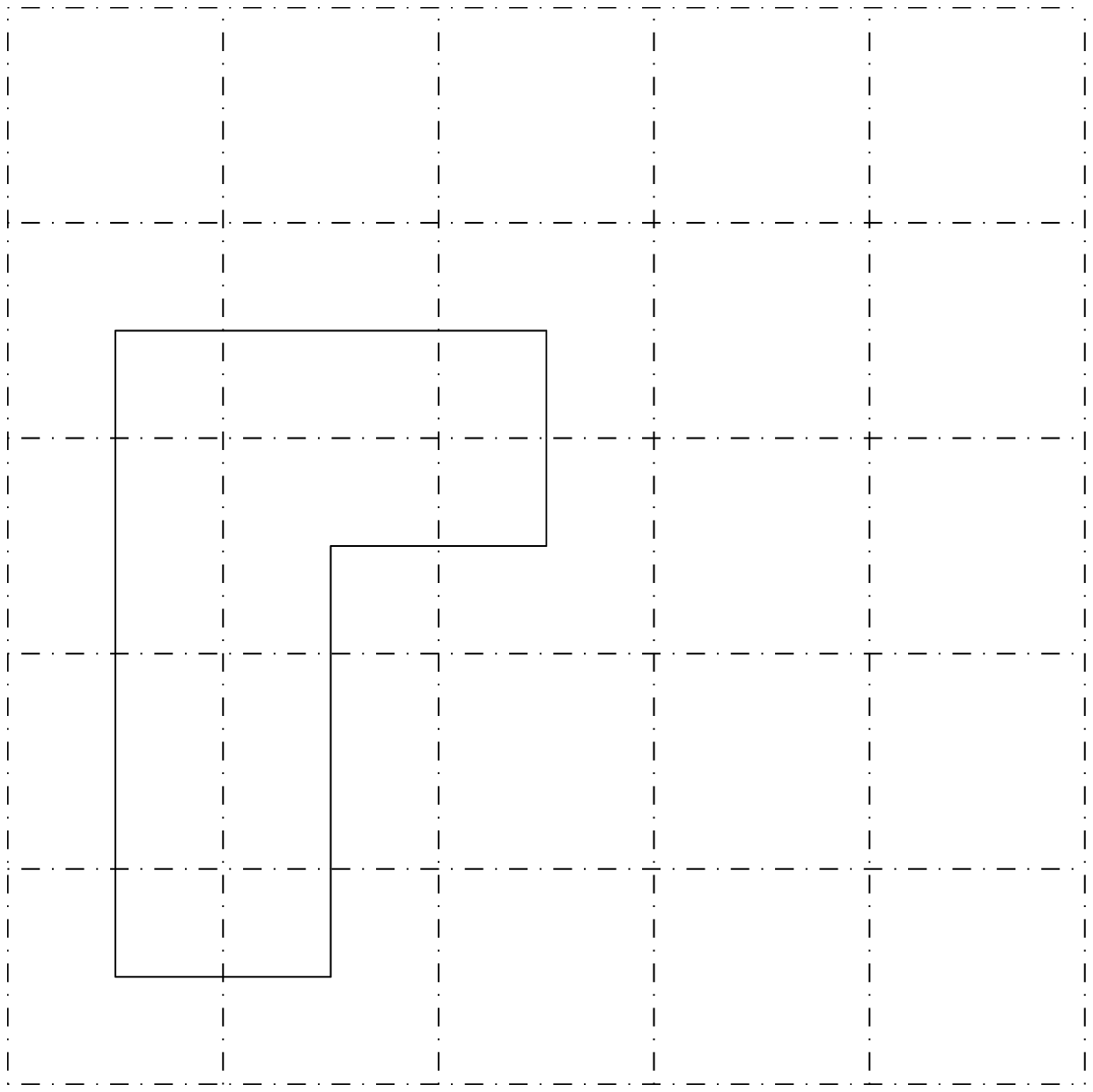, width=3.5 cm}}
\qquad
%%% sous figure 4
\subfigure[\label{circuits_complets_10_piecmax_ev3brio_new_saw4}]
{\epsfig{file=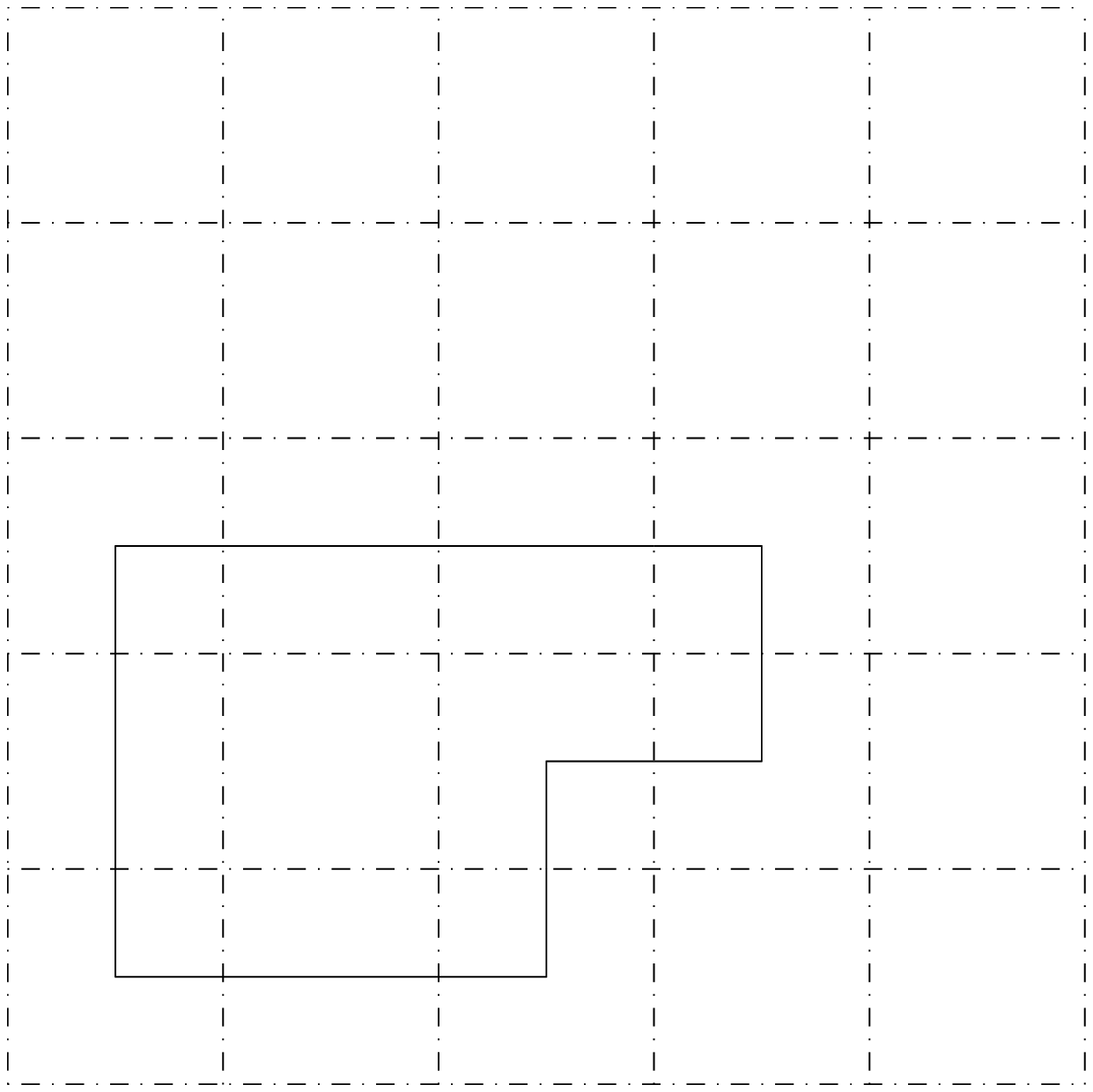, width=3.5 cm}}
\qquad
%%% sous figure 5
\subfigure[\label{circuits_complets_10_piecmax_ev3brio_new_saw5}]
{\epsfig{file=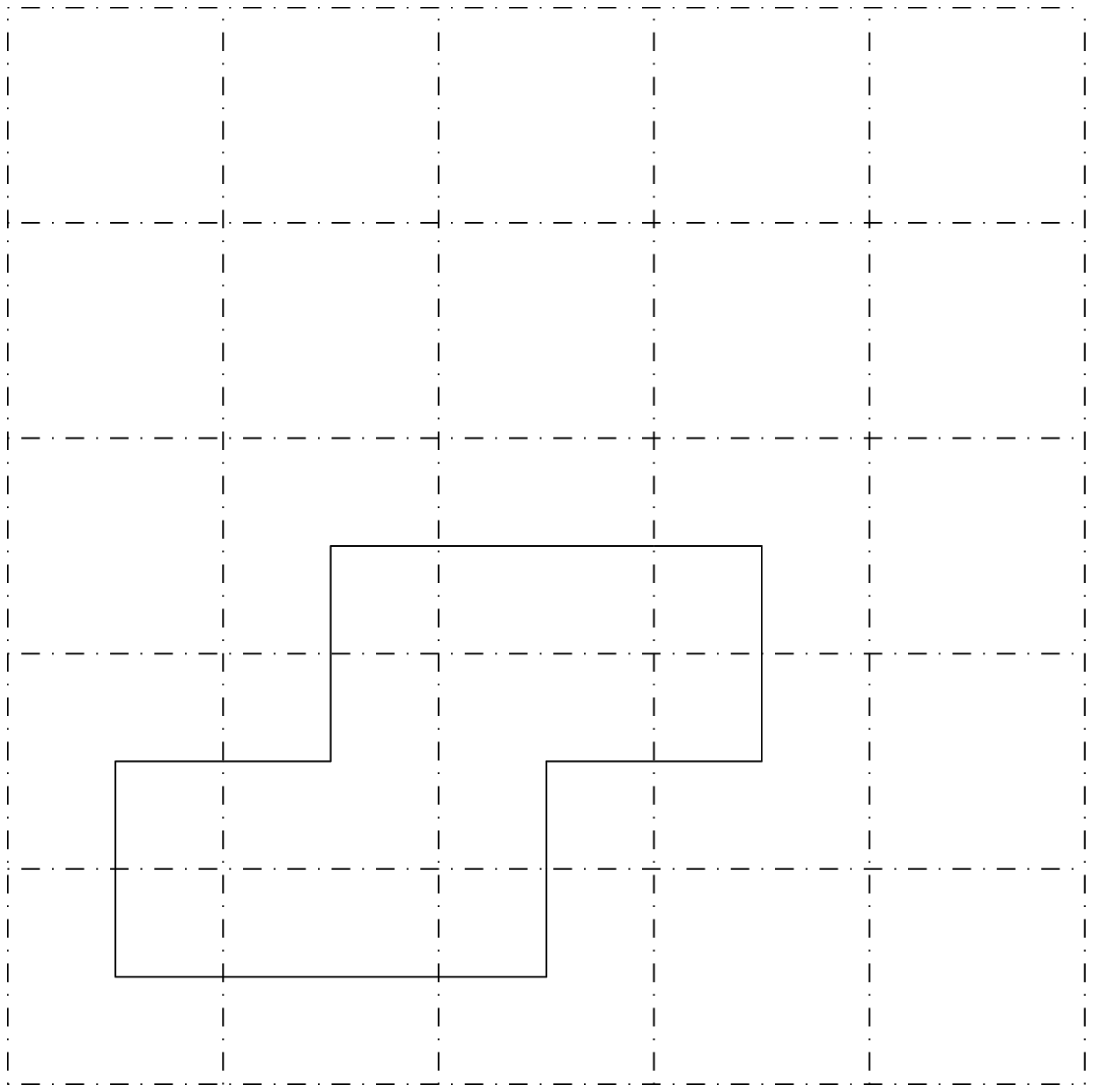, width=3.5 cm}}
\qquad
%%% sous figure 6
\subfigure[\label{circuits_complets_10_piecmax_ev3brio_new_saw6}]
{\epsfig{file=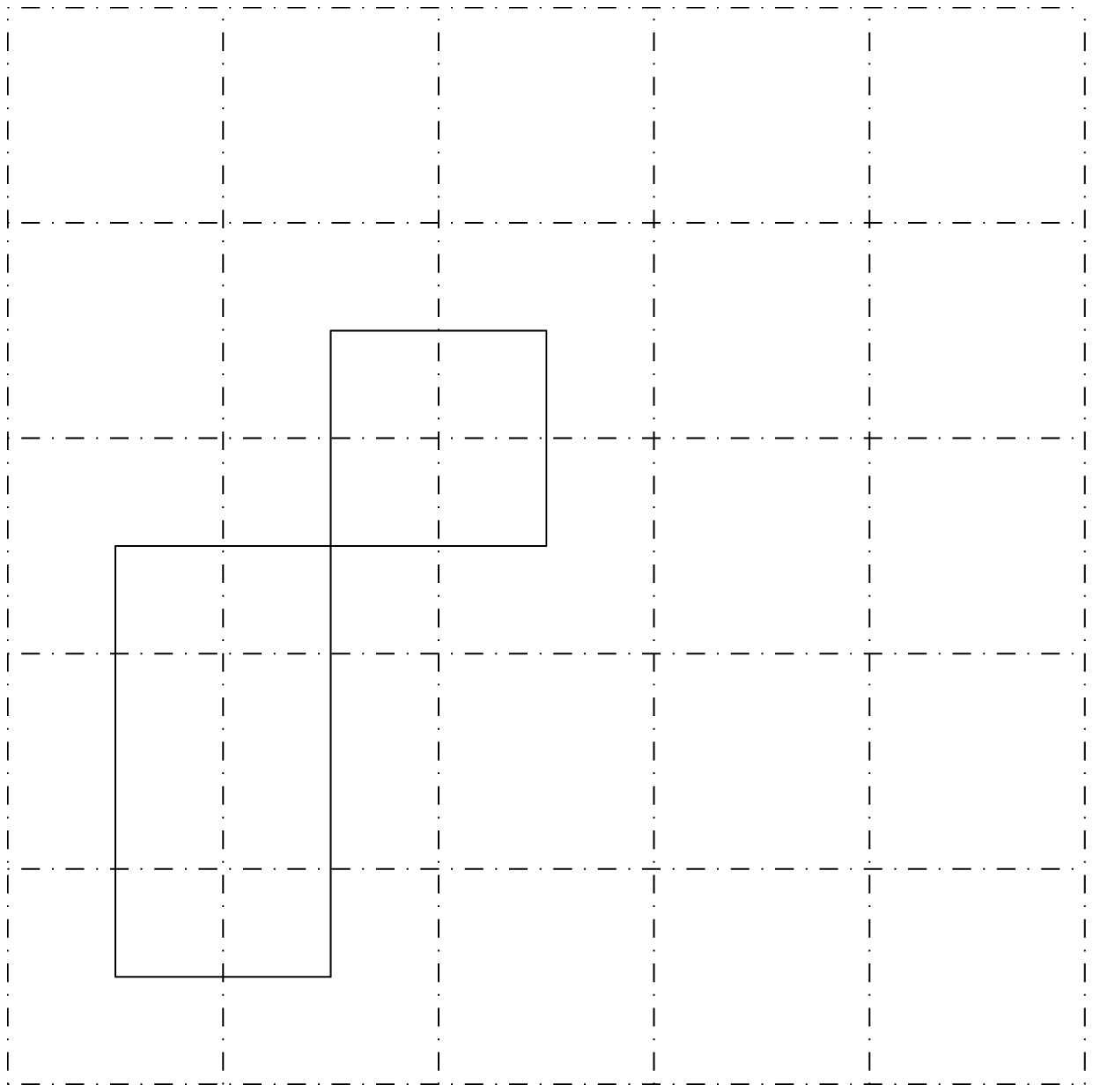, width=3.5 cm}}
\qquad
%%% sous figure 7
\subfigure[\label{circuits_complets_10_piecmax_ev3brio_new_saw7}]
{\epsfig{file=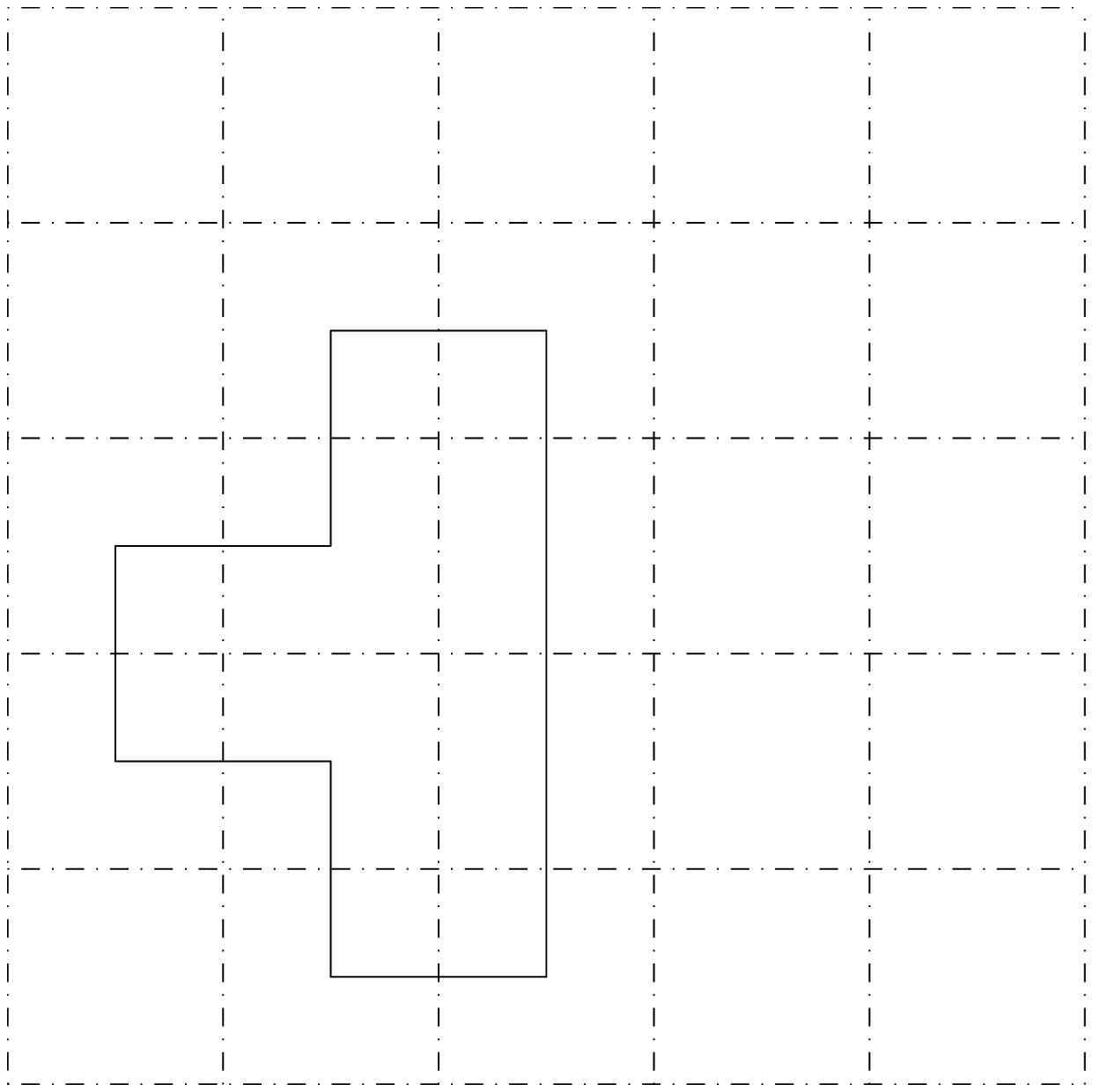, width=3.5 cm}}
\caption{\label{circuits_complets_10_piecmax_ev3brio_new_saw}\iflanguage{french}{Les 7 circuits à 10 pièces  tracés sous forme de polygone}{The 7 circuits 
with 10  pieces,  
drawn under polygons form}.}
\end{figure}

\iflanguage{french}{%
Traçons sous formes de polygones fermés
les circuits obtenus dans les exemples 
\ref{examplesimulation600},
\ref{examplesimulation601},
\ref{examplesimulation602} et 
\ref{examplesimulation603}.
Voir les figures respectives
\ref{circuits_complets_4_piecmax_ev3brio_new_saw},
\ref{circuits_complets_6_piecmax_ev3brio_new_saw},
\ref{circuits_complets_8_piecmax_ev3brio_new_saw} et 
\ref{circuits_complets_10_piecmax_ev3brio_new_saw}.

\begin{itemize}
\item
Pour $N=4$ (voir figure 
\ref{circuits_complets_4_piecmax_ev3brio_new_saw}),
on obtient dans les deux cas, 1 circuit.
\item
Pour $N=6$ (voir figure 
\ref{circuits_complets_6_piecmax_ev3brio_new_saw}), 
on obtient 1 circuit traditionnel et 2 polygones autoévitants, de formes rectangulaires. Ces deux derniers 
correspondent au circuit traditionnel et son image par une rotation d'angle $\pi/2$.
\item
Pour $N=8$ (voir figure 
\ref{circuits_complets_8_piecmax_ev3brio_new_saw}), 
on obtient 4 circuits traditionnels et 7  polygones autoévitants. 
Comme écrit dans 
\url{http://oeis.org/A002931},
les 7 polygones autoévitants correspondent aux 
1, 2  et 4 rotations (d'angle $\pi/2$) respectives des circuits des circuits des figures
\ref{circuits_complets_8_piecmax_ev3brio_new_saw1}
\ref{circuits_complets_8_piecmax_ev3brio_new_saw2} et 
\ref{circuits_complets_8_piecmax_ev3brio_new_saw3}. 
Le circuit de la figure 
\ref{circuits_complets_8_piecmax_ev3brio_new_saw4} ne correspond à aucun polygone autoévitant,
puisqu'un des carrés est occupé par deux pièces.
On a donc bien retrouvé $7=1+2+4$.
\item
Pour $N=10$ (voir figure 
\ref{circuits_complets_10_piecmax_ev3brio_new_saw}), 
on obtient 7 circuits traditionnels et 28  polygones autoévitants. 
En effet, 
les 2 circuits des figures 
\ref{circuits_complets_10_piecmax_ev3brio_new_saw1} et 
\ref{circuits_complets_10_piecmax_ev3brio_new_saw2}
fournissent chacun par 2 rotations (d'angle $\pi/2$) 2 polygones autoévitants.
Les 2 circuits des figures 
\ref{circuits_complets_10_piecmax_ev3brio_new_saw3} et 
\ref{circuits_complets_10_piecmax_ev3brio_new_saw4} 
fournissent chacun par 4 rotations (d'angle $\pi/2$) et une réflexion  8 polygones autoévitants.
Le circuit de la figure  
\ref{circuits_complets_10_piecmax_ev3brio_new_saw5}
fournit par une rotation (d'angle $\pi/2$) et  une réflexion, 4 polygones autoévitants..
Le  circuit de  la figures
\ref{circuits_complets_10_piecmax_ev3brio_new_saw7},
fournit   par 4 rotations (d'angle $\pi/2$) 4  polygones autoévitants.
Le circuit de la figure 
\ref{circuits_complets_10_piecmax_ev3brio_new_saw6} ne correspond à aucun polygone autoévitant,
puisqu'un des carrés est occupé par deux pièces.
On a donc bien retrouvé $28=2\times 2+2\times 2 \times 4+4 +4$ polygones autoévitants.
\end{itemize}%
}%
{%

We draw the circuits obtained in examples \ref{examplesimulation600},
\ref{examplesimulation601},
\ref{examplesimulation602} and 
\ref{examplesimulation603} in the form of closed polygons. See Figures
\ref{circuits_complets_4_piecmax_ev3brio_new_saw},
\ref{circuits_complets_6_piecmax_ev3brio_new_saw},
\ref{circuits_complets_8_piecmax_ev3brio_new_saw} and 
\ref{circuits_complets_10_piecmax_ev3brio_new_saw} respectively.

\begin{itemize}
\item
For $N=4$ (see Figure 
\ref{circuits_complets_4_piecmax_ev3brio_new_saw}),
we obtain in both cases, 1 circuit.
\item
For $N=6$ (see Figure 
\ref{circuits_complets_6_piecmax_ev3brio_new_saw}), 
we obtain 1 traditional circuit and two self-avoiding polygons, rectangular in shape. The latter two correspond to the traditional circuit and its image under a rotation by angle $\pi/2$.
\item
For $N=8$ (see Figure 
\ref{circuits_complets_8_piecmax_ev3brio_new_saw}), 
we obtain 4 traditional circuits and 7 self-avoiding polygons. 
As noted in
\url{http://oeis.org/A002931},
the 7 self-avoiding polygons correspond to the 1, 2 and 4 rotations (by angle $\pi/2$) of the circuits in Figures 
\ref{circuits_complets_8_piecmax_ev3brio_new_saw1}
\ref{circuits_complets_8_piecmax_ev3brio_new_saw2} and 
\ref{circuits_complets_8_piecmax_ev3brio_new_saw3} respectively. 
The circuit in Figure 
\ref{circuits_complets_8_piecmax_ev3brio_new_saw4} does not correspond to any self-avoiding polygon, since one of the squares is occupied by two pieces.
We have thus indeed found $7=1+2+4$.
\item
For $N=10$ (see Figure 
\ref{circuits_complets_10_piecmax_ev3brio_new_saw}), 
we obtain 7 traditional circuits and 28 self-avoiding polygons. 
In fact, the two circuits in Figures
\ref{circuits_complets_10_piecmax_ev3brio_new_saw1} and 
\ref{circuits_complets_10_piecmax_ev3brio_new_saw2}
each provide, by 2 rotations (by angle $\pi/2$), 2 self-avoiding polygons.
The 2 circuits in Figures
\ref{circuits_complets_10_piecmax_ev3brio_new_saw3} and 
\ref{circuits_complets_10_piecmax_ev3brio_new_saw4} 
each provide, by 4 rotations (by angle $\pi/2$) and one reflection, 8 self-avoiding polygons.
The circuit in Figure  
\ref{circuits_complets_10_piecmax_ev3brio_new_saw5}
provides, by one rotation (by angle $\pi/2$) and one reflection, 4 self-avoiding polygons.
The circuit in Figure
\ref{circuits_complets_10_piecmax_ev3brio_new_saw7},
provides, by 4 rotations (by angle $\pi/2$), 4 self-avoiding polygons.
The circuit in Figure 
\ref{circuits_complets_10_piecmax_ev3brio_new_saw6} does not correspond to any self-avoiding polygon, since one of the squares is occupied by two pieces.
We have thus indeed found $28=2\times 2+2\times 2 \times 4+4 +4$ self-avoiding polygons.
\end{itemize}%
}

%%%%%%%%%%%%%%%%%%%%%%%%%%%%%%%%%%%%%%%%%%%%%%%%%%%

If one now compares the last column of Table \ref{enumeration_construction_circuit_tab_new03} with the first two, one notices that the comparison is then in favor of the \textit{Easyloop} system.%

The number of operations is
\begin{equation}
\label{nombreoperation}
\mathcal{O}\left(5^N\right)
\end{equation}
Thus, passing from $N$ to $N+1$ multiplies the number of operations by $5$. All of the results in this section were realized
\ifcase 0
 in a time greater than   $23$ hours.
Replacing $11$ by 
$12$
therefore demands more than $4$  days of calculation!%
\or
\textbf{ATTENTION  !!!!!!!!!!!!!!!!!}
Durée inférieure à une heure !!%
\fi
}

\iflanguage{french}{%

%%%%%%%%%%%%%%%%%%%%%%%%%%%%%%%%%%%%%%%%%%%%%%%%%%%%%%%%%%%%%
%%%%%%%%%%%%%%%%%%%%%%%%%%%%%%%%%%%%%%%%%%%%%%%%%%%%%%%%%%%%%
\section{Estimation du nombre de circuits à nombre de pièces important}
\label{estimation}%
}{%

%%%%%%%%%%%%%%%%%%%%%%%%%%%%%%%%%%%%%%%%%%%%%%%%%%%%%%%%%%%%%
%%%%%%%%%%%%%%%%%%%%%%%%%%%%%%%%%%%%%%%%%%%%%%%%%%%%%%%%%%%%%
\section{Estimation of the number of circuits with a significant number of pieces}
\label{estimation}%
}

%%% Vérifier passage de 11 à 12

\iflanguage{french}{%
Dans le cas où $N$ est supérieur à  $11$,
les calculs sont trop longs et il n'est pas possible d'utiliser l'énumération et le dénombrement 
des circuits.
Nous utilisons, avec un grand abus, l'estimation donnée par 
\eqref{estime_nombre_circuit}, qui provient de  
\cite{%
MR2883859,%
MR1985492,%
MR1718791,%
Guttmann2012,%
Guttmann2012b,%
MR2902304}.
Nous allons utiliser les nombres de circuits
donnés dans la section 
\ref{enumerationpetit}
 et nous nous en servirons pour évaluer les valeurs des constantes 
$A$, $\mu$ et $\gamma$
de la formule
\eqref{estime_nombre_circuit}. Cette évaluation est remplacée par une égalité et les coefficients $A$, $\mu$ et $\gamma$
sont déterminés grâce  à un système aux moindres carrés
qui devient linéaire si l'on en prend le logarithme.%
}{%

In the case where $N$ is greater than $11$, the calculations take too long, and it is not possible to use the enumeration of the circuits. We use, with much impropriety, the estimation given by \eqref{estime_nombre_circuit}, which comes from  
\cite{%
MR2883859,%
MR1985492,%
MR1718791,%
Guttmann2012,%
Guttmann2012b,%
MR2902304}.
We will take the number of circuits given in Section~\ref{enumerationpetit}, and we will make use of it to evaluate the constants $A$, $\mu$ and $\gamma$ in formula~\eqref{estime_nombre_circuit}. This evaluation is replaced by an equality and the coefficients  $A$, $\mu$ and $\gamma$ are determined by solving  a least squares system, which becomes linear when we take the logarithm.%
}

% table crée par stocke_tableau_tex le 29-Mar-2016 14:34:39
\begin{table}[h]
\begin{center}
\begin{tabular}{|l||l|l|l|l|}
\hline
&looping circuits&direct isom.&isometries&constructible
\\
\hline
\hline
$\gamma-1$ &$-4.41060$ & $-10.08120$ & $-9.74233$ & $-8.75998$
\\
\hline
$\mu$ &$7.68523$ & $13.10739$ & $11.22176$ & $9.13739$
\\
\hline
\end{tabular}
\vspace{1 cm}
\caption{\label{enumeration_construction_circuit_tab_new10}
The values of $\gamma-1$ and $\mu$ obtained using Table~\ref{enumeration_construction_circuit_tab_new01} (corresponding to $N_j=+\infty$).}
\end{center}
\end{table}

%\iflanguage{french}{\input{enumeration_construction_circuit_tab_new10}}{\input{enumeration_construction_circuit_eng_tab_new10}}

\iflanguage{french}{%
Utilisons les différents résultats du tableau 
\ref{enumeration_construction_circuit_tab_new01} (correspondant à $N_j=+\infty$) pour déterminer les valeurs de $\gamma-1$ et de $\mu$.
Voir tableau~\ref{enumeration_construction_circuit_tab_new10}.
Les valeurs obtenues sont différentes naturellement des valeurs données
par \eqref{valeurmugamma}.
Cela est normal puisque l'estimation \eqref{estime_nombre_circuit} a été remplacée par une égalité et naturellement, rien ne valide \textit{a priori} cette égalité.
Notons que les signes de ces coefficients sont cohérents avec ceux de la littérature.%
}{%
Let us now estimate $\gamma-1$ and $\mu$ using the different results from Table~\ref{enumeration_construction_circuit_tab_new01} (corresponding to $N_j=+\infty$). See Table~\ref{enumeration_construction_circuit_tab_new10}. 
The obtained values are naturally different from the values given in \eqref{valeurmugamma}. This is normal since the estimation \eqref{estime_nombre_circuit} has been replaced by an equality, and naturally, nothing \textit{a priori} validates this equality.
We note that the signs of the coefficients are consistent with those in the literature.
%It is interesting to note that the values of $\mu$ are smaller than the factor of $5$ appearing in~\eqref{nombreoperation}, which is normal since the number of circuits retained is less than the number of circuits examined. On the other hand, when the constraints increase, the number $\mu$ approaches the value given by~\eqref{valeurmu}, being very close to it in the case where only constructible circuits are retained, up to an isometry. This proximity in value is troubling in this latter case. Subsequent use of formula~\eqref{estime_nombre_circuit}, coming from the study of self-avoiding walks, is not at all justified in theory in the context of  \textit{Easyloop} circuits, but the proximity of the values of $\mu$ with those in the literature will be the only, \textit{a posteriori}, justification of this formula.%
}

% table crée par stocke_tableau_tex le 29-Mar-2016 14:34:39
\begin{table}[h]
\begin{center}
\begin{tabular}{|l||l|l|l|l|}
\hline
&looping circuits&direct isom.&isometries&constructible
\\
\hline
\hline
$\gamma-1$ &$-3.93476$ & $-9.65660$ & $-9.76067$ & $-8.69817$
\\
\hline
$\mu$ &$6.80151$ & $11.68005$ & $10.61788$ & $8.69023$
\\
\hline
\end{tabular}
\vspace{1 cm}
\caption{\label{enumeration_construction_circuit_tab_new11}
The values of $\gamma-1$ and $\mu$ obtained using Table~\ref{enumeration_construction_circuit_tab_new01} (corresponding to $N_j=4$).}
\end{center}
\end{table}

%\iflanguage{french}{\input{enumeration_construction_circuit_tab_new11}}{\input{enumeration_construction_circuit_eng_tab_new11}}

\iflanguage{french}{%
Si l'on utilise maintenant l'estimation de $\mu$ en utilisant les différents résultats du tableau 
\ref{enumeration_construction_circuit_tab_new02} (correspondant à $N_j=4$), 
on obtient les valeurs données dans le tableau \ref{enumeration_construction_circuit_tab_new11}, avec des observations
similaires aux précédentes.}{%
If we now do the estimation of $\mu$ using the different results in Table~\ref{enumeration_construction_circuit_tab_new02} (corresponding to $N_j=4$), we obtain the values given in Table~\ref{enumeration_construction_circuit_tab_new11}, with some observations similar to the previous ones.%
}

\iflanguage{french}{%
Pour valider notre méthode, nous pouvons donner les 
valeurs retenues des coefficients $A$, $\mu$ et $\gamma-1$ de la formule \eqref{estime_nombre_circuit}
en utilisant les valeurs issues de l'identification obtenue en prenant la première colonne du tableau
\ref{enumeration_construction_circuit_tab_new03}.
Elles sont données par 
\begin{equation*}
A=2.1356,\quad
\gamma-1=-3.87223,\quad
\mu=3.15940,
\end{equation*}
qui sont de l'ordre des valeurs données par \eqref{valeurmugamma}.
Pour $N=24$, on obtient une valeurs du nombre de \slws\ donnée par 
$102070280$, au lieu de $2$ $521$ $270$, donnée par la littérature.%
}{%
To validate our method, we can give the retained values of the coefficients $A$, $\mu$ and $\gamma-1$ in formula \eqref{estime_nombre_circuit}
by using the values from the identification obtained by taking the first column of Table \ref{enumeration_construction_circuit_tab_new03}.
They are given by
\begin{equation*}
A=2.1356,\quad
\gamma-1=-3.87223,\quad
\mu=3.15940,
\end{equation*}
which are on the order of the values obtained by \eqref{valeurmugamma}.
For $N=24$, we obtain a value of $102070280$ for the number of \slws\ , instead of the value $2$ $521$ $270$, given in the literature.%
}

\iflanguage{french}{%
Les valeurs retenues des coefficients $A$, $\mu$ et $\gamma-1$ de la formule \eqref{estime_nombre_circuit}
dans le cas des boîtes \textit{Easyloop}, correspondant à $N_j=4$
correspondent donc à la dernière case du tableau \ref{enumeration_construction_circuit_tab_new11} et sont données
par    
\begin{equation}
\label{refalpha}
A=4.5900 \, 10^{1},\quad
\gamma-1=-8.69817,\quad
\mu=8.69023.
\end{equation}%
}{%
The retained values of the coefficients $A$, $\mu$ and $\gamma-1$ in formula~\eqref{estime_nombre_circuit} in the case of the \textit{Easyloop} boxes, corresponding to $N_j=4$, therefore correspond to the last case in Table~\ref{enumeration_construction_circuit_tab_new11} and are given by
\begin{equation}
\label{refalpha}
A=4.5900 \, 10^{1},\quad
\gamma-1=-8.69817,\quad
\mu=8.69023.
\end{equation}%
}

\iflanguage{french}{%
Dans le cas \eqref{refalpha}, 
les nombres estimés de circuits sont alors
 $\left\{0,0,0,2,2,3,8,21,65,226,857\right\}$ 
proches des nombres exacts de circuits ($\left\{0,0,0,2,1,5,6,28,63,244,753\right\}$).
Pour $N=24$, nous obtenons 
}{%
In case~\eqref{refalpha}, 
the estimated numbers of circuits are then $\left\{0,0,0,2,2,3,8,21,65,226,857\right\}$, which is close to the exact numbers of circuits ($\left\{0,0,0,2,1,5,6,28,63,244,753\right\}$). For $N=24$, we obtain%
}
\begin{equation}
\label{estimationeq02}
q(24)\approx
1560511691458.
\end{equation}%
\begin{figure}
\centering
%%% sous figure 1
\subfigure%[\label{fig:02a}]
{\iflanguage{french}{\epsfig{file=estime_nombre_circuit_fig_easy_loop_new01A.eps, width=12cm}}{\epsfig{file=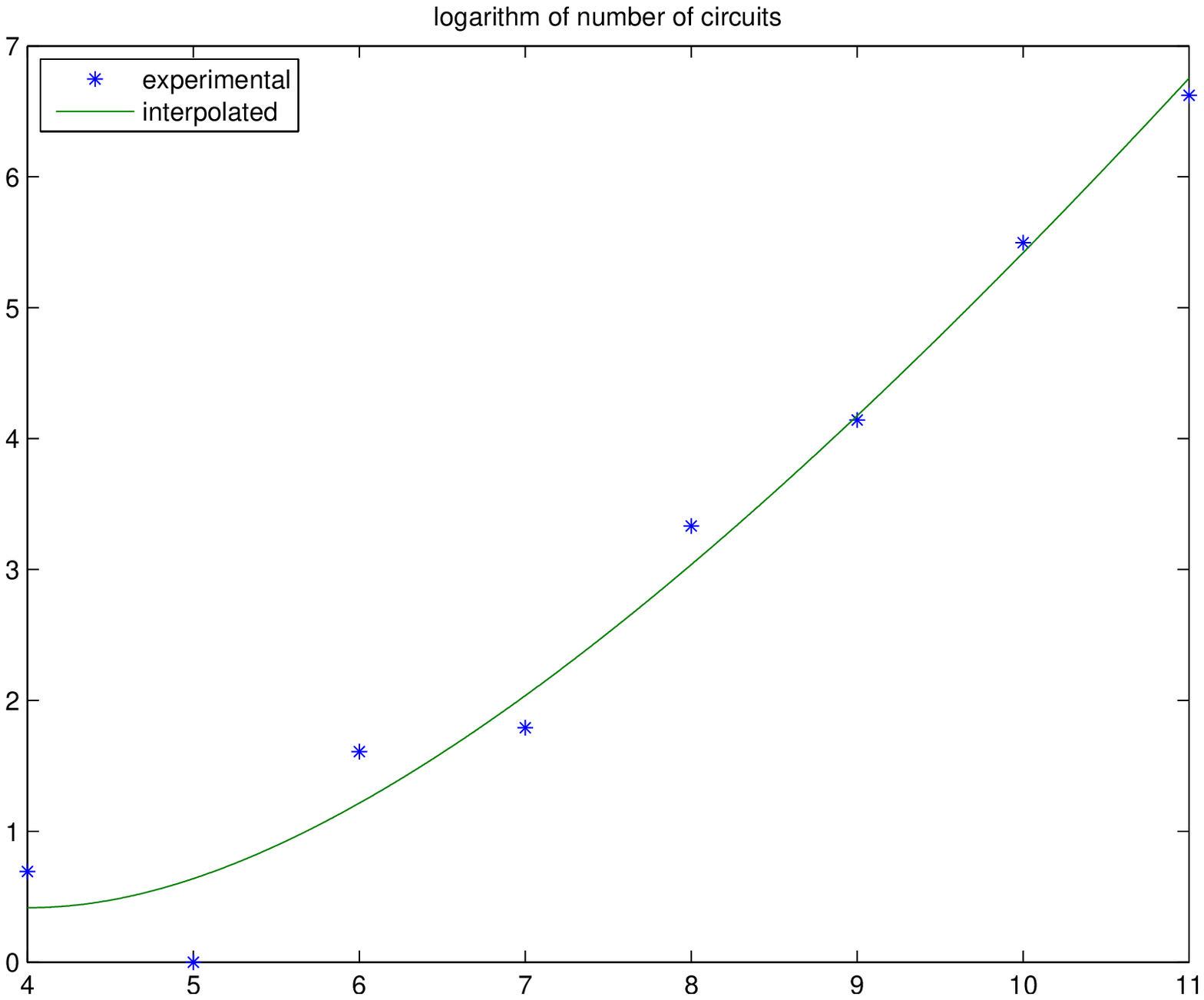, width=12cm}}}
\qquad
%%% sous figure 2
\subfigure%[\label{fig:02b}]
{\iflanguage{french}{\epsfig{file=estime_nombre_circuit_fig_easy_loop_new01B.eps, width=12cm}}{\epsfig{file=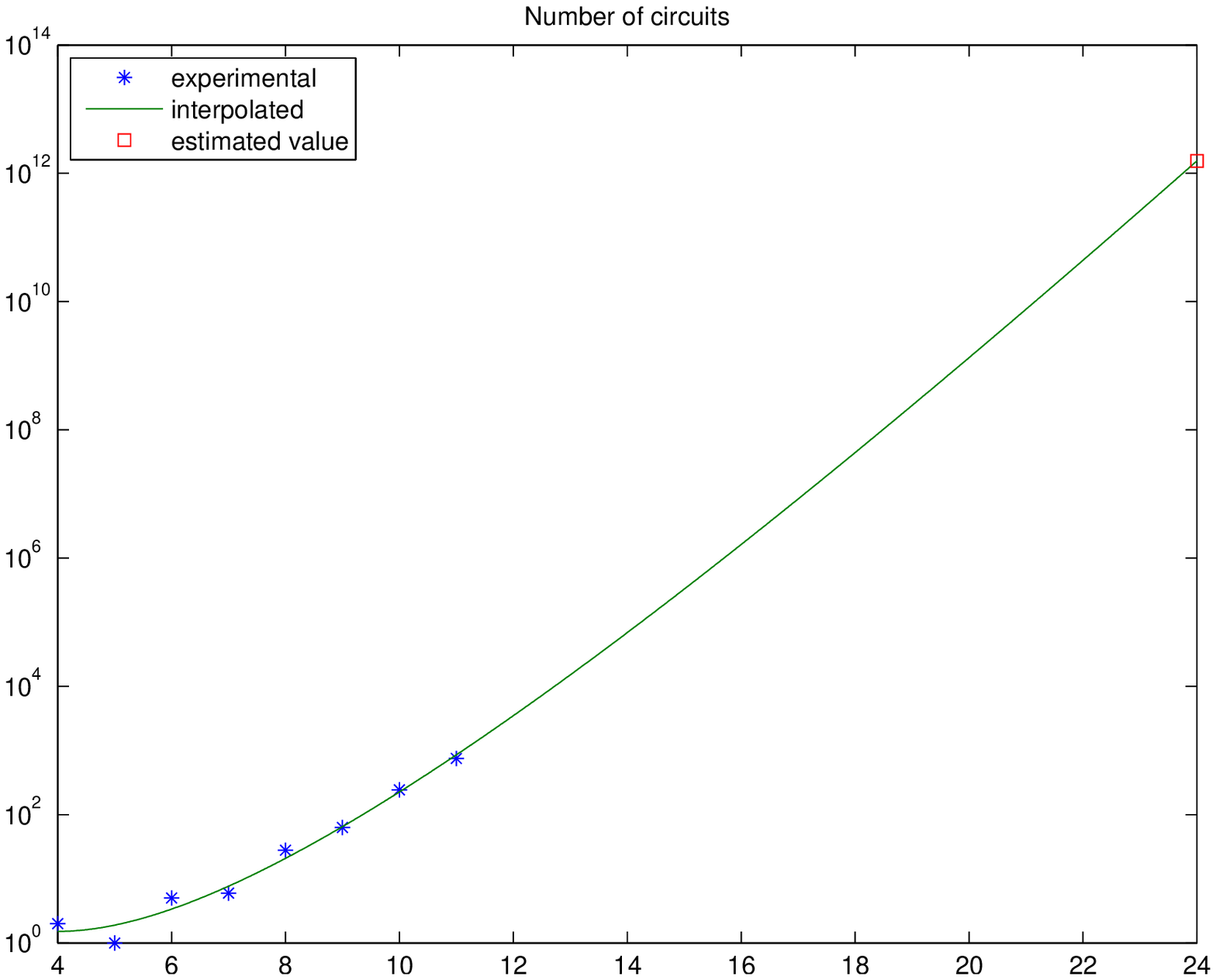, width=12cm}}}
\caption{\label{estime_nombre_circuit_fig_easy_loop_new01tot}\iflanguage{french}{Estimation du nombre de circuits dans le cas $N_j=4$.}{Estimation of the number of circuits in the case $N_j=4$.}}
 \end{figure}
\iflanguage{french}{%
On obtient alors
les courbes
montrées en figure \ref{estime_nombre_circuit_fig_easy_loop_new01tot}.}{We then obtain the curves shown in Figure~\ref{estime_nombre_circuit_fig_easy_loop_new01tot}.}

\iflanguage{french}{%
Pour comparer les circuits \textit{Easyloop} avec les systèmes traditionnels,  on peut de même utiliser
les nombres de la dernière colonne 
du tableau \ref{enumeration_construction_circuit_tab_new03}, ce qui fournit
\begin{equation}
\label{eqmu}
\mu=2.52246, 
\end{equation}
$\gamma-1=-3.71496$ et un nombre de circuits estimé
égal à 
\begin{equation}
\label{estimationeq03}
q(24)\approx
130229,
\end{equation}
ce qui reste  beaucoup plus faible que \eqref{estimationeq02}.%
Notons que le cas \eqref{eqmu}, 
le plus proche de celui des polygones autoévitants, fournit des valeurs relativement proches de celles données par \eqref{valeurmu}.%
}{%
To compare the \textit{Easyloop} circuits with traditional systems, one may likewise use the numbers in the last column of Table~\ref{enumeration_construction_circuit_tab_new03}, which gives 
\begin{equation}
\label{eqmu}
\mu=2.52246, 
\end{equation}
$\gamma-1=-3.71496$ and an estimated number of circuits equal to 
\begin{equation}
\label{estimationeq03}
q(24)\approx
130229,
\end{equation}
which is  still much smaller than~\eqref{estimationeq02}.% 
We note that the  case \eqref{eqmu}, the closest to that of self-avoiding polygons, provides values relatively close to those given by \eqref{valeurmu}.}

\iflanguage{french}{%

%%%%%%%%%%%%%%%%%%%%%%%%%%%%%%%%%%%%%%%%%%%%%%%%%%%%%%%%%%%%%
%%%%%%%%%%%%%%%%%%%%%%%%%%%%%%%%%%%%%%%%%%%%%%%%%%%%%%%%%%%%%
\section{Construction aléatoire de circuits à nombre de pièces important}
\label{constructionaleat}

Nous sommes capables pour les valeurs de $N$ inférieures
à  $11$, 
d'obtenir tous les circuits constructibles et en particulier de les exhiber. 
Un autre objectif d'un industriel serait de proposer un catalogue de 
plans de circuits de train pouvant contenir des circuit avec $N$ quelconque.
Malheureusement, au delà de $N=11$, 
cela n'est plus envisageable.
Une conception manuelle\footnote{voir les exécutables distribués pour windows, page \pageref{execatalogueweb}.}
est possible, mais fastidieuse et non programmable.
Nous proposons ainsi dans cette section un moyen de générer automatiquement des 
circuits pour $N$ et $N_j$ donnés pour des valeurs plus grandes que 
$N=11$, sans avoir à créer
tous les circuits possibles, comme proposé en section 
\ref{enumerationtouscircuit}, en s'appuyant sur une méthode aléatoire.

Pour $N\geq 1$ donné, on considère $r,s\in \En^*$ tels que $N=r+s$.
On est capable de déterminer tous les circuits à $r$ pièces partant de l'origine 
en décrivant un produit cartésien d'ensemble fini. Pour éviter cette longue
étape, nous nous contentons de choisir aléatoirement $q$ circuits
en choisissant les paramètres dans ce   produit cartésien. 
Pour chacun de ces circuits, le dernier carré occupé par la dernière pièce, n'est pas nécessairement égal à l'origine.
% faux ci-dessous !!
%, même si en moyenne, ce résultat est asymptotiquement vrai. 
On se donne $R\in \En^*$ et on ne conserve de ces circuits
que ceux pour lesquels la valeur absolue de l'abscisse et de l'ordonnée est inférieure à $R$.
Pour chacun de ces circuits conservés, on est capable de déterminer tous les circuits à $s$ pièces partant du dernière carré
et retournant à l'origine. On considère donc tous les circuits obtenus par les concaténation des circuits
à $r$ pièces allant de l'origine à un carré quelconque et des circuits à $s$ pièces retournant à l'origine.
Enfin, sur ces circuits, on ne conserve que ceux donc les types des pièces sont inférieurs à $N_j$.
On applique aussi la sélection des isométries et des contraintes de constructibilités.
On a donc obtenu un certain nombre de circuits constructibles à $N$ pièces, sans avoir eu à construire 
le produit cartésien des paramètres des circuits déterminant tous les circuits possibles, dont le cardinal est trop important.
Naturellement, pour augmenter les chances de réussites, 
on doit choisir 
$r$, $s$, $q$ et $R$ les plus grands possibles. Informatiquement, il ne faut pas que ces nombres soient trop importants.
La détermination aléatoire consistera donc à choisir convenablement ces paramètres.
On peut créer soi-même de tels circuits en utilisant 
les exécutables distribués pour windows, cités page \pageref{execatalogueweb}.%
}{%

%%%%%%%%%%%%%%%%%%%%%%%%%%%%%%%%%%%%%%%%%%%%%%%%%%%%%%%%%%%%%
%%%%%%%%%%%%%%%%%%%%%%%%%%%%%%%%%%%%%%%%%%%%%%%%%%%%%%%%%%%%%
\section{Random construction of circuits with a large number of pieces}
\label{constructionaleat}

For values of $N$ smaller than $11$, we are capable of obtaining all of the constructible circuits, and in particular to show them. Another objective of a manufacturer would be to offer a catalogue of train-track designs which may contain circuits with any $N$. Unfortunately, beyond $N=11$, this is no longer conceivable. Manual designs\footnote{see the executables distributed for Windows on page \pageref{execatalogueweb}.} are possible, but tiresome, and non-programmable. We thus propose in this section a way to automatically generate circuits for given $N$ and $N_j$ for values larger than $N=11$, without having to create all of the possible circuits as proposed in Section~\ref{enumerationtouscircuit}, by relying on a random method.

For a given $N\geq 1$, we consider $r,s\in \En^*$ such that $N=r+s$. We are capable of determining all of the circuits with $r$ pieces starting from the origin by describing a Cartesian product of finite sets. To avoid this long stage, we will simply randomly choose $q$ circuits by choosing the parameters in this Cartesian product. For each of these circuits, the last square occupied by the last piece, is not necessarily equal to the origin.
% faux ci-dessous !!
%, even if on average this result is asymptotically true. 
We take $R\in \En^*$ and we keep from these circuits only those for which the absolute value of the abscissa and the ordinate is less than $R$. For each of these retained circuits, we are capable of determining all of the circuits with $s$ pieces starting from the last square and returning to the origin. We therefore consider all of the circuits obtained by the concatenation of the circuits with $r$ pieces going from the origin to any square with the circuits with $s$ pieces returning to the origin. Finally, of these circuits, we keep only those for which they types of pieces are less than $N_j$. We also apply the selection of the isometries and the constructibility constraints. One has hence obtained a certain number of circuits constructible with $N$ pieces, without having had to to construct the Cartesian product of the circuit parameters determining all of the possible circuits, whose cardinality is too great. Naturally, to increase the chances of success, one must choose $r$, $s$, $q$ and $R$ as large as possible. Computationally, it is necessary that these numbers not be too great. The random determination will therefore consist of choosing these parameters appropriately. One may create such circuits oneself using the executables distributed for Windows, quoted on page~\pageref{execatalogueweb}.%

}

\ifcase \choarticle

\begin{example}\
\label{examplesimulation700}

%%%%%%%%%%%%%%%%%%%%%%%%%%%%%%%%%%%%%%%%%%
%\input{simulations_circuit/simulation700}
% fichier tex crée par MaTeXBuild02 le 04-Sep-2015 09:20:42
% à compiler avec 
% MaTeXBuild02('simulation700',0)
% après le fichier 'enumeration_construction_circuit.matex'
% Attention, compilation aléatoire

\begin{figure}[h] 
\begin{center} 
\epsfig{file=circuitaleatoire01.eps, width=11cm}
\end{center} 
\caption{\label{circuitaleatoire01}\iflanguage{french}{Un circuit aléatoire à $13$ pièces.}{A random circuit with $13$ pieces.}} 
\end{figure}

\iflanguage{french}{%
On choisit $N_j=4$ et les paramètres suivants 
\begin{equation*}
r=9,\quad 
s=4,\quad 
q=10,\quad 
R=4.
\end{equation*}
On obtient le circuit aléatoire à $13$ pièces donné en figure \ref{circuitaleatoire01}.%
}{%
We choose $N_j=4$ and the following parameters
\begin{equation*}
r=9,\quad 
s=4,\quad 
q=10,\quad 
R=4.
\end{equation*}
We obtain the random circuits with $13$ pieces given in Figure~\ref{circuitaleatoire01}.%
}
%%%%%%%%%%%%%%%%%%%%%%%%%%%%%%%%%%%%%%%%%%
\end{example}

\or
\fi

\begin{example}\
\label{examplesimulation710}

%%%%%%%%%%%%%%%%%%%%%%%%%%%%%%%%%%%%%%%%%%
%\input{simulations_circuit/simulation710}
% fichier tex crée par MaTeXBuild02 le 04-Sep-2015 09:12:43
% à compiler avec 
% MaTeXBuild02('simulation710',0)
% après le fichier 'enumeration_construction_circuit.matex'
% Attention, compilation aléatoire

\begin{figure}[h] 
\begin{center} 
\epsfig{file=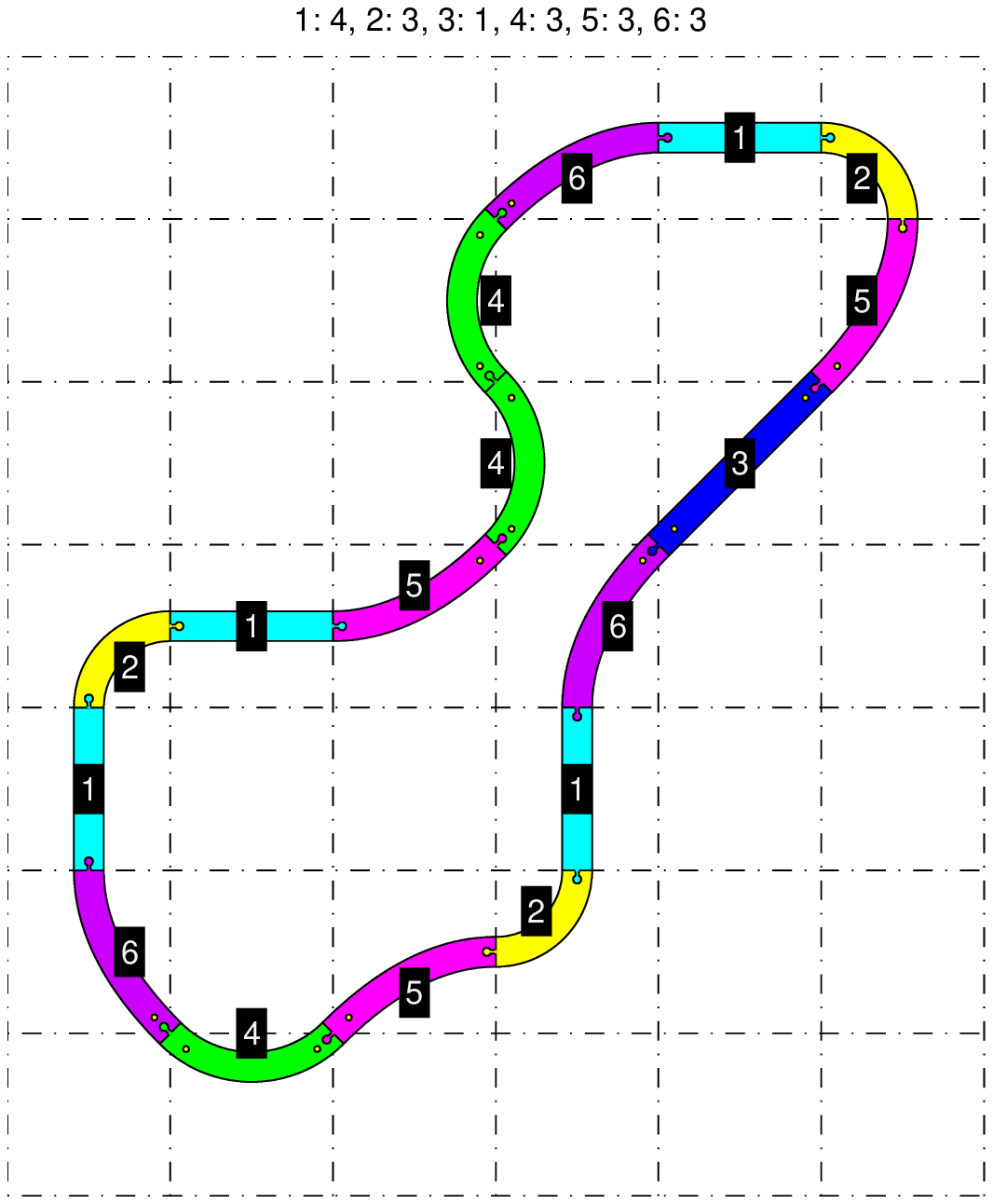, width=11cm}
\end{center} 
\caption{\label{circuitaleatoire02}\iflanguage{french}{Un circuit aléatoire à $17$ pièces.}{A random circuit with $17$ pieces.}} 
\end{figure}

\iflanguage{french}{%
On choisit $N_j=4$ et les paramètres suivants 
\begin{equation*}
r=12,\quad 
s=5,\quad 
q=18,\quad 
R=8.
\end{equation*}
On obtient le circuit aléatoire à $17$ pièces donné en figure \ref{circuitaleatoire02}.%
}{%
We choose $N_j=4$ and the following parameters
\begin{equation*}
r=12,\quad 
s=5,\quad 
q=18,\quad 
R=8.
\end{equation*}
We obtain the random circuit with $17$ pieces given in Figure~\ref{circuitaleatoire02}.%
}

%%%%%%%%%%%%%%%%%%%%%%%%%%%%%%%%%%%%%%%%%%
\end{example}

\iflanguage{french}{%
Nous avons obtenu des circuits de façon aléatoire en pouvant prendre des valeurs de  $N$ strictement plus grandes que  $11$,
et, finalement, en un temps beaucoup plus bref.

D'autres circuits créés aléatoirement avec cette méthode sont disponibles sur Internet
(voir page \pageref{execatalogueweb}).%
}{%

We have obtained some circuits in a random way, being able to take values of $N$ strictly larger than $11$, and, finally, in a much shorter time.

Other circuits created randomly with this method are available on the internet (see page \pageref{execatalogueweb}).%

}
%%% par la suite, essayer d'automatiser l'augmentation du nombre de pièces !

\iflanguage{french}{%

%%%%%%%%%%%%%%%%%%%%%%%%%%%%%%%%%%%%%%%%%%%%%%%%%%%%%%%%%%%%%
%%%%%%%%%%%%%%%%%%%%%%%%%%%%%%%%%%%%%%%%%%%%%%%%%%%%%%%%%%%%%
\section{Généralisations}
\label{generalisation}

%%%%%%%%%%%%%%%%%%%%%%%%%%%%%%%%%%%%%%%%%%%%%%%%%%%%%%%%%%%%%
\subsection{Forme des tuiles du pavage}
\label{formetuile_generalisation}

Nous avons vu que pour un pavage carré, le nombre de courbes nécessaires pour relier chaque point de ${\mathcal{H}}_i$
à tous les autres points distincts de ${\mathcal{H}}_i$ était égal à 5 ou 6, selon que l'on prenne les pièces numérotées
\pieces\ et \piecesb\ ou non.
Ce nombre dépend intrinsèquement du nombre de points de ${\mathcal{H}}_i$ et du cardinal du groupe des isométries laissant le carré invariant.

Se pose la question de savoir si la méthode de construction des rails des circuits étudiés %\textit{Easyloop} 
peut s'appliquer à d'autres types de pavages
que le carré et si l'on est capable de déterminer le nombre de courbes de base, ici égal à 5 ou 6, uniquement à partir du pavage
et des points ${\mathcal{H}}_i$ considérés.
Cette généralisation est aussi évoquée pour des chemins autoévitants dans \cite{MR2104301},
dans un cas plus simple, puisque les circuits ne peuvent passer que par les bords des tuiles constituant le pavage du plan.%

}{%

%%%%%%%%%%%%%%%%%%%%%%%%%%%%%%%%%%%%%%%%%%%%%%%%%%%%%%%%%%%%%
%%%%%%%%%%%%%%%%%%%%%%%%%%%%%%%%%%%%%%%%%%%%%%%%%%%%%%%%%%%%%
\section{Generalizations}
\label{generalisation}

%%%%%%%%%%%%%%%%%%%%%%%%%%%%%%%%%%%%%%%%%%%%%%%%%%%%%%%%%%%%%
\subsection{Shape of the tiling}
\label{formetuile_generalisation}

We have seen that for a square tiling, the number of curves necessary to connect each point of ${\mathcal{H}}_i$ to every other distinct point in ${\mathcal{H}}_i$  was equal to 5 or 6, according to whether or not one takes the pieces numbered \pieces\ and \piecesb. This number depends intrinsically on the number of points in ${\mathcal{H}}_i$ and on the cardinality of the group of the isometries leaving the square invariant.

The question arises whether or not the method of constructing the % \textit{Easyloop} 
rails of the studied circuits
can be applied to types of tiling other than the square, and if one is capable of determining the number of basic curves, here equal to 5 or 6, uniquely from the tiling and the points ${\mathcal{H}}_i$ considered. This generalization is also mentioned for self-avoiding walks in \cite{MR2104301}, which is a simpler case since the circuits may only follow the edges of the tiles constituting the tiling of the plane.%

}

\begin{figure}
\centering
%%% sous figure 1
\subfigure[\label{autrepavagetrianglea}\iflanguage{french}{ : avec les 3 milieux}{: with the 3 middles}]
{\epsfig{file=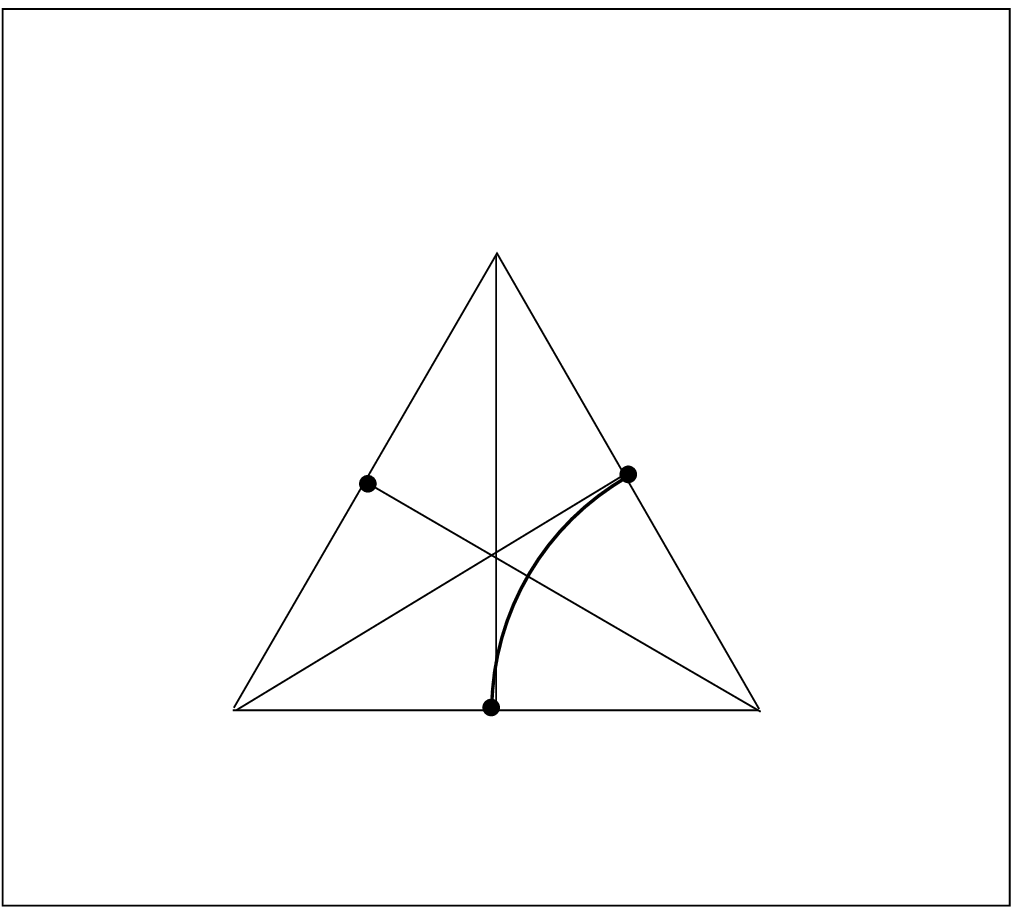, width=7cm}}
\qquad
%%% sous figure 2
\subfigure[\label{autrepavagetriangleb}\iflanguage{french}{ : avec les 3 milieux et les 3 sommets}{: with the 3 middles and the 3 vertices}]
{\epsfig{file=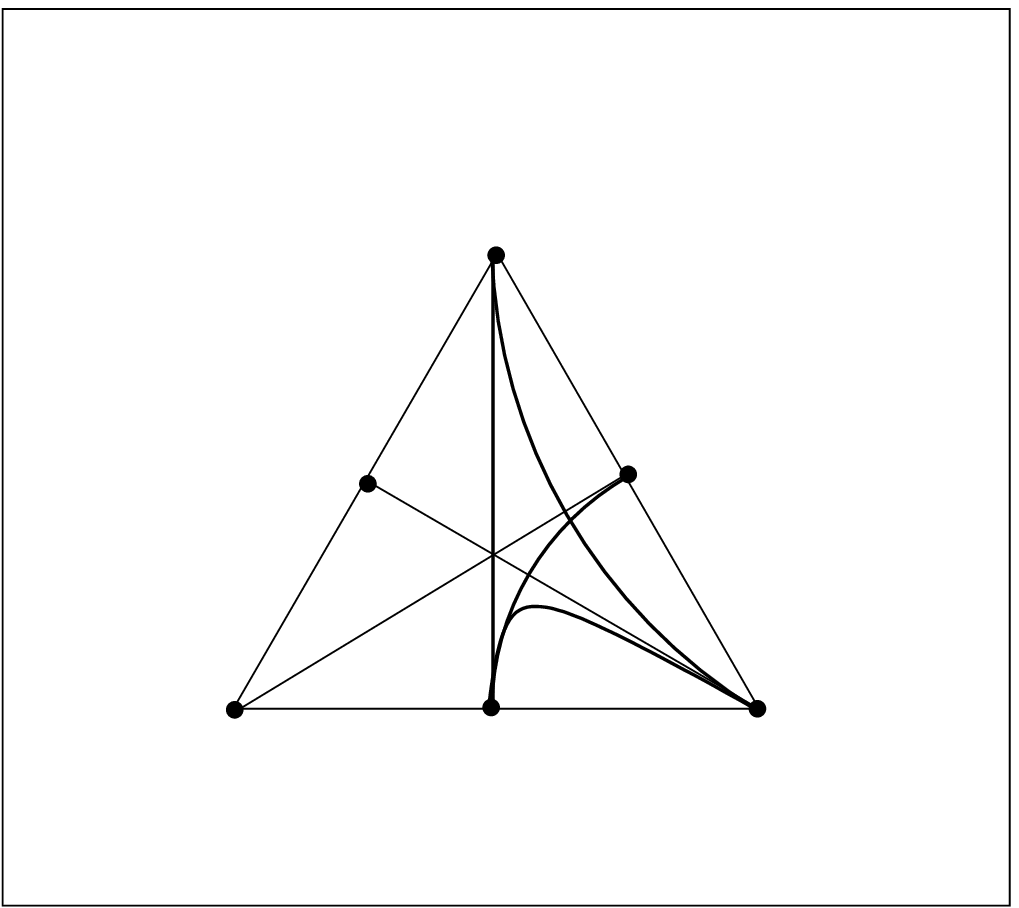, width=7cm}}
\caption{\label{autrepavagetriangle}\iflanguage{french}{Autres pavages possibles : le triangle équilatéral}{Other possible tilings: the equilateral triangle}}
\end{figure}
\iflanguage{french}{%
Par exemple, on peut considérer des pavages du plan
par des triangles équilatéraux, en ne prenant que les 3 milieux des 3 côtés (voir figure \ref{autrepavagetrianglea}) ou les 3 milieux et les 3  sommets  du triangle
(voir figure \ref{autrepavagetriangleb}).
On impose que la courbe passe par deux points distincts de cet ensemble ${\mathcal{H}}_i$, en étant tangente 
à la droite reliant ce point avec le centre du triangle.
Dans le premier cas, une seule courbe est nécessaire, dans le second, 3 le sont.
D'autres solutions peuvent être envisagées avec d'autres types de pavages possibles.

Le nombre de courbes nécessaires peut-il s'exprimer en fonction du polygone de pavage et de la nature des points ${\mathcal{H}}_i$ ?
L'énumération de circuits, construits sur ces méthodes, semble être de nouveau, un problème ouvert.%
}{%
For example, one may consider tiling the plane with equilateral triangles, taking only the 3 middles of the 3 sides (see Figure~\ref{autrepavagetrianglea}) or the 3 middles and the 3 vertices of the triangle (see Figure~\ref{autrepavagetriangleb}). We impose that the curve passes through two distinct points of this set ${\mathcal{H}}_i$, while being tangent to the straight line connecting this point with the center of the triangle. In the first case, a single curve is necessary, while in the second, 3 are. Other solutions can be envisaged, with other types of possible tilings.

Can the number of necessary curves be expressed as a function of the tiling polygon and the nature of the points ${\mathcal{H}}_i$?
The enumeration of circuits, constructed using these methods, seems once again to be an open problem.%
}

\iflanguage{french}{%

%%%%%%%%%%%%%%%%%%%%%%%%%%%%%%%%%%%%%%%%%%%%%%%%%%%%%%%%%%%%%
\subsection{Présence  d'aiguillages}
\label{geneaiguillage}

Restons dans le cas du pavage carré.
Comme décrit plus haut, les pièces des  circuits constructibles 
ne doivent pas avoir d'extrémités %correspondant à des milieux de côtés de carrés 
en commun, si 
la présence d'aiguillages n'est pas prévue. On peut, au contraire, laisser certaines
extrémités %correspondant à des milieux de côtés de carrés 
être en commun à plusieurs pièces, ce qui
correspond simplement à autoriser des aiguillages.
\input{simulations_circuit/simulation403}
Dans ce cas, les circuits orientés peuvent contenir plusieurs boucles, ce qui en fait des graphes orientés dont
l'énumération est un problème beaucoup plus ardu, à cause de la multiplicité des types d'aiguillages possibles
et donc des types de graphes possibles.
Voir la figure \ref{aiguillage_boucles} qui montre un exemple d'un circuit avec plusieurs 
\ifcase \choarticle
aiguillages ainsi que la figure \ref{aiguillage_bouclesb} qui illustre l'une des branches à 6 pièces
de ce circuit.
\or
aiguillages.
\fi
}{%

%%%%%%%%%%%%%%%%%%%%%%%%%%%%%%%%%%%%%%%%%%%%%%%%%%%%%%%%%%%%%
\subsection{Presence of switches}
\label{geneaiguillage}

We keep to the case of square tiling. As described earlier, if the presence of switches is not anticipated, the pieces of constructible circuits must not have extremities 
%corresponding to the middles of sides of squares  
in common. One can, %on the other hand, 
on the contrary, 
allow certain extremities 
%corresponding to the middles of sides of squares 
be common to several pieces, which simply corresponds to allowing switches.
%%%%%%%%%%%%%%%%%%%%%%%%%%%%%%%%%%%%%%%%%%%%%%%%%
%\input{simulations_circuit/simulation403}
% fichier tex crée par MaTeXBuild02 le 04-Sep-2015 05:56:02
% A compiler avec MaTeXBuild02('simulation403',0)
% après le fichier 'enumeration_construction_circuit.matex'
% Sortie statique, pas très pertinent avec matex !!

% pour faire la figure : exemple_aiguillage.eps 
% Lancer  
% global impnb;impnb=2;
% naiguillage=3;
%[ntyp,ctot,ptot,ittot,rottot,symettot,nfemtot,casex,casret]=...
%    exhaustif_circuit_aiguillage(naiguillage,9,[],[],[],5,1,1)
% et faire eps en cours d'execution 
% Attendre "Avec 6 pièces"

\ifcase \choarticle

\begin{figure}
\centering
%%% sous figure 1
\subfigure[\label{aiguillage_boucles}\iflanguage{french}{Un exemple de plan de circuit avec
2 aiguillages}{An example track design with 2 switches}]
{\epsfig{file=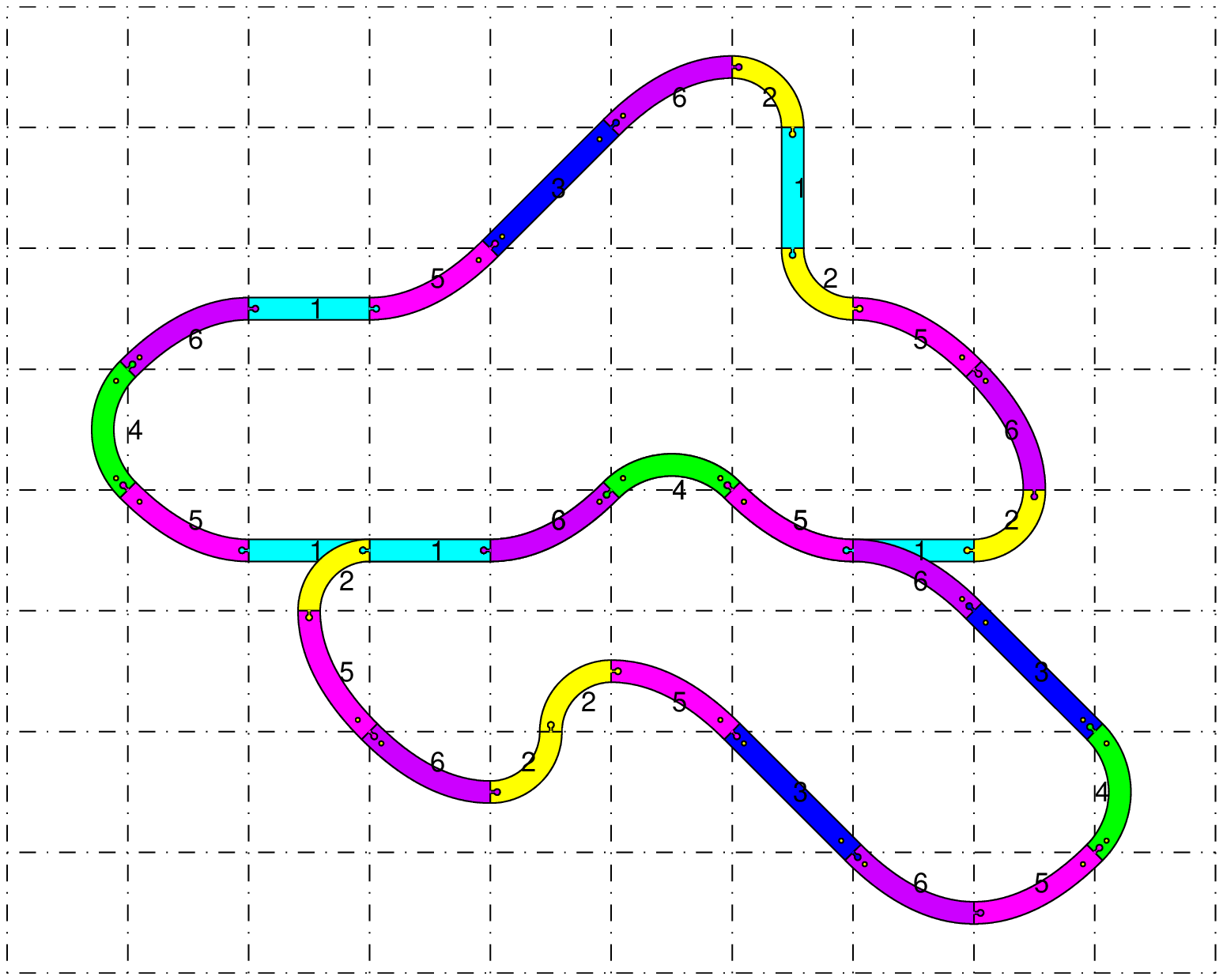, width=12cm}}
\qquad
%%% sous figure 2
\subfigure[\label{aiguillage_bouclesb}\iflanguage{french}{Circuit avec
$3$ aiguillages : l'ensemble de l'une des branches à 6 pièces}{Circuit with  $3$ switches: the set of branches with 6 pieces}]
{\epsfig{file=exemple_aiguillage.eps, width=12cm}}
\caption{\label{aiguillage_bouclestot}\iflanguage{french}{Aiguillage.}{Switches.}}
\end{figure}

\or

\begin{figure}[h] 
\begin{center} 
\epsfig{file=aiguillage_boucles.eps, width=12cm}
\end{center} 
\caption{\label{aiguillage_boucles}\iflanguage{french}{Un exemple de plan de circuit avec
2 aiguillages.}{An example track design with 2 switches}}
\end{figure}

\fi

%%%%%%%%%%%%%%%%%%%%%%%%%%%%%%%%%%%%%%%%%%%%%%%%%
In this case, the oriented circuits may contain several loops, which makes them directed graphs, whose enumeration is a much more arduous problem due to the multiplicity of the types of possible switches and therefore the types of possible graphs. See Figure~\ref{aiguillage_boucles}, which shows an example of a circuit with 
\ifcase \choarticle
switches, as well as Figure~\ref{aiguillage_bouclesb} which illustrates for this circuit the set branches with 6 pieces.
\or
switches.
\fi
}

%%%%%%%%%%%%%%%%%%%%%%%%%%%%%%%%%%%%%%%%%%%%%%%%%%%%%%%%%%%%%
%%%%%%%%%%%%%%%%%%%%%%%%%%%%%%%%%%%%%%%%%%%%%%%%%%%%%%%%%%%%%
\section{Conclusion}
\label{conclusion}%

\iflanguage{french}{%
La question \questionindustrielle, simple à exprimer est plus difficile à résoudre. 
Nous y sommes parvenus, dans cet article : l'énumération et la construction de tels
circuit sont possibles et ont été implémentées informatiquement jusqu'à  $N=11$.
Au-delà, une estimation  du nombre de circuits possibles à été fournie (voir \eqref{estimationeq02}, correspondant 
au cas $N_j=4$  et \eqref{estimationeq03}).
Notons que ces nombres correspondent  aux nombres de circuits
qui contiennent exactement $24$  pièces.
Si l'on veut dénombrer tous les circuits qui contiennent au plus 
$24$  pièces, il suffit de faire la somme de la dernière colonne du tableau 
\ref{enumeration_construction_circuit_tab_new02}, ce qui donne 
$1102$ circuits, puis d'appliquer l'estimation 
pour $N$ variant de  $12$ à $24$,
avec les paramètres estimés donnés par \eqref{refalpha},
ce qui donne au total
\begin{equation}
\label{estimationeq02final}
1873804310490,
\end{equation}%
}{%
The question \questionindustrielle\ 
is simple to express, and more difficult to resolve. We have succeeded in this article: the enumeration and the construction of such circuits is possible and have been implemented computationally up to $N=11$. Beyond that, an estimation of the number of possible circuits has been provided (see~\eqref{estimationeq02}, corresponding to the case $N_j=4$, and~\eqref{estimationeq03}). We note that these numbers correspond to the number of circuits which contain exactly $24$ pieces. If we want to tally all of the circuits which contain at most $24$ pieces, it is sufficient to sum the last column in Table~\ref{enumeration_construction_circuit_tab_new02}, which gives $1102$ circuits, then to apply the estimation for $N$ varying from $12$ to $24$, with the estimated parameters given by~\eqref{refalpha}, which gives in total
\begin{equation}
\label{estimationeq02final}
1873804310490,
\end{equation}%
}
%\iflanguage{french}{}{\textbf{ERREUR !!! Traduire cela dans la langue choisie !!!!!}}
\iflanguage{french}{%
soit un total de plus de 
\begin{equation}
\label{estimationeq02finaltlter}
\text{$\numberstringnum{1}$ billion   de  circuits réalisables
avec $24$ pièces.}
\end{equation}%
}{%
that is, a total of more than
\begin{equation}
\label{estimationeq02finaltlter}
\text{$\numberstringnum{1}$ trillion  feasible circuits 
with $24$ pieces.}
\end{equation}%
}
\iflanguage{french}{%
De plus, une construction aléatoire de circuit a été proposée permettant d'obtenir des valeurs de $N$
strictement plus grandes que $N=11$.
Des exécutables et un catalogue de circuits sont disponibles sur internet.
\input{simulations_circuit/simulation402}
Notons enfin, qu'à la main, 
\ifcase \choarticle
quatre
\or
deux 
\fi
circuits constructibles correspondant à 
$N_j=4$ contenant
exactement le nombre maximum de pièces ($24$)
ont être créés.
\ifcase \choarticle
Voir figures \ref{exemple_circuit_plan_24max01}--\ref{exemple_circuit_plan_24max04}.
\or
Voir figures \ref{exemple_circuit_plan_24max02}
et \ref{exemple_circuit_plan_24max03}.
\fi

Il est intéressant de constater que la tradionnelle théorie des chemins autoévitants correspond presque
aux circuits de trains déjà existants (cf. section \ref{compartradi}) tandis que le jeu breveté 
correspond à la notion de \slw.

Il resterait à améliorer les algorithmes d'énumération de circuits pour obtenir des valeurs
de $N$ plus élevées, dans le cas déterministe, en tentant par exemple d'éviter
la très longue énumération de circuits possibles ; est-ce qu'une construction
directe des circuits constructibles sans passer par cette énumération est possible ?
Une application des techniques parallélisables proposées par 
G. Slade, 
I. Jensen, 
ou 
A. J. Guttmann 
pourrait être essayée sur les circuits afin d'augmenter les valeurs de $N$ pour lesquelles les énumérations
de circuits sont exactes.

Il serait intéressant de démontrer si l'estimation 
\eqref{estime_nombre_circuit}
est valable, avec un éventuel calcul des constantes $A$, $\mu$ et $\gamma$.
La généralisation évoquée en section \ref{generalisation} permettrait de créer d'autres types de circuits, mais aussi d'essayer de comprendre
la nature algébrique du système proposé avec des carrés.%
}{%
In addition, a random circuit construction has been offered, allowing values of $N$ strictly larger than $N=11$ to be obtained. Some executables and a catalogue of circuits are available online.
%%%%%%%%%%%%%%%%%%%%%%%%%%%%%%%%%%%%%%
%\input{simulations_circuit/simulation402}
% fichier tex crée par MaTeXBuild02 le 13-Sep-2015 11:04:09
% A compiler avec MaTeXBuild02('simulation402',0)
% après le fichier 'enumeration_construction_circuit.matex'
% Sortie statique, pas très pertinent avec matex !!

\ifcase \choarticle
\begin{figure}[h] 
\begin{center} 
\epsfig{file=exemple_circuit_plan_24max01.eps, width=11cm}
\end{center} 
\caption{\label{exemple_circuit_plan_24max01}\iflanguage{french}{Un exemple de plan de circuit avec
$24$ pièces.}{An example track design with $24$ pieces.}}
\end{figure}
\or
\fi

\begin{figure}[h] 
\begin{center} 
\epsfig{file=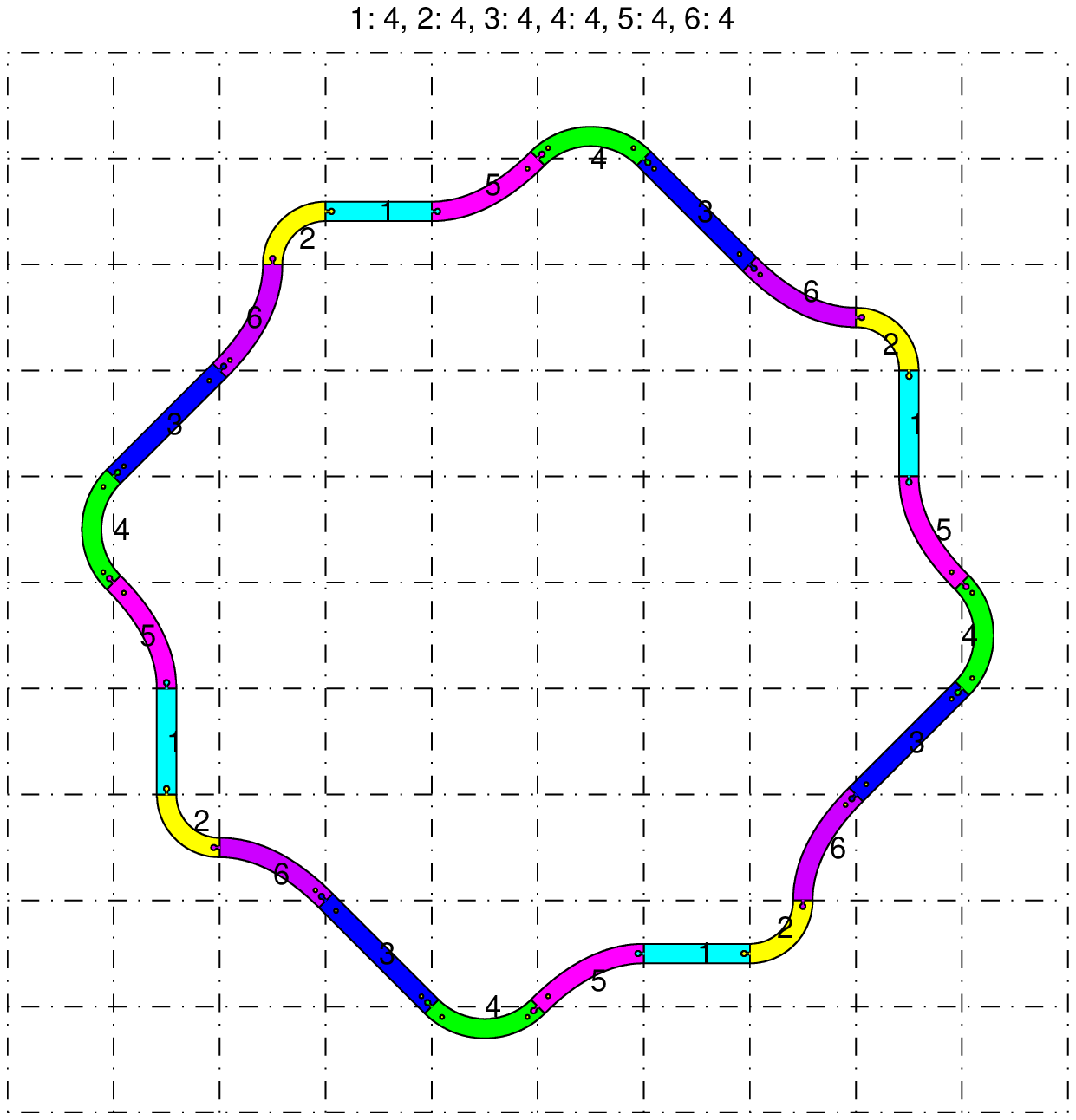, width=11cm}
\end{center} 
\caption{\label{exemple_circuit_plan_24max02}\iflanguage{french}{Un exemple de plan de circuit avec
$24$ pièces.}{An example track design with $24$ pieces.}}
\end{figure}

\begin{figure}[h] 
\begin{center} 
\epsfig{file=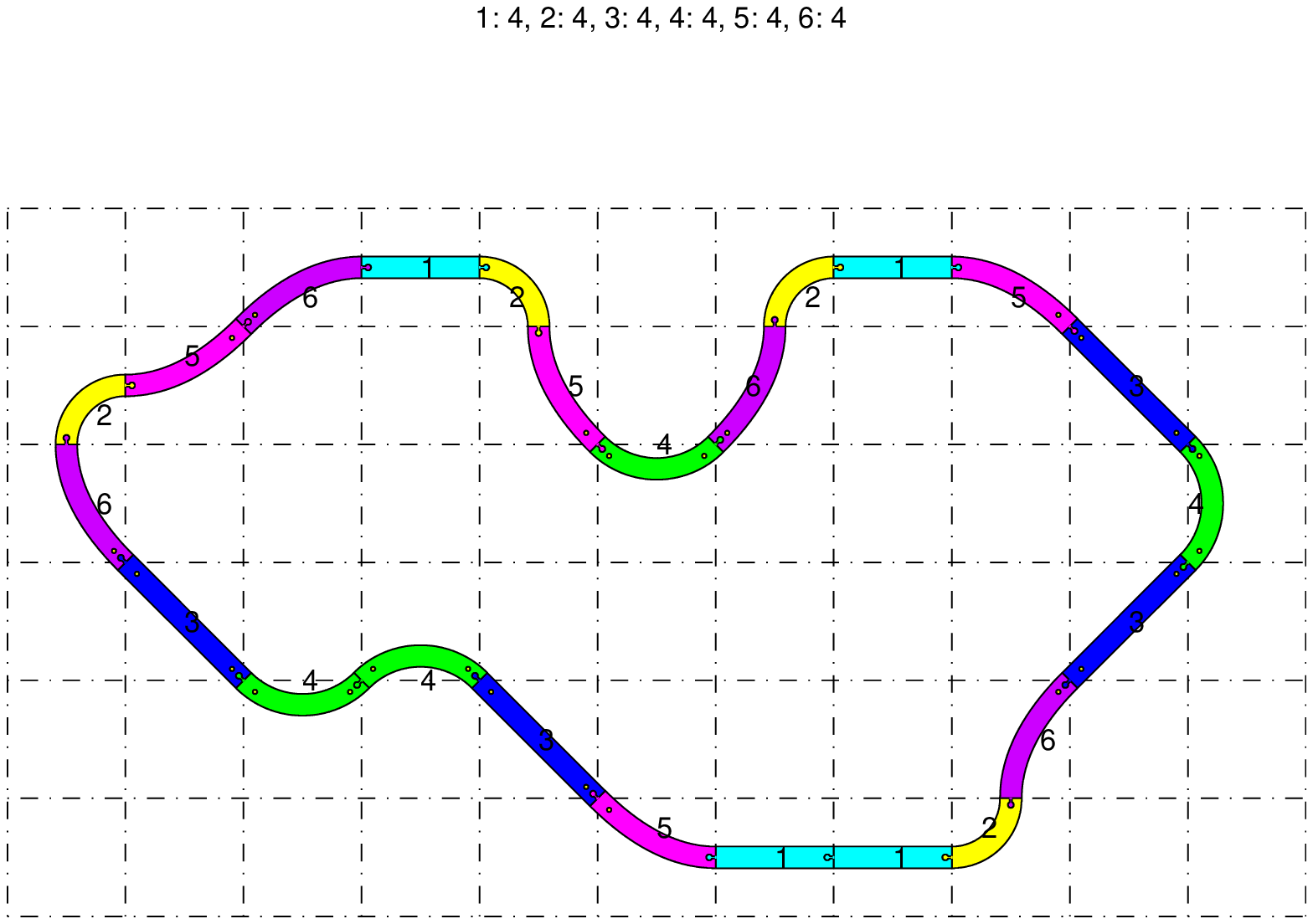, width=14cm}
\end{center} 
\caption{\label{exemple_circuit_plan_24max03}\iflanguage{french}{Un exemple de plan de circuit avec
$24$ pièces.}{An example track design with $24$ pieces.}}
\end{figure}

\ifcase \choarticle
\begin{figure}[h] 
\begin{center} 
\epsfig{file=exemple_circuit_plan_24max04.eps, width=11cm}
\end{center} 
\caption{\label{exemple_circuit_plan_24max04}\iflanguage{french}{Un exemple de plan de circuit avec
$24$ pièces.}{An example track design with $24$ pieces.}}
\end{figure}
\or
\fi
%%%%%%%%%%%%%%%%%%%%%%%%%%%%%%%%%%%%%%
Finally, we note that 
\ifcase \choarticle
four
\or
two 
\fi
constructible circuits corresponding to $N_j=4$, and containing exactly the maximum number of pieces ($24$) have been created by hand. 
See 
\ifcase \choarticle
Figures~\ref{exemple_circuit_plan_24max01}--\ref{exemple_circuit_plan_24max04}.
\or
Figures~\ref{exemple_circuit_plan_24max02} and \ref{exemple_circuit_plan_24max03}.
\fi

% Attention, anglais ci-dessous non traduit par ml 
It is interesting that the traditional theory of self-avoiding paths corresponds almost to existing trains circuits (see section \ref{compartradi}) while the patend studied system 
corresponds to the notion of \slw.

It remains to improve the circuit enumeration algorithms to obtain higher values of $N$, in the deterministic case, by trying, for example, to avoid the very long enumeration of possible circuits; is a direct construction of constructible circuits possible, without going through this enumeration? 
An application 
%A application 
of the parallelizable techniques proposed by
G. Slade, 
I. Jensen, 
or 
A. J. Guttmann 
might be tried on the circuits in order to increase the values of $N$ for which the circuit enumerations are exact.

It would be interesting to prove if estimation \eqref{estime_nombre_circuit} is valid, with an eventual calculation of the constants $A$, $\mu$ and $\gamma$. The generalization raised in Section~\ref{generalisation} allows the creation of other types of circuits, but also an attempt to understand the algebraic nature of the system proposed with squares.%
}

%%%%%%%%%%%%%%%%%%%%%%%%%%%%%%%%%%%%%%%%%%%%%%%%%%%%%%%%%%%%
%%%%%%%%%%%%%%%%%%%%%%%%%%%%%%%%%%%%%%%%%%%%%%%%%%%%%%%%%%%%%
% Bibliograhie  faite avec biblatex
\printbibliography

\end{document}